\documentclass[11pt,twoside]{book}

\usepackage[T1]{fontenc} 
\usepackage{fix-cm} 
\usepackage{a4wide}
\usepackage{stmaryrd}
\usepackage{mathrsfs}
\usepackage{amsmath}
\usepackage{amssymb}
\usepackage{mathtools}
\usepackage{bm}
\usepackage[normalem]{ulem}
\usepackage{graphicx}
\usepackage{tikz}
\usetikzlibrary{cd}
\usepackage{tikz-cd}
\usepackage{shadow}
\usepackage{fancybox}
\usepackage{scalerel,stackengine}
\usepackage{eso-pic}

\usepackage{setspace}
\usepackage{fancyhdr}
\usepackage[nottoc]{tocbibind}
\usepackage{alltt}
\usepackage{amsthm}

\usepackage{hyperref}
\usepackage[nameinlink]{cleveref}

\usepackage{epigraph}
\usepackage{makeidx}
\usepackage{aliascnt}
\usepackage{multicol}
\usepackage{dcolumn}

\stackMath

\hypersetup{
	unicode=false,          
	pdftoolbar=true,        
	pdfmenubar=true,        
	pdffitwindow=false,     
	pdfstartview={FitH},    
	pdftitle={My title},    
	pdfauthor={Author},     
	pdfsubject={Subject},   
	pdfcreator={Creator},   
	pdfproducer={Producer}, 
	pdfkeywords={keywords}, 
	pdfnewwindow=true,      
	colorlinks=true,       
	linkcolor=blue,          
	citecolor=blue,        
	filecolor=magenta,      
	urlcolor=green           
}

\newcommand\reallywidehat[1]{%
	\savestack{\tmpbox}{\stretchto{%
			\scaleto{%
				\scalerel*[\widthof{\ensuremath{#1}}]{\kern-.6pt\bigwedge\kern-.6pt}%
				{\rule[-\textheight/2]{1ex}{\textheight}}
			}{\textheight}%
		}{0.5ex}}%
	\stackon[1pt]{#1}{\tmpbox}%
}

\theoremstyle{definition}
\newtheorem{theorem}{Theorem}[section]
\newtheorem{lemma}[theorem]{Lemma}
\newtheorem{proposition}[theorem]{Proposition}
\newtheorem{corollary}[theorem]{Corollary}
\newtheorem{definition}[theorem]{Definition}
\newtheorem{example}[theorem]{Example}
\newtheorem{remark}[theorem]{Remark}

\DeclareMathOperator{\pr}{pr}
\DeclareMathOperator{\op}{op}

\makeatletter
\newcommand{\doublewidetilde}[1]{{%
		\mathpalette\double@widetilde{#1}%
}}
\newcommand{\double@widetilde}[2]{%
	\sbox\z@{$\m@th#1\widetilde{#2}$}%
	\ht\z@=.9\ht\z@
	\widetilde{\box\z@}%
}
\makeatother

\newcommand{\mc}{\mathcal}
\newcommand{\mf}{\mathfrak}
\newcommand{\mb}{\mathbb}
\newcommand{\xra}{\xrightarrow}
\newcommand{\ra}{\rightarrow}
\newcommand{\rra}{\rightrightarrows}

\newcommand {\apgt} {\ {\raise-.5ex\hbox{$\buildrel>\over\sim$}}\ }
\newcommand {\aplt} {\ {\raise-.5ex\hbox{$\buildrel<\over\sim$}}\ }

\newcommand{\addsymbol}{}
\def\addsymbol #1: #2#3{$#1$ \> \parbox{5in}{#2 \dotfill \pageref{#3}}\\}

\makeatletter
\newcommand\myfnt{\ifx\protect\@typeset@protect\expandafter\footnote\else\expandafter\@gobble\fi}
\makeatother

\makeindex

\allowdisplaybreaks
\numberwithin{equation}{section}
\makeatletter
\newcommand*{\relrelbarsep}{.386ex}
\newcommand*{\relrelbar}{%
	\mathrel{%
		\mathpalette\@relrelbar\relrelbarsep
	}%
}
\newcommand*{\@relrelbar}[2]{%
	\raise#2\hbox to 0pt{$\m@th#1\relbar$\hss}%
	\lower#2\hbox{$\m@th#1\relbar$}%
}
\providecommand*{\rightrightarrowsfill@}{%
	\arrowfill@\relrelbar\relrelbar\rightrightarrows
}
\providecommand*{\leftleftarrowsfill@}{%
	\arrowfill@\leftleftarrows\relrelbar\relrelbar
}
\providecommand*{\xrightrightarrows}[2][]{%
	\ext@arrow 0359\rightrightarrowsfill@{#1}{#2}%
}
\providecommand*{\xleftleftarrows}[2][]{%
	\ext@arrow 3095\leftleftarrowsfill@{#1}{#2}%
}
\makeatother

\begin{document}
	
	\pagenumbering{roman}

\thispagestyle{empty}
\begin{center}{\textbf{\LARGE{Geometric structures on Lie groupoids and differentiable stacks}}}
\end{center}{ \par}

\vspace*{\fill}

\bigskip{}
\begin{center}{
{\it {\large A thesis submitted for the degree of}}\\
\vspace {0.5cm}
{{\Large\bf Doctor of Philosophy}}\\
\vspace {0.5cm}
{\it {\large in the}}\\
\vspace {0.5cm} {\bf {\Large School of Mathematics}}}
\end{center}
\vspace*{\fill}
\begin{center}
\end{center}
\vspace*{\fill}
\begin{center}{\textbf{\Large {\bf Praphulla Koushik}}}\end{center}{\Large \par}

\vspace*{\fill}

\begin{center}
\includegraphics[scale=0.5]{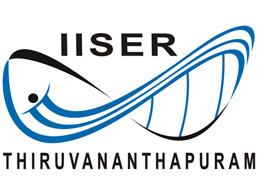}{\small}\\
{\small SCHOOL OF MATHEMATICS}\\
{\small INDIAN INSTITUTE OF SCIENCE EDUCATION AND RESEARCH}\\
{\small THIRUVANANTHAPURAM - 695551}\\
{\small INDIA}\end{center}{\small
\par}

\begin{center}{\small October 2021}\end{center}{\small \par}
\clearpage
\newpage
\thispagestyle{empty}
\bigskip{}
\vspace{1cm} \cleardoublepage
\thispagestyle{plain}
\begin{center}\textbf{\huge Declaration}\end{center}{\huge \par}
\vspace{1cm} I hereby declare that the work reported in this thesis
is original and was carried out by me during my tenure as a PhD
student at the \textbf{School of Mathematics, Indian Institute of
Science Education and Research Thiruvananthapuram}. Such materials that
have been obtained from other sources are duly acknowledged in the thesis.\\

\vspace{1in} \noindent Place: IISER Thiruvananthapuram\hfill {\bf
Praphulla Koushik}\\\smallskip Date\ : October 2021 \hfill
\\\smallskip

\vspace{0.05cm}

\noindent
{\bf Certified} that the work incorporated in the thesis

\vspace{0.1cm}

\begin{center} {\bf ``Geometric structures on Lie groupoids and differentiable stacks"}
\end{center}

\vspace{0.05cm}

\noindent submitted by {\bf Praphulla Koushik} was carried out by the
candidate under my supervision. This thesis has not formed the basis
for the award of any degree, of any university or institution.

\vspace{1in} \noindent Place: IISER Thiruvananthapuram\hfill {\bf
Dr. Saikat Chatterjee}\\ \smallskip Date\ : October 2021 \hfill(Thesis supervisor)\\
\smallskip

	\newpage
	\thispagestyle{empty}
	\bigskip{}
	\vspace{1cm} \cleardoublepage

	\newpage
	
	\thispagestyle{plain}
	\begin{center}\textbf{\Large  Acknowledgments}\end{center}{\huge \par}
	\bigskip{}
	
	First of all, I would like to thank my Ph.D. thesis supervisor Saikat Chatterjee for his support. I want to thank Indranil Biswas and Frank Neumann for giving me the oppurtunity to work  with them in collaboration.

	I want to thank all the members of the school of mathematics, IISER Trivandrum, in particular, my Ph.D. thesis committee members Viji Zac Thomas, Sarbeswar pal, and all the former \& present heads of SOM (namely Uptal Manna, P. Devaraj, Viji Zac Thomas) for their continuous support and encouragement. Special thanks to supporting staff (maintenance, mess, and security) of IISER Trivandrum. 
	
	I want to thank all members of the Department of Mathematics, University of Hyderabad (2008-2013), in particular Rajat Tandon, TKS Moothathu, S. Kumaresan, T. Amarnath for their support.  Special thanks to the staff of Indira Gandhi Memorial Library University of Hyderabad. I want to thank all members of the Department of Mathematics, IIT Bombay (2014-2016), in particular, U. K. Anandavardhanan, Anant R. Shastri, and S.G. Dani  for their support.
	
	I want to thank the community of \href{https://mathoverflow.net/}{Mathoverflow}, in particular, David Roberts, Ben McKay, Todd Trimble, Dmitri Pavlov, Mamuka Jibladze, for their questions, comments, and answers.
	
	This acknowledgment would be incomplete without offering my thanks to 
	my friends who made my stay in University of Hyderabad, IIT Bombay, IISER Trivandrum memorable, in particular, I mention the names of, Chiranjeevi, Sampath, Ambika, Bhargav, Gobinda, Sudeshna, Vivek, Pramod, Debasmita, Ankush, Adittya, Mahendra, Soham, Sunil, and Yogesh.
	\quad \\ \quad \\ \indent\qquad\qquad\qquad\qquad\qquad\qquad\qquad
	\qquad\qquad\qquad\qquad\qquad\quad Praphulla Koushik
	
	\cleardoublepage
	
	\pagenumbering{arabic}
	
	\tableofcontents
	
	\chapter{Introduction}

In this thesis, we have studied certain geometric structures over Lie groupoids and differentiable stacks.

A Lie groupoid is a groupoid $\mc{G}=[\mc{G}_1\rra \mc{G}_0]$, where $\mc{G}_0,\mc{G}_1$ are smooth manifolds, the maps $s,t:\mc{G}_1\ra \mc{G}_0$ are submersions, and all other structure maps are smooth maps. The notion of a Lie groupoid covers many interesting geometric structures over a manifold. For example,  \begin{itemize}
	\item a Lie group action on a manifold $M$ assigns a Lie groupoid whose base is $M$,
	\item a principal $G$-bundle over a manifold $M$ assigns a Lie groupoid whose base is $M$,
	\item a vector bundle over a manifold $M$ assigns a Lie groupoid whose base is $M$,
	\item a foliation on a manifold $M$ assigns a Lie groupoid whose base is $M$.
\end{itemize} 
Some standard ``textbook references'' for the notion of Lie groupoids are \cite{MR2012261, MR896907, MR2157566}. 

The geometric structures over Lie groupoids that we study in this thesis are Lie groupoid extensions and principal bundles over Lie groupoids (\Cref{Chap.2,Chap.3}). 

A Lie groupoid extension is a morphism of Lie groupoids $(F,1_M):[\mc{G}_1\rra M]\ra [\mc{H}_1\rra M]$, such that, $F:\mc{G}_1\ra \mc{H}_1$ is (at least) a surjective submersion. Various author require different conditions on the map $F$. We ask $F$ to be a surjective submersion. The notion of (and other geometric structures on) Lie groupoid extensions is discussed, for example, in \cite{Moerdijk2,MR2493616,MR3480061}.

Let $G$ be a Lie group. The principal $G$-bundles over a manifold have a rich geometric structure, and many interesting questions arise out of it. Motivated by the notion of principal bundle over a manifold, we introduce the notion of a principal $G$-bundle over a Lie groupoid. Let $[X_1\rra X_0]$ be a Lie groupoid. A principal $G$-bundle over $[X_1\rra X_0]$ is given by a principal $G$-bundle $E_G\ra X_0$, over the manifold $X_0$, with an action of the Lie groupoid $[X_1\rra X_0]$ on $E_G$, such that the action of $[X_1\rra X_0]$ and $G$ on $E_G$ are compatible. The notion of principal $G$-bundles over a Lie groupoid cover the following classical cases of principal bundles:
\begin{itemize}
	\item a principal $G$-bundle over a Lie groupoid $[M\rra M]$ is a principal $G$-bundle over the manifold $M$,
	\item a principal $G$-bundle over a Lie groupoid $[H\times M\rra M]$ is an $H$-equivariant principal $G$-bundle over the manifold $M$. 
\end{itemize}
Some of the papers that discuss (geometric structures on) principal bundle over Lie groupoids are \cite{MR3150770,MR2270285,MR1950948,MR2119241}.

Let $\mc{C}$ be a category, and $\mc{J}$ be a Grothendieck topology on $\mc{C}$. We call the pair $(\mc{C},\mc{J})$ to be a site. A stack over the site $(\mc{C},\mc{J})$ can be thought of as a ``generalized object'' of the category $\mc{C}$.  To explain the terminology of ``generalized object'' of the category $\mc{C}$, we need to recall the Yoneda lemma (Yoneda embedding) of category theory.

Let $\mc{C}$ be a (locally-small) category. Let $A$ be an object of $\mc{C}$. Consider the functor of points $h_A:\mc{C}^{\op}\ra \text{Set}$. This assigns a functor $\mc{C}\ra [\mc{C}^{\op},\text{Set}]$. The Yoneda embedding says that this functor $\mc{C}\ra [\mc{C}^{\op},\text{Set}]$ is an embedding of categories. So, we see the objects of $[\mc{C}^{\op},\text{Set}]$ as ``generalized objects of $\mc{C}$''. 

Once we fix a Grothendieck topology on $\mc{C}$, we can talk about the notion of covering of an object $U$ of $\mc{C}$. Using this notion of a covering, we can restrict our attention to functors $\mc{C}^{\op}\ra \text{Set}$ that are ``locally determined''. A functor $\mc{C}^{\op}\ra \text{Set}$ that is locally determined is called a sheaf on the site $(\mc{C},\mc{J})$. So, sheaves on $\mc{C}$ should be seen as a generalized object of $\mc{C}$. We want to ensure that this Grothendieck topology is such that the functors $h_A$, for each object $A$, are ``locally determined''. A Grothendieck topology on $\mc{C}$ for which $h_A$ is locally determined for each object $A$ of $\mc{C}$ is called a subcanonical site. 

Let us further generalize the notion of sheaf by considering  ``maps'' from $\mc{C}$ that takes values in categories instead of sets. By ``maps'', we actually mean pseudo-functors on $\mc{C}$. A pseudo-functor on a category consists of the following data:
\begin{itemize}
	\item a category for each object $U$ of $\mc{C}$,
	\item a functor for each morphism $U\ra V$ of $\mc{C}$,
\end{itemize}
satisfying certain compatibility conditions. Fixing a Grothendieck topology $\mc{J}$ on $\mc{C}$, we can also declare the notion for a pseudo-functor to be ``locally determined''. We call a locally determined pseudo-functor to be a stack on the site $(\mc{C},\mc{J})$. Some standard references for the notion of stacks are \cite{MR2223406,Johan,laumon2018champs,MR0344253}.

In this thesis, we are interested in stacks on the category $\mc{C}=\text{Man}$, equipped with the Grothendieck topology $\mc{J}$ defined by the open cover topology. One can define ``many'' stacks over the site $(\text{Man},\mc{J})$. We have mentioned that, for each object of the category $\text{Man}$, we have a stack over the site $(\text{Man},\mc{J})$. To work with stacks coming from manifolds would be uninteresting, and to work with arbitrary stacks over $(\text{Man},\mc{J})$ would be too general. So, in this thesis, we settle with stacks that come from Lie groupoids. The stacks that come from Lie groupoids are called differentiable stacks. 

In particular, the geometric structures over differentiable stacks that we study in this thesis are gerbes over differentiable stacks and principal bundles over differentiable stacks (\Cref{Chap.2,Chap.4}).

A gerbe over a topological space $X$ is a sheaf of groupoids on $X$ that is ``locally connected'' and ``locally non-empty'' (\cite[Definition $3.1$]{Moerdijk2}). The notion of a gerbe over a stack is a generalization of the notion of a gerbe over a topological space. A gerbe over a differentiable stack is a morphism of differentiable stacks satisfying some conditions similar to ``locally non-emptiness'' and ``locally connectedness'' (\cite[Definition 
$4.7$]{MR2817778}). Some papers that discuss (geometric structures on) gerbes over differentiable stacks are \cite{MR2817778,MR1969572,MR3480061}. In \cite{MR2817778} the authors discuss connections and curvature for groupoid central $S^1$-extensions. They also describe a prequantization result for principal $S^1$-bundles and $S^1$-gerbes. In \cite{MR3480061} the authors give a one-one correspondence between principal $[\text{Aut}(G)\rtimes G\rra G]$-bundles over differentiable stacks and $G$-gerbes over differentiable stacks.  There is a different notion of gerbe, introduced by M.K. Murray, which are called \textit{bundle gerbes} (\cite[Section $3$]{MR1405064}). In \cite{MR2681698} the author gives an exposition of the notion of bundle gerbe and relates them to the Hitchin-Chatterjee gerbes (\cite{MR1876068,chatterjee1998gerbs}). In \cite{MR2117631} the authors discuss a higher version of principal bundles, which they call the nonabelian bundle gerbes. They also study the connection, curvature, gauge transformations of such higher principal bundles.  

As mentioned earlier, a differentiable stack is the same as a Morita equivalence class of a Lie groupoid. The notion of principal bundle over Lie groupoid is ``Morita invariant'', thus giving the notion of principal bundles over differentiable stacks. In (the lecture notes) \cite{MR2206877}, the author gives a detailed exposition to differentiable stacks and their cohomology.  In \cite{MR2923408} the authors determine a criterion for a principal bundle (over a root stack) to have an algebraic connection. The notion of parallel transport on principal bundles over differentiable stacks is studied in \cite{MR3521476}. A slight generalization of the notion of principal $G$-bundle over a differentiable stack, namely principal $[H/G]$-bundle over a differentiable stack, for a Lie crossed module $(G,H,\tau,\alpha)$ is studied in \cite{Ginot}. The principal $S^1$-bundles over differentiable stacks are the main content of the study in \cite{MR1969572}.

In \Cref{Chap.5}, we study some topological structures over topological stacks. A topological stack is a stack over the category $\mc{C}=\text{Top}$, with the Grothendieck topology being the open cover topology. 
In \cite{noohi2005foundations,noohiquick} the author give an exposition to the notion of topological stacks. In his thesis \cite{Carchedi}, David Carchedi gives a detailed primer on topological and differentiable stacks. In  \cite{metzler2003topological}, the author offers a detailed comparison between topological and differentiable stacks. In \cite{MR2900442}, the author discusses topological stacks arising out of Cartesian closed categories. Algebraic topological properties of topological stacks, for example, singular chains, fibrations, Homotopy types, mapping stacks, and fundamental groups, are discussed in \cite{MR3552548,MR3144243,MR2927363,MR2719557, MR2443246,ebert2009homotopy}. In this thesis, we introduce the notion of topological groupoid extension and study some topological properties of such extensions. 
\section{Contents of the thesis}

\subsection{Preliminaries} In \Cref{Chap.Preliminaries}, we have recalled some standard material from different backgrounds of category theory and differential geometry. 

In \Cref{Section:Categorytheoryand2-categorytheory}, we have recalled some background from category theory and $2$-category theory. Some references on which this section is based are \cite{MR1291599, MR1438546, MR2858226}

In \Cref{Section:Fiberedcategoriespseudofunctorsandstacks}, we have recalled the notions of fibered categories, pseudo-functors, and stacks. This section is mainly based on Angelo Vistoli's ``Grothendieck topologies, fibered categories, and descent
theory'' \cite{MR2223406}.

In \Cref{Section:LiegroupoidsMoritaequivalencedifferentiablestacks}, we have recalled the necessary background of Lie groupoids, Morita equivalence, and differentiable stacks. We have used \cite{MR3089760, MR2778793, MR2012261} as references for this section.  

In \Cref{Section:Atiyahforsmoothmanifold}, we have recalled the notion of Atiyah sequence and Chern-Weil theory for principal bundles over the smooth manifold. This section is mainly based on \cite{MR896907, MR1393941}.

\subsection{On two notions of a gerbe over a stack}
In \Cref{Chap.2}, based on our paper \cite{MR4124773}, we have given a correspondence between two notions of a gerbe over a stack. 

The notion of gerbes was introduced by Giraud in \cite{MR0344253} to study nonabelian cohomology. An algebraic geometry point of view of gerbes are studied by Breen in \cite{MR1301844, MR2664623}, and by Breen and Messing in \cite{MR2183393}. Nonabelian gerbes over manifolds are studied in \cite{MR2764890,MR2520993,MR3084724}. The abelian bundle gerbe is introduced by Murray (\cite{MR1405064,MR2681698}).  The notion of a nonabelian bundle gerbe is introduced by Aschieri, Jurco, Cantini in \cite{MR2117631}. The paper \cite{MR3089401} discusses four equivalent versions of nonabelian gerbes.


The outline of this chapter is as follows:

In \Cref{Section:Twonotionsofgerbe}, we have recalled the two existing notions of a gerbe over a stack (one is a morphism of stacks \Cref{Definition:GerbeoverstackasinMR2817778}, the other is a morphism of Lie groupoids \Cref{Definition:GerbeoverstackVersion1}). We give some examples of a gerbe over a differentiable stack. 

In \Cref{Section:gerbegivingLiegroupoid}, we start with the definition of a gerbe over a stack. Assuming an extra condition on a gerbe over a stack, we associate a Lie groupoid extension. This is done in multiple steps. Given a gerbe over a stack $F:\mc{D}\ra \mc{C}$, we first, prove that, there exists a special kind of atlas $q:\underline{X}\ra\mc{C}$ for $\mc{C}$, and a morphism of stacks $p:\underline{X}\ra \mc{D}$ with some compatibility conditions. Using the property that the diagonal morphism associated to the morphism $F:\mc{D}\ra \mc{C}$ is an epimorphism, we prove that $p:\underline{X}\ra \mc{D}$ is an epimorphism of stacks. We then use our assumption that the diagonal morphism is a representable surjective submersion to prove that $p:\underline{X}\ra \mc{D}$ is an atlas. These atlases give Lie groupoids for $\mc{D}, \mc{C}$. We prove that this morphism of stacks $\mc{D}\ra \mc{C}$ thus associates a morphism of Lie groupoids (that is, a Lie groupoid extension).

In \Cref{Section:GoSassociatedtoALiegrupoidExtension} we have proved that the morphism of stacks associated to a Lie groupoid extension is a gerbe over a stack. 
Firstly, we have associated a $\mc{G}-\mc{H}$-bibundle for a morphism of Lie groupoids $\mc{G}\ra \mc{H}$ (\Cref{Subsection:bibundleAssociatedtoMorphismofLiegroupoids}). We have then assigned a morphism of stacks $B\mc{G}\ra B\mc{H}$ for a $\mc{G}-\mc{H}$-bibundle (\Cref{SubSubsection:MorphismofStacksassociatedtobibundle}). Combining these two constructions, we have a description of assigning a morphism of stacks $B\mc{G}\ra B\mc{H}$ for a morphism of Lie groupoids $\mc{G}\ra \mc{H}$. Finally, in \Cref{Subsection:LieGpdExtngivesGoS}, we proved that the morphism of stacks $B\mc{G}\ra B\mc{H}$, associated to a Lie groupoid extension $\mc{G}\ra \mc{H}$, is a gerbe over the stack $B\mc{H}$.

\subsection{Chern-Weil theory for principal bundles over Lie groupoids} In \Cref{Chap.3}, we have assigned a Chern-Weil map for a principal bundle over Lie groupoid using the Atiyah sequence approach. 

Let $G$ be a Lie group, $M$ a smooth manifold, and $P(M, G)$ a principal bundle. The standard way to introduce the notion of connection on $P(M,G)$ is in terms of differential forms (and distributions). The notion of connection can also be seen as splitting of the Atiyah sequence of the $G$-bundle $P(M, G)$. The Atiyah sequence is a short exact sequence 
\[0\ra (P\times \mf{g})/G\ra (TP)/G\ra TM\ra 0,\]
of vector bundles over the manifold $M$. 

In this chapter, we introduce the notion of a connection on a principal bundle over a Lie groupoid using the notion of the Atiyah sequence. 

Let $G$ be a Lie group, $[X_1\rra X_0]$ a Lie groupoid. A principal $G$-bundle over $[X_1\rra X_0]$ is a principal $G$-bundle $E_G\ra X_0$ over the manifold $X_0$, with an action of $[X_1\rra X_0]$ on $E_G$, that is compatible with the action of $G$ on $E_G$. This notion of principal bundles $(E_G\ra X_0,[X_1\rra X_0])$ include the example of principal bundle over a manifold, equivariant principal bundles. 

In \cite{MR2270285} the authors introduced the Chern-Weil map from a different point of view. In our paper, we use the Atiyah sequence approach. It turns out that an extra structure of a connection on $[X_1\rra X_0]$ is necessary to get a short exact sequence of vector bundles over $[X_1\rra X_0]$, which we call the Atiyah sequence associated to the principal bundle $(E_G\ra X_0,[X_1\rra X_0])$. A vector bundle over a Lie groupoid $[X_1\rra X_0]$ is a vector bundle $E\ra X_0$ over the manifold $X_0$ with an action of $[X_1\rra X_0]$ on $E$ such that the induced maps on fibers of $E\ra X_0$ are linear maps.


The outline of this chapter is as follows:

In \Cref{Section:principalbundnleoverLiegroupoids}, we introduce the notion of principal bundles over Lie groupoids. Some standard references for this notion are \cite{MR3150770,MR2270285,MR1950948,MR2119241}. In \Cref{Definition:connectiononLiegroupoid}, we introduce the notion of a connection on a Lie groupoid $[X_1\rra X_0]$. Using this notion of connection, we turn the tangent bundle $TX_0\ra X_0$ into a vector bundle over the Lie groupoid $[X_1\rra X_0]$ (\Cref{Section:TXXisavectorbundle}). 

In \Cref{Section:AtiyahforbundleoverLiegroupoid}, we introduce the notion of Atiyah sequence associated to a principal bundle $(E_G\ra X_0,[X_1\rra X_0])$, as a short exact sequence of vector bundles over the Lie groupoid $[X_1\rra X_0]$. 

In \Cref{Section:differentialformsassociatedtoconnectiononLiegroupoid}, we introduce the notion of differential form on a Lie groupoid $[X_1\rra X_0]$, equipped with an integrable connection $[X_1\rra X_0]$. We also prove some results about differential forms on Lie groupoids (\Cref{Lemma:differentialisdifferentialform}, \Cref{Proposition:pullbackisdifferrentialform}, \Cref{Proposition:pullbackofdifferrentialform}).

In \Cref{Section:differentialformsassociatedtoconnection}, we assign differential $1$-form (on Lie groupoid $[s^*E_G\rra E_G]$) associated to a connection $\mc{D}$ on the principal bundle $(E_G\ra X_0,[X_1\rra X_0])$. In \Cref{Section:differential2formsassociatedtoconnection}, we associate a differential $2$-form (on the Lie groupoid $[X_1\rra X_0]$) associated to a connection $\mc{D}$ on $(E_G\ra X_0,[X_1\rra X_0])$. 

In \Cref{Section:ChernWeilforEGX0X1X0}, we first introduce the notion of de Rham cohomology of the pair $(\mb{X},\mc{H})$ (\Cref{Section:deRhamcohomologyofLiegroupoidwithConnection}). Imitating the construction of the Chern-Weil map for a principal bundle over a manifold, we introduce the notion of Chern-Weil map for a principal bundle over a Lie groupoid (\Cref{Section:ChernWeilMap}).

\subsection{Connections on principal bundles over differentiable stacks} In \Cref{Chap.4}, we have introduced the notion of connection on principal bundles over differentiable stacks. 

Let $G$ be a Lie group, and $[X_1\rra X_0]$ a Lie groupoid. There is a notion of pullback of a principal bundle along a morphism of Lie groupoids, and more generally, along a ``generalized morphism of Lie groupoids''. It is a standard result that, under this notion of pullback of a principal bundle along generalized morphism of Lie groupoids, for two Morita equivalent Lie groupoids $\mb{X},\mb{Y}$, there is an equivalence of categories between category of principal bundles over $\mb{X}$ and category of principal bundles over $\mb{Y}$ (\cite[Corollary $2.12$]{MR2270285}). As a differentiable stack is a Morita equivalence class of a Lie groupoid, the notion of principal bundle over Lie groupoid induces the notion of principal bundle over a differentiable stack. The notion of parallel transport on principal bundles over differentiable stacks is studied in \cite{MR3521476}. 

The outline of this chapter is as follows:

In \Cref{Section:principalbundleoverdifferentiablestacks}, we recall the notion of principal bundles over differentiable stacks. This notion of principal bundle over differentiable stacks can be found in \cite{MR3521476,MR2206877,MR2817778}.

In \Cref{Definition:principalbundleoverstacks}, we introduce the notion of connection on a principal bundle over differentiable stacks. In \Cref{Section:Atiyahsequenceforbundleoverstacks}, we introduce the Atiyah sequence for a principal bundle over a Deligne-Mumford stack. We also mention some results about connection on principal bundle over differentiable stacks.

\subsection{Extensions of topological groupoids} In \Cref{Chap.5}, we introduce the notion of topological groupoid extensions and study their correspondence with gerbes over topological stacks. 

In \Cref{Section:top groupoids and top stacks}, we recall the notion of topological groupoids and topological stacks. These notions can be found in \cite{Carchedi, metzler2003topological, noohi2005foundations}.

In \Cref{Section:Topogrpdextensionmorstacks}, we introduce the notion of topological groupoid extension. We prove that the morphism of stacks associated to a topological groupoid extension is a gerbe over a topological stack. In \Cref{Section:Serre-hurewiczfibrations}, we recall the notions of Serre/Hurewicz stacks and see some results about gerbes over Serre/Hurewicz stacks. 

\subsection{Future directions of research} There are several possible directions of enquiries may emerge out of this thesis. Here we mention some of them.  One can think of extending the correspondence discussed in \Cref{Chap.2} for a gerbe over stack with a connection. It would be interesting to relate our Chern-Weil theory developed in \Cref{Chap.3,Chap.4} with   simplicial homotopy theoritic ideas introduced in  \cite{MR3049871}. Another more ambitious project would be to apply our constructions to study the notion of  gerbes over an algebraic stack.
	\chapter{Preliminaries} \label{Chap.Preliminaries}
In this chapter, we recall some preliminaries needed for the thesis. This chapter is split into four sections:
\begin{enumerate}
	\item Category theory and $2$-category theory (\Cref{Section:Categorytheoryand2-categorytheory})
	\item Fibered categories, pseudo-functors, and stacks (\Cref{Section:Fiberedcategoriespseudofunctorsandstacks}).
	\item Lie groupoids, Morita equivalence, and differentiable stacks (\Cref{Section:LiegroupoidsMoritaequivalencedifferentiablestacks}).
	\item Atiyah sequence and Chern-Weil theory for principal bundles over smooth manifold (\Cref{Section:Atiyahforsmoothmanifold}).
\end{enumerate}

\section{Category theory and $2$-category theory}\label{Section:Categorytheoryand2-categorytheory}
In this section, we recall some of the standard terminology, examples, results from the category theory and the $2$-category theory. The content of this section is mostly standard, and can be found in books, such as  \cite{MR1291599, MR1438546, MR2858226}.

\begin{definition}[Category]\label{Definition:category}
	A \textit{category, denoted as $\mc{C}$}, consists of the following data:
	\begin{itemize}
		\item a class $\mc{C}_0$, whose elements are called \textit{objects of $\mc{C}$},
		\item a class $\mc{C}_1$, whose elements are called \textit{morphisms of $\mc{C}$},
		\item a map $s_{\mc{C}}:\mc{C}_1\ra \mc{C}_0$, called \textit{the source map of $\mc{C}$},
		\item a map $t_{\mc{C}}:\mc{C}_1\ra \mc{C}_0$, called \textit{the target map of $\mc{C}$},
		\item a map $m_{\mc{C}}:\mc{C}_1\times_{t_{\mc{C}},\mc{C}_0,s_{\mc{C}}}\mc{C}_1\ra \mc{C}_1$, called \textit{the composition map of $\mc{C}$},
		\item a map $e_{\mc{C}}:\mc{C}_0\ra \mc{C}_1$, called \textit{the identity map of $\mc{C}$},
	\end{itemize}	
	satisfying the following conditions:
	\begin{itemize}
		\item $s_{\mc{C}}\circ m_{\mc{C}}=
		s_{\mc{C}}\circ \pr_1$, and $t_{\mc{C}}\circ m_{\mc{C}}=t_{\mc{C}}\circ \pr_2$,
		\item $s_{\mc{C}}\circ e_{\mc{C}}=1_{\mc{C}_0}$, and $t_{\mc{C}}\circ e_{\mc{C}}=1_{\mc{C}_0}$,
		\item $m_{\mc{C}}\circ (m_{\mc{C}},1_{\mc{C}_1})=m_{\mc{C}}\circ (1_{\mc{C}_1},m_{\mc{C}})$,
		\item $m_{\mc{C}}(f,e_{\mc{C}}(a))=f$ for all $a\in C_0$, and $f\in \mc{C}_1$ with $t_{\mc{C}}(f)=a$,
		\item $m_{\mc{C}}(e_{\mc{C}}(a),f)=f$ for all $a\in \mc{C}_0$, and $f\in \mc{C}_1$ with $s_{\mc{C}}(f)=a$.
	\end{itemize}	
	A category, with all the associated maps, can be expressed by the following diagram,
	\[\mc{C}_1\times_{\mc{C}_0}\mc{C}_1\xrightarrow{m_{\mc{C}}}
	\mc{C}_1\xrightrightarrows[t_{\mc{C}}]{s_{\mc{C}}} \mc{C}_0\xrightarrow{e_{\mc{C}}}\mc{C}_1.\]
	However, often, we use a less elaborate notation and denote a category as $[\mc{C}_1\rightrightarrows \mc{C}_0]$. We call the maps $s_{\mc{C}},t_{\mc{C}},m_{\mc{C}}, e_{\mc{C}}$ to be the structure maps of the category $\mc{C}$. When there is no confusion, we ignore the subscript $\mc{C}$ in structure maps and write them as $s,t,m,e$. The classes $\mc{C}_0$, and $\mc{C}_1$ will be denoted alternatively as $\text{Obj}(\mc{C})$, and $\text{Mor}(\mc{C})$ respectively.
\end{definition}
\begin{remark}[Notation]\label{Remark:notationofcategory}
	For $(\gamma_1,\gamma_2)\in \mc{C}_1\times_{t_{\mc{C}},\mc{C}_0,s_{\mc{C}}}\mc{C}_1$, we denote the element $m_{\mc{C}}(\gamma_1,\gamma_2)$ as $\gamma_2\circ \gamma_1$.
	For $a\in \mc{C}_0$, we denote the element $e_{\mc{C}}(a)$ by $1_a$.  With this notations, the conditions in the definition of a category can be seen as the following:
	\begin{itemize}
		\item for all $(\gamma_1,\gamma_2)\in \mc{C}_1\times_{\mathcal{C}_0}\mc{C}_1$, we have $s_{\mc{C}}(\gamma_2\circ \gamma_1)=s_{\mc{C}}(\gamma_1)$, and $t_{\mc{C}}(\gamma_2\circ \gamma_1)=t_{\mc{C}}(\gamma_2)$,
		\item for all $a\in \mc{C}_0$, we have $s_{\mc{C}}(1_a)=a$ and $t_{\mc{C}}(1_a)=a$,
		\item for all $(\gamma_1,\gamma_2,\gamma_3)\in \mc{C}_1\times_{\mathcal{C}_0}\mc{C}_1\times_{\mathcal{C}_0}\mc{C}_1$, we have $\gamma_3\circ (\gamma_2\circ \gamma_1)=(\gamma_3\circ\gamma_2)\circ \gamma_1$,
		\item for all $a\in \mc{C}_0, \gamma\in \mc{C}_1$, with $s_{\mc{C}}(\gamma)=a$, we have  $\gamma\circ 1_a=\gamma$,
		\item for all $a\in \mc{C}_0, \gamma\in \mc{C}_1$, with $t_{\mc{C}}(\gamma)=a$, we have $1_a\circ \gamma=\gamma$.
	\end{itemize}
	For $\gamma\in \mc{C}_1$ with $s_{\mc{C}}(\gamma)=a$ and $t_{\mc{C}}(\gamma)=b$, we say that, $\gamma$ is a \textit{morphism from $a$ to $b$}. For $a,b\in \mc{C}_0$, we denote the collection $\{\gamma\in \mc{C}_1: s_{\mathcal{C}}(\gamma)=a, t_{\mathcal{C}}(\gamma)=b\}$, as $\mc{C}(a,b)$. An element $\gamma$ of $\mc{C}(a,b)$ is denoted as $\gamma:a\ra b$.
	
\end{remark}

\begin{remark}
	Depending on the ``size of the class $\mc{C}_1$'', we call a category $\mc{C}$ to be a small, or, a locally-small, or, a large category. A category $\mc{C}$ in which $\mc{C}_1$ is a set, is called a \textit{small category}. A category $\mc{C}$, in which $\mc{C}(a,b)$ is a set for each $a,b\in \mc{C}_0$, is called a \textit{locally-small category}. A category $\mc{C}$, in which $\mc{C}(a,b)$ is not a set for some $a,b\in \mc{C}_0$, is called a \textit{large category}. In this thesis we will not discuss such foundational issues in much detail.
\end{remark}
In this thesis, our primary interest would be a special kind of small categories, namely, Lie groupoids. 

\begin{example}\label{Example:MMcategory}
	Let $M$ be a set. Consider the category $[M\rra M]$, whose object class is the set $M$, whose morphism class is the set $M$. All the structure maps are identity maps. We call $[M\rra M]$ to be the \textit{discrete category associated to the set $M$}. In particular, when $M$ has the structure of a smooth manifold, we get a category $[M\rra M]$, which we call \textit{the category associated to the smooth manifold $M$}.
\end{example}
\begin{example}\label{Example:G*category}
	Let $G$ be a group. Consider the category $[G\rra *]$, whose object class is the singleton set $\{*\}$, whose morphism class is the set $G$. The source, target maps are given by the constant map that takes every element of $G$ to the element $*$. The composition map $m: G\times G\ra G$ is the group operation. The identity map $*\ra G$ is the map that takes $*$ to the identity element $e_G$. We call $[G\rra *]$ to be the \textit{category associated to the group $G$}.  In particular, when $G$ has the structure of a Lie group, we get a category $[G\rra *]$, which we call \textit{the category associated to the Lie group $G$}.
\end{example}
\begin{example}\label{Example:actioncategory}
	Let $G$ be a group, $M$ a set, and $\mu:M\times G\ra M$ an action of the group $G$ on the set $M$. Consider the category $[M\times G\rra M]$ with the following description:
	\begin{itemize}
		\item the source map is given by the first projection,
		\item the target map is given by the action map,
		\item the composition map is given by $((m_1,g_1),(m_2,g_2))\mapsto (m_1,g_1g_2)$,
		\item the identity map is given by $m\mapsto (m,e_G)$.
	\end{itemize}
	We call $[M\times G\rra M]$ to be the \textit{category associated to the action of $G$ on $M$}. In particular, when $G$ is a Lie group, $M$ a smooth manifold, and $\mu$ an action of Lie group on manifold, we call $[M\times G\rra M]$ to be the \textit{category associated to the Lie group action of $G$ on $M$}.
\end{example}
\begin{definition}[Groupoid]\label{Definition:Groupoid}
	A \textit{groupoid} is a category $\mc{C}$ (with structure maps $(s_{\mc{C}},t_{\mc{C}},m_{\mc{C}},e_{\mc{C}})$) along with a map $i_{\mc{C}}:\mc{C}_1\ra \mc{C}_1$ such that, the following conditions are satisfied:
	\begin{itemize}
		\item $s_{\mc{C}}\circ i_{\mc{C}}=t_{\mc{C}}$, and $t_{\mc{C}}\circ i_{\mc{C}}=s_{\mc{C}}$,
		\item $m_{\mc{C}}\circ (1_{\mc{C}_1},i_{\mc{C}})=e_{\mc{C}}\circ s_{\mc{C}}$, and $m_{\mc{C}}\circ (i_{\mc{C}},1_{\mc{C}_1})=e_{\mc{C}}\circ t_{\mc{C}}$.
	\end{itemize}
\end{definition}
\begin{remark}[Notation]
	In a groupoid $\mc{C}$, for an element $\gamma\in \mc{C}_1$, we call $i(\gamma)$ as, \textit{the inverse of $\gamma$} and denote by $\gamma^{-1}$. The first condition in the above definition says that $s(\gamma)=t(\gamma^{-1})$, and $t(\gamma)=s(\gamma^{-1})$. The second condition in the above definition says that, for $\gamma\in \mc{C}_1$, we have, 
	$ \gamma\circ \gamma^{-1}=1_{t(\gamma)}$, and $\gamma^{-1}\circ \gamma=1_{s(\gamma)}$.	
	
\end{remark}
\begin{example}
	All the three examples $[M\rra M], [G\rra *]$, and $[M\times G\rra M]$, mentioned as 
	\Cref{Example:MMcategory,Example:G*category,Example:actioncategory} are groupoids. The maps $i_{\mc{C}}:\mc{C}_1\ra \mc{C}_1$ for these categories have the following description:
	\begin{itemize}
		\item the map $M\ra M$, defined as $m\mapsto m$ for each $m\in M$,
		\item the map $G\ra G$, defined as $g\mapsto g^{-1}$ for each $g\in G$,
		\item the map $M\times G\ra M\times G$, defined as $(m,g)\mapsto (mg,g^{-1})$ for each $(m,g)\in M\times G$.
	\end{itemize}
\end{example}

\begin{definition}[Functor]\label{Definition:functor}
	Let $\mc{C}=[\mc{C}_1\rra \mc{C}_0]$, and $\mc{D}=[\mc{D}_1\rra \mc{D}_0]$ be categories with structure maps given by $(s_{\mc{C}},t_{\mc{C}},m_{\mc{C}},e_{\mc{C}})$, and 
	$(s_{\mc{D}},t_{\mc{D}},m_{\mc{D}},e_{\mc{D}})$ respectively. A \textit{functor from $\mc{C}$ to $\mc{D}$}, denoted as $F:\mc{C}\ra \mc{D}$, is given by a pair of maps 
	\[( F_1:\mc{C}_1\ra \mc{D}_1,F_0:\mc{C}_0\ra \mc{D}_0),\] such that, the following conditions are satisfied: 
	\begin{itemize}
		\item $F_0\circ s_{\mc{C}}=s_{\mc{D}}\circ F_1 $, and $F_0\circ t_{\mc{C}}=t_{\mc{D}}\circ F_1$,
		\item $F_1\circ m_{\mc{C}}=m_{\mc{D}}\circ (F_1,F_1)$,
		\item $F_1\circ e_{\mc{C}}=e_{\mc{D}}\circ F_0$.
	\end{itemize}
	One can express a functor $F:\mc{C}\ra \mc{D}$ as the following diagram,
	\[
	\begin{tikzcd}
		\mc{C}_1\times_{\mc{C}_0}\mc{C}_1 \arrow[d, "m_{\mc{C}}"'] \arrow[rr, "{(F_1,F_1)}"] & & \mc{D}_1\times_{\mc{D}_0}\mc{D}_1 \arrow[d, "m_{\mc{D}}"] \\
		\mc{C}_1 \arrow[d,xshift=0.75ex,"t_{\mc{C}}"]\arrow[d,xshift=-0.75ex,"s_{\mc{C}}"'] \arrow[rr, "F_1"]  & & \mc{D}_1 \arrow[d,xshift=0.75ex,"t_{\mc{D}}"]\arrow[d,xshift=-0.75ex,"s_{\mc{D}}"']     \\
		\mc{C}_0 \arrow[d, "e_{\mc{C}}"'] \arrow[rr, "F_0"]     & & \mc{D}_0 \arrow[d, "e_{\mc{D}}"]    \\
		\mc{C}_1 \arrow[rr, "F_1"]        & & \mc{D}_1       
	\end{tikzcd}.\]
	We simplify it by expressing a functor as the middle square in the above diagram. We denote a functor as $(F_1,F_0):[\mc{C}_1\rra \mc{C}_0]\ra[\mc{D}_1\rra \mc{D}_0]$.
\end{definition}

\begin{remark}
	In the definition of a functor, the condition $F_0\circ s_{\mc{C}}=s_{\mc{D}}\circ F_1 $, and $F_0\circ t_{\mc{C}}=t_{\mc{D}}\circ F_1$, says that, for each $\gamma\in \mc{C}_1$, the element $F_1(\gamma)$ is a morphism from $F_0(s_{\mc{C}}(\gamma))$ to $F_0(t_{\mc{C}}(\gamma))$. In other words, for $\gamma:a\ra b$ in $\mc{C}_1$, we have $F_1(\gamma):F_0(a)\ra F_0(b)$. The condition $F_1\circ m_{\mc{C}}=m_{\mc{D}}\circ(F_1,F_1)$ says that, for each $(\gamma,\gamma')\in \mc{C}_1\times_{\mathcal{C}_0}\mc{C}_1$, we have $F_1(\gamma'\circ \gamma)=F_1(\gamma')\circ F_1(\gamma)$.
	The condition $F_1\circ e_{\mc{C}}=e_{\mc{D}}\circ F_0$ says that $F_1(1_a)=1_{F_0(a)}$ for all $a\in \mc{C}_0$.  Thus, for each $a,b\in \mc{C}_0$, we have a map $\mc{C}(a,b)\ra \mc{D}(F_0(a),F_0(b))$.
\end{remark}

Based on the properties satisfied by the map $\mc{C}(a,b)\ra \mc{D}(F_0(a),F_0(b))$, for each $a,b\in \mc{C}_0$, we name the functor $(F_1,F_0):[\mc{C}_1\rra \mc{C}_0]\ra[\mc{D}_1\rra \mc{D}_0]$.
\begin{definition}
	A functor $(F_1,F_0):[\mc{C}_1\rra \mc{C}_0]\ra[\mc{D}_1\rra \mc{D}_0]$ is said to be a \textit{faithful (or full) functor} if, for each $a,b$ in $\mc{C}_0$, the map
	$\mc{C}(a,b)\ra \mc{D}(F_0(a),F_0(b))$ is  an injective (or a surjective) map.
	
\end{definition}

\begin{example}\label{Example:functorMMNN}
	Let $M,N$ be smooth manifolds, and $f:M\ra N$ a smooth map. Then, $(f,f):[M\rra M]\ra [N\rra N]$ is a functor.
\end{example}
\begin{example}\label{Example:functorG*H*}
	Let $G,H$ be Lie groups, and $f:G\ra H$ a Lie group homomorphism. Then, 
	$(f,*):[G\rra *]\ra [H\rra *]$ is a functor.
\end{example}
\begin{example}\label{Example:functorMGMNGN}
	Let $G$ be a Lie group, acting on smooth manifolds $M,N$. Let $f:M\ra N$ be a $G$-equivariant map, that is $f(mg)=f(m)g$ for all $(m,g)\in M\times G$. Then, \[(f\times 1, f):[M\times G\rra M]\ra [N\times G\rra N]\] is a functor.
\end{example}

\begin{definition}[composition of functors]\label{Definition:compositionoffunctors}
	Let $\mc{C},\mc{D}$, and $\mc{E}$ be categories. Consider the functors $F:\mc{C}\ra \mc{D}$ and $G:\mc{D} \ra \mc{E}$ that are given by the maps $(F_1:\mc{C}_1\ra \mc{D}_1, F_0:\mc{C}_0\ra \mc{D}_0)$, and $(G_1:\mc{D}_1\ra \mc{E}_1, G_0:\mc{D}_0\ra \mc{E}_0)$ respectively. Then, the maps \[(G_1\circ F_1:\mc{C}_1\ra \mc{E}_1, G_0\circ F_0:\mc{C}_0\ra \mc{E}_0),\] gives a functor from $\mc{C}$ to $\mc{E}$, which we call \textit{the composition of $F, G$}, and denote by $G\circ F$.
\end{definition}

\begin{definition}[Natural transformation]\label{Definition:Naturaltransformation}
	Let $\mc{C}=[\mc{C}_1\rra \mc{C}_0]$ and $\mc{D}=[\mc{D}_1\rra \mc{D}_0]$ be categories with structure maps given by $(s_{\mc{C}},t_{\mc{C}},m_{\mc{C}},e_{\mc{C}})$ and 
	$(s_{\mc{D}},t_{\mc{D}},m_{\mc{D}},e_{\mc{D}})$ respectively. Let $(F_1,F_0)$, and $(G_1,G_0)$ be functors from $\mc{C}$ to $\mc{D}$. A \textit{natural transformation $\eta$ from $F$ to $G$}, denoted by $\eta:F\Rightarrow G$, is a map $\eta:\mc{C}_0\ra \mc{D}_1$, such that, the following conditions are satisfied:
	\begin{itemize}
		\item $s_{\mc{D}}\circ \eta=F_0$, and $t_{\mc{D}}\circ \eta=G_0$,
		\item $s_{\mc{D}}\circ \eta\circ t_{\mc{C}}=t_{\mc{D}}\circ F_1$, and $t_{\mc{D}}\circ \eta\circ s_{\mc{C}}=s_{\mc{D}}\circ G_1$,
		\item $m_{\mc{D}}\circ (F_1,\eta\circ t_{\mc{C}})=m_{\mc{D}}\circ (\eta\circ s_{\mc{C}},G_1)$ as maps $\mc{C}_1\ra \mc{D}_1\times_{\mathcal{D}_0}\mc{D}_1\ra \mc{D}_1$.
	\end{itemize}We see a natural transformation as the following diagram,
	\[
	\begin{tikzcd}[sep=large]
		\mc{C}_1 \arrow[dd,xshift=0.75ex,"t_{\mc{C}}"]\arrow[dd,xshift=-0.75ex,"s_{\mc{C}}"'] \arrow[rr, "F_1", bend left] \arrow[rr, "G_1", bend right]  & & \mc{D}_1 \arrow[dd,xshift=0.75ex,"t_{\mc{D}}"]\arrow[dd,xshift=-0.75ex,"s_{\mc{D}}"'] \\
		& &   \\
		\mc{C}_0 \arrow[rr, "F_0", bend left] \arrow[rruu, "\eta"] \arrow[rr, "G_0", bend right] & & \mc{D}_0  
	\end{tikzcd}.\]
\end{definition}
\begin{remark}
	In the \Cref{Definition:Naturaltransformation}, the condition $s_{\mc{D}}\circ \eta=F_0$, and $t_{\mc{D}}\circ \eta=G_0$ says that, for each $a\in \mc{C}_0$, we have $\eta(a):F_0(a)\ra G_0(a)$. The second condition says that, for each $\gamma\in \mc{C}_1$, the elements $\eta(t_{\mc{C}}(\gamma)), F_1(\gamma)$ are composable, giving an element 
	$\eta(t_{\mc{C}}(\gamma))\circ F_1(\gamma)\in \mc{D}_1$. Similarly, the elements 
	$G_1(\gamma),\eta(s_{\mc{C}}(\gamma))$ being composable, gives 
	$G_1(\gamma)\circ \eta(s_{\mc{C}}(\gamma))\in \mc{D}_1$. The third condition says that, for each $\gamma\in \mc{C}_1$, we have 
	$\eta(t_{\mc{C}}(\gamma))\circ F_1(\gamma)=G_1(\gamma)\circ \eta(s_{\mc{C}}(\gamma))$. If $\gamma\in\mc{C}(a, a')$, then, the third condition says that the following diagram is commutative,
	\[
	\begin{tikzcd}
		F_0(a) \arrow[dd, "F_1(\gamma)"'] \arrow[rr, "\eta(a)"] & & G_0(a) \arrow[dd, "G_1(\gamma)"] \\
		& &     \\
		F_0(a') \arrow[rr, "\eta(a')"]    & & G_0(a')    
	\end{tikzcd}.\]
\end{remark}
\begin{definition}[Natural isomorphism]
	Let $\mc{C}=[\mc{C}_1\rra \mc{C}_0],\mc{D}=[\mc{D}_1\rra \mc{D}_0]$ be categories, and $(F_1,F_0),(G_1,G_0):[\mc{C}_1\rra \mc{C}_0]\ra[\mc{D}_1\rra \mc{D}_0]$ functors from $\mc{C}$ to $\mc{D}$. A natural transformation $\eta:\mc{C}_0\ra \mc{D}_1$, from $F$ to $G$, in which $\eta(a)$ is an isomorphism for each $a\in \mc{C}_0$, is called a \textit{natural isomorphism from $(F_1,F_0)$ to $(G_1,G_0)$}. 
\end{definition}
\begin{definition}[Equivalence of categories]
	Let $\mc{C}=[\mc{C}_1\rra \mc{C}_0],\mc{D}=[\mc{D}_1\rra \mc{D}_0]$ be categories. A functor $F:\mc{C}\ra \mc{D}$ is said to be \textit{an equivalence of categories}, if, there exists a functor $G:\mc{D}\ra \mc{C}$, such that, there are natural isomorphism of functors $F\circ G\Rightarrow 1_{\mc{D}}$, and $G\circ F\Rightarrow1_{\mc{C}}$. An equivalence of categories is said to be an isomorphism of categories if the natural isomorphisms are given by identity maps. Two categories $\mc{C}$ and $\mc{D}$ are said to be \textit{equivalent (or isomorphic) categories} if there is an equivalence (or an isomrophism) of categories $F:\mc{C}\ra \mc{D}$. 
\end{definition}


\begin{theorem}[{\cite[Theorem $4.4.1$]{MR1712872}}]
	The following properties of a functor $F:\mc{D}\ra \mc{C}$ are logically equivalent:
	\begin{enumerate}
		\item $F:\mc{D}\ra \mc{C}$ is an equivalence of categories.
		\item $F$ is full and faithful, and each object $c\in \mc{C}_0$ is isomorphic to $F(d)$ for some object $d\in \mc{D}_0$.
	\end{enumerate}
\end{theorem}

To motivate the definition of a $2$-category, 
using \Cref{Remark:notationofcategory}, we rewrite the definition of a category (\Cref{Definition:category}).
\begin{definition}\label{Definition:categoryEquivalentdefinition}
	A \textit{category}, denoted by $\mathcal{C}$, consists of the following data, 
	\begin{enumerate}
		\item a class $\mathcal{C}_0$, whose elements are called \textit{objects of $\mc{C}$},
		\item a class $\mc{C}(x,y)$ for each $x,y\in \mc{C}_0$,
		\item a map $\{1\}\ra \mc{C}(x,x)$ for each $x\in \mc{C}_0$,
		\item a map $\circ :\mc{C}(x,y)\times\mc{C}(y,z)\ra \mc{C}(x,z)$, with $(f,g)\mapsto g\circ f$, for each $x,y,z\in \mc{C}_0$,
	\end{enumerate}
	satisfying the following conditions:
	\begin{enumerate}
		\item for each $x,y\in \mc{C}_0$, the following diagram is commutative,
		\begin{equation} \label{Diagram:definitionofcategoryidentity}
				\begin{tikzcd}[sep=large]
				{\mc{C}(x,y)} \arrow[dd, "{(1_x,1)}"'] \arrow[rrr,"{(1,1_y)}"] \arrow[rrrdd, "1"] & & & {\mc{C}(x,y)\times \mc{C}(y,y)} \arrow[dd, "\circ"] \\
				& & &       \\
				{\mc{C}(x,x)\times\mc{C}(x,y)} \arrow[rrr,"\circ"']  & & & {\mc{C}(x,y)}     
			\end{tikzcd},
		\end{equation}
		\item for each $x,y,z,w\in \mc{C}_0$, the following diagram is commutative,
	\begin{equation}\label{Diagram:definitionofcategorycomposition}
			\begin{tikzcd}[sep=large]
				{\mc{C}(x,y)\times\mc{C}(y,z)\times\mc{C}(z,w)} \arrow[rr, "{(\circ,1)}"] \arrow[dd, "{(1,\circ)}"'] & & \mc{C}(x,z)\times \mc{C}(z,w) \arrow[dd, "\circ"] \\
				& &        \\
				{\mc{C}(x,y)\times\mc{C}(y,w)} \arrow[rr, "\circ"]       & & {\mc{C}(x,w)}      
			\end{tikzcd}.
		\end{equation}
	\end{enumerate}
\end{definition}

The Diagram \ref{Diagram:definitionofcategoryidentity} is a fancy way of saying that $1_y\circ f=f$ and $f\circ 1_x=f$ for all $f\in \mc{C}(x,y)$.
The Diagram \ref{Diagram:definitionofcategorycomposition} says that $h\circ (g\circ f)=(h\circ g)\circ f$ for all $(f,g,h)\in 
\mc{C}(x,y)\times\mc{C}(y,z)\times\mc{C}(z,w)$.

\begin{example}\label{Example:OXcategory}
	Let $X$ be a topological space. Consider a category $\mc{O}(X)$, with the following description:
	\begin{itemize}
		\item objects are open subsets of $X$,
		\item for open sets $U, V$ of $X$, we declare 
		\begin{equation*}
			\mc{O}(X)(U,V)= 
			\begin{cases}
				\emptyset & \text{ if } U\nsubseteq V \\
				\{i_U:U\hookrightarrow V\} & \text{ if } U\subseteq V \\
			\end{cases}
		\end{equation*}
	\end{itemize}
	We call $\mc{O}(X)$ to be the \textit{category of open subsets of $X$}.
\end{example}
\begin{example}\label{Example:SetCategory}
	Let $\rm{Set}$ be the category described as follows:
	\begin{itemize}
		\item objects are sets,
		\item for sets $A,B$, the morphism class $\text{Set}(A,B)$ is the collection of functions from $A$ to $B$,
		\item the composition map is given by composition of functions.
	\end{itemize}
	We call $\rm{Set}$ to be the category of sets.
\end{example}
\begin{example}\label{Example:TopCategory}
	Let $\text{Top}$ be the category with the following description:
	\begin{itemize}
		\item objects are topological spaces,
		\item for topological spaces $X,Y$, the morphism class $\text{Top}(X,Y)$ is the set of all continuous maps from $X$ to $Y$,
		\item the composition map is given by composition of functions. 
	\end{itemize}
	We call $\text{Top}$ to be the category of topological spaces.
\end{example}
\begin{example}\label{Example:ManCategory}
	Let $\text{Man}$ be the category with the following description:
	\begin{itemize}
		\item objects are smooth manifolds,
		\item for smooth manifolds $M,N$, the morphism class $\text{Man}(M,N)$ is the set of all smooth maps from $M$ to $N$,
		\item the composition map is given by composition of functions. 
	\end{itemize}
	We call $\text{Man}$ to be the category of smooth manifolds.
\end{example}

We have given some examples of functors 
(\Cref{Example:functorMMNN,Example:functorG*H*,Example:functorMGMNGN}). Those examples are between small categories. Now we give some more examples of functors (between arbitrary categories).

\begin{example}\label{Example:forgetfulfunctor}
	Consider the categories $\rm{Top}$ and $\rm{Set}$. Consider the assignment $\text{Top}_0\ra \text{Set}_0$, that maps a topological space $(A,\tau)$ to the underlying set $A$, and the assignment $\text{Top}(A,B)\ra \text{Set}(A,B)$ that maps a continuous map to the underlying set map. This gives a functor $\text{Top}\ra \text{Set}$, which we call the forgetful functor from $\text{Top}$ to $\text{Set}$.
\end{example}

\begin{example}\label{Example:covarianthomfunctor}
	Let $\mc{C}$ be a (locally small) category. Fix an object $A$ in $\mc{C}$. The assignment $B\mapsto \mc{C}(A,B)$ for objects $B$ of $\mc{C}$, and the assignment $(B\ra C)\mapsto (\mc{C}(A,B)\ra \mc{C}(A,C))$ for morphisms $B\ra C$ in $\mc{C}$, gives a functor $h^A:\mc{C}\ra \text{Set}$. We call $h^A$ to be \textit{the (covariant) functor of points associated to the object $A$}. 
\end{example} 
\begin{definition}[Opposite category of a category]
	Let $\mc{C}$ be a category. The \textit{opposite category of $\mc{C}$}, denoted by $\mc{C}^{\op}$, is a category with the following description:
	\begin{itemize}
		\item the class $\mc{C}^{\op}_0$ is equal to the class $\mc{C}_0$,
		\item for each $x,y\in \mc{C}^{\op}_0$, we have $\mc{C}^{\op}(x,y)=\mc{C}(y,x)$,
		\item for each triple $(x,y,z)$ in $\mc{C}^{\op}_0$, the map $\circ:\mc{C}^{\op} (x,y)\times \mc{C}^{\op}(y,z)\ra \mc{C}^{\op}(x,z)$ is the map 
		$\circ:\mc{C}(z,y)\times \mc{C}(y,x)\ra \mc{C}(z,x)$,
		\item the map $\{1\}\ra \mc{C}^{\op}(x,x)$ is the map $\{1\}\ra \mc{C}(x,x)$.
	\end{itemize}
\end{definition}


\begin{remark}
	A functor on the category $\mc{C}^{\op}$ is called a \textit{contravariant functor on the category $\mc{C}$}. Let $F:\mc{C}^{\op}\ra \mc{D}$ be a functor.
	Consider a morphism $\gamma:a\ra b$ in $\mc{C}$. This can be seen as $\gamma:b\ra a$ in $\mc{C}^{\op}$. 
	The condition $s_{\mc{D}}\circ F_1=F_0\circ s_{\mc{C}^{\op}}$ says that $s_{\mc{D}}(F_1(\gamma))
	=F_0(s_{\mc{C}^{\op}}(\gamma))=F_0(t_{\mc{C}}(\gamma))$. Similarly, the condition $t_{\mc{D}}\circ F_1=F_0\circ t_{\mc{C}^{\op}}$ says that $t_{\mc{D}}(F_1(\gamma))
	=F_0(t_{\mc{C}^{\op}}(\gamma))=F_0(s_{\mc{C}}(\gamma))$. Thus, we have $F_1(\gamma):F_0(t_{\mc{C}}(\gamma))\ra F_0(s_{\mc{C}}(\gamma))$. So, for $\gamma:a\ra b$ in $\mc{C}$, we have $F_1(\gamma):F_0(b)\ra F_0(a)$. Similarly, the compatibility with composition map says that, $F_1(\gamma\circ \gamma')=F_1(\gamma')\circ F_1(\gamma)$ for all $(\gamma',\gamma)\in \mc{C}_1\times_{\mc{C}_0}\mc{C}_1$.
\end{remark}

In the \Cref{Example:covarianthomfunctor}, we have mentioned the (covariant) functor of points. Now, we give another example of a functor, the contravariant functor of points. 

\begin{example}\label{Example:contravarianthomfunctor}
	Let $\mc{C}$ be a (locally small) category. Fix an object $A$ in $\mc{C}$. The assignment $B\mapsto \mc{C}(B,A)$ for objects $B$ of $\mc{C}$, and the assignment $(B\ra C)\mapsto (\mc{C}(C,A)\ra \mc{C}(B,A))$ for morphisms $B\ra C$ in $\mc{C}$, gives a functor $h_A:\mc{C}^{\op}\ra \text{Set}$. We call $h_A$ to be \textit{the (contravariant) functor of points associated to the object $A$}. 
\end{example}
Let $\mc{C}$ be a (locally-small) category and  $\mc{F}:\mc{C}^{\op}\ra \text{Set}$ a functor. 
Let $\text{Hom}(h_U,\mc{F})$ be the set of all natural transformations from the functor $h_U$ to the functor $\mc{F}$, and $\eta:h_U\Rightarrow \mc{F}$ a natural transformation. Consider the morphism $\eta(U):h_U(U)\ra \mc{F}(U)$. The identity element $1_U\in h_U(U)=\mc{C}(U,U)$ gives an element $\eta(U)(1_U)\in \mc{F}(U)$. Thus, we have a morphism of sets $\text{Hom}(h_U,\mc{F})\ra \mc{F}(U)$. 
The following result says that, this morphism of sets $\text{Hom}(h_U,\mc{F})\ra \mc{F}(U)$ is a bijection.

\begin{lemma}[Yoneda Lemma {\cite[Section $2.1.2$]{MR2223406}}]
	Let $\mc{C}$ be a (locally-small) category and  $\mc{F}:\mc{C}^{\op}\ra \text{Set}$ a functor. For an object $U$ of $\mc{C}$, the morphism of sets $\text{Hom}(h_U,\mc{F})\ra \mc{F}(U)$ is a bijection.
\end{lemma}
Given a category $\mc{C}$, we have the functor category $[\mc{C}^{\op},\text{Set}]$ with the following description:
\begin{itemize}
	\item objects are functors $\mc{C}^{\op}\ra \text{Set}$,
	\item morphisms are natural transformations. 
\end{itemize}
For a locally small category $\mc{C}$, we have the functor $\mc{C}\ra [\mc{C}^{\op}, \text{Set}]$ with the following description:
\begin{itemize}
	\item an object $U$ of $\mc{C}$ is mapped to the functor of points $h_U:\mc{C}^{\op}\ra \text{Set}$,
	\item a morphism $f:U\ra V$ of $\mc{C}$ is mapped to the natural transformation $f^*:h_U\Rightarrow h_V$, obtained by composition.
\end{itemize}
As a corollary of the Yoneda lemma, we have the following result, which says that the above functor is an embedding of categories.
\begin{lemma}[Yoneda embedding {\cite[Corolalrry $2.2.8$]{riehl2017category}}]\label{Lemma:YonedaEmbedding}
	For a (locally small) category $\mc{C}$, the functor $\mc{C}\ra [\mc{C}^{\op},\text{Set}]$ is an embedding of categories. 
\end{lemma}

Similar to the notion of the composition of functors, we can talk about the composition of two natural transformations. There are two ways we can compose; one is that of a horizontal composition, and the other is that of the vertical composition.

\textbf{Vertical composition of natural transformations} :
Let $\mc{C}=[\mc{C}_1\rra \mc{C}_0]$ and $\mc{D}=[\mc{D}_1\rra \mc{D}_0]$ be categories. Let \[(F_1,F_0),(G_1,G_0), (H_1,H_0):[\mc{C}_1\rra \mc{C}_0]\ra [\mc{D}_1\rra \mc{D}_0]\] be functors, $\alpha:\mc{C}_0\ra \mc{D}_1$ a natural transformation from $F$ to $G$, and
$\beta:\mc{C}_0\ra \mc{D}_1$ a natural transformation from $G$ to $H$, as in the following diagram,
\[
\begin{tikzcd}[sep=large]
	\mathcal{C}_1 \arrow[dd,xshift=0.75ex,"t_{\mc{C}}"]\arrow[dd,xshift=-0.75ex,"s_{\mc{C}}"'] \arrow[rrr, "F_1", bend left, shift right] \arrow[rrr, "H_1", bend right, shift left] \arrow[rrr, "G_1"]   & & & \mc{D}_1 \arrow[dd,xshift=0.75ex,"t_{\mc{D}}"]\arrow[dd,xshift=-0.75ex,"s_{\mc{D}}"'] \\
	& & &   \\
	\mathcal{C}_0 \arrow[rrr, "F_0"', bend left, shift right] \arrow[rrr, "H_0"', bend right, shift left] \arrow[rrr, "G_0"'] \arrow[rrruu, "{\alpha, \beta}"] & & & \mc{D}_0  
\end{tikzcd}.\] 
The morphisms $\alpha:\mc{C}_0\ra \mc{D}_1, \beta:\mc{C}_0\ra \mc{D}_1$ combine to give a map $(\alpha,\beta):\mc{C}_0\ra \mc{D}_1\times_{\mathcal{D}_0}\mc{D}_1$. Composing the map $(\alpha,\beta)$ with the composition map $m_{\mc{D}}:\mc{D}_1\times_{\mathcal{D}_0}\mc{D}_1\ra \mc{D}_1$ gives the map $m_{\mc{D}}\circ (\alpha,\beta):\mc{C}_0\ra \mc{D}_1$. It is easy to check that this map $m_{\mc{D}}\circ (\alpha,\beta):\mc{C}_0\ra \mc{D}_1$ gives a natural transformation from the functor $F$ to the functor $H$ 
(\Cref{Definition:Naturaltransformation}). We call $m_{\mc{D}}\circ (\alpha,\beta)$ to be the \textit{vertical composition of the natural transformations $\alpha, \beta$}, and denote by $\beta\circ_{\textit{ver}}\alpha$. So, for $a\in \mc{C}_0$, we have 
\begin{align}
	(\beta\circ_{\textit{ver}}\alpha)(a)=\beta(a)\circ \alpha(a). 
\end{align}

\textbf{Horizontal composition of natural transformations} : 
Let $\mc{C}=[\mc{C}_1\rra \mc{C}_0], \mc{D}=[\mc{D}_1\rra \mc{D}_0]$, and $\mc{E}=[\mc{E}_1\rra \mc{E}_0]$ be categories. Let \[(F_1,F_0),(G_1,G_0):[\mc{C}_1\rra \mc{C}_0]\ra [\mc{D}_1\rra \mc{D}_0],\] \[(F'_1,F'_0),(G'_1,G'_0):[\mc{D}_1\rra \mc{D}_0]\ra [\mc{E}_1\rra \mc{E}_0],\]
be functors, and $\alpha:\mc{C}_0\ra \mc{D}_1$ and $\beta:\mc{D}_0\ra \mc{E}_1$ are natural transformations, as in the following diagram \[
\begin{tikzcd}[sep=large]
	\mathcal{C}_1 \arrow[dd,xshift=0.75ex,"t_{\mc{C}}"]\arrow[dd,xshift=-0.75ex,"s_{\mc{C}}"'] \arrow[rr, "F_1", bend left] \arrow[rr, "G_1", bend right]  & & \mathcal{D}_1 \arrow[rr, "F_1'", bend left] \arrow[dd,xshift=0.75ex,"t_{\mc{D}}"]\arrow[dd,xshift=-0.75ex,"s_{\mc{D}}"'] \arrow[rr, "G_1'", bend right]  & & \mathcal{E}_1 \arrow[dd,xshift=0.75ex,"t_{\mc{E}}"]\arrow[dd,xshift=-0.75ex,"s_{\mc{E}}"'] \\
	& &             & &    \\
	\mathcal{C}_0 \arrow[rr, "F_0"', bend left] \arrow[rruu, "\alpha"] \arrow[rr, "G_0"', bend right] & & \mathcal{D}_0 \arrow[rr, "F_0'"', bend left] \arrow[rruu, "\beta"] \arrow[rr, "G_0'"', bend right] & & \mathcal{E}_0  
\end{tikzcd}.\] Similar to the notion of vertical composition, we would like to define a natural transformation from the functor $(F_1'\circ F_1,F_0'\circ F_0)$ to the functor $(G_1'\circ G_1,G_0'\circ G_0)$. 

So, given $\alpha:\mc{C}_0\ra \mc{D}_1$ and $\beta:\mc{D}_0\ra \mc{E}_1$, we need to construct a map $\mc{C}_0\ra \mc{E}_1$. Motivated by the vertical composition, we look for a map $\mc{C}_0\ra \mc{E}_1\times_{\mc{E}_0}\mc{E}_1$ and then using the composition map $\mc{E}_1\times_{\mc{E}_0}\mc{E}_1\ra \mc{E}_1$, we get a map $\mc{C}_0\ra \mc{E}_1\times_{\mc{E}_0}\mc{E}_1\ra \mc{E}_1$. 
Given an element $a\in \mc{C}_0$, there are $4$ ways to get an element in $\mc{E}_1$, namely $(F_1'\circ\alpha)(a),(G_1'\circ\alpha)(a), (\beta\circ F_0)(a), (\beta\circ G_0)(a)$. Looking at the source, target of these $4$ morphisms, we see that 
\[ ((G_1'\circ \alpha)(a), (\beta\circ F_0)(a))\in \mc{E}_1\times_{\mc{E}_0}\mc{E}_1,\]
\[ ((\beta\circ G_0)(a),(F_1'\circ \alpha)(a))\in \mc{E}_1\times_{\mc{E}_0}\mc{E}_1.\]
Thus, there are two maps $\mc{C}_0\ra \mc{E}_1\times_{\mc{E}_0}\mc{E}_1$ given by $(G_1'\circ \alpha, \beta\circ F_0)$ and $(\beta\circ G_0, F_1'\circ \alpha)$. But, as $\beta:\mc{D}_0\ra \mc{E}_1$ is a natural transformation, we see that, $m_{\mc{E}}\circ (G_1'\circ \alpha, \beta\circ F_0)=m_{\mc{E}}\circ (\beta\circ G_0, F_1'\circ \alpha)$. Thus, we have ``unique'' map $\mc{C}_0\ra \mc{E}_1$, given by $m_{\mc{E}}\circ (G_1'\circ \alpha, \beta\circ F_0)$ or $m_{\mc{E}}\circ (\beta\circ G_0, F_1'\circ \alpha)$. It is easy to check that, this map $\mc{C}_0\ra \mc{E}_1$ satisfies the condition for a map to be a natural transformation 
(\Cref{Definition:Naturaltransformation}). We call $m_{\mc{E}}\circ (G_1'\circ \alpha, \beta\circ F_0)$ or $m_{\mc{E}}\circ (\beta\circ G_0, F_1'\circ \alpha)$ to be \textit{the horizontal composition of the natural transformations $\alpha, \beta$}, and denote by $\beta\circ_{\textit{hor}}\alpha$. So, for $a\in \mc{C}_0$, we have 
\begin{align}
	(\beta\circ_{\textit{hor}}\alpha)(a)=\beta (F_0(a))\circ G_1'( \alpha(a))= F_1'(\alpha(a))\circ \beta (G_0(a)).
\end{align}

\subsection{$2$-category theory notions} Now, we introduce some terminology from the $2$-category theory. 

\begin{definition}[$2$-category {\cite[Definition $7.1.1$]{MR1291599}}]
	A \textit{$2$-category}, denoted by $\mc{C}$, consists of the following data,
	\begin{enumerate}
		\item a collection $\mc{C}_0$, whose elements are called \textit{objects of $\mc{C}$},
		\item a small category $\mc{C}(x,y)$ for each pair of objects $(x,y)$ in $\mc{C}$,
		\item a functor $\{1\}\ra \mc{C}(x,x)$ for each object $x$ in $\mc{C} $,
		\item a functor $\mc{C}_{xyz}:\mc{C}(x,y)\times\mc{C}(y,z)\ra \mc{C}(x,z)$ for each triple $(x,y,z)$ of objects in $\mc{C}$,
	\end{enumerate}
	satisfying the following conditions:
	\begin{enumerate}
		\item for each $x,y\in \mc{C}_0$, the following diagram is commutative,\[
		\begin{tikzcd}
			& & {\mc{C}(x,y)} \arrow[dd, "1"'] \arrow[rrd, "{(1,1_y)}"] \arrow[lld, "{(1_x,1)}"'] & &        \\
			{\mc{C}(x,x)\times\mc{C}(x,y)} \arrow[rrd, "\mc{C}_{xxy}"'] & &           & & {\mc{C}(x,y)\times \mc{C}(y,y)} \arrow[lld, "\mc{C}_{xyy}"] \\
			& & {\mc{C}(x,y)}         & &        
		\end{tikzcd},\]
		\item for each $x,y,z,w\in \mc{C}_0$, the following diagram is commutative,
		\[
		\begin{tikzcd}
			& {\mc{C}(x,y)\times\mc{C}(y,z)\times\mc{C}(z,w)} \arrow[ld, "{(1,\mc{C}_{yzw})}"'] \arrow[rd, "{(\mc{C}_{xyz},1)}"] &        \\
			{\mc{C}(x,y)\times\mc{C}(y,w)} \arrow[rd, "\mc{C}_{xyw}"'] &               & {\mc{C}(x,z)\times \mc{C}(z,w)} \arrow[ld, "\mc{C}_{xzw}"] \\
			& {\mc{C}(x,w)}            &        
		\end{tikzcd}.\] 
	\end{enumerate}
\end{definition}
Before we proceed further, we fix some notation.
\begin{remark}[Notation]
	A \textit{$1$-morphism in a $2$-category $\mc{C}$} is an object in the category $\mc{C}(x,y)$ for some $x,y$ in $\mc{C}_0$. A \textit{$2$-morphism in a $2$-category $\mc{C}$} is a morphism in the category $\mc{C}(x,y)$ for some $x,y$ in $\mc{C}_0$.
\end{remark}

The word ``composition'' occurs in $3$ places in a $2$-category:
\begin{enumerate}
	\item Let $f,g$ be $1$-morphisms in a $2$-category $\mc{C}$. If $(f,g)\in \mc{C}(x,y)_0\times \mc{C}(y,z)_0$ for some $x,y,z\in \mc{C}_0$, then, the functor $\mc{C}_{xyz}$ gives an element $\mc{C}_{xyz}(f,g)\in \mc{C}(x,z)_0$. We denote this element $\mc{C}_{xyz}(f,g)$ by $g\circ f$.
	
	\item Let $\alpha, \beta$ be $2$-morphisms in a $2$-category $\mc{C}$. There are two possibilities here: either $\alpha,\beta\in \mc{C}(x,y)_1$ for some $x,y\in \mc{C}_0$, or $(\alpha,\beta)\in \mc{C}(x,y)_1\times \mc{C}(y,z)_1$ for some $x,y,z\in \mc{C}_0$.
	
	\begin{itemize}
		\item Suppose that $\alpha,\beta\in \mc{C}(x,y)_1$. Further, if $t(\alpha)=s(\beta)$, then, the composition map in the category $\mc{C}(x,y)$ gives an element $m_{\mc{C}(x,y)}(\alpha,\beta)$. We denote this element by $\beta\odot\alpha$, and call it \textit{the vertical composition of the $2$-morphisms $\alpha, \beta$}. We see it as the following diagram,
		
		\[\begin{tikzcd}[column sep=2cm]
			\vphantom{f}x
			\arrow[r, "f", ""{name=F, below}, bend left=49]
			\arrow[r, "h"', ""{name=H, above}, bend right=49]
			\arrow[r, "g" description, ""{name=GA,above}, ""{name=GB,below}]
			& \vphantom{f}y
			\arrow[Rightarrow, from=F, to=GA, "\alpha", shorten <=-3pt]
			\arrow[Rightarrow, from=GB, to=H, "\beta", shorten >=-3pt]
		\end{tikzcd} \text{ gives } 
		\begin{tikzcd}
			x
			\arrow[rr, "f"{name=F}, bend left=49]
			\arrow[rr, "h"'{name=H}, bend right=49]
			& & y
			\arrow[Rightarrow, from=F, to=H, "\beta\odot \alpha"]
		\end{tikzcd}.\]
		\item Suppose that $(\alpha,\beta)\in \mc{C}(x,y)_1\times \mc{C}(y,z)_1$ for some $x,y,z\in \mc{C}_0$. The morphism level map of the functor
		$\mc{C}_{xyz}:\mc{C}(x,y)\times \mc{C}(y,z)\ra \mc{C}(x,z)$ gives an element $\mc{C}_{xyz}(\alpha,\beta)\in \mc{C}(x,z)_1$. We denote this element $\mc{C}_{xyz}(\alpha,\beta)\in \mc{C}(x,z)_1$ by $\beta\star\alpha$, and call it the horizontal composition of the $2$-morphisms $\alpha,\beta$. We see it as the following diagram, 
		\[\begin{tikzcd}
			x \arrow[rr, "f"{name=F}, bend left=49]
			\arrow[rr, "\mathclap{g}"'{name=G}, bend right=49] \arrow[Rightarrow, "\alpha", from=F, to=G]
			&& y \arrow[rr, "h"{name=H}, bend left=49]
			\arrow[rr, "\mathclap{l}"'{name=L}, bend right=49]
			\arrow[Rightarrow, "\beta", from=H, to=L] && z 
		\end{tikzcd}\text{ gives } 
		\begin{tikzcd}
			x \arrow[rr, "h\circ f"{name=HF}, bend left=49] \arrow[rr, "l\circ g"'{name=LG}, bend right=49] \arrow[Rightarrow, "\beta\star \alpha", from=HF, to=LG]& & z
		\end{tikzcd}. \]
	\end{itemize}	
\end{enumerate}

\begin{definition}[{\cite[Definition A.12.]{MR3762701}}]\label{Definition:2-isomorphism}
	A $2$-morphism in a $2$-category is a \textit{$2$-isomorphism} if it is invertible with respect to vertical composition of $2$-morphisms.
\end{definition}
\begin{definition}[2-commutative diagram {\cite[Appendix A]{MR3762701}}]\label{Definition:2-commutativediagram}
	Let $\mc{C}$ be a $2$-category. A diagram
	\[
	\begin{tikzcd}
		A \arrow[rr, "f"] \arrow[dd, "k"'] & & B \arrow[dd, "g"] \\
		& &   \\
		C \arrow[rr, "h"] \arrow[Rightarrow, shorten >=20pt, shorten <=20pt, uurr, "\alpha"]   & & D  
	\end{tikzcd},\]
	in $\mc{C}$ is said to be a \textit{$2$-commutative diagram}, if $\alpha:g\circ f\Rightarrow h\circ k$ is a $2$-isomorphism in the $2$-category.
\end{definition}

We end this section with some miscellaneous concepts from category theory. 
\subsection{Pull-back in a category} One notion from category theory we come across in this thesis quite often is that of a pull-back in a category. 

\begin{definition}[pull-back in a category]
	Let $\mc{C}$ be a category, and $\gamma:a\ra b, \delta:c\ra b$ morphisms in $\mc{C}$. A \textit{pull-back of the pair $(\gamma,\delta)$} is given by an object $d$ of $\mc{C}$, along with morphisms ${\pr}_2:d\ra c, {\pr}_1:d\ra a$, such that, $\delta\circ {\pr}_2=\gamma\circ {\pr}_1$, and for any pair of morphisms $(\tau:e\ra c, \sigma:e\ra a)$ with $\delta\circ \tau=\gamma\circ \sigma$, there exists a unique morphism $\Phi:e\ra d$ with ${\pr}_2\circ \Phi=\tau$, and ${\pr}_1\circ \Phi=\sigma$. We see a pull-back as the following diagram, 
	\[
	\begin{tikzcd}[sep=large]
		e \arrow[rd, "\Phi", dashed] \arrow[rrrd, "\tau"] \arrow[rddd, "\sigma"'] &      & &   \\
		& d \arrow[dd, "{\pr}_1"'] \arrow[rr, "{\pr}_2"] & & c \arrow[dd, "\delta"] \\
		&      & &   \\
		& a \arrow[rr, "\gamma"]   & & b   
	\end{tikzcd}.\]\end{definition}
Alternatively, we use the terminology ``fiber product" to mean the pull-back. 
\begin{remark}
	Not all categories have all pull-backs. If a pull-back of a pair $(\gamma,\delta)$ exists in a category, then it is unique up to a unique isomorphism. 
\end{remark}
\begin{example}
	Let ${\rm Set}$ be the category of sets. Let $f:A\ra B$ and $g:C\ra B$ be set maps. Then, the pull-back of the pair $(f,g)$ is given by the set $\{(a,c)\in A\times C:f(a)=g(c)\}$. We denote this set as $A\times_BC$. 
\end{example}
\begin{example}
	Let ${\rm Top}$ be the category of topological space. Let $f:A\ra B$ and $g:C\ra B$ be continuous maps. Then, the pull-back of the pair $(f,g)$ is given by the set $\{(a,c):f(a)=g(c)\}$, with subspace topology induced from the product topology on $A\times C$. We denote this space as $A\times_BC$. 
\end{example}
\begin{example}
	Let $\text{Man}$ be the category of smooth manifolds. Let $f:M\ra N$ and $g:M'\ra N$ be smooth maps. Then, the set-theoretic pull-back $M\times_NM'$ does not always have a nice smooth structure to become a pull-back in the category of smooth manifolds. When either $f:M\ra N$ or $g: M'\ra N$ is a submersion, then $M\times_NM'$ becomes the pull-back of the pair $(f,g)$, in the category of smooth manifolds. 
\end{example}

\subsection{Epimophisms in a category}
We come across a special kind of morphism in a category, namely an epimorphism. One of the good references to read more about epimorphisms is \cite[$1.8$ Epimorphisms]{MR1291599}.
\begin{definition}[epimorphism in a category]
	Let $\mc{C}$ be a category. A morphism $\theta:X\ra Y$ in $\mc{C}$ is said to be \textit{an epimorphism}, if for any $g_1,g_2:Y\ra Z$ with $g_1\circ \theta=g_2\circ \theta$ implies $g_1=g_2$.
\end{definition}

Consider the category $\text{Set}$, whose objects are sets and whose morphisms are functions. Then, epimorphisms in this category are precisely the surjective functions in it \cite[Example  $1.8.5 (a)$]{MR1291599}.

Let $\mc{C}$ be a category whose objects are sets with extra structure, whose morphisms are functions with extra conditions (such categories, up to some technicality, are called concrete categories). 

Let $\theta:X\ra Y$ be a morphism in $\mc{C}$. Ignoring the extra conditions on $\theta$, this map $\theta$ can be seen as a morphism in $\text{Set}$. 
If $\theta$ is an epimorphism in $\text{Set}$, then $g_1\circ \theta=g_2\circ \theta$ implies $g_1=g_2$ for all appropriate $g_1,g_2$ in $\text{Set}$. In particular, $g_1\circ \theta=g_2\circ \theta$ implies $g_1=g_2$ for all appropriate $g_1,g_2$ in $\mc{C}$. Thus, every morphism in $\mc{C}$, whose underlying set map is an epimorphism in $\text{Set}$, will be an epimorphism in $\mc{C}$. In other words, every morphism in $\mc{C}$, whose underlying set map is surjective, is an epimorphism in $\mc{C}$.

Suppose $\theta$ is an epimorphism in $\mc{C}$, then, it only means that $g_1\circ \theta=g_2\circ \theta$ implies $g_1=g_2$ for all appropriate $g_1,g_2$ in $\mc{C}$. This does not say that $g_1\circ \theta=g_2\circ \theta$ implies $g_1=g_2$ for all appropriate $g_1,g_2$ in $\text{Set}$, as the collection of appropriate  morphism in $\text{Set}$ may be bigger than the collection in $\mc{C}$. Thus, not every epimorphism in $\mc{C}$ needs to be a surjective map on the underlying sets.


\section{Fibered categories, pseudo-functors, and stacks}\label{Section:Fiberedcategoriespseudofunctorsandstacks}
In this section, we recall the notions of fibered categories, pseudo-functors and stacks. Most of this section is based on Angelo Vistoli's ``Notes on Grothendieck topologies, fibered categories and descent theory'' \cite{MR2223406}.
\subsection{Sheaf on a topological space} Let $X$ be a topological space and $\mc{O}(X)$ the category of open subsets of $X$ (\Cref{Example:OXcategory}). A functor $\mc{F}\colon \mathcal{O}(X)^{\op}\ra\text{Set}$ is called a \textit{presheaf (of sets) on the topological space $X$}. For an object $U$ of $\mc{O}(X)$, we call the elements of $\mc{F}(U)$ as the sections of $U$. 
For a morphism $f:U\ra V$ in $\mc{O}(X)$, we denote the map $\mc{F}(f):\mc{F}(V)\ra \mc{F}(U)$ as $f^*:\mc{F}(V)\ra \mc{F}(U)$. 

Let $U$ be an open subset of $X$ and $\{U_i\ra U\}_{i\in \Lambda}$ an open cover of $U$. For indices $i,j\in \Lambda$, let $U_{ij}$ denote the pull-back of the inclusion maps $\sigma_i:U_i\ra U$ and $\sigma_j:U_j\ra U$; that is, $U_{ij}=U_i\cap U_j$. We denote the corresponding (inclusion) maps $U_{ij}\ra U_i$ and $U_{ij}\ra U_j$ as ${\pr}_1$ and ${\pr}_2$ respectively, as in the following (commutative) diagram,
\[
\begin{tikzcd}
	U_{ij}=U_i\cap U_j \arrow[dd, "{\pr}_2"'] \arrow[rr, "{\pr}_1"] & & U_i \arrow[dd, "\sigma_i"] \\
	& &     \\
	U_j \arrow[rr, "\sigma_j"']      & & U    
\end{tikzcd}.\]
The functor $\mathcal{F}:\mc{O}(X)^{\op}\ra \text{Set}$, maps the above commutative diagram, to the following commutative diagram,
\[
\begin{tikzcd}
	\mathcal{F}(U) \arrow[dd, "\sigma_j^*"'] \arrow[rr, "\sigma_i^*"] & & \mathcal{F}(U_i) \arrow[dd, "{\pr}_1^*"] \\
	& &      \\
	\mathcal{F}(U_j) \arrow[rr, "{\pr}_2^*"']   & & \mathcal{F}(U_i\cap U_j)  
\end{tikzcd}.\]
Consider the set $\mc{F}(\{U_i\ra U\})$ with the following description:
\begin{itemize} 
	\item an element of $\mc{F}(\{U_i\ra U\})$ is a collection $\{s_{i}\}$, where $s_{i}\in \mc{F}(U_i)$ for each index $i\in \Lambda$ and ${\pr}_1^*(s_{i}) 
	={\pr}_2^*(s_{j})$ in $\mc{F}(U_{ij})$ for each pair of indices $i,j\in \Lambda$.
\end{itemize}

Given an open set $U$ of $X$ and an open cover $\{\sigma_i:U_i\ra U\}$ of $U$, we have the function $\mc{F}(U)\ra \mc{F}(\{U_i\ra U\})$ given by $s\mapsto\{\sigma_i^*(s)\}_{i \in \Lambda}$. 

\begin{definition}[Sheaf on a topological space] \label{Definition:SheafOnTopspace} Let $X$ be a topological space. Let $\mc{F}:\mc{O}(X)^{\op}\ra \text{Set}$ be a functor. We say that $\mc{F}$ is a \textit{locally determined functor}, if, for each open set $U$ of $X$ and an open cover $\{U_i\ra U\}$ of $U$, the function $\mc{F}(U)\ra \mc{F}(\{U_i\ra U\})$ is a bijection. A locally determined functor $\mc{F}\colon \mc{O}(X)^{\op}\ra \text{Set}$ is said to be \textit{a sheaf (of sets) on the topological space $X$}.
\end{definition} 

\begin{example}\label{Example:SheafContinuousMap}
	Consider the functor $\mc{F}\colon \mc{O}(X)^{\op}\ra \text{Set}$ with the following description:
	\begin{itemize}
		\item for an object $U$ of $\mc{O}(X)$, we have $\mc{F}(U)=\{\text{continuous maps } U\ra \mb{R}\}$,
		\item for a morphism $f:U\ra V$ of $\mc{O}(X)$, the morphism 
		$\mc{F}(f):\mc{F}(V)\ra \mc{F}(U)$ is given by restriction.	\end{itemize}
	This functor $\mc{F}$ is locally determined, thus gives a sheaf on $X$.
\end{example}

\subsection{Sheaf on a site} Let $\mc{C}$ be a category. A functor $\mc{F}\colon \mc{C}^{\op}\ra\text{Set}$ is called a \textit{presheaf (of sets) on the category $\mc{C}$}. 

The notion of   ``a functor $\mc{F}:\mc{O}(X)^{\op}\ra \text{Set}$ being locally determined'' can be extended to a functor $\mc{F}:\mc{C}^{\op}\ra \text{Set}$, for a general category $\mc{C}$, when $\mc{C}$ is equipped with the so called Grothendieck topology, which we recall below.
\begin{definition}[Grothendieck topology 
	{\cite[Definition $2.24$]{MR2223406}}] \label{Definition:GrothendieckTopology}
	Let $\mathcal{C}$ be a category. For an object $U$ of $\mc{C}$, a \textit{covering of $U$} is a collection of morphisms in $\mc{C}$ whose target is $U$. Consider a collection $\mc{J}=\{\mc{J}_U:U\in \text{Obj}(\mc{C})\}$, where $\mc{J}_U$ be a collection of coverings of $U$.  We say that $\mc{J}$ is a \textit{Grothendieck topology on $\mc{C}$} if
	the following conditions are satisfied:
	\begin{enumerate}
		\item for every object $U$ of $\mc{C}$ and an isomorphism $V\ra U$ of $\mc{C}$, the covering $\{V\ra U\}$ is in $\mc{J}_U$,
		\item for every object $U$ of $\mc{C}$, a covering $\{U_i\rightarrow U\}_{i\in \Lambda}$ in $\mc{J}_U$, and a covering $\{U_{ij}\rightarrow U_i\}_{j\in \Lambda_i}$ in $\mc{J}_{U_i}$, for each $i\in \Lambda$, the composition covering $\{U_{ij}\rightarrow U\}_{i\in \Lambda, j\in \Lambda_i}$ is in $\mc{J}_U$, 
		\item for every object $U$ of $\mc{C}$, a covering $\{U_i\ra U\}_{i\in \Lambda}$ in $\mc{J}_U$ and a morphism $V\ra U$ in $\mc{C}$, the fiber product $U_i\times_UV$ exists for each $i\in \Lambda$, and the covering $\{{\pr}_2\colon U_i\times_U V\ra V\}_{i\in \Lambda}$ is in $\mc{J}_V$.
	\end{enumerate}
	A category $\mc{C}$ equipped with a Grothendieck topology $\mc{J}$ is called a \textit{site}.  We denote a site as a pair $(\mc{C},\mc{J})$, or simply as $\mc{C}$ when it is clear what $\mc{J}$ is. 
\end{definition}
\begin{remark}
	Only for the convenience of defining the notion of a Grothendieck topology, we used the phrase ``covering of $U$'' for an arbitrary collection of morphisms $\{U_i\rightarrow U\}_{i\in \Lambda}$. Once we fix a Grothendieck topology $\mc{J}$ on $\mc{C}$, by a covering of $U$, we always mean a collection $\{U_i\ra U\}_{i\in \Lambda}$ with $\{U_i\ra U\}_{i\in \Lambda}\in \mc{J}_U$. 
\end{remark}

\begin{example}[big \'etale topology on $\text{Man}$]\label{Example:big\'etaletopology}
	Consider the category $\text{Man}$ of smooth manifolds. For an object $M$ of the category $\text{Man}$, define a covering of $M$ to be a collection of \'etale morphisms $\{\sigma_i:U_i\ra M\}_{i\in \Lambda}$ such that $\bigcup_{i\in \Lambda}\sigma_i(U_i)=M$. This gives a Grothendieck topology on $\text{Man}$, which we call the big \'etale topology on $\text{Man}$. We denote this site by $(\text{Man}, \mc{J}_{{e}tale})$.
\end{example}

Before we proceed further, we fix some notations. Let $(\mc{C},\mc{J})$ be a site. Let $U$ be an object of $\mc{C}$ and $\{U_i\ra U\}_{i\in \Lambda}$ a covering of $U$. 
For indices $i,j\in \Lambda$, let $U_{ij}$ denote the fiber product of morphisms $\sigma_i :U_i\ra U$ and $\sigma_j:U_j\ra U$, as in the following diagram,
\begin{equation}
	\label{Diagram:2pull-back}
	\begin{tikzcd}
		U_{ij}=U_i\times_UU_j \arrow[dd, "{\pr}_1"] \arrow[rr, "{\pr}_2"] & & U_j \arrow[dd,"\sigma_j"] \\
		& &   \\
		U_i \arrow[rr,"\sigma_i"]        & & U   
	\end{tikzcd}.
\end{equation}

For indices $i,j,k\in \Lambda$, let $U_{ijk}$ be the pullback of the projections
$U_{ij}\ra U_j$, and $U_{jk}\ra U_j$. The universal property of a pullback  gives a unique morphism $U_{ijk}\ra U_{ik}$, which we denote by ${\rm pr}_{13}$. It then follows that, every square in the below diagram commutes,
\begin{equation}\label{Diagram:3pull-back}
	\begin{tikzcd}
		& U_{ijk} \arrow[dd, " {\pr}_{13}"] \arrow[rr, "{\pr}_{23}"] \arrow[ld, "{\pr}_{12}"'] &   & U_{jk} \arrow[dd] \arrow[ld] \\
		U_{ij} \arrow[dd] \arrow[rr] &            & U_j \arrow[dd,"\sigma_j"] &     \\
		& U_{ik} \arrow[rr] \arrow[ld]       &   & U_k \arrow[ld,"\sigma_k"]   \\
		U_i \arrow[rr,"\sigma_i"]   &            & U   &     
	\end{tikzcd}.
\end{equation}

Let $\mc{F}:\mc{C}^{\op}\ra \text{Set}$ be a functor. Let $U$ be an object of $\mc{C}$ and $\{U_i\ra U\}_{i\in \Lambda}$ a covering of $U$. Consider the set $\mc{F}(\{U_i\ra U\})$ with the following description:
\begin{itemize}
	\item an element of $\mc{F}(\{U_i\ra U\})$ is of the form $\{s_i\}_{i\in \Lambda}$, where $s_i\in \mc{F}(U_i)$ for each $i\in \Lambda$ and ${\pr}_2^*(s_j)={\pr}_1^*(s_i)$ in $\mc{F}(U_{ij})$ for each pair of indices $i,j\in \Lambda$.
\end{itemize}
Let $U$ be an object of $\mc{C}$ and $\{U_i\ra U\}$ a covering of $U$. Then, we have the function $\mc{F}(U)\ra \mc{F}(\{U_i\ra U\})$ defined by $s\mapsto \{ \sigma_i^*(s)\}_{i\in \Lambda}$.

\begin{definition}[Sheaf on a site] \label{Definition:Sheafonasite} Let $(\mc{C},\mc{J})$ be a site.	A functor $\mc{F}:\mc{C}^{\op}\ra \text{Set}$ is said to be a \textit{locally determined functor}, if, for each object $U$ of $\mc{C}$ and a covering $\{U_i\ra U\}$ of $U$, the function $\mc{F}(U)\ra \mc{F}(\{U_i\ra U\})$ is a bijection. A locally determined functor $\mc{F}:\mc{C}^{\op}\ra \text{Set}$ is said to be a \textit{sheaf on the site $(\mc{C},\mc{J})$}.
\end{definition}

\begin{remark}\label{Remark:subcanonicalsite}
	A site $(\mc{C},\mc{J})$ is called a \textit{subcanonical site} if for each object $U$ of $\mc{C}$, the functor of points $h_U:\mc{C}^{\op}\ra \text{Set}$ is a sheaf on the site $(\mc{C},\mc{J})$. In this thesis, all our sites are assumed to be subcanonical.
\end{remark}

\subsection{Stack on a site} Let $(\mc{C},\mc{J})$ be a site. 
The notion of ``a functor $\mc{F}:\mc{C}^{\op}\ra \text{Set}$ being locally determined'' can be extended to a pseudo-functor $\mc{F}:\mc{C}^{\op}\ra \text{Cat}$.
%
%
Firstly, we recall the notion of a pseudo-functor on a category $\mc{C}$.

\begin{definition}[{\cite[Definition $3.10$]{MR2223406}}]\label{Definition:pseudofunctor} Let $\mc{C}$ be a category.
	A \textit{pseudo-functor on $\mc{C}$}, denoted by $\mc{F}\colon \mc{C}^{\op}\ra \text{Cat}$, consists of the following data:
	\begin{enumerate}
		\item a category $\mc{F}(U)$ for each object $U$ of $\mc{C}$,
		\item a functor $f^*\colon \mc{F}(V)\ra \mc{F}(U)$ for each morphism $f\colon U\ra V$ of $\mc{C}$,
		\item a natural isomorphism of functors $\epsilon_U\colon (1_U)^*\Rightarrow 1_{\mc{F}(U)}\colon \mc{F}(U)\ra \mc{F}(U)$ for each object $U$ of $\mc{C}$,
		\item a natural isomorphism of functors $\alpha_{f,g}\colon f^*g^* \Rightarrow (g\circ f)^*\colon \mc{F}(W)\ra \mc{F}(U)$ for each pair of morphisms $U\xra{f}V\xra{g}W$ of $\mc{C}$,
	\end{enumerate}
	satisfying the following conditions:
	\begin{enumerate}
		\item for every morphism $f\colon U\ra V$ in $\mc{C}$, and an object $\eta$ of $\mc{F}(V)$, we have 
		\[\alpha_{f,1_V}(\eta)=f^*(\epsilon_V\eta), \text{~and~} \alpha_{1_U,f}(\eta)=\epsilon_U(f^*\eta),\]
		\item for every triple of composable morphisms $U\xra{f}V\xra{g}W\xra{h}T$ in $\mc{C}$, and an object $\eta$ of $\mc{F}(T)$, we have 
		\[\alpha_{gf,h}(\eta)\circ \alpha_{f,g}(h^*\eta)=
		\alpha_{f,hg}(\eta)\circ f^*\alpha_{g,h}(\eta).\]
	\end{enumerate}
\end{definition}

Let $(\mc{C},\mc{J})$ be a site. Consider a pseudo-functor $\mc{F}:\mc{C}^{\op}\ra \text{Cat}$. Let $U$ be an object of $\mc{C}$ and $\{U_i\ra U\}_{i\in \Lambda}$ a covering of $U$. Consider
the category $\mc{F}(\{U_i\ra U\})$, with the following description:
\begin{itemize}
	\item an object in $\mc{F}(\{U_i\ra U\})$ is a collection $\big\{\{s_i\}_{i\in \Lambda},\{\phi_{ij}\}_{i,j\in \Lambda}\big\}$, where
	\begin{enumerate}
		\item $s_i$ is an object of $\mc{F}(U_i)$ for each $i\in \Lambda$,
		\item $\phi_{ij}\colon {\pr}_2^*s_j\ra {\pr}_1^*s_i$ is an isomorphism in $\mc{F}(U_{ij})$ for each $i,j\in \Lambda$,
	\end{enumerate}
	satisfying the condition $ {\pr}_{13}^*(\phi_{ik})={\pr}_{12}^*(\phi_{ij})\circ {\pr}_{23}^*(\phi_{jk})$ in the category $\mc{F}(U_{ijk})$ for each $i,j,k\in \Lambda$,
	\item a morphism from an object $\big\{\{s_i\},\{\phi_{ij}\}\big\}$ to an object $\big\{\{t_i\},\{\psi_{ij}\}\big\}$ is given by a collection $\{\theta_i\colon s_i\ra t_i \}_{i\in \Lambda}$, here, $\theta_i$ is a morphism in $\mc{F}(U_i)$ for each $i\in \Lambda$, satisfying the condition ${\pr}_1^*(\theta_i)\circ \phi_{ij}=\psi_{ij}\circ {\pr}_2^*(\theta_j)$ in $\mc{F}(U_{ij})$ for each $i,j\in \Lambda$.
\end{itemize}
We call the category $\mc{F}(U_i\ra U)$ to be \textit{the descent category of the covering $\{U_i\ra U\}$}.

Let $U$ be an object of $\mc{C}$ and $\{\sigma_i:U_i\ra U\}$ a covering of $U$. Consider the functor $\mc{F}(U)\ra \mc{F}(\{U_i\ra U\})$ with the following description:
\begin{itemize}
	\item an object $s$ of $\mc{F}(U)$ is mapped to the object $\{\{\sigma_i^*(s)\},\{\phi_{ij}\}\}$ of $\mc{F}(\{U_i\ra U\})$, where $\phi_{ij}:{\pr}_2^*(\sigma_j^*(s))\ra 	{\pr}_1^*(\sigma_i^*(s))$ is the isomorphism given by composing the morphisms $\alpha_{{\pr}_1,\sigma_i}(s):{\pr}_1^*(\sigma_i^*(s))\ra (\sigma_i\circ {\pr}_1)^*(s)$ and $ \alpha_{{\pr}_2,\sigma_j}(s):{\pr}_2^*(\sigma_j^*(s))\ra (\sigma_j\circ {\pr}_2)^*(s)$,
	\item a morphism $\theta:s\ra t$ in $\mc{F}(U)$ is mapped to the morphism $\{\sigma_i^*(\theta):\sigma_i^*(s)\ra \sigma_i^*(t)\}_{i\in \Lambda}$ in $\mc{F}(\{U_i\ra U\})$.
\end{itemize}

\begin{definition}[stack on a site]\label{Definition:StackFiberedCategory}
	Let $(\mc{C},\mc{J})$ be a site. A pseudo-functor $\mc{F}:\mc{C}^{\op}\ra \text{Cat}$ is said to be a \textit{locally determined pseudo-functor}, if, for each object $U$ of $\mc{C}$ and a covering $\{U_i\ra U\}$ the functor $\mc{F}(U)\ra \mc{F}(\{U_i\ra U\})$ is an equivalence of categories. A locally determined pseudo-functor on a site $(\mc{C},\mc{J})$ is said to be a \textit{stack on the site $(\mc{C},\mc{J})$}.
\end{definition}

\subsection{Fibered categories}It is sometimes useful to think of a pseudo-functor $\mc{F}:\mc{C}^{\op}\ra \text{Cat}$ as a fibered category $\pi_{ \mc{D}}:\mc{D}\ra \mc{C}$. In this subsection, we recall the notion of a fibered category, give an outline of correspondence between pseudo-functors and fibered categories, and introduce the notion of a stack as a fibered category.

\begin{definition}[{\cite[Definition $3.1$]{MR2223406}}]\label{Definition:Cartesianarrow}
	Let $\mc{C}$ be a category and $\pi_{\mc{D}}:\mc{D}\ra\mc{C}$ a functor. A morphism $\theta\colon \xi\ra \eta$ in $\mc{D}$ is said to be a \textit{Cartesian arrow in $\mc{D}$} if, for every morphism $\theta'\colon \xi'\ra \eta$ in $\mc{D}$ and a morphism $h\colon \pi_{\mc{D}}(\xi')\ra \pi_{\mc{D}}(\xi)$ in $\mc{C}$ with $\pi_{\mc{D}}(\theta)\circ h=\pi_{\mc{D}}(\theta')$, there exists a unique morphism $\Phi\colon \xi'\ra \xi$ in $\mc{D}$ such that, $\theta\circ \Phi=\theta'$ in $\mc{D}$ and $\pi_{\mc{D}}(\Phi)=h$ in $\mc{C}$. We visualize a Cartesian arrow as the morphism $\theta$ in the following diagram,
	\begin{equation}\label{Diagram:Cartesianarrow}
		\begin{tikzcd}
			\xi' \arrow[rd,"\theta'"] \arrow[dd, dotted,"\Phi"'] \arrow[rrr,maps to] & & & \pi_{\mc{D}}(\xi') \arrow[rd,"\pi_{\mc{D}}(\theta')"] \arrow[dd,"h"'] & \\
			& \eta \arrow[rrr, maps to] & & & \pi_{\mc{D}}(\eta) \\
			\xi \arrow[ru,"\theta"'] \arrow[rrr, maps to] & & & \pi_{\mc{D}}(\xi) \arrow[ru,"\pi_{\mc{D}}(\theta)"'] & 
		\end{tikzcd}.
	\end{equation}
\end{definition}

\begin{definition}[{\cite[Definition $3.5$]{MR2223406}}] \label{Definition:fiberedcategory}
	Let $\mc{C}$ be a category and $\pi_{\mc{D}}:\mc{D}\ra\mc{C}$ a functor. We say that $(\mc{D},\pi_{\mc{D}},\mc{C})$ is a \textit{fibered category over $\mc{C}$}, if, for every
	$(\eta,f)\in \mc{D}_0\times_{\mc{C}_0,t}\mc{C}_1$, there exists a Cartesian arrow $\theta\colon \xi\ra \eta$ with $\pi_{\mc{D}}(\theta)=f$. We call such $\xi$ to be \textit{a pull-back of $\eta$ along $f$}.
\end{definition}
A pull-back of an object of $\mc{D}$ along a morphism of $\mc{C}$, if it exists, is unique up to a unique isomorphism.

Following are some examples of fibered categories that we come across in this thesis often.
\begin{example}\label{Example:fiberedcategoryCU}
	Let $\mc{C}$ be a category, and $U$ an object of $\mc{C}$. Let $\mc{C}/U$ be the category with the following description:
	\begin{itemize}
		\item an object of $\mc{C}/U$ is a morphism in $\mc{C}$ whose target is $U$; that is of the form $f:V\ra U$ for some object $V$ of $\mc{C}$,
		\item a morphism in $\mc{C}/U$ from an object $(f:V\ra U)$ to an other object $(f':V'\ra U)$ is a morphism $g:V\ra V'$ in $\mc{C}$ with $f'\circ g=f$.
	\end{itemize}
	Consider the functor $\pi_U:\mc{C}/U\ra \mc{C}$ with the following description:
	\begin{itemize}
		\item an object $(f:V\ra U)$ of $\mc{C}/U$ is mapped to the object $V$ of $\mc{C}$,
		\item a morphism $g:(f:V\ra U)\ra (f':V'\ra U)$ of $\mc{C}/U$ is mapped to the morphism $g:V\ra V'$ of $\mc{C}$.
	\end{itemize}
	The functor $\pi_U:\mc{C}/U\ra \mc{C}$ is a fibered category over $\mc{C}$. We denote this fibered category by $\underline{U}\ra \mc{C}$ or just by $\underline{U}$.
\end{example}
\begin{example}\label{Example:FiberedCategoryBG}
	Let $G$ be a Lie group. Consider the category $BG$, whose objects are principal $G$-bundles and morphisms are morphisms of principal $G$-bundles. Consider the functor $\pi_G:BG\ra \text{Man}$ with the following description:
	\begin{itemize}
		\item an object $P(M,G)$ of $BG$ is mapped to the object $M$ of $\text{Man}$,
		\item a morphism $(F,f):P(M,G)\ra P'(M',G)$ of $BG$ is mapped to the morphism $f:M\ra M'$ of $\text{Man}$.
	\end{itemize}
	The functor $\pi_G:BG\ra \text{Man}$ is a fibered category. 
\end{example}
\begin{definition}[{\cite[Definition $3.6$]{MR2223406}}]\label{Definition:MorphismofFiberedcategories}
	Let $\mc{C}$ be a category. Let $(\mc{D},\pi_{\mc{D}},\mc{C})$ and $(\mc{D}',\pi_{\mc{D}'},\mc{C})$ be fibered categories. A \textit{morphism of fibered categories} from 
	$(\mc{D},\pi_{\mc{D}},\mc{C})$ to $(\mc{D}',\pi_{\mc{D}'},\mc{C})$ is given by a functor $F\colon \mc{D}\ra \mc{D}'$ such that,
	\begin{enumerate}
		\item the functor $F$ is compatible with the functors $\pi_{\mc{D}}$ and $\pi_{\mc{D}'}$; that is, $\pi_{\mc{D}'}\circ F=\pi_{\mc{D}}$,
		\item the functor $F$ maps Cartesian arrows in $\mc{D}$ to Cartesian arrows in $\mc{D}'$.
	\end{enumerate}
\end{definition}

\begin{definition}\label{Definition:isomorphismoffiberedcategories}
	A morphism of fibered categories $F:(\mc{D},\pi_{\mc{D}},\mc{C})\ra (\mc{D}',\pi_{\mc{D}'},\mc{C})$ is called an \textit{isomorphism of fibered categories}, if, the underlying functor $F:\mc{D}\ra \mc{D}'$ is an equivalence of categories.
\end{definition}
\begin{remark}
	As mentioned in \cite[Section $1.1$]{noohiquick}, ``It is not uncommon in the literature to call two equivalent fibered categories isomorphic''.
\end{remark}

Let $A,B$ be sets and $f:A\ra B$ a set map. For $b\in B$, the fiber of $b$ is defined to be the set $\{a\in A:f(a)=b\}$. Similar to this notion of the fiber of an element in a set map, we define the notion of fiber of an object of $\mc{C}$ in a fibered category $(\mc{D},\pi_{\mc{D}},\mc{C})$. For obvious reasons, fiber of an object $U$ of $\mc{C}$ will be a category.
\begin{definition}[{\cite[Definition $3.8$]{MR2223406}}]\label{Definition:Fiberofanobject}
	Let $\mc{C}$ be a category and $(\mc{D},\pi_{\mc{D}},\mc{C})$ a fibered category. For an object $U$ of $\mc{C}$, \textit{the fiber of $U$}, denoted by $\mc{D}(U)$, is a category with the following description:
	\begin{align*}
		\text{Obj}(\mc{D}(U))&=\{\eta\in \text{Obj}(\mc{D})| \pi_{\mc{D}}(\eta)=U \},\\
		\text{Mor}_{\mc{D}(U)}(\eta,\eta')&=\{f\in \text{Mor}(\mc{D})|\pi_{\mc{D}}(f)=1_U\}.
	\end{align*}
\end{definition}
\begin{definition}\label{Definition:Cleavage} Let $\mc{C}$ be a category and $(\mc{D},\pi_{\mc{D}},\mc{C})$ a fibered category. A \textit{cleavage of the fibered category} $(\mc{D},\pi_{\mc{D}},\mc{C})$ is, a collection $\mc{A}$, of Cartesian arrows in $\mc{D}$, such that, for each $(\eta,f)\in \mc{D}_0\times_{\mc{C}_0,t}\mc{C}_1$, there exists exactly one element $\theta\in \mc{A}$, with $t(\theta)=\eta$ and $\pi_{\mc{D}}(\theta)=f$.
\end{definition}
Let $\mc{C}$ be a category and $(\mc{D},\pi_{\mc{D}},\mc{C})$ a fibered category. Fix a cleavage $\mc{A}$ of $(\mc{D},\pi_{\mc{D}},\mc{C})$. 
Let $f\colon U\ra V$ be a morphism in $\mc{C}$, and $\eta$ an object of $\mc{D}(V)$. Then, there exists unique element $\theta_{\eta}\in \mc{A}$ such that $t(\theta_{\eta})=\eta$
and $\pi_{\mc{D}}(\theta_\eta)=f$. Denote $s(\theta_{\eta})$ by $f^*\eta$. This gives a map $\text{Obj}(\mc{D}(V))\ra \text{Obj}(\mc{D}(U))$ defined as $\eta\mapsto f^*\eta$. Let $\tau\colon \eta\ra \eta'$ be a morphism in $\mc{D}(V)$. Then, by the definition of a fibered category, there exists a unique morphism $\Psi\colon f^*\eta\ra f^*\eta'$ in $\mc{D}(U)$ such that $\tau\circ \theta_{\eta}=\theta_{\eta'}\circ \Psi$. We denote the morphism $\Psi$ by $f^*\tau\colon f^*\eta\ra f^*\eta'$. This gives a map $\text{Mor}(\mc{D}(V))\ra \text{Mor}(\mc{D}(U))$, defined by $\tau\mapsto f^*\tau$. The maps $\text{Obj}(\mc{D}(V))\ra \text{Obj}(\mc{D}(U))$ and 
$\text{Mor}(\mc{D}(V))\ra \text{Mor}(\mc{D}(U))$ combine to give a functor $ \mc{D}(V)\ra \mc{D}(U)$, which we denote by $f^*\colon \mc{D}(V)\ra \mc{D}(U)$.

Let $\mc{C}$ be a category and $(\mc{D},\pi_{\mc{D}},\mc{C})$ a fibered category. Fixing a cleavage $\mc{A}$ on the fibered category $(\mc{D},\pi_{\mc{D}},\mc{C})$ we get the collection of functors $\{f^*\colon \mc{D}(V)\ra \mc{D}(U)\}_{f\in \text{Mor}(\mc{C})}$. This collection $\{f^*\colon \mc{D}(V)\ra \mc{D}(U)\}_{f\in \text{Mor}(\mc{C})}$ gives a pseudo-functor on the category $\mc{C}$ (\Cref{Definition:pseudofunctor}).
Similarly, given a pseudo-functor on a category $\mc{C}$, we can associate a fibered category over $\mc{C}$ (along with a cleavage of the fibered category). We refer to Chapter $3$ of \cite{MR2223406} for a detailed correspondence between pseudo-functors and fibered categories. 

\begin{example}\label{Example:PseudofunctorFiberedcategory}
	Let $\mc{C}$ be a category. Let $U$ be an object of $\mc{C}$. Consider the functor (of points) $h_U:\mc{C}^{\op}\ra \text{Set}$. Treating a set as a discrete category, we see $h_U$ as a pseudo-functor $h_U:\mc{C}^{\op}\ra \text{Cat}$. Then, the associated fibered category for this pseudo-functor is the fibered category $\mc{C}/U\ra \mc{C}$ mentioned in \Cref{Example:fiberedcategoryCU}. 
\end{example}



The correspondence between pseudo-functors and fibered categories allows us to see a stack over a site $(\mc{C},\mc{J})$ as a fibered category over $\mc{C}$. The notion of descent category of a covering $\{U_\alpha\ra U\}$ of an object $U$ of $\mc{C}$, for a pseudo-functor is transferred to the set up of fibered categories to give the notion of descent category of the covering $\{U_\alpha\ra U\}$ for a fibered category over $\mc{C}$. 

\subsubsection{Descent Category associated to a covering $\{U_i\ra U\}$}
Let $(\mathcal{C},\mc{J})$ be a site, and $(\mc{D},\pi_{\mc{D}},\mc{C})$ a fibered category.
Fix a cleavage $\mc{A}$ of the fibered category $(\mc{D},\pi_{\mc{D}},\mc{C})$. Let $U$ be an object of $\mathcal{C}$ and $\{\sigma_i\colon U_i\rightarrow U\}$ a covering of $U$. 
For this covering $\{\sigma_i\colon U_i\rightarrow U\}$, we associate a category $\mathcal{D}(\{U_i\rightarrow U\})$. For indices $i,j$, let $U_{ij}$ denote the pull-back in the Diagram \ref{Diagram:2pull-back}, and for indices $i,j,k$, let $U_{ijk}$ denote the pull-back in the Diagram \ref{Diagram:3pull-back}.

An object in the category $\mathcal{D}(\{U_i\rightarrow U\})$ is of the form $(\{\xi_i\},\{\phi_{ij}\})$, where,
\begin{enumerate}
	\item $\xi_i$ is an object of $\mathcal{D}(U_i)$ for each index $i\in\Lambda$,
	\item $\phi_{ij}\colon {\pr}_2^*(\xi_j)\rightarrow {\pr}_1^*(\xi_i)$ is an isomorphism in $\mathcal{D}(U_{ij})$ for each pair of indices $i,j\in\Lambda$ such that $ {\pr}_{13}^*(\phi_{ik})={\pr}_{12}^*(\phi_{ij})\circ {\pr}_{23}^*(\phi_{jk})$ in $\mathcal{D}(U_{ijk})$ for three indices $i,j,k\in \Lambda$.
\end{enumerate}

A morphism from $(\{\xi_i\},\{\phi_{ij}\}) $ to $(\{\eta_i\},\{\psi_{ij}\})$ is given by a collection $\{\alpha_i\colon \xi_i\rightarrow \eta_i\}_{i\in \Lambda}$. Here, $\alpha_i:\xi_i\ra \eta_i$ is a morphism in $\mc{D}(U_i)$ for each $i\in \Lambda$ satisfying the condition $\psi_{ij}\circ ({\pr}_2^*\alpha_j)=({\pr}_1^*\alpha_i)\circ \phi_{ij}$ in the category $\mc{D}(U_{ij})$ for each pair of indices $i,j\in \Lambda$.
We call $\mathcal{D}(\{U_i\rightarrow U\})$ to be the \textit{descent category associated to the covering $\{U_i\rightarrow U\}$}. 

Let $\xi$ be an object of $\mc{D}(U)$. For each index $i\in \Lambda$, the functor $\sigma_i^*\colon \mc{D}(U)\ra \mc{D}(U_i)$ (associated to the morphism $\sigma_i\colon U_i\ra U$) gives an object $\sigma_i^*(\xi)$ in $\mc{D}(U_i)$. For indices $i,j$, the uniqueness of pull-back gives an isomorphism $\phi_{ij}\colon {\pr}_2^*(\sigma_j^*(\xi))\ra {\pr}_1^*(\sigma_i^*(\xi))$ in $\mc{D}(U_{ij})$. It is easy to see that, for indices $i,j,k$ the morphisms $\phi_{ij},\phi_{jk}, \phi_{ik}$ satisfies the co-cycle condition mentioned above. Thus, we get an object $\{\{\sigma_i^*(\xi)\},\{\phi_{ij}\}\}$ in $\mc{D}(\{U_i\ra U\})$. This gives a map $\text{Obj}(\mc{D}(U))\ra \text{Obj}(\mc{D}(\{U_i\ra U\}))$. Similarly, the pull-back defines a map $\text{Mor}(\mc{D}(U))\ra \text{Mor}(\mc{D}(\{U_i\ra U\}))$. Thus, we have a functor $\mc{D}(U)\ra \mc{D}(\{U_i\ra U\})$.
\begin{definition}[Stack over a site] \label{Definition:stack}
	Let $(\mc{C},\mc{J})$ be a site.
	A fibered category $(\mc{D},\pi_{\mc{D}},\mc{C})$ is said to be \textit{a stack over the site $(\mc{C},\mc{J})$} if, for each object $U$ of $\mc{C}$ and a covering $\{U_i\ra U\}$ of $U$, the functor $\mc{D}(U)\rightarrow \mc{D}(\{U_i\ra U\})$ is an equivalence of categories.
\end{definition}

In this thesis, we are interested in fibered categories $(\mc{D},\pi_{\mc{D}},\mc{C})$ that are fibered in groupoids, that is, the fiber $\mc{D}(U)$ is a groupoid for each object $U$ of $\mc{C}$. We call them \textit{categories fibered in groupoids over $\mc{C}$}. By a stack, we always mean a stack fibered in groupoids.

\begin{definition}[morphism of stacks]\label{Definition:morphismofstacks}
	Let $(\mc{C},\mc{J})$ be a site, $(\mc{D},\pi_{\mc{D}},\mc{C})$ and $(\mc{D}',\pi_{\mc{D}'},\mc{C})$ are stacks. A \textit{morphism of stacks from $(\mc{D},\pi_{\mc{D}},\mc{C})$ to $(\mc{D}',\pi_{\mc{D}'},\mc{C})$} is a morphism of fibered categories from $(\mc{D},\pi_{\mc{D}},\mc{C})$ to $(\mc{D}',\pi_{\mc{D}'},\mc{C})$. A morphism of stacks $F:(\mc{D},\pi_{\mc{D}},\mc{C})\ra (\mc{D}',\pi_{\mc{D}'},\mc{C})$ is called an \textit{isomorphism of stacks} if it is an isomorphism of fibered categories over $\mc{C}$ (\Cref{Definition:isomorphismoffiberedcategories}).
\end{definition}
\begin{example}\label{Example:subcanonicalstack}
	Let $\mc{S}$ be a category, and $\mc{J}$ a subcanonical topology on $\mc{S}$. For an object $U$ of $\mc{S}$, the fibered category $\pi_U:\mc{S}/U\ra \mc{S}$ is a stack over the site $(\mc{S},\mc{J})$. 
\end{example}
\begin{definition}[representable stack] \label{Definition:representablestack}
	Let $(\mc{S},\mc{J})$ be a subcanonical site. A stack $(\mc{D},\pi_{\mc{D}},\mc{S})$ over the site $(\mc{S},\mc{J})$, is said to be a \textit{representable stack}, if, there exists an object $U$ of $\mc{S}$ such that 
	$(\mc{D},\pi_{\mc{D}},\mc{S})$ is isomorphic to $(\underline{U},\pi_{U},\mc{S})$ as stacks over $(\mc{S},\mc{J})$.
\end{definition}
We are interested in morphisms of stacks that are representable by objects of the site. To define the notion of representable morphisms we need to introduce the notion of $2$-fiber product of fibered categories (in particular, $2$-fiber product of stacks).
\begin{definition}[$2$-fiber product of fibered categories]\label{Definition:2fiberproduct}Let $\mc{S}$ be a category, $\pi_{\mc{D}}\colon \mc{D}\rightarrow \mc{S}$, $\pi_{\mc{D}'}\colon \mc{D}'\rightarrow \mc{S}$ and $\pi_{\mc{C}}\colon \mc{C}\rightarrow \mc{S}$
	are fibered categories. 
	Let
	$F\colon \mc{D}\rightarrow \mc{C}, F'\colon \mc{D}'\rightarrow \mc{C}$ 
	be a pair of morphisms of fibered categories.
	Consider the category $\mc{D}\times_{\mc{C}}\mc{D}'$ with 
	the following description:
	\begin{itemize}
		\item an object of $\mc{D}\times_{\mc{C}}\mc{D}'$ is a triple $(d,d',\alpha)$, where $d\in \text{Obj}(\mc{D}), d'\in \text{Obj}(\mc{D}')$ with $\pi_{\mc{D}}(d)=\pi_{\mc{D}'}(d')$ and $\alpha\in \text{Mor}(\mc{C})$ with $s(\alpha)=F(d), t(\alpha)=F'(d')$,
		\item a morphism of $\mc{D}\times_{\mc{C}}\mc{D}'$ from $(d_1,d_1',\alpha_1)$ to $(d_2,d_2',\alpha_2)$ is given by a pair $(u,v)$ where $u\colon d_1\ra d_2$ is in $\text{Mor}(\mc{D})$
		and $v\colon d_1'\ra d_2'$ is in $\text{Mor}(\mc{D}')$ such that 
		$\pi_{\mc{D}}(u)=\pi_{\mc{D}'}(v)$ and $F'(v)\circ \alpha_1=\alpha_2\circ F(u)$.
	\end{itemize}	
	Consider the functor $\Phi\colon \mc{D}\times_{\mc{C}}\mc{D}'\ra \mc{S}$ with the following description 
	\begin{itemize}
		\item an object $(d,d',\alpha)$ of $\mc{D}\times_{\mc{C}}\mc{D}'$ is mapped to the object $\pi_{\mc{D}}(d)=\pi_{\mc{D}'}(d')$ of $\mc{S}$,
		\item a morphism $(u,v)$ of $\mc{D}\times_{\mc{C}}\mc{D}'$ is mapped to the morphism $\pi_{\mc{D}}(u)=\pi_{\mc{D}'}(v)$ of $\mc{S}$.
	\end{itemize} Then, $\Phi\colon \mc{D}\times_{\mc{C}}\mc{D}'\ra \mc{S}$ is a fibered category over $\mc{S}$, which we call \textit{the $2$-fiber product of $\mc{D}$ and $\mc{D}'$ with respect to the morphisms
		$F\colon \mc{D}\rightarrow \mc{C}$ and $F'\colon \mc{D}'\rightarrow \mc{C}$}. We see the $2$-fiber product as the following diagram,
	\[
	\begin{tikzcd}
		\mc{D}\times_{\mc{C}}\mc{D}' \arrow[rr, "{\pr}_2"] \arrow[dd, "{\pr}_1"'] & & \mc{D}' \arrow[dd, "F'"] \arrow[rddd, "\pi_{\mc{D}'}"] & \\
		& &       & \\
		\mc{D} \arrow[rr, "F"] \arrow[rrrd, "\pi_{\mc{D}}"']   & & \mc{C} \arrow[rd, "\pi_{\mc{C}}"]   & \\
		& &       & \mc{S}
	\end{tikzcd}.\]
\end{definition}

\begin{remark}\label{Remark:2-fiberproduct}
	If $(\mc{D},\pi_{\mc{D}},\mc{S}), (\mc{D}',\pi_{\mc{D}'},\mc{S})$ and $(\mc{C},\pi_{\mc{C}},\mc{S})$ are categories fibered in groupoids over $\mc{S}$, then $\mc{D}\times _{\mc{C}}\mc{D}'\ra \mc{S}$ is a category fibered in groupoids over $\mc{S}$. If $\mc{J}$ is a Grothendieck topology on $\mc{S}$ such that $(\mc{D},\pi_{\mc{D}},\mc{S}), (\mc{D}',\pi_{\mc{D}'},\mc{S})$ and $(\mc{C},\pi_{\mc{C}},\mc{S})$ are stacks over the site $(\mc{S},\mc{J})$, then $\mc{D}\times _{\mc{C}}\mc{D}'\ra \mc{S}$ is a stack over the site $(\mc{S},\mc{J})$.
\end{remark}

\begin{definition}[representable morphism of stacks]
	\label{Definition:representablemorphism}
	Let $(\mc{S},\mc{J})$ be a subcanonical site. Let $(\mc{D},\pi_{\mc{D}},\mc{S})$ and $(\mc{C},\pi_{\mc{C}},\mc{S})$ be stacks. A morphism of stacks $F:\mc{D}\ra \mc{C}$ is said to be a \textit{representable morphism of stacks }if for every object $U$ of $\mc{S}$ and a morphism of stacks $\underline{U}\ra \mc{C}$, the $2$-fiber product $\mc{D}\times_{\mc{C}}\underline{U}$ is a representable stack (\Cref{Definition:representablestack}).
\end{definition}
We are also interested in another special type of morphism of stacks, namely an epimorphism of stacks. 
\begin{definition}[epimorphism of stacks \cite{MR4139032}]\label{Definition:epimorphismofstacks}
	Let $(\mc{S},\mc{J})$ be a subcanonical site, and $(\mc{D},\pi_{\mc{D}},\mc{S}), (\mc{C},\pi_{\mc{C}},\mc{S})$ are stacks. A morphism of stacks $F:\mc{D}\ra \mc{C}$ is said to be an \textit{epimorphism of stacks} if, for every object $U$ of $\mc{S}$ and a morphism of stacks $q:\underline{U}\ra \mc{C}$, there exists a covering $\{\sigma_\alpha:U_\alpha\ra U\}_{\alpha\in \Lambda}$ of $U$ and a morphism of stacks $p_\alpha:\underline{U_\alpha}\ra \mc{D}$ with a $2$-commutative diagram,
	\[
	\begin{tikzcd}
		\underline{U_\alpha} \arrow[rr, "\sigma_\alpha"] \arrow[dd, "p_\alpha"'] & & \underline{U} \arrow[dd, "q"] \\
		& &    \\
		\mc{D} \arrow[rr, "F" ] \arrow[Rightarrow, shorten >=20pt, shorten <=20pt, uurr, "\alpha"]       & & \mc{C}   
	\end{tikzcd}\] 
	for each $\alpha\in \Lambda$.
\end{definition}
\subsection{$2$-Yoneda lemma}\label{Subsection:2Yoneda} 

Let $\mc{C}$ be a category. For a functor $\mc{F}:\mc{C}^{\op}\ra \text{Set}$, the Yoneda lemma gives a bijection $\text{Hom} (h_U,\mc{F})\ra \mc{F}(U)$ for objects $U$ of $\mc{C}$. 

 Let $\pi_{\mc{F}}:\mc{F}\ra \mc{C}$ be the fibered category associated to the functor $\mc{F}:\mc{C}^{\op}\ra \text{Set}$. As mentioned before, $\pi_U:\mc{C}/U\ra \mc{C}$ is the fibered category associated to the functor $h_U:\mc{C}^{\op}\ra \text{Set}$. Then, the Yoneda lemma says that the map $\text{Hom}(\pi_U,\pi_{\mc{F}})\ra \mc{F}(U)$ is a bijection, where $\mc{F}(U)$ is fiber of the object $U$ in the fibered category $(\mc{F},\pi_{\mc{F}},\mc{C})$. The $2$-Yoneda lemma is a higher version of the Yoneda lemma where the pseudo-functors do not necessarily takes values in Sets.

Let $\mc{F}:\mc{C}^{\op}\ra \text{Cat}$ be a pseudo-functor, and $\mc{F}\ra \mc{C}$ the associated fibered category. Let $\text{Hom}(\mc{C}/U,\mc{F})$ be the category whose objects are morphisms of fibered categories over $\mc{C}$ from $\mc{C}/U$ to $\mc{F}$ and morphisms are natural transformations. Consider the functor $\text{Hom}(\mc{C}/U,\mc{F})\ra \mc{F}(U)$ with the following description
\begin{itemize}
	\item an object $\Phi$ of $\text{Hom}(\mc{C}/U,\mc{F})$ is mapped to the object $\Phi(U)(1_U)$,
	\item a morphism $\alpha$ of $\text{Hom}(\mc{C}/U,\mc{F})$ is mapped to the morphism $\alpha(1_U)$.
\end{itemize}

\begin{lemma}[$2$-Yoneda lemma {\cite[Lemma $4.19$]{MR2778793}}]\label{Lemma:2-yoneda}
	The functor $\text{Hom}(\mc{C}/U,\mc{F})\ra \mc{F}(U)$, mentioned above, is an equivalence of categories. 
\end{lemma}
The following result is an application of the $2$-Yoneda lemma, giving a characterization for epimorphisms of stacks. We will see the proof of this result in \Cref{Chap.2}, as part of the \Cref{Lemma:equivalentnotionfordifferentiablestack}.

\begin{lemma} Let $\mc{D}\ra \text{Man}$ and $\mc{C}\ra \text{Man}$ be stacks. 
	A morphism of stacks $F\colon \mc{D}\rightarrow \mc{C}$ is an epimorphism of stacks if and only if the following conditions holds:
	for each smooth manifold $U$, and an object $a$ of $ \mc{C}(U)$, there exists an open cover $\{U_i\rightarrow U\}_{i\in \Lambda}$ of $U$ and objects $x_i$ of $\mc{D}(U_i)$ with an isomorphism $F(x_i)\rightarrow a|_{U_i}$ in $\mc{C}(U_i)$ for each $i\in \Lambda$. 
\end{lemma}


\section{Lie groupoids, Morita equivalence, and differentiable stacks}\label{Section:LiegroupoidsMoritaequivalencedifferentiablestacks}
In this thesis, we are interested in stacks of a special kind, namely differentiable stacks. A characterization of differentiable stacks is that, it is ``representable by a Lie groupoid''. In this section, we recall the notion of a Lie groupoid, Morita equivalence of Lie groupoids, and differentiable stacks. A major part of this section is based on \cite{MR3089760, MR2778793, MR2012261}.

\begin{definition}\label{Definition:Liegroupoid}
	A \textit{Lie groupoid} is a groupoid $\mc{G}=[\mc{G}_1\rra \mc{G}_0]$, where $\mc{G}_0,\mc{G}_1$ are smooth manifolds, the structure maps $s,t:\mc{G}_1\ra \mc{G}_0$ are submersions, the other structure maps $m:\mc{G}_1\times_{\mc{G}_0}\mc{G}_1\ra \mc{G}_1, e:\mc{G}_0\ra \mc{G}_1, u:\mc{G}_1\ra \mc{G}_1$ are smooth maps. 
\end{definition}
\begin{example}
	Let $M$ be a smooth manifold. The groupoid $[M\rra M]$, associated to the smooth manifold $M$, mentioned in \Cref{Example:MMcategory}, is a Lie groupoid.
\end{example}
\begin{example}
	Let $G$ be a Lie group. The groupoid $[G\rra *]$, associated to the Lie group $G$, mentioned in \Cref{Example:G*category} is a Lie groupoid.
\end{example}
\begin{example}
	Let $G$ be a Lie group, $M$ be a smooth manifold.
	Let $\mu: M\times G\ra M$ be an action of the Lie group $G$ on the smooth manifold $M$. The groupoid $[M\times G\rra M]$, associated to a Lie group action, mentioned in \Cref{Example:actioncategory} is a Lie groupoid.
\end{example}

\begin{example}
	Let $G$ be a Lie group, $M$ a smooth manifold, and $P\ra M$ a principal $G$-bundle over $M$. Consider the category $[(P\times P)/G\rra M]$ with the following description:
	\begin{itemize}
		\item objects are the points of the smooth manifold $M$,
		\item morphisms are the points of the quotient space $(P\times P)/G$,
		\item the source map is given by $[(x_1,x_2)]\mapsto \pi(x_1)$, for all $[(x_1,x_2)]\in (P\times P)/G$
		\item the target map is given by 
		$[(x_1,x_2)]\mapsto \pi(x_2)$, for all $(x_1,x_2)\in (P\times P)/G$,
		\item the composition map is given by $([(x_1,x_2)],[(x_3,x_4)])\mapsto [(x_1,x_4)]$, for all 
		$[(x_1,x_2)]$, and $[(x_3,x_4)]$ in $(P\times P)/G$ with $\pi(x_2)=\pi(x_3)$.
	\end{itemize}
	We call $[(P\times P)/G\rra M]$ to be the \textit{Gauge groupoid of the principal bundle $\pi:P\ra M$}. 
\end{example}
\begin{example}
	Let $M$ be a smooth manifold, and $E\ra M$ a vector bundle over the smooth manifold $M$. Consider the category $[{\rm GL}(E)\rra M]$, with the following description:
	\begin{itemize}
		\item objects are the points of the smooth manifold $M$,
		\item a morphism from an element $x$ to an element $y$ of $M$ is given by a linear isomorphism $E_x\ra E_y$,
		\item the source map is given by $(\gamma:E_x\ra E_y)\mapsto x$,
		\item the target map is given by $(\gamma:E_x\ra E_y)\mapsto y$,
		\item the composition map is given by the usual composition of functions. 
	\end{itemize}
	We call $[{\rm GL}(E)\rra M]$ to be the \textit{general linear groupoid of the vector bundle $E\ra M$}.	
\end{example}

\begin{definition}[morphism of Lie groupoids {\cite[Section $5.1$]{MR2012261}}]
	Let $\mc{G}=[\mc{G}_1\rra \mc{G}_0]$ and $\mc{H}=[\mc{H}_1\rra \mc{H}_0]$ be Lie groupoids. A \textit{morphism of Lie groupoids from $\mc{G}=[\mc{G}_1\rra \mc{G}_0]$ to $\mc{H}=[\mc{H}_1\rra \mc{H}_0]$} is a functor $(F_1,F_0):[\mc{G}_1\rra \mc{G}_0]\ra [\mc{H}_1\rra \mc{H}_0]$ such that $F_0:\mc{G}_0\ra \mc{H}_0$ and $F_1:\mc{G}_1\ra \mc{H}_1$ are smooth maps.
\end{definition}
We will express a morphism of Lie groupoids by the following diagram, 
\begin{equation}\label{Diagram:MorphismofLiegroupoids} \begin{tikzcd}
		\mc{G}_1 \arrow[dd,xshift=0.75ex,"t_{\mc{G}}"]\arrow[dd,xshift=-0.75ex,"s_{\mc{G}}"'] \arrow[rr, "F_1"] & & \mc{H}_1 \arrow[dd,xshift=0.75ex,"t_{\mc{H}}"]\arrow[dd,xshift=-0.75ex,"s_{\mc{H}}"'] \\
		& & \\
		\mc{G}_0 \arrow[rr, "F_0"] & & \mc{H}_0
	\end{tikzcd}.
\end{equation}
\begin{example}
	Let $G$ and $H$ be Lie groups acting on smooth manifolds $M$ and $N$ respectively. Let $\psi:G\ra H$ be a morphism of Lie groups, and $f:M\ra N$ a smooth map preserving the action of Lie groups; that is, $f(mg)=f(m)\psi(g)$ for all $m\in M, g\in G$. Then, $((f\times \psi),f):[M\times G\rra M]\ra [N\times H\rra N]$ is a morphism of Lie groupoids.
\end{example}
\begin{example}
	Let $G$ be a Lie group, $P\ra M, P'\ra M'$ are principal $G$-bundles. A morphism of principal bundles $(F,f):P(M,G)\ra P'(M',G)$ induce a morphism of Lie groupoids $((F,F)^{/G},f):[(P\times P)/G\rra M]\ra [(P'\times P')/G\rra M']$.
\end{example}
\begin{example}
	Let $M$ be a smooth manifold. Consider vector bundles $E\ra M$, and $E'\ra M$. A morphism of vector bundles $(F,f):(E\ra M)\ra (E'\ra M')$ induce a morphism of Lie groupoids $[{\rm GL}(E)\rra M]\ra [{\rm GL}(E')\rra M']$.
\end{example}

For each smooth manifold $M$, using \Cref{Example:PseudofunctorFiberedcategory} we can associate a stack $\underline{M}:\text{Man}\ra {\rm Gpd}$. Also for a Lie group $G$, we have associated a stack $BG:{\rm Man}\ra \text{Gpd}$. As smooth manifolds, Lie groups are special cases of Lie groupoids; we expect a similar construction of stack for Lie groupoids.

For a Lie group $G$, the stack $BG$ assigns for each manifold $M$, a category $BG(M)$, of principal $G$-bundles over $M$. Inspired by the notion of a principal $G$-bundle over a smooth manifold, we introduce the notion of a principal $\mc{G}$-bundle over a smooth manifold for a Lie groupoid $\mc{G}$. For that, firstly, we introduce the notion of action of a Lie groupoid on a smooth manifold. 

\begin{definition}[action of a Lie groupoid on a Manifold]\label{actionOfLiegroupoid} Let $\mc{G}=[\mc{G}_1\rra \mc{G}_0]$ be a Lie groupoid and $P$ a smooth manifold. 
	
	A \textit{right action of $\mc{G}$ on $P$} consists of a pair of smooth maps 
	$(a:P\ra \mc{G}_0,\mu:P\times_{a,\mc{G}_0,t}\mc{G}_1\ra P)$ satisfying the following conditions:
	\begin{enumerate}
		\item $\mu(p,1_{a(p)})=p$ for all $p\in P$,
		\item $a(\mu(p,\gamma))=s(\gamma)$ for all $(p,\gamma)\in P\times_{a,\mc{G}_0,t}\mc{G}_1$,
		and
		\item $\mu(\mu(p,\gamma),\gamma')=\mu(p,\gamma\circ \gamma')$ for all $(p,\gamma,\gamma')\in P\times_{a,\mc{G}_0,t}\mc{G}_1\times_{s,\mc{G}_0,t}\mc{G}_1$.
	\end{enumerate}
	
	A \textit{left action of $\mc{G}$ on $P$} consists of a pair of smooth maps $(a\colon P\ra \mc{G}_0, \mu\colon \mc{G}_1\times_{s,\mc{G}_0,a}P\ra P)$ satisfying the following conditions:
	\begin{enumerate}
		\item $\mu(1_{a(p)},p)=p$ for all $p\in P$,
		\item $a(\mu(\gamma,p))=t(\gamma)$ for all $(\gamma,p)\in \mc{G}_1\times_{s,\mc{G}_0,a}P$, and
		\item $\mu(\gamma,\mu(\gamma',p))=\mu(\gamma\circ \gamma',p)$ for all $(\gamma,\gamma',p)\in \mc{G}_1\times_{s,\mc{G}_0,t}\mc{G}_1\times_{s,\mc{G}_0,a}P$.
	\end{enumerate}
\end{definition}

\begin{example}
	The action of a Lie groupoid $[M\rra M]$ on a smooth manifold $N$ consists of a smooth map $f: N\ra M$.
\end{example}

\begin{example}
	An action of a Lie groupoid $[G\rra *]$ on a smooth manifold $M$ is given by an action $\mu:M\times G\ra M$ of the Lie group $G$ on the manifold $M$.
\end{example}	
\begin{example}
	Let $\mc{G}=[\mc{G}_1\rra \mc{G}_0]$ be a Lie groupoid. Consider the composition map
	\[m:\mc{G}_1\times_{t,\mc{G}_0,s} \mc{G}_1\ra \mc{G}_1.\] 
	The pair $(t:\mc{G}_1\ra \mc{G}_0, m:\mc{G}_1\times_{t,\mc{G}_0,s} \mc{G}_1\ra \mc{G}_1)$ gives a left action of $[\mc{G}_1\rra \mc{G}_0]$ on $\mc{G}_1$. The pair $(s:\mc{G}_1\ra \mc{G}_0, m:\mc{G}_1\times_{t,\mc{G}_0,s} \mc{G}_1\ra \mc{G}_1)$ gives a right action of $[\mc{G}_1\rra \mc{G}_0]$ on $\mc{G}_1$.
\end{example}

\begin{definition}[{\cite[Definition $3.17$]{MR2778793}}]
	Let $\mc{G}$ be a Lie groupoid and $M$ a smooth manifold. A \textit{(right) principal $\mc{G}$-bundle over $M$} consists of 
	\begin{enumerate}
		\item a smooth manifold $P$,
		\item a surjective submersion $\pi:P\ra M$,  
		\item an action of $\mc{G}$ on $P$ given by a pair of maps $(a_{\mc{G}}:P\ra \mc{G}_0, \mu:P\times_{a_{\mc{G}},\mc{G}_0,t}\mc{G}_1\ra P)$, 
	\end{enumerate}
	satisfying the following conditions:
	\begin{enumerate}
		\item the map $\pi:P\ra M$ is such that, $\pi(p\cdot \gamma)=\pi(p)$ for all $(p,\gamma)\in P\times_{a_{\mc{G}},\mc{G}_0,t}\mc{G}_1$,
		\item the map $P\times_{a_{\mc{G}},\mc{G}_0,t}\mc{G}_1\ra P\times_MP$ given by $(p,\gamma)\mapsto (p\cdot\gamma,p)$ is a diffeomorphism.
	\end{enumerate}
	A \textit{left principal $\mc{G}$-bundle over $M$} is defined similarly using the left action of $\mc{G}$.
\end{definition}

We see a (right) principal $\mc{G}$-bundle over $M$ as the following diagram,
\begin{equation}\label{Diagram:principalLiegroupoidbundle}
\begin{tikzcd}
	&      & \mathcal{G}_1 \arrow[dd,xshift=0.75ex,"t"]\arrow[dd,xshift=-0.75ex,"s"'] \\
	& P \arrow[ld, "\pi"'] \arrow[rd, "a_{\mc{G}}"] &    \\
	M &      & \mathcal{G}_0  
\end{tikzcd}. 
\end{equation}

We see a (left) principal $\mc{G}$-bundle over $M$ as the following diagram,\[
\begin{tikzcd}
	\mc{G}_1 \arrow[dd,xshift=0.75ex,"t"]\arrow[dd,xshift=-0.75ex,"s"'] &      & \\
	& P \arrow[ld, "a_{\mc{G}}"'] \arrow[rd, "\pi"] & \\
	\mc{G}_0  &      & M
\end{tikzcd}.\]

%
%

\begin{example}\label{Example:targetmapisLGbundle}
	Let $\mc{G}=[\mc{G}_1\rra \mc{G}_0]$ be a Lie groupoid. The target map $t:\mc{G}_1\ra \mc{G}_0$ is considered as a right principal $\mc{G}$-bundle, with the action of $[\mc{G}_1\rra \mc{G}_0]$ on $\mc{G}_1$ is given by the pair of maps $(s:\mc{G}_1\ra \mc{G}_0, m:\mc{G}_1\times_{s,\mc{G}_0,t} \mc{G}_1\ra \mc{G}_1)$. 
We see this as the following diagram,
\begin{equation}\label{Diagram:targetmapisLGbundle}
\begin{tikzcd}
	&      & \mathcal{G}_1 \arrow[dd,xshift=0.75ex,"t"]\arrow[dd,xshift=-0.75ex,"s"'] \\
	& \mc{G}_1 \arrow[ld, "t"'] \arrow[rd, "s"] &    \\
	\mc{G}_0 &      & \mathcal{G}_0  
\end{tikzcd}.
\end{equation} 
The source map $s:\mc{G}_1\ra \mc{G}_0$ is considered as a left principal $\mc{G}$-bundle, with the action of $[\mc{G}_1\rra \mc{G}_0]$ on $\mc{G}_1$ is given by $(t:\mc{G}_1\ra \mc{G}_0, m:\mc{G}_1\times_{s,\mc{G}_0,t} \mc{G}_1\ra \mc{G}_1)$. 
\end{example}

\begin{example}\label{Example:P(M,G)isLGbundle}
Let $M$ be a smooth manifold, $G$ a Lie group, and $\pi:P\ra M$ a principal $G$ bundle. We see $\pi:P\ra M$ as a principal $[G\rra *]$-bundle over $M$ in the following diagram,
\begin{equation}\label{Diagram:P(M,G)isLGbundle}
	\begin{tikzcd}
		&      & G \arrow[dd,xshift=0.75ex,"t"]\arrow[dd,xshift=-0.75ex,"s"'] \\
		& P \arrow[ld, "\pi"'] \arrow[rd] &    \\
		M &      & * 
	\end{tikzcd}.
\end{equation} 
\end{example}

\begin{definition}
	Let $\mc{G}$ be a Lie groupoid, $P(M,\mc{G})$ and $P'(M',\mc{G})$ are principal $\mc{G}$-bundles. A \textit{morphism of principal $\mc{G}$-bundles} from $P(M,\mc{G})$ to $P'(M',\mc{G})$ consists of a pair of smooth maps $(F:P\ra P', f:M\ra M')$ such that, $\pi'\circ F=f\circ \pi$, $a_{\mc{G}}'\circ F=a_\mc{G}$, and $\mu'\circ (F,1)=F\circ \mu$. 
\end{definition}

We see a morphism of principal $\mc{G}$-bundles as the following commutative diagram,
\[
\begin{tikzcd}
	& P\times \mathcal{G}_1 \arrow[dd, "\mu"'] \arrow[rr, "{(F,1)}"] &  & P'\times \mathcal{G}_1 \arrow[dd, "\mu'"]  & \\
	&         &  &       & \\
	& P \arrow[rd, "a_{\mathcal{G}}"'] \arrow[rr, "F"] \arrow[ld, "\pi"'] &  & P' \arrow[ld, "a_{\mathcal{G}}'"] \arrow[rd, "\pi'"] & \\
	M \arrow[rrrr, "f"', bend right] &         & \mathcal{G}_0 &       & M'
\end{tikzcd}.\]

Similar to the notion of pull-back of a geometric structure on a manifold $M$ along a smooth map $f:N\ra M$, we have the notion of pull-back of a $\mc{G}$-bundle over a manifold $M$ along a smooth map $f:N\ra M$.
\subsection{Pullback of a principal $\mc{G}$-bundle}
\label{Subsection:PullbackofPrincipalLiegroupoidbundle}
Let $\mc{G}=[\mc{G}_1\rra \mc{G}_0]$ be a Lie groupoid and $\pi\colon P\rightarrow M$ a principal $\mc{G}$-bundle. Let $f\colon N\rightarrow M$ be a smooth map. As $\pi\colon P\rightarrow M$ is a submersion, the 
pull-back $f^*P=N\times_MP=\{(n,p)\colon f(n)=\pi(p)\}$ is a manifold (an embedded submanifold of $N\times P$).

We will express the pull-back by the following diagram,
\begin{equation}\label{Diagram:PullbackofPrincipalLiegroupoidbundle}
	\begin{tikzcd}
		N\times_MP \arrow[dd,"{\pr}_1"'] \arrow[rr,"{\pr}_2"] & & P \arrow[dd, "\pi"] \\
		& & \\
		N \arrow[rr, "f"] & & M
	\end{tikzcd}.\end{equation} 
Adjoining the above pull-back Diagram \ref{Diagram:PullbackofPrincipalLiegroupoidbundle} to the principal bundle Diagram \ref{Diagram:principalLiegroupoidbundle} gives the following diagram,
\begin{equation}\begin{tikzcd}
		& N\times_MP \arrow[ld,"{\pr}_1"'] \arrow[rd, "{\pr}_2"] & & \mc{G}_1 \arrow[dd,xshift=0.75ex,"t"]\arrow[dd,xshift=-0.75ex,"s"'] \\
		N \arrow[rd, "f"] & & P \arrow[ld, "\pi"'] \arrow[rd, "a"] & \\
		& M & & \mc{G}_0
	\end{tikzcd}.\end{equation} 
Consider the map $\mu\colon (N\times_M P)\times \mc{G}_1\ra N\times_MP$, defined as  $\big((n,p),g\big)\mapsto (n,pg)$. 
The pair of maps 
\[(a\circ {\pr}_2\colon N\times_BP\ra \mc{G}_0, \mu\colon (N\times_M P)\times \mc{G}_1\ra N\times_MP),\]
gives a (right) action of the Lie groupoid $\mc{G}$ on the manifold $N\times_MP$. The projection map ${\pr}_1:N\times_MP\ra N$ will then be a principal $\mc{G}$-bundle over the manifold $N$.  We call $(N\times_BP, {\pr}_1,N)$ to be  \textit{the pull-back of the principal $\mc{G}$-bundle $\pi\colon P\ra M$ along $f\colon N\ra M$}.

Similar to the example of a fibered category coming from a Lie group, we have a fibered category 
(\Cref{Definition:fiberedcategory}) coming from a Lie groupoid $\mc{G}$.

\subsection{A stack associated to a Lie groupoid}\label{Section:associatingastackforaLiegroupoid}
Let $\mc{G}$ be a Lie groupoid. Let $B\mc{G}$ be the category, whose objects are principal $\mc{G}$-bundles and morphisms are morphisms of principal $\mc{G}$-bundles. Consider the functor $\pi_{\mc{G}}:B\mc{G}\ra\text{Man}$ with the following description:
\begin{itemize}
	\item an object $P(M,\mc{G})$ of the category $B\mc{G}$ is mapped to the object $M$ of the category $\text{Man}$, 
	\item a morphism $(F,f):P(M,\mc{G})\ra P'(M',\mc{G})$ of the category $B\mc{G}$ is mapped to the morphism $f:M\ra M'$ of the category $\text{Man}$.
\end{itemize}
The functor $\pi_{\mc{G}}:B\mc{G}\ra \text{Man}$ is a fibered category. Consider the big \'etale topology on the category $\text{Man}$ 
(\Cref{Example:big\'etaletopology}). Then, the functor $\pi_{\mc{G}}:B\mc{G}\ra \text{Man}$ becomes a stack (\Cref{Definition:stack}). We call the stack $\pi_{\mc{G}}:B\mc{G}\ra \text{Man}$ to be the \textit{classifying stack of the Lie groupoid $\mc{G}$}.

%

The stack $B\mc{G}$ comes with an extra property. We have a morphism of stacks $\underline{\mc{G}_0}\ra B\mc{G}$ defined by the pull-back of the trivial $\mc{G}$-bundle $t:\mc{G}_1\ra \mc{G}_0$ along a smooth map $M\ra \mc{G}_0$ for a smooth manifold $M$. This morphism is a representable morphism of stacks (\cite[Example $4.24$]{MR2778793}). This map is actually a ``surjective submersion''.
To make sense of a (representable) morphism of stacks being surjective submersion, we need to recall an outcome of the Yoneda lemma. 

Let $\mc{C}$ be a category. 
Let $A$ be an object of $\mc{C}$. 
Then, we have the functor (of points) $h_A:\mc{C}^{\op}\rightarrow \text{Set}$.
This gives a functor $\mc{C}\ra [C^{\op},\text{Set}]$. 
The Yoneda lemma says that this functor is an embedding of categories. 
Consider a functor $\mc{F}:\mc{C}^{\op}\ra \text{Set}$.
Seeing a set as a discrete category, we can see the functor $\mc{F}:\mc{C}^{\op}\ra \text{Set}$ as a pseudo-functor. 
As mentioned before, for every pseudo-functor on $\mc{C}$, we can associate a fibered category over $\mc{C}$. 
The same would be true for the category $\text{Man}$.
Thus, we have an embedding of categories $\text{Man}\rightarrow \text{Fib}_{\text{Man}}$. 
For an object $M$ of $\text{Man}$, the fibered category associated to the functor $h_M:\text{Man}^{\op}\ra \text{Set}$ is given by $\underline{M}\ra \text{Man}$. As mentioned before, these fibered categories $\underline{M}\ra \text{Man}$ are stacks. So, we have an embedding of categories $\text{Man}\rightarrow \text{Stacks}_{\text{Man}}$. In particular, for smooth manifolds $M, N$, we have a bijection $\text{Hom}_{\text{Man}}(M,N)\rightarrow \text{Hom}_{\text{Stacks}}(\underline{M},\underline{N})$. We have the following lemma.
\begin{lemma}[{\cite[Corollary $4.16$]{MR2778793}}]
	The functor $\text{Man}\rightarrow \text{Stacks}_{\text{Man}}$ defined above is an embedding of categories. 
\end{lemma}
Thus, for smooth manifolds $M, N$, any morphism of stacks $\underline{M}\rightarrow \underline{N}$ comes from a unique smooth map $M\rightarrow N$. Further, if a stack $\mc{D}\ra \text{Man}$ is representable by a smooth manifold, then there  exists, up to a diffeomorphism, exactly one smooth manifold representing the stack $\mc{D}$.  

\begin{definition}
	Let $\mc{D}\ra\text{Man}$ and $\mc{C}\ra \text{Man}$ be stacks. 
	A representable morphism of stacks $\mc{D}\rightarrow \mc{C}$ is said to be a \textit{surjective submersion} if for each smooth manifold $M$ and a morphism of stacks $\underline{M}\rightarrow\mc{C}$, the projection morphism ${\pr}_2:\mc{D}\times_{\mc{C}}\underline{M}\rightarrow \underline{M}$ induce a surjective submersion of the smooth manifolds $X\ra M$, where $\underline{X}=\mc{D}\times_{\mc{C}}\underline{M}$.
\end{definition}
With this definition of a representable morphism being a surjective submersion, the morphism of stacks $\underline{\mc{G}_0}\ra B\mc{G}$ is a representable surjective submersion (\cite[Example $4.24$]{MR2778793}). 
The existence of a representable surjective submersion $M\ra \mc{D}$ for a given stack $\mc{D}\ra \text{Man}$ is a very special property, and there is a specific name to that.

\begin{definition}[an atlas for the stack {\cite[Definition $4.27$]{MR2778793}}]\label{Definition:Atlasofastack}
	Let $\mc{D}\rightarrow \text{Man}$ be a stack. An \textit{atlas for the stack $\mc{D}\rightarrow {\rm Man}$} consists of a smooth manifold $X$ and a representable morphism of stacks $\underline{X}\ra \mc{D}$ which is a surjective submersion.
\end{definition}

\begin{definition}[differentiable stack {\cite[Definition $4.27$]{MR2778793}}]\label{Definition:differentiablestack}
	A stack $\mc{D}\ra \text{Man}$ is said to be a \textit{differentiable stack} if there exists an atlas for the stack $\mc{D}$.
\end{definition}

\begin{remark}\label{Remark:compositionisatlas}
	Given an atlas, $X\ra \mc{D}$ and a surjective submersion $Y\ra X$, the composition $Y\ra \mc{D}$ would also be an atlas for the stack $\mc{D}$.
\end{remark}

\subsection{A Lie groupoid associated to a differentiable stack}
\label{Subsection:Liegroupoidassociatedtoatlas}
In the last section, we have introduced the notion of a differentiable stack as a stack, which shares a common property with the $B\mc{G}$, namely both admits an atlas. In this section, we associate a Lie groupoid for a differentiable stack. More details can be found at \cite[Proposition $4.31$]{MR2778793}.

Let $\pi_{\mc{D}}\colon \mc{D}\rightarrow \text{Man}$ be a differentiable stack, and $r\colon \underline{X}\rightarrow \mc{D}$ an atlas 
for the stack $\mc{D}$. 
Let $\underline{X}\times_{\mc{D}}\underline{X}$ be the $2$-fiber product (\Cref{Definition:2fiberproduct}) expressed in the following diagram,
\begin{equation}\label{Diagram:XXDLiegroupoid}
	\begin{tikzcd}
		\underline{X}\times_{\mc{D}}\underline{X} \arrow[dd, "{\pr}_2"'] \arrow[rr, "{\pr}_1"] & & \underline{X} \arrow[dd, "r"] \\
		& & \\
		\underline{X} 
		\arrow[Rightarrow, shorten >=30pt, shorten <=30pt, uurr]
		\arrow[rr, "r"] & & \mc{D} 
	\end{tikzcd}.\end{equation}

As $r\colon \underline{X}\rightarrow \mc{D}$ is an atlas for the stack $\mc{D}$, the $2$-fiber product $\underline{X}\times_{\mc{D}}\underline{X} $ is representable by a smooth manifold, which we denote by $X\times_{\mc{D}}X$. 
Since, stacks are fibered in groupoids, we can define the Lie groupoid   $\mc{G}_r=[X\times_{\mc{D}}X\rightrightarrows X]$, with the following structure maps:
\begin{enumerate}
	\item the source map $s:X\times_{\mc{D}}X\ra X$ is given by the morphism of smooth manifolds associated to the morphism of stacks ${\pr}_1\colon \underline{X}\times_{\mc{D}}\underline{X}\rightarrow \underline{X}$,
	\item the target map $t:X\times_{\mc{D}}X\ra X$ is given by the morphism of smooth manifolds associated to the morphism of stacks ${\pr}_2\colon \underline{X}\times_{\mc{D}}\underline{X}\rightarrow \underline{X}$,
	\item the composition is given by the morphism of stacks \[m\colon (\underline{X}\times_{\mc{D}}\underline{X})\times_{\underline{X}}(\underline{X}\times_{\mc{D}}\underline{X})\ra \underline{X}\times_{\mc{D}}\underline{X},\] defined as \[m\big((a,b,\alpha\colon r(a)\ra r(b)),(b,c,\beta\colon r(b)\ra r(c))\big)=\big(a,c,\beta\circ \alpha\colon r(a)\ra r(c)\big),\]
	\item the unit map is given by the morphism of stacks
	\[u\colon \underline{X}\ra \underline{X}\times_{\mc{D}}\underline{X},\] defined as \[u(a)= \big(a,a, \text{Id}\colon r(a)\ra r(a)\big),\]
	\item the inverse map is given by the morphism of stacks \[i\colon \underline{X}\times_{\mc{D}}\underline{X}\ra \underline{X}\times_{\mc{D}}\underline{X},\] defined as \[i\big((a,b,\alpha\colon r(a)\ra r(b))\big)=\big(b,a,\alpha^{-1}\colon r(b)\ra r(a)\big).\]
\end{enumerate}

It turns out that there is an isomorphism of stacks $\mc{D}\cong B\mc{G}_r$. 

\begin{lemma}[{\cite[Proposition $4.31$]{MR2778793}}]\label{Lemma:Liegroupoidrepresentingstack} Let $\pi_{\mc{D}}\colon \mc{D}\rightarrow \text{Man}$ be a differentiable stack and $r\colon \underline{X}\rightarrow \mc{D}$ an atlas for $\mc{D}$. 
	Then,
	there exists a Lie groupoid $\mc{G}$ with an isomorphism of stacks $\mc{D}\cong B\mc{G}$. Moreover, we may take $\mc{G}_0=X$ and $\mc{G}_1=X\times_{\mc{D}}X$, where $X\times_{\mc{D}}X$ is the smooth manifold representing the $2$-fiber product $\underline{X}\times_{\mc{D}}\underline{X}$. 
\end{lemma}

\subsection{Morita equivalence of Lie groupoids}
In the previous section, we have associated a Lie groupoid for a differentiable stack after making a choice of an atlas for the stack. We have also mentioned that, for every surjective submersion $Y\ra X$, and an atlas $X\ra \mc{D}$, the composition $Y\ra \mc{D}$ is an atlas for $\mc{D}$ (\Cref{Remark:compositionisatlas}).
In this way, we get a large collection of Lie groupoids associated to a single differentiable stack. We can ask the following question: how are these Lie groupoids related?

Let $r_X:X\ra \mc{D}$ and $r_Y:Y\ra \mc{D}$ be atlases with the property that, $r_Y\circ f=r_X$ for some surjective submersion $f:X\ra Y$. The $2$-fiber product diagrams (as in Diagram \ref{Diagram:XXDLiegroupoid}) for atlases $r_X, r_Y$ give the following diagram,
\[
\begin{tikzcd}
	&  &                                                    & \underline{Y}\times_{\mc{D}} \underline{Y} \arrow[dd, "{\pr}_1"] \arrow[rr, "{\pr}_2"] &  & \underline{Y} \arrow[dd, "r_Y"] \\
	&  & \underline{X} \arrow[rrrd, "r_X"]                  &                                                                                            &  &                                 \\
	\underline{X}\times_{\mc{C}}\underline{X} \arrow[rrd, "{\pr}_1"'] \arrow[rru, "{\pr}_2"] &  &                                                    & \underline{Y} \arrow[rr, "r_Y"]                                                            &  & \mc{D}                          \\
	&  & \underline{X} \arrow[rrru, "r_X"'] \arrow[ru, "f"] &                                                                                            &  &                                
\end{tikzcd}.\]
We expect some relation between the Lie groupoids $[X\times_{\mc{D}}X\rra X]$ and $[Y\times_{\mathcal{D}}Y\rra Y]$, given by the surjective submersion $f:X\ra Y$. It turns out that the ``pull-back of $[Y\times_{\mathcal{D}}Y\rra Y]$ along $f:X\ra Y$'' is isomorphic to the Lie groupoid $[X\times_{\mathcal{D}}Y\rra X]$. 

\subsubsection{pull-back of a Lie groupoid along a surjective submersion} Let $\Gamma=(\Gamma_1\rightrightarrows\Gamma_0)$ be a Lie groupoid and
$J\colon P_0\rightarrow \Gamma_0$ a surjective submersion. 

Consider the  pull-back of the source map $s\colon \Gamma_1\rightarrow \Gamma_0$ along $J\colon P_0\rightarrow \Gamma_0$ to obtain $P_0\times_{\Gamma_0}\Gamma_1$. We then pull-back the map $J\colon P_0\rightarrow \Gamma_0$ along $t\circ {\pr}_2\colon P_0\times_{\Gamma_0}\Gamma_1\rightarrow \Gamma_0$ to obtain $(P_0\times_{\Gamma_0}\Gamma_1)\times_{\Gamma_0}P_0$. Above pull-backs can be expressed by the following diagram,
\begin{equation}\label{Diagram:PullbackLiegroupoid}
	\begin{tikzcd}
		( P_0\times_{\Gamma_0}\Gamma_1)\times_{\Gamma_0} P_0 \arrow[dd] \arrow[rrrr] & & & & P_0 \arrow[dd, "J"] \\
		& & & & \\
		P_0\times_{\Gamma_0}\Gamma_1 \arrow[dd,"{\pr}_1"'] \arrow[rr,"{\pr}_2"] & & \Gamma_1 \arrow[dd, "s"] \arrow[rr, "t"] & & \Gamma_0 \\
		& & & & \\
		P_0 \arrow[rr, "J"] & & \Gamma_0 & & 
	\end{tikzcd}.
\end{equation}
Denote the smooth manifold $ (P_0\times_{\Gamma_0}\Gamma_1)\times_{\Gamma_0} P_0$ by $P_1$. The smooth manifold $P_1$ along with $P_0$ gives a Lie groupoid $[P_1\rightrightarrows P_0]$, with the following description \begin{enumerate}
	\item the source map $s\colon P_1\ra P_0$ is given by $(p,x,q)\mapsto p$,
	\item the target map $t\colon P_1\ra P_0$ is given by $(p,x,q)\mapsto q$,
	\item the composition map $m\colon P_1\times_{P_0}P_1\ra P_1$ is given by $\big((p,x,q),(q,y,r)\big)\mapsto (p,x\circ y,r)$,
	\item the identity map $u\colon P_0\ra P_1$ is given by $a\mapsto (a,1_{J(a)},a)$, 
	\item the inverse map $i\colon P_1\ra P_1$ is given by $(a,\gamma,b)\mapsto (b,\gamma^{-1},a)$.
\end{enumerate}

\begin{definition}
	Let $\Gamma=[\Gamma_1\rightrightarrows\Gamma_0]$ be a Lie groupoid and
	$J\colon P_0\rightarrow \Gamma_0$ a surjective submersion. The Lie groupoid $[P_1\rightrightarrows P_0]$ mentioned above (in the Diagram \ref{Diagram:PullbackLiegroupoid}) is called the \textit{pull-back groupoid  of the Lie groupoid $[\Gamma_1\rightrightarrows \Gamma_0]$ along the map $J\colon P_0\rightarrow \Gamma_0$}. 
\end{definition}
Under this notion of a pull-back of a Lie groupoid along a surjective submersion, the Lie groupoids $[X\times_{\mc{D}}X\rra X]$ and $[Y\times_{\mathcal{D}}Y\rra Y]$ associated to atlases $X\ra \mc{D}$ and $Y\ra \mc{D}$ (related by a surjective submersion $f:X\ra Y$) are related. More precisely, we have the following result.


\begin{proposition}\label{Proposition:relatedatlases}
	Let $\mc{D}$ be a differentiable stack, and $r_X:X\ra \mc{D},r_Y:Y\ra \mc{D}$ are atlases for the stack $\mc{D}$. Suppose that, there exists a surjective submersion $f:X\ra Y$, such that, $r_Y\circ f=r_X$. Then, the induced morphism of Lie groupoids $[X\times_{\mc{D}}X\rra X]\ra [Y\times_{\mc{D}}Y\rra Y]$ is with the property that $[X\times_{\mc{D}}X\rra X]$ is isomorphic to the pull-back groupoid of $[Y\times_{\mc{D}}Y\rra Y]$ along the surjective submersion. 
\end{proposition}

The above property has a name associated to it, a \textit{Morita morphism of Lie groupoids}.

\begin{definition}[Morita morphism of Lie groupoids {\cite[Section $2.2$]{MR2493616}}]\label{Definition:MoritamorphismofLiegroupoids}
	Let $[X_1\rra X_0],[Y_1\rra Y_0]$ be Lie groupoids. We say that a morphism of Lie groupoids $(\phi_1,\phi_0):[X_1\rra X_0]\ra [Y_1\rra Y_0]$ is a \textit{Morita morphism of Lie groupoids}, if, $\phi_0:X_0\ra Y_0$ is a surjective submersion, and $[X_1\rra X_0]$ is isomorphic to the pull-back of the Lie groupoid $[Y_1\rra Y_0]$ along the surjective submersion $\phi_0:X_0\ra Y_0$.
\end{definition}

Suppose that $r_X:X\ra \mc{D}$ and $r_Y:Y\ra \mc{D}$ are atlases. Consider the following $2$-fiber product diagram,
\[
\begin{tikzcd}
	\underline{X}\times_{\mc{D}}\underline{Y} \arrow[dd, "{\pr}_1"'] \arrow[rr, "{\pr}_2"] &  & \underline{Y} \arrow[dd, "r_Y"] \\
	&  &                                 \\
	\underline{X} \arrow[rr, "r_X"]                                                            &  & \mathcal{D}                    
\end{tikzcd}.\]
As $r_X$ is an atlas, the map ${\pr}_2:X\times_{\mc{D}}Y\ra Y$ is a surjective submersion. 
As $r_Y$ is an atlas, the map ${\pr}_1:X\times_{\mathcal{D}}Y\ra X$ is a surjective submersion. As $r_Y$ is an atlas for $\mc{D}$, and ${\pr}_2$ is a surjective submersion, the composition $X\times_{\mc{D}}Y\ra \mc{D}$ is an atlas for $\mc{D}$. This gives a morphism of Lie groupoids $[(X\times_{\mc{D}}Y)\times_{\mathcal{D}}(X\times_{\mc{D}}Y)\rra (X\times_{\mc{D}}Y)]\rra [Y\times_{\mc{D}}Y\rra Y]$. Similarly, treating $r_X$ as an atlas, for the surjective submersion ${\pr}_1$ as surjective submersion, we have morphism of Lie groupoids
$[(X\times_{\mc{D}}Y)\times_{\mathcal{D}}(X\times_{\mc{D}}Y)\rra (X\times_{\mc{D}}Y)]\rra [X\times_{\mc{D}}X\rra X]$.

As a consequence of \Cref{Proposition:relatedatlases}, we have the following result.
\begin{proposition}
	Let $\mc{D}$ be a differentiable stack, and $r_X:X\ra \mc{D}, r_Y:Y\ra \mc{D}$ atlases for the stack $\mc{D}$. Then, there exists a Lie groupoid $[(X\times_{\mc{D}}Y)\times_{\mathcal{D}}(X\times_{\mc{D}}Y)\rra (X\times_{\mc{D}}Y)]$ and a pair of morphism of Lie groupoids 
	\begin{align*}
		[(X\times_{\mc{D}}Y)\times_{\mathcal{D}}(X\times_{\mc{D}}Y)\rra (X\times_{\mc{D}}Y)]&\ra[X\times_{\mc{D}}X\rra X],\\
		[(X\times_{\mc{D}}Y)\times_{\mathcal{D}}(X\times_{\mc{D}}Y)\rra (X\times_{\mc{D}}Y)]&\ra[Y\times_{\mc{D}}Y\rra Y],
	\end{align*}
	such that, both these morphisms are Morita morphisms of Lie groupoids. 
\end{proposition}
Two Lie groupoids sharing the property mentioned above are called \textit{Morita equivalent Lie groupoids}, which we define below. 

\begin{definition}[Morita equivalent Lie groupoids {\cite[Section $2.2$]{MR2493616}}]\label{Definition:MoritaequivalentLiegroupoids}
	Let $[X_1\rra X_0]$ and $[Y_1\rra Y_0]$ be Lie groupoids. We say that $[X_1\rra X_0]$ and $[Y_1\rra Y_0]$ are \textit{Morita equivalent Lie groupoids}, if, there exists a Lie groupoid $[Z_1\rra Z_0]$ and Morita morphisms of Lie groupoids $[Z_1\rra Z_0]\ra [X_1\rra X_0]$, and
	$[Z_1\rra Z_0]\ra [Y_1\rra Y_0]$.
\end{definition}
Thus, we have the following result.
\begin{lemma}\label{Lemma:LiegroupoidsareMoritaequivalent}
	Let $\mc{D}\ra \text{Man}$ be a differentiable stack. If $X\ra \mc{D}$, and $Y\ra \mc{D}$ are atlases for the stack $\mc{D}$, then, the Lie groupoids $[X\times_{\mathcal{D}}X\rra X]$ and $[Y\times_{\mathcal{D}}Y\rra Y]$ are Morita equivalent.
\end{lemma}

\subsection{Bibundles}
Let $[X_1\rra X_0]$ and $[Y_1\rra Y_0]$ be Morita equivalent Lie groupoids. Thus, there exists a Lie groupoid $[Z_1\rra Z_0]$ and Morita morphisms $[Z_1\rra Z_0]\ra [X_1\rra X_0]$, and $[Z_1\rra Z_0]\ra [Y_1\rra Y_0]$. In particular, we have surjective submersions $Z_0\ra X_0, Z_0
\ra Y_0$. As $[Z_1\rra Z_0]$ is actually determined by the maps $Z_0\ra X_0, Z_0\ra Y_0$, we can see the notion of Morita equivalent Lie groupoids as the following. 

\begin{lemma}\label{Lemma:alternateforMoritaequivalence}
	Suppose that $[X_1\rra X_0]$ and $[Y_1\rra Y_0]$ be Morita equivalent Lie groupoids. Then, there exists a smooth manifold $P$ and surjective submersions $P\ra X_0, P\ra Y_0$, along with an action of $[X_1\rra X_0], [Y_1\rra Y_0]$ on $P$, such that $P\ra X_0$ is a principal $[Y_1\rra Y_0]$-bundle and $P\ra Y_0$ is a principal $[X_1\rra X_0]$-bundle. 
\end{lemma}
Relaxing the condition of $P\ra Y_0$ being a principal $[Y_1\rra Y_0]$-bundle in the \Cref{Lemma:alternateforMoritaequivalence}, gives an interesting notion of a ``generalized morphism from $[X_1\rra X_0]$ to $[Y_1\rra Y_0]$''. This structure goes by the name of an $\mb{X}-\mb{Y}$-bibundle.  
\begin{definition}[$\mc{G}-\mc{H}$ bibundle]\label{Definition:GHbibundle}
	Let $\mc{G},\mc{H}$ be Lie groupoids. A \textit{$\mc{G}-\mc{H}$ bibundle} consists of,
	\begin{enumerate}
		\item a smooth manifold $P$,
		\item a left action of $\mc{G}$ on $P$, with anchor map 
		$a_{\mc{G}}\colon P\rightarrow \mc{G}_0$,
		\item a right action of $\mc{H}$ on $P$, with anchor map $a_{\mc{H}}\colon P\rightarrow \mc{H}_0$,
	\end{enumerate}
	such that,
	\begin{enumerate}
		\item the anchor map $a_{\mc{G}}\colon P\rightarrow \mc{G}_0$ is a principal $\mc{H}$-bundle,
		\item the anchor map $a_{\mc{H}}\colon P\rightarrow \mc{H}_0$ is a $\mc{G}$-invariant map; that is, $a_{\mc{H}}(gp)=a_{\mc{H}}(p)$ for $p\in P$ and $g\in \mc{G}_1$ with $s(g)=a_{\mc{G}}(p)$,
		\item the action of $\mc{G}$ on $P$ is compatible with the action of $\mc{H}$ on $P$; that is, $(gp)h=g(ph)$ for $g\in \mc{G}_1, p\in P$ and $h\in \mc{H}_1$ with $s(g)=a_{\mc{G}}(p)$ and $t(h)=a_{\mc{H}}(p)$.
	\end{enumerate}
	We will express a $\mc{G}-\mc{H}$ bibundle by the following diagram,
	\begin{equation} \begin{tikzcd}
			\mc{G}_1 
			\arrow[dd,xshift=0.75ex,"t"]
			\arrow[dd,xshift=-0.75ex,"s"'] & & \mc{H}_1 
			\arrow[dd,xshift=0.75ex,"t"]
			\arrow[dd,xshift=-0.75ex,"s"'] \\
			& P \arrow[rd, "a_{\mc{H}}"] \arrow[ld, "a_{\mc{G}}"'] & \\
			\mc{G}_0 & & \mc{H}_0
		\end{tikzcd}.\end{equation}
\end{definition}
We denote a $\mc{G}-\mc{H}$ bibundle by $P\colon \mc{G}\ra \mc{H}$.
\begin{remark}\label{Remark:GPrincipalbibunde}
	If the anchor map $a_{\mc{H}}\colon P\ra \mc{H}_0$, in a $\mc{G}-\mc{H}$ bibundle $P\colon \mc{G}\ra \mc{H}$, is a principal $\mc{G}$-bundle, then, we call $P\colon \mc{G}\ra \mc{H}$ to be \textit{a $\mc{G}$-principal bibundle}.
\end{remark}

\begin{example}\label{Example:P(M,G)isbibundle}
Let $G$ be a Lie group, $M$ a smooth manifold, and $\pi:P\ra M$ a principal $G$ bundle. As mentioned before, we see the Lie group $G$ as the Lie groupoid $[G\rra *]$, and the smooth manifold $M$ as the Lie groupoid $[M\rra M]$. The action of $G$ on $P$ gives an action of the Lie groupoid $[G\rra *]$ on the manifold $P$. The smooth map $\pi:P\ra M$ gives an action of the Lie groupoid $[M\rra M]$ on the manifold $P$. These actions are seen to be compatible, giving a $[M\rra M]-[G\rra *]$-bibundle, as in the following diagram,
\[\begin{tikzcd}
	M
	\arrow[dd,xshift=0.75ex,"t"]
	\arrow[dd,xshift=-0.75ex,"s"'] & & G
	\arrow[dd,xshift=0.75ex,"t"]
	\arrow[dd,xshift=-0.75ex,"s"'] \\
	& P \arrow[rd] \arrow[ld, "\pi"'] & \\
	M & & *
\end{tikzcd}.\]
\end{example}
\begin{example}\label{Example:actionofLiegroupsisbibundle}
Let $G,H$ be Lie groups, seen as Lie groupoids $[G\rra *]$ and $[H\rra *]$ respectively. An action of $G$ on $H$, can be considered as the $[G\rra *]-[H\rra *]$-bibundle in the following diagram,
\[\begin{tikzcd}
	G
	\arrow[dd,xshift=0.75ex,"t"]
	\arrow[dd,xshift=-0.75ex,"s"'] & & H
	\arrow[dd,xshift=0.75ex,"t"]
	\arrow[dd,xshift=-0.75ex,"s"'] \\
	& H \arrow[rd] \arrow[ld] & \\
	* & & *
\end{tikzcd}.\]
Let $G,H$ be Lie groups and $\phi:G\ra H$ be a Lie group morphism. Then, we have an action of $G$ on $H$ given by $(g,h)\mapsto \phi(g)h$. Thus, as a special case, we have a $[G\rra *]-[H\rra *]$-bibundle. 
\end{example}
\begin{remark}\label{Remark:generalizedmorphismofLiegroupoids}
The \Cref{Example:actionofLiegroupsisbibundle} says that, for a morphism of Lie groups $\phi:G\ra H$, we can associate a $[G\rra *]-[H\rra *]$-bibundle. 
This would give a hope to associate a $\mc{G}-\mc{H}$-bibundle for a morphism of Lie groupoids $\phi:\mc{G}\ra \mc{H}$. We postpone this construction of associating a $\mc{G}-\mc{H}$-bibundle to a morphism of Lie groupoids $\mc{G}\ra \mc{H}$ to \Cref{Subsection:bibundleAssociatedtoMorphismofLiegroupoids}. This justifies the alternate terminology of a \textit{generalized morphism from $\mc{G}$ to $\mc{H}$} to a $\mc{G}-\mc{H}$-bibundle. 
 \end{remark}
\begin{example}\label{Example:targetmapisbibundle}
Let $[\mc{G}_1\rra \mc{G}_0]$ be a Lie groupoid. We have seen 
in \Cref{Example:targetmapisLGbundle} that the target map $t:\mc{G}_1\ra \mc{G}_0$ is considered as a principal $\mc{G}$-bundle (see Diagram \ref{Diagram:targetmapisLGbundle}). The map $t:\mc{G}_1\ra \mc{G}_0$ gives an action of the Lie groupoid $[\mc{G}_0\rra \mc{G}_0]$ on the manifold $\mc{G}_1$. These actions of $[\mc{G}_0\rra \mc{G}_0]$ from left side, $[\mc{G}_1\rra \mc{G}_0]$ from right side, on the manifold $\mc{G}_1$ are compatible. This gives the following $[\mc{G}_0\rra \mc{G}_0]-[\mc{G}_1\rra \mc{G}_0]$-bibundle,
\[\begin{tikzcd}
	\mc{G}_0
	\arrow[dd,xshift=0.75ex,"t"]
	\arrow[dd,xshift=-0.75ex,"s"'] & & \mc{G}_1
	\arrow[dd,xshift=0.75ex,"t"]
	\arrow[dd,xshift=-0.75ex,"s"'] \\
	& \mc{G}_1 \arrow[rd, "s"] \arrow[ld, "t"'] & \\
	\mc{G}_0 & & \mc{G}_0
\end{tikzcd}.\]
\end{example}

\begin{remark}\label{Remark:MoritaequivalentLiegroupoidsgiveisomorphicstacks}
	In  \Cref{Lemma:LiegroupoidsareMoritaequivalent}, we have seen that, for a differentiable stack $\mc{D}$, the associated Lie groupoids coming from any two atlases are Morita equivalent. In other words, if two differentiable stacks are isomorphic, then they can be represented by Morita equivalent Lie groupoids. It is also true that if two Lie groupoids $\mc{G}$, and $\mc{H}$ are Morita equivalent, then the associated differentiable stacks $B\mc{G}$, and $B\mc{H}$ are isomorphic. We see the proof of this result in  \Cref{Chap.2}.
\end{remark}

We end this section with a result about transitive Lie groupoids, which we will use in  \Cref{Chap.2}.

\begin{lemma}\label{Lemma:TransitiveLiegroupoidisMEtoLiegroup}
	Any transitive Lie groupoid $\mc{G}$ is Morita equivalent to the Lie group $\mc{G}_x$ for any $x\in \mc{G}_0$, that is, 
	the Lie groupoid $[\mc{G}_1\rightrightarrows \mc{G}_0]$ is Morita equivalent to the Lie groupoid $[\mc{G}_x\rightrightarrows *]$.
	\begin{proof}
	Let $x$ be an element in $\mc{G}_0$. Consider   a morphism of Lie groupoids $\psi\colon [\mc{G}_x\rightrightarrows *]\rightarrow [\mc{G}_1\rightrightarrows \mc{G}_0]$ with the following description:	\begin{itemize}
\item the element $*$ is mapped to the element $x$ in $\mc{G}_0$,
\item an element $g$ in $\mc{G}_x$ is mapped to the element $g$ in $\mc{G}_1$.
	\end{itemize}
We see the above morphism as the following diagram,
		\[\begin{tikzcd}
			\mc{G}_x \arrow[dd,xshift=0.75ex,"t"]\arrow[dd,xshift=-0.75ex,"s"'] \arrow[rr,"\psi_1"] & & \mc{G}_1 \arrow[dd,xshift=0.75ex,"t"]\arrow[dd,xshift=-0.75ex,"s"'] \\
			& & \\
			* \arrow[rr,"\psi_0"] & & \mc{G}_0 
		\end{tikzcd}.\]
		
		Observe that the morphism set of the pull-back groupoid is
		\begin{align*}
			*\times_{\psi_0,\mc{G}_0,s}\mc{G}_1\times_{t\circ {\pr}_2,\mc{G}_0,\psi_0}*&=\{(*,g,*)|\psi_0(*)=s(g), (t\circ {\pr}_2)(*,g)=\psi_0(*)\}\\
			&=\{ (*,g,*)|s(g)=x, t(g)=x \}\\
			&=\{(*,g,x)|g\in s^{-1}(x)\bigcap t^{-1}(x)\}\\
			&=\{(*,g,*)|g\in \mc{G}_x\}=\mc{G}_x.
		\end{align*}
		So, the pull-back groupoid of  $[\mc{G}_1\rightrightarrows \mc{G}_0]$ along the  map $*\rightarrow \mc{G}_0$ is the Lie groupoid $ [\mc{G}_x\rightrightarrows *]$. Thus, we conlude that  $\psi\colon [\mc{G}_x\rightrightarrows *]\rightarrow [\mc{G}_1\rightrightarrows \mc{G}_0]$ is a Morita morphism of Lie groupoids. So, $[\mc{G}_1\rightrightarrows \mc{G}_0]$ is Morita equivalent to the Lie groupoid $[\mc{G}_x\rightrightarrows *]$.
	\end{proof}
\end{lemma}

\section{Atiyah sequence and Chern-Weil theory for principal bundles over smooth manifold}\label{Section:Atiyahforsmoothmanifold}
This section recalls the differential geometry of principal bundles over smooth manifolds, their connections, curvature, and other results. We also recall the notion of the Atiyah sequence associated to a principal bundle and use it to introduce the notion of a connection on a principal bundle.

Throughout this section we will assume our smooth manifolds to be second countable and Hausdorff, $G$ will denote a Lie group and $\mf{g}$ its Lie algebra. For $g\in G$, let $\text{Ad}(g):G\ra G$ be the adjoint map defined as $\text{Ad}(g)(h)=ghg^{-1}$ for all $h\in G$. Let ${\rm ad}(g):\mf{g}\ra \mf{g}$ be the differential of $\text{Ad}(g):G\ra G$ at the identity. 

Let $P$ be a smooth manifold with an action of $G$. For $g\in G$, let $\delta_g:P\ra P$ be the right translation map, defined as $p\mapsto pg$ for $p\in P$. For $p\in P$, let $\delta_p:G\ra P$ be the map defined as $g\mapsto pg$ for $g\in G$. Given an element $A\in \mf{g}$, let us denote by $A^*:P\ra TP$, the vector field on $P$, defined as
$A^*(p)=(\delta_p)_{*,e}(A)$ for all $p\in P$.
We call $A^*$ to be the fundamental vector field on $P$ associated to the element $A$ of the Lie algebra $\mf{g}$.

\begin{definition}[connection structure on a principal bundle]\label{Definition:ConnectionOnP(M,G)} Let $M$ be a smooth manifold, $G$ a Lie group and $P(M,G)$ a principal bundle. 
	A \textit{connection structure on $P(M,G)$} is given by a $\mf{g}$-valued differential $1$-form $\omega:P\ra \Lambda^1_{\mf{g}}T^*P$, satisfying the following conditions:
	\begin{enumerate}
		\item $\omega(p)(A^*(p))=A$ for all $p\in P$ and $A\in \mf{g}$,
		\item $(\delta_g^*\omega)(p)(v)=(
		{\rm ad}(g^{-1}))(\omega(p)(v))$ for all $g\in G, p\in P$ and $v\in T_pP$.
	\end{enumerate}
\end{definition}
\begin{example}Let $M$ be a smooth manifold, and $G$ a Lie group. Let $\theta:G\ra \Lambda^1_{\mf{g}}T^*G$ be the Maurer-Cartan form of $G$. Let ${\pr}_1:M\times G\ra M$ be the first projection map, considered as a principal $G$-bundle. Let ${\pr}_2:M\times G\ra G$ be the second projection map. The pull-back ${\pr}_2^*\theta:M\times G\ra \Lambda^1_{\mf{g}}T^*(M\times G)$ is a connection structure on the principal bundle ${\pr}_1:M\times G\ra M$.
\end{example}
\begin{remark}
	Using the partition of unity argument, we can conclude that a connection exists on any principal bundle $P(M, G)$. Later in this section, we see that giving a connection structure on $P(M, G)$ is the same as providing a splitting of some short exact sequence of vector bundles over $M$. As any short exact sequence of vector bundles over a smooth manifold splits (again using the partition of unity argument), we can say that connection exists on any principal bundle $P(M, G)$.
\end{remark}

The notion of connection (and curvature) on a principal bundle over a smooth manifold in terms of a differential form (and distribution) is classical and is discussed in many books. In \Cref{Chap.3}, we define a connection on a principal bundle over a Lie groupoid using the Atiyah sequence. 
The language of Atiyah sequence and its splitting to define a connection on principal bundle over a manifold is relatively less commonly found in standard textbooks. So, here we will discuss this topic in bit more details.
 The only place where we could find it is the appendix of \cite{MR896907}.
%
We use \cite{MR896907} as a reference for this section. We first introduce the construction of quotient vector bundles. 

\subsection{Quotient vector bundles} 
\begin{proposition}[{\cite[Appendix A, Proposition $2.1$]{MR896907}}]\label{Proposition:quotientbundleproposition}
	Let $\pi:P\ra M$ be a principal $G$-bundle, and $p^E\colon E\rightarrow P$ a vector bundle. Suppose that $G$ acts on $E$, from the right side such that the following two conditions are satisfied: 
	\begin{enumerate}
		\item for each $g\in G$, the right translation maps $\delta_g^E:E\ra E, \delta_g^P:P\ra P$ induce an isomorphism of vector bundles $(\delta^E_g,\delta^P_g)\colon (E,p^E,P)\rightarrow (E,p^E,P)$,
		\item there exists an open cover $\{U_\alpha\}_{\alpha\in \Lambda}$ of $P$ with a trivialization $\{\psi_\alpha:U_\alpha\times V\ra (p^E)^{-1}(U_\alpha)\}$ of the vector bundle $p^{E}:E\ra P$, satisfying the following conditions:
		\begin{enumerate}
			\item each $U_\alpha$ is a $\pi$-saturated open set\footnote{that is $\pi^{-1}(\pi(U_\alpha))=U_\alpha$.} of the form $\pi^{-1}(V_\alpha)$ for some open $V_\alpha\subseteq M$, for each $\alpha\in \Lambda$,
			\item the trivialization map $\psi_\alpha:U_\alpha\times V\ra (p^E)^{-1}(U_\alpha)$ is $G$-equivariant, in the sense that 
			$\psi_\alpha(ug,v)=\psi(u,v)g$ for all $u\in U_\alpha, v\in V, g\in G$ and $\alpha\in \Lambda$.
		\end{enumerate}
	\end{enumerate} 
	Then, the orbit set $E/G$ has a unique vector bundle structure over $M$ such that the quotient map $q^E\colon E\rightarrow E/G$ is a surjective submersion, and $(q^E,\pi)\colon (E,p^E,P)\ra (E/G,p^{E/G},M)$ is a morphism of vector bundles. Further, 
 	\begin{equation}\label{PullbackDiagram}
		\begin{tikzcd}
			E \arrow[dd, "p^E"'] \arrow[rr, "q^E"] & & E/G \arrow[dd, "p^{E/G}"] \\
			& &   \\
			P \arrow[rr, "\pi"]   & & M  \end{tikzcd} ,
	\end{equation}
	is a pull-back diagram. In particular, for each $u\in P$, the map $q^E_u\colon E_u\rightarrow (E/G)_{\pi(u)}$ is an isomorphism of vector spaces. 
	
	We call $(E/G,p^{E/G},M)$ to be the \textit{quotient vector bundle} of $(E,p^E,P)$ by the action of $G$.
\end{proposition} 
As an application of the above proposition we construct the adjoint bundle and the Atiyah bundle associated to a principal bundle $P(M,G)$.
\begin{example}Let $M$ be a smooth manifold, $G$ a Lie group and $P(M,G)$ a principal bundle. Consider the tangent bundle $TP\ra P$, which is a vector bundle over the smooth manifold $P$. Consider the action of $G$ on $TP$, given by $((p,v),g)\mapsto (pg,(\delta_g)_{*,p}(v))$ for $g \in G, p\in P, v\in T_pP$. This action of $G$ on $TP$ satisfies the conditions mentioned in \Cref{Proposition:quotientbundleproposition}, giving the quotient vector bundle $(TP)/G\ra M$, which we call \textit{the Atiyah bundle associated to the principal bundle $P(M,G)$}. We see this construction as the following diagram of vector bundles,
	\begin{equation}\label{Tangentbundlequotient}
		\begin{tikzcd}
			TP \arrow[dd,"p^{TP}"'] \arrow[rr, "q^{TP}"] & & (TP)/G \arrow[dd,"p^{(TP)/G}"] \\
			& &   \\
			P \arrow[rr, "\pi"]  & & M  
		\end{tikzcd}.\end{equation}
\end{example}

\begin{example}
	Let $M$ be a smooth manifold, $G$ a Lie group and $P(M,G)$ a principal bundle. Consider the trivial bundle $p^{P\times \mf{g}}\colon P\times \mathfrak{g}\rightarrow P$, a vector bundle over the smooth manifold $P$. Consider the action of $G$ on $P\times \mf{g}$ given by $((p,A),g)\mapsto (pg,
	{\rm ad}(g^{-1})(A))$.
	This action of $G$ on $P\times \mf{g}$ satisfies the conditions mentioned in the \Cref{Proposition:quotientbundleproposition}, giving the quotient vector bundle $(P\times \mf{g})/G\ra M$, which we call \textit{the adjoint bundle associated to the principal bundle $P(M,G)$}. We see this construction as the following diagram of vector bundles,
	\begin{equation}\label{Trivialbundlequotient}
		\begin{tikzcd}
			P\times \mathfrak{g} \arrow[dd,"p^{P\times \mf{g}}"'] \arrow[rr, "q^{P\times \mathfrak{g}}"] & & (P\times \mathfrak{g})/G \arrow[dd,"p^{(P\times \mf{g})/G}"] \\
			& &       \\
			P \arrow[rr, "\pi"]       & & M      
		\end{tikzcd}.
	\end{equation}
\end{example}
\subsection{Morphisms of quotient vector bundles} Let $M,M'$ be smooth manifolds, $G,G'$ are Lie groups and $P(M,G),P'(M',G')$ are principal bundles. Consider vector bundles $(E,p^E,P)$ and $(E',p^{E'},P')$ over $P(M,G)$ and $P'(M',G')$ respectively satisfying the conditions mentioned in \Cref{Proposition:quotientbundleproposition}. This gives quotient vector bundles $(E/G,p^{E/G},M)$ and $(E'/G',p^{E'/G'},M')$. The following result gives a criterion for a morphism of vector bundles $(E,p^E,P)\ra (E',p^{E'},P')$ to induce a morphism of vector bundles 
$(E/G,p^{E/G},M)\ra (E'/G,p^{E'/G'},M')$.
\begin{proposition}[{\cite[Appendix $A$, Proposition $2.2(i)$]{MR896907}}]\label{Proposition:morphismofquotientbundles}
	Let $P(M,G), P'(M',G')$ be principal bundles, and 
	$(E,p^E,P), (E',p^{E'},P')$ are vector bundles over $P(M,G)$ and $P'(M',G')$ respectively, satisfying the conditions in the \Cref{Proposition:quotientbundleproposition}. Let $F(f,\phi)\colon P(M,G)\ra P'(M',G')$ be a morphism of principal bundles, and $(\Phi,F)\colon (E,p^E,P)\ra (E',p^{E'},P')$ a morphism of vector bundles, as in the following diagram, 
	\[\begin{tikzcd}
		E \arrow[dd, "p^E"] \arrow[rr, "\Phi"] \arrow[rrrrd, "q^E"] & & E' \arrow[dd, "p^{E'}"] \arrow[rrrrd, "q^{E'}"] & &       & &    \\
		& &       & & E/G \arrow[dd, "p^{E/G}"] & & E'/G' \arrow[dd, "p^{E'/G'}"] \\
		P \arrow[rr, "F"] \arrow[rrrrd, "\pi"]   & & P' \arrow[rrrrd, "\pi'"]   & &       & &    \\
		& & & & M \arrow[rr, "f"] & & M'    
	\end{tikzcd}.\]
	If $\Phi(\xi g)=\Phi(\xi) \phi(g)$ for all $\xi \in E, g\in G$, then, there is a unique morphism of vector bundles $(\Phi^{/G},f)\colon (E/G,p^{E/G},M)\ra (E'/G',p^{E'/G'},M')$ giving following diagram,
	\[\begin{tikzcd}
		E \arrow[dd, "p^E"] \arrow[rr, "\Phi"] \arrow[rrrrd, "q^E"] & & E' \arrow[dd, "p^{E'}"] \arrow[rrrrd, "q^{E'}"] & &       & &    \\
		& &       & & E/G \arrow[rr, "\Phi^{/G}"] \arrow[dd, "p^{E/G}"] & & E'/G' \arrow[dd, "p^{E'/G'}"] \\
		P \arrow[rr, "F"] \arrow[rrrrd, "\pi"]   & & P' \arrow[rrrrd, "\pi'"]   & &       & &    \\
		& &       & & M \arrow[rr, "f"]     & & M'    
	\end{tikzcd}.\]
\end{proposition}
\begin{example}\label{Example:morphismofquotientbundles}
	The \Cref{Proposition:morphismofquotientbundles}, allows us to  construct a morphism of vector bundles from the adjoint bundle of $P(M,G)$ to the Atiyah bundle of $P(M,G)$, which would be the first half of the Atiyah sequence associated to the principal bundle $P(M,G)$. 
	
	Consider the morphism of vector bundles $(j,1_P)\colon (P\times \mf{g},p^{P\times\mf{g}},P)\ra (TP,p^{TP},P)$ given by $j(p,A)= A^*(p)$, where $A^*\colon P\ra TP$ is the fundamental vector field on $P$ associated to the element $A\in \mf{g}$. It turns out that, this morphism of vector bundles satisfies the condition in the \Cref{Proposition:morphismofquotientbundles} giving the morphism of quotient vector bundles,
	\[
	\begin{tikzcd}
		(P\times \mf{g})/G \arrow[dd, "p^{(P\times \mf{g})/G}"'] \arrow[rr, "j^{/G}"] & & (TP)/G \arrow[dd, "p^{(TP)/G}"] \\
		& &     \\
		M \arrow[rr, "1_M"]        & & M    
	\end{tikzcd}.\] 
	This morphism $(j^{/G},1_M)\colon ((P\times \mf{g})/G,p^{(P\times \mf{g})/G},M)\ra ((TP)/G, p^{(TP)/G}, M)$ of quotient vector bundles is the first half of the Atiyah sequence of vector bundles over $M$.
\end{example}
\subsection{Morphism from a quotient vector bundle} Let $G$ be a Lie group, and $M$ a smooth manifold. Consider a principal bundle $P(M,G)$ and a vector bundle $(E,p^E,P)$ over $P(M,G)$. Suppose further that there is an action of $G$ on $E$ satisfying the conditions in \Cref{Proposition:quotientbundleproposition}. This gives the quotient vector bundle $(E/G,p^{E/G},M)$ associated to the vector bundle $(E,p^E,P)$. Let $(E',p^{E'},M')$ be another vector bundle. The following result gives a criterion for a morphism of vector bundles $(E,p^E,P)\ra (E',p^{E'},M')$ to induce a morphism of vector bundles 
$(E/G,p^{E/G},M)\ra (E',p^{E'},M')$. 
\begin{proposition}[{\cite[Appendix $A$, Proposition $2.2(ii)$]{MR896907}}]\label{Proposition:morphismfromquotientbundle}Let $P(M,G)$ be a principal bundle, and $(E,p^E,P)$ a vector bundle over $P(M,G)$ satisfying the conditions as in the \Cref{Proposition:quotientbundleproposition}. Let $(E',p^{E'},M')$ be a vector bundle, and $(\Phi,F)\colon (E,p^E,P)\ra (E',p^{E'},M')$ a morphism of vector bundles such that $\Phi(\xi g)=\Phi(\xi)$ for all $\xi\in E, g\in G$ and $F(ug)=F(u)$ for all $u\in P, g\in G$.
	Then, there exists a unique vector bundle morphism \[(\Phi^{/G},F^{/G})\colon (E/G,p^{E/G},M)\ra (E',p^{E'},M'),\] as in the following diagram,
	\[ \begin{tikzcd}
		E \arrow[dd, "p^E"'] \arrow[rr, "\Phi"] \arrow[rrrrd, "q^E"'] & & E' \arrow[dd, "p^{E'}"] & &         \\
		& &    & & E/G \arrow[ddd, "p^{E/G}"] \arrow[llu, "\Phi^{/G}"', bend right] \\
		P \arrow[rr, "F"] \arrow[rrrrdd, "\pi"']   & & M'   & &         \\
		& &    & &         \\
		& &    & & M \arrow[lluu, "F^{/G}"', bend right]    
	\end{tikzcd}.\]
	In particular, we have $F=F^{/G}\circ \pi$ and $\Phi = \Phi^{/G}\circ q^E$.
\end{proposition}
\begin{example}\label{Example:morphismfromquotientbundle}
The \Cref{Proposition:morphismfromquotientbundle}, allows us to construct the other half of the Atiyah sequence associated with $P(M, G)$.
	
	Consider the morphism of vector bundles $(\pi_*,\pi)\colon (TP,p^{TP},P)\ra (TM,p^{TM},M)$. It turns out that, the action of $G$ on $TP$ satisfies the conditions in the \Cref{Proposition:morphismfromquotientbundle}, giving the morphism of vector bundles as in the following diagram,
	\[
	\begin{tikzcd}
		(TP)/G \arrow[dd, "p^{(TP)/G}"'] \arrow[rr, "\pi_*^{/G}"] & & TM \arrow[dd, "p^{TM}"] \\
		& &    \\
		M \arrow[rr, "1_M"]     & & M   
	\end{tikzcd}.\] 
	This morphism of vector bundles $(\pi_*^{/G},1_M)\colon ((TP)/G,p^{(TP)/G},M)\ra (TM,p^{TM},M)$, is the second half of the Atiyah sequence of vector bundles over $M$.
\end{example}
\subsection{Sections of quotient vector bundle} Let $P(M,G)$ be a principal bundle. Let $(E,p^E,P)$ be a vector bundle over $P(M,G)$ satisfying conditions as in the \Cref{Proposition:quotientbundleproposition} giving the following quotient vector bundle,
\begin{equation}
	\begin{tikzcd}
		E \arrow[dd, "p^E"'] \arrow[rr, "q^E"] & & E/G \arrow[dd, "p^{E/G}"] \\
		& &    \\
		P \arrow[rr, "\pi"]   & & M  \end{tikzcd}. 
\end{equation}
As the above diagram is a pull-back diagram, it turns out that there is a bijective correspondence between the set of $G$-invariant sections of the vector bundle $(E,p^E, P)$ and the set of sections of the quotient vector bundle $(E/G,p^{E/G}, M)$.

Let $X\colon M\rightarrow E/G$ be a section of the vector bundle $p^{E/G}\colon E/G\rightarrow M$. As the above diagram is a pull-back diagram, by the uniqueness of pull-back, the following diagram 
\begin{equation}\label{PtoEmap}\begin{tikzcd}
		P \arrow[rrdddd, "1_P"', bend right] \arrow[rrrrdd, "X\circ \pi", bend left] & &     & &    \\
		& &     & &    \\
		& & E \arrow[dd, "p^E"'] \arrow[rr, "q^E"] & & E/G \arrow[dd, "p^{E/G}"] \\
		& &     & &    \\
		& & P \arrow[rr, "\pi"]   & & M   
	\end{tikzcd},\end{equation}
gives a unique smooth map $\overline{X}\colon P\ra E$ such that $p^E\circ \overline{X}=1_P$ and $q^E\circ \overline{X}=X\circ \pi$. 
The property $p^E\circ \overline{X}=1_P$ says that $\overline{X}\colon P\ra E$ is a section of the vector bundle $p^E\colon E\ra P$. The property $q^E\circ \overline{X}=X\circ \pi$ says that $q^E(\overline{X}(u))=X(\pi(u))$ for all $u\in P$.
Moreover, the section $\overline{X}:P
\ra E$ is $G$-equivariant. Thus, given a section $X\colon M\ra E/G$ of the vector bundle $p^{E/G}\colon E/G\ra M$, we have constructed a $G$-equivariant section $\overline{X}\colon P\ra E$ of the vector bundle $p^E\colon E\ra P$.

Conversely, let $X\colon P\rightarrow E$ be a section of the vector bundle $p^E\colon E\rightarrow P$; that is $p^E\circ X=1_P$. 
We impose the condition that $X\colon P\ra E$ is $G$-equivariant; 
that is, $X(ug)=X(u)g$ for all $u\in P$ and $g\in G$. Define $\underline{X}\colon M\ra E/G$ as $\underline{X}(m)=(q^E\circ X)(p)$, for each $m\in M$, where $p\in P$ is such that $\pi(p)=m$. 
To see this is well defined, take $p_1,p_2\in P$ such that $\pi(p_1)=\pi(p_2)$. 
Then, there exists $g\in G$ such that $p_2=p_1g$.
Then, $q^E(X(p_2))=q^E(X(p_1g))=q^E(X(p_1)g)=q^E(X(p_1))$.
Thus, $\underline{X}\colon M\ra E/G$
is well-defined. The condition $\underline{X}\circ \pi=q^E\circ X$ implies that $p^{E/G}\circ \underline{X}=1_M$.
This can be seen as follows: as $\underline{X}\circ \pi=q^E\circ X$,
we have \[p^{E/G}\circ \underline{X}\circ \pi=p^{E/G}\circ q^E\circ X=\pi\circ p^E\circ X=\pi;\] that is, $p^{E/G}\circ \underline{X}\circ \pi=\pi$. 
As $\pi\colon P\ra M$ is surjective, the condition $p^{E/G}\circ \underline{X}\circ \pi=\pi$ imply that $p^{E/G}\circ \underline{X}=1_M$; that is $\underline{X}\colon M\ra E/G$ is a section of the vector bundle $p^{E/G}\colon E/G\ra M$.

We have $\underline{X}\circ \pi=q\circ X$. As $\pi$ is a surjective submersion, the condition that the composition $\underline{X}\circ \pi$ is smooth implies that $\underline{X}$ is smooth. Thus, given a section $X\colon P\rightarrow E$ of $p^E\colon E\rightarrow P$ that is $G$-equivariant; in the sense that $X(pg)=X(p)g$, we get a section $\underline{X}\colon M\rightarrow E/G$ of $p^{E/G}\colon E/G\rightarrow M$. 

The conditions $q\circ \overline{X}=X\circ \pi$ for all $X\in \Gamma(E/G)$ and $\underline{X}\circ \pi=q\circ X$ for each $G$-equivariant $X\in \Gamma(E)$ says that the constructions $X\ra \overline{X}$ and $X\ra \underline{X}$ are  inverses of each other. Thus, we have the following result:
\begin{proposition}\label{sections}
	Let $G$	be a Lie group and $P(M,G)$ a principal bundle. Let $(E,p^E,P)$ be a vector bundle over $P(M,G)$ as in \Cref{Proposition:quotientbundleproposition}. Denote by $\Gamma^G(E)$ the set of sections of $E$ which are $G$-equivariant. Then, there is an isomorphism of $C^{\infty}(M)$- modules $\Gamma^G(E)\rightarrow \Gamma(E/G)$.
\end{proposition}

\subsection{Atiyah sequence of vector bundles} 
Consider the principal bundle $P(M,G)$. In  \Cref{Example:morphismofquotientbundles} and \Cref{Example:morphismfromquotientbundle} we have found morphisms of vector bundles 
\begin{align*}
	(j^{/G},1_M)&\colon ((P\times \mf{g})/G,p^{(P\times \mf{g})/G},M)\ra ((TP)/G, p^{(TP)/G}, M),\\
	(\pi_*^{/G},1_M)&\colon ((TP)/G,p^{(TP)/G},M)\ra (TM,p^{TM},M).
\end{align*}
These two morphisms of vector bundles combine to give a short exact sequence
\begin{equation}\label{Equation:AtiyahSequence}
	0\ra (P\times \mf{g})/G \xra{j^{/G}} (TP)/G\xra{\pi_*^{/G}} TM\ra 0,
\end{equation}
of vector bundles over the smooth manifold $M$. We call this short exact sequence to be the \textit{Atiyah sequence associated to the principal bundle $P(M,G)$}. It turns out that this  short exact sequence of vector bundles over $M$ induce a short exact sequence of (vector spaces) of sections, 
\[0\rightarrow \Gamma(M,(P\times \mf{g})/G) \ra \Gamma(M,(TP)/G) \ra \Gamma(M,TM)\ra 0.\]
Observe that the vector space $\Gamma(M, TM)$, of vector fields on $M$, has a Lie algebra structure. Using the Lie algebra structure of $\mf{g}$, we introduce a Lie algebra structure on the vector space $\Gamma(M,(P\times \mf{g})/G)$. Similarly, the Lie algebra structure on the vector space $\Gamma(P, TP)$ induces a Lie algebra structure on the vector space $\Gamma(M,(TP)/G)$.
\subsubsection{Lie algebra structure on $\Gamma(M,(P\times \mf{g})/G)$}
Let $X\colon M\ra (P\times \mf{g})/G$ be an element of $\Gamma(M,(P\times \mf{g})/G)$; 
that is, a section of the vector bundle $p^{(P\times \mf{g})/G}\colon (P\times \mf{g})/G\ra M$. By \Cref{sections}, this section $X:M\ra (P\times \mf{g})/G$ gives a $G$-equivariant section $\overline{X}:P\ra P\times \mf{g}$ of the trivial bundle $p^{(P\times \mf{g})}:P\times \mf{g}\ra P$. 
We can consider $\overline{X}:P\ra P\times \mf{g}$ as a $G$-equivariant map $\overline{X}\colon P\ra \mf{g}$. 
Let $Y\colon M\ra (P\times \mf{g})/G$ be an element of $\Gamma(M,(P\times \mf{g})/G)$; a section of the vector bundle $p^{(P\times \mf{g})/G}\colon (P\times \mf{g})/G\ra M$. This gives a $G$-equivariant map $\overline{Y}\colon P\ra \mf{g}$. 
For elements $X, Y$ of $\Gamma(M,(P\times \mf{g})/G)$, define $[\overline{X},\overline{Y}]\colon P\ra \mf{g}$ as $p\mapsto [\overline{X}(p),\overline{Y}(p)]$ for $p\in P$. This map $[\overline{X},\overline{Y}]\colon P\ra \mf{g}$ is in fact $G$-equivariant. The map $[\overline{X},\overline{Y}]:P\ra \mf{g}$ can be seen as a $G$-equivariant section $[\overline{X},\overline{Y}]\colon P\ra P\times \mf{g}$ of the trivial bundle $p^{(P\times \mf{g})}:P\times \mf{g}\ra P$. By \Cref{sections}, this section $[\overline{X},\overline{Y}]\colon P\ra P\times \mf{g}$ gives a section $[X,Y]\colon M\ra (P\times \mf{g})/G$ of the adjoint bundle $p^{(P\times \mf{g})/G}:(P\times \mf{g})/G\ra M$. It is easy to see that the description $(X,Y)\mapsto [X,Y]$ of the map 
$\Gamma(M,(P\times \mf{g})/G)\times \Gamma(M,(P\times \mf{g})/G)\ra \Gamma(M,(P\times \mf{g})/G)$ satisfies the condition for a binary operation to be a Lie bracket. With this Lie bracket structure, we consider $\Gamma(M, (P\times \mf{g})/G)$ as a Lie algebra.


\subsubsection{Lie algebra structure on $\Gamma(M,(TP)/G)$} Let $X\colon M\ra (TP)/G$ be an element of $\Gamma(M,(TP)/G)$; that is a section of the vector bundle $p^{(TP)/G}\colon (TP)/G\ra M$. By \Cref{sections}, this section $X\colon M\ra (TP)/G$ gives a $G$-equivariant section $\overline{X}\colon P\ra TP$. Let $Y\colon M\ra (TP)/G$ be an element of $\Gamma(M,(TP)/G)$. This gives a $G$-equivariant section $\overline{Y}\colon P\ra TP$. For elements $X,Y$ of $\Gamma(M,(TP)/G)$, consider the Lie bracket $[\overline{X},\overline{Y}]\colon P\ra TP$. The Lie bracket $[\overline{X},\overline{Y}]\colon P\ra TP$ is a $G$-equivariant map, thus producing a section $[X,Y]\colon M\ra (TP)/G$ of the Atiyah bundle $(TP)/G\ra M$. Moreover, the assignment  $(X,Y)\mapsto [X,Y]$ induces a Lie algebra structure on $\Gamma(M,(TP)/G)$. With this Lie bracket structure, we consider $\Gamma(M,(TP)/G)$ as a Lie algebra. 


Given a principal bundle $P(M,G)$, we have the associated Atiyah sequence 
\[0\ra (P\times \mf{g})/G\xra{j^{/G}} (TP)/G\xra{\pi_*^{/G}} TM\ra 0.\]
This produces a short exact sequence of vector spaces 
\begin{equation}\label{Eqaution:sectionsofAtiyah}
0\ra \Gamma(M,(P\times \mf{g}))/G\xra{j^{/G}} \Gamma(M,(TP)/G)
\xra{\pi_*^{/G}} \Gamma(M,TM)\ra 0.
	\end{equation}
The Lie bracket structures on $\Gamma(M,(P\times \mf{g}))/G$, $\Gamma(M,(TP)/G)$ and $\Gamma(M, TM)$ turn the Diagram \ref{Eqaution:sectionsofAtiyah} into a short exact sequence of Lie algebras. 

\subsection{Connection on $P(M,G)$} We have covered the necessary background to introduce the notion of connection on a principal bundle as a splitting of the Atiyah sequence. The following definition appears for example in 
\cite[Section $2$, page $188$]{MR86359}, \cite[Definition $1.21$]{cohen1998topology}, and 
\cite[Section $2$]{MR2748598}

\begin{definition}[connection on $P(M,G)$]
	Let $P(M,G)$ be a principal bundle. A \textit{connection on $P(M,G)$} is defined to be a splitting of the Atiyah sequence of vector bundles over $M$ (Diagaram \ref{Equation:AtiyahSequence}).
\end{definition}
\begin{remark}
	Recall that any exact sequence of vector bundles admits a splitting. So, for every principal bundle, there exists a connection.
\end{remark}
A splitting of a short exact sequence $0\ra P\xra{f} Q\xra{g} R\ra 0$ constructs a retract of the morphism $f\colon P\ra Q$; that is a morphism $\theta:Q\ra P$ such that $\theta\circ f=1_P$ or, equivalently, a section of the morphism $g\colon Q\ra R$; that is a morphism $\psi:R\ra Q$ such that $g\circ \psi=1_R$. 
\begin{definition}[infinitesimal connection]
	Let $G$ be a Lie group and $P(M,G)$ a principal bundle. An \textit{infinitesimal connection on $P(M,G)$} is defined to be a section $\gamma\colon TM\ra (TP)/G$ of the morphism of vector bundles $\pi_*^{/G}\colon (TP)/G\ra TM$.
\end{definition}

\begin{definition}[back-connection]
	Let $G$ be a Lie group and $P(M,G)$ a principal bundle. A \textit{back-connection on $P(M,G)$} is defined to be a retract $\mc{D}\colon (TP)/G\ra (P\times\mf{g})/G$ of the morphism of vector bundles $j_*^{/G}\colon (P\times \mf{g})/G\ra (TP)/G $.
\end{definition}

Let $P(M,G)$ be a principal bundle, and $\mc{D}:(TP)/G\ra (P\times \mf{g})/G$ a back connection  of $P(M,G)$. Consider the following diagram, obtained by combining the Diagrams \ref{Tangentbundlequotient}, \ref{Trivialbundlequotient}, and adding the maps $1_P:P\ra P, 1_M:M\ra M, \mc{D}:(TP)/G\ra (P\times \mf{g})/G$,
\[\begin{tikzcd} 
	TP \arrow[dd] \arrow[rr, "q^{TP}"] &                                                                      & (TP)/G \arrow[dd] \arrow[rd, "\mc{D}"] &                                     \\
	& P\times\mathfrak{g} \arrow[dd] \arrow[rr, "q^{P\times\mathfrak{g}}"] &                                        & (P\times \mathfrak{g})/G \arrow[dd] \\
	P \arrow[rr] \arrow[rd, "1_P"]            &                                                                      & M \arrow[rd, "1_M"]                           &                                     \\
	& P \arrow[rr]                                                         &                                        & M                                  
\end{tikzcd}.\]
Rewriting the above diagram gives the following diagram,
\[
\begin{tikzcd} 
TP \arrow[rrrd, "\mc{D}\circ q^{TP}"] \arrow[rddd, "p^{TP}"'] &                                                                      &   &                                     \\
& P\times\mathfrak{g} \arrow[dd, "p^{P\times \mf{g}}"] \arrow[rr, "q^{P\times\mathfrak{g}}"] &   & (P\times \mathfrak{g})/G \arrow[dd, "p^{(P\times \mf{g})/G}"] \\
&                                                                      &  &                                     \\
& P \arrow[rr]                                                         &   & M                                  
\end{tikzcd}\]

As the square is a pull-back diagram, there exists unique map $\widetilde{\mc{D}}:TP\ra P\times \mf{g}$, such that 
$q^{P\times \mf{g}}\circ \widetilde{\mc{D}}=\mc{D}\circ q^{TP}$. This map $\widetilde{\mc{D}}:TP\ra P\times \mf{g}$ gives  a $\mf{g}$-valued differential $1$-form $\omega:P\ra \Lambda^1_{\mf{g}}T^*P$, with the property that,
$\widetilde{\mc{D}}(v)=(p,\omega(p)(v))$ for all $p\in P$, and $v\in T_pP$.
\begin{lemma}\label{Lemma:propertiesofomega}
The differential $1$-form $\omega:P\ra \Lambda^1_{\mf{g}}T^*P$ associated to the back connection $\mc{D}$ satisfies the following conditions:
\begin{enumerate}
\item $\omega(p)(A^*(p))=A$ for all $p\in P$ and $A\in \mf{g}$,
\item $\omega(pg)((\delta_g)_{*,p}(v))={\rm ad}(g^{-1})(\omega(p)(v))$ for all $g\in G, p\in P, v\in T_pP$.
\end{enumerate}
\begin{proof}
We first observe that $\widetilde{\mc{D}}:TP\ra P\times \mf{g}$ is a retract of the map $j:P\times \mf{g}\ra TP$.  Observe that, the smooth maps $j:P\times \mf{g}\ra TP$ and $j^{/G}:(P\times \mf{g})/G\ra (TP)/G$ are related as $q^{TP}\circ j=j^{/G}\circ q^{P\times \mf{g}}$. Thus, we have 
\[\mc{D}\circ q^{TP}\circ j=\mc{D}\circ j^{/G}\circ q^{P\times \mf{g}}.\]
As $\mc{D}$ is a retract of $j^{/G}$, we have $\mc{D}\circ j^{/G}=1$. Further, $q^{P\times \mf{g}}\circ \widetilde{\mc{D}}=\mc{D}\circ q^{TP}$. So, 
\[q^{P\times \mf{g}}\circ \widetilde{\mc{D}}\circ j=q^{P\times \mf{g}}.\]
Consider the following diagram,
\[
\begin{tikzcd}
	P\times \mathfrak{g} \arrow[rddd, "p^{P\times \mathfrak{g}}"', bend right] \arrow[rrrd, "q^{P\times \mathfrak{g}}", bend left] &                                                                                                     &  &                                                                     \\
	& P\times \mathfrak{g} \arrow[dd, "p^{P\times \mathfrak{g}}"'] \arrow[rr, "q^{P\times \mathfrak{g}}"] &  & (P\times \mathfrak{g})/G \arrow[dd, "p^{(P\times \mathfrak{g})/G}"] \\
	&                                                                                                     &  &                                                                     \\
	& P \arrow[rr, "\pi"]                                                                                 &  & M                                                                  
\end{tikzcd}.\]
As the square in above diagram is a pull-back diagram, there exists unique map $P\times \mf{g}\ra P\times \mf{g}$ that fill the above diagram. As the maps $\widetilde{\mc{D}}\circ j, 1_{P\times \mf{g}}$ both fill the above diagram, we have $\widetilde{\mc{D}}\circ j= 1_{P\times \mf{g}}$. Thus, $\widetilde{\mc{D}}:TP\ra P\times \mf{g}$ is a retract of the map $j:P\times \mf{g}\ra TP$.

Let $(p,A)\in P\times \mf{g}$. We have $\widetilde{\mc{D}}(j(p,A))=(p,A)$; that is, $\widetilde{\mc{D}}(A^*(p))=(p,A)$. By the definition of $\omega$, we have 
$\widetilde{\mc{D}}(A^*(p))=(p, \omega(p)(A^*(p)))$. Thus, the condition 
$(p,A)=\widetilde{\mc{D}}(A^*(p))=(p, \omega(p)(A^*(p)))$ imply that 
$A=\omega(p)(A^*(p))$.

Let $g\in G, p\in P, v\in T_pP$. Consider $(\delta_g)_{*,p}(v)\in T_{pg}P$. 
We have 
\begin{align*}
q^{P\times \mf{g}} (\widetilde{\mc{D}}((\delta_g)_{*,p}(v)))&=(q^{P\times \mf{g}}\circ \widetilde{\mc{D}})((\delta_g)_{*,p}(v))\\
&=(\mc{D}\circ q^{TP})((\delta_g)_{*,p}(v))\\
&=\mc{D} (q^{TP}((\delta_g)_{*,p}(v))\\
&=\mc{D} (q^{TP}(v))\\
&=q^{P\times \mf{g}}(\widetilde{\mc{D}}(v)).
\end{align*}
Thus, there exists unique $h\in G$ such that 
$\widetilde{\mc{D}}((\delta_g)_{*,p}(v))=
h \widetilde{\mc{D}}(v)$. Rewriting this in terms of $\omega$, we have 
$(pg,\omega(pg)((\delta_g)_{*,p}(v)))=(ph, 
{\rm ad}(h^{-1})(\omega(p)(v)))$. As the action is free, we have $h=g$. So, 
$\omega(pg)((\delta_g)_{*,p}(v))=
{\rm ad}(g^{-1})(\omega(p)(v))$.
\end{proof}
\end{lemma}

 Similarly, given a $\mf{g}$-valued differential $1$-form  $\omega:P\ra \Lambda^1_{\mf{g}}T^*P$, satisfying the conditions in \Cref{Lemma:propertiesofomega}, one can associate a back-connection $\mc{D}$ on $P(M,G)$. For $\omega:P\ra \Lambda^1_{\mf{g}}T^*P$, define $\widetilde{\mc{D}}:TP\ra P\times \mf{g}$ as $v\mapsto (a,\omega(a)(v))$ for all $v\in T_aP$ and $a\in P$. The $G$-equivariance of $\omega$ implies that this map $\widetilde{\mc{D}}:TP\ra P\times \mf{g}$ induce a map 
 $ \mc{D}:(TP)/G\ra (P\times \mf{g})/G$. The condition $\omega(p)(A^*(p))=A$ for all $p\in P$ and $A\in \mf{g}$ implies that $ \mc{D}:(TP)/G\ra (P\times \mf{g})/G$ is a retract of the map $j^{/G}:(P\times \mf{g})/G\ra (TP)/G$. Thus, we get a back connection $ \mc{D}:(TP)/G\ra (P\times \mf{g})/G$ on the principal bundle $P(M,G)$.

Let $\gamma\colon TM\ra (TP)/G$ be an infinitesimal connection on $P(M,G)$. Consider the image $\gamma(TM)\subseteq (TP)/G$. Let $p\in P$. Consider the inverse image $(q^{TP}_p)^{-1}(\gamma(T_{\pi(p)}M))\subseteq T_pP$. This defines a distribution $\mc{H}\subseteq TP$. Similarly, given a distribution $\mc{H}\subseteq TP$ satisfying the conditions of definition of a connection as in \cite{MR1393940}, we get an infinitesimal connection on $P(M,G)$.

\subsection{Curvature form of the connection on the principal bundle $P(M,G)$}\label{CurvatureassociatedtoSplitting} 

Let $G$ be a Lie group and $P(M,G)$ a principal bundle. Let $\mc{D}\colon TM\ra (TP)/G$ be a connection of the principal bundle $P(M,G)$; that is, a section of the morphism $\pi_*^{/G}\colon (TP)/G\ra TM$.

Let $X\colon M\ra TM$ be a vector field on the smooth manifold $M$; that is, a section of the vector bundle $TM\ra M$. Then, the composition $\mc{D}\circ X\colon M\ra (TP)/G$ is a section of the quotient vector bundle $(TP)/G\ra M$. By \Cref{sections}, this section $\mc{D}\circ X\colon M\ra (TP)/G$ of the quotient vector bundle $q^{(TP)/G}:(TP)/G\ra M$ gives a section $\overline{\mc{D}\circ X}\colon P\ra TP$ of the vector bundle $TP\ra P$; that is, $\overline{\mc{D}\circ X}\colon P\ra TP$ is a vector field on $P$.
We denote the vector field $\overline{\mc{D}\circ X}$ as $\Phi(X)$.
Thus, a connection $\mc{D}$ gives a linear map $\Phi\colon \mf{X}(M)\ra \mf{X}(P)$, assigning to each vector field $X\colon M\ra TM$ on $M$, a vector field $\Phi(X)\colon P\ra TP$ on $P$. 
As the vector spaces $\mf{X}(M)$ and $\mf{X}(P)$ has the structure of a Lie algebra, one can ask if the map $\Phi$ is a homomorphism of Lie algebras. The failure of the map $\Phi$ to be a homomorphism of Lie algebras is measured as the curvature of the connection $\mc{D}$. 

Let $X,Y : M\ra TM$ be vector fields on $M$.
Then, $\Phi([X,Y])-[\Phi(X),\Phi(Y)]\colon P\ra TP$ is a vector field on $P$; that is a section of the vector bundle $TP\ra P$. From \Cref{sections}, this section $\Phi([X,Y])-[\Phi(X),\Phi(Y)]\colon P\ra TP$ of the vector bundle $TP\ra P$ gives a section $\underline{\Phi([X,Y])-[\Phi(X),\Phi(Y)]}\colon M\ra (TP)/G$ of the quotient vector bundle $(TP)/G\ra M$. Observe that \[\underline{\Phi([X,Y])-[\Phi(X),\Phi(Y)]}=\underline{\Phi([X,Y])}
-[\underline{\Phi(X)},\underline{\Phi(Y)}]=\underline{\overline{\mc{D}\circ [X,Y]}}-[\underline{\overline{\mc{D}\circ X}},\underline{\overline{\mc{D}\circ Y}}].\]
Observe from \Cref{sections} that the operations $X\mapsto \overline{X}$ and $Y\mapsto \underline{Y}$ are inverses to each other. Thus,
$\underline{\overline{\mc{D}\circ [X,Y]}}-[\underline{\overline{\mc{D}\circ X}},\underline{\overline{\mc{D}\circ Y}}]=\mc{D}\circ [X,Y]-[\mc{D}\circ X,\mc{D}\circ Y]$. Consider the composition \[\pi_*^{/G}\circ\underline{\Phi([X,Y])-[\Phi(X),\Phi(Y)]}\colon M\ra TM;\]
that is 
\[\pi_*^{/G}\circ(\mc{D}\circ [X,Y]-[\mc{D}\circ X,\mc{D}\circ Y])\colon M\ra TM\]
As $\pi_{*}^{/G}$ induce a Lie algebra homomorphism on sections, we have \[\pi_*^{/G}\circ[\mc{D}\circ X,\mc{D}\circ Y]=[\pi_*^{/G}\circ \mc{D}\circ X,\pi_*^{/G}\circ \mc{D}\circ Y].\] So, for the section $\underline{\Phi([X,Y])-[\Phi(X),\Phi(Y)]}\colon M\ra (TP)/G$ of the quotient vector bundle $(TP)/G\ra M$ we have 
\begin{align*}
	\pi_*^{/G}\circ \underline{\Phi([X,Y])-[\Phi(X),\Phi(Y)]}
	&=\pi_*^{/G}\circ(\mc{D}\circ [X,Y]-[\mc{D}\circ X,\mc{D}\circ Y])\\
	&=\pi_*^{/G}\circ(\mc{D}\circ [X,Y])-\pi_*^{/G}\circ([\mc{D}\circ X,\mc{D}\circ Y])\\
	&=(\pi_*^{/G}\circ \mc{D})\circ [X,Y]-[(\pi_*^{/G}\circ \mc{D})\circ X,(\pi_*^{/G}\circ \mc{D})\circ Y]\\
	&=[X,Y]-[X,Y]=0.
\end{align*}
As $\underline{\Phi([X,Y])-[\Phi(X),\Phi(Y)]}\colon M\ra (TP)/G$ is a section of the vector bundle $(TP)/G\ra M$ such that $\pi_*^{/G}\circ \underline{\Phi([X,Y])-[\Phi(X),\Phi(Y)]}\colon M\ra TM$ is the zero section of the vector bundle $TM\ra M$, there exists a section $\mc{K}_{\mc{D}}(X,Y)\colon M\ra (P\times \mf{g})/G$ of the vector bundle $(P\times \mf{g})/G\ra M$ such that $j^{/G}\circ \mc{K}_{\mc{D}}(X,Y)=\underline{\Phi([X,Y])-[\Phi(X),\Phi(Y)]}$. Thus, given vector fields $X,Y\colon M\ra TM$ of the smooth manifold $M$, we have a section $\mc{K}_{\mc{D}}(X,Y)\colon M\ra (P\times \mf{g})/G$ of the adjoint bundle $(P\times \mf{g})/G\ra M$. Observe that, for each $m\in M$, the element $\mc{K}_{\mc{D}}(X,Y)(m)$ belongs to the fiber of $m\in M$ in the vector bundle $(P\times \mf{g})/G$.
As mentioned before, the fiber of $p\in P$ in the vector bundle $P\times \mf{g}\ra P$ is isomorphic to the fiber of $m=\pi(p)\in M$ in the vector bundle $(P\times \mf{g})/G\ra M$; that is, $\mf{g}$ is the fiber of $m\in M$ in the vector bundle $(P\times \mf{g})/G\ra M$. 
So, we can see the section $\mc{K}_{\mc{D}}(X,Y)\colon M\ra (P\times \mf{g})/G$ as a smooth map $\mc{K}_{\mc{D}}(X,Y)\colon M\ra \mf{g}$. 
Thus, for vector fields $X,Y\colon M\ra TM$ of $M$, we have a smooth map $\mc{K}_{\mc{D}}(X,Y)\colon M\ra \mf{g}$. 
Thus, a connection $\mc{D}\colon TM\ra (TP)/G$, gives a $\mf{g}$-valued differential $2$-form on $M$, denoted by 
\begin{equation}\label{Equation:curvatureKDforclassicalprincipalbundle}
\mc{K}_{\mc{D}}\colon M\ra \Lambda^2_{\mf{g}}T^*M,
\end{equation} defined as $\mc{K}_{\mc{D}}(m)(v_1,v_2)=\mc{K}_{\mc{D}}(X,Y)(m)$, 
for all $m\in M$ and $v_1,v_2\in T_mM$, 
where, $X,Y$ are vector fields on $M$ such that $X(m)=v_1$ and $Y(m)=v_2$. This $\mf{g}$-valued $2$-form on $M$ is called \textit{the curvature form associated to the connection $\mc{D}$} of the principal bundle $P(M,G)$.

\begin{proposition}[{\cite[Proposition $4.17$]{MR896907}}]\label{Proposition:propertyofKD}
Let $P(M,G)$ be a principal bundle. Let $\mc{D}$ be a connection on $P(M,G)$, and $\omega:P\ra \Lambda^1_{\mf{g}}T^*P, \mc{K}_{\mc{D}}:M\ra \Lambda^2_{\mf{g}}T^*M$ the differential forms associated to $\mc{D}$. Then, we have \[\pi^*\mc{K}_{\mc{D}}=d\omega-[\omega,\omega].\]
\end{proposition}
\begin{remark}
Given a principal bundle $P(M,G)$, the \textit{curvature form associated to connection $\mc{D}$} is usually a declared to be the $2$-form on $P$, defined to be the pull-back $\pi^*\mc{K}_{\mc{D}}:P\ra \Lambda^2_{\mf{g}}T^*P$.
\end{remark}
\subsection{Chern-Weil morphism associated to a principal bundle}\label{ChernWeiloversmooth manifold} Let $G$ be a Lie group and $P(M,G)$ a principal bundle. In this section, we outline the construction of Chern-Weil morphism and recall the fact that the construction is independent of the choice of connection.

Let $\text{sym}_k(\mf{g})^G$ denote the set of all symmetric multi-linear mappings $f\colon \underbrace{\mathfrak{g}\times\cdots\times\mathfrak{g}}_{k ~\text{times}}\ra \mb{R}$ such that 
$f({\rm ad}(a)t_1,\cdots,{\rm ad}(a)t_k)=f(t_1,\cdots,t_k)$ for all $a\in G$ and $t_i\in \mf{g}$ for $1\leq i\leq k$. This set $\text{sym}_k(\mf{g})^G$ is an $\mb{R}$-vector space for each $k\geq 1$. The direct sum \[\text{sym}(\mf{g})^G=\bigoplus_{k\geq 0}\text{sym}_k(\mf{g})^G,\]is then made into an $\mb{R}$-algebra by defining the multiplication of $f\in \text{sym}_k(\mf{g})^G$ and $g\in \text{sym}_l(\mf{g})^G$ as 
\[fg(t_1,\cdots,t_{k+l})=\frac{1}{(k+l)!}\sum_{\sigma\in S_{k+l}}f(t_{\sigma(1)},\cdots,t_{\sigma(k)})g(t_{\sigma(k+1)},\cdots,t_{\sigma(k+l)}),\]
and extending it linearly to all of $\text{sym}(\mf{g})^G$. Thus, we have an $\mb{R}$-algebra $\text{sym}(\mf{g})^G$.

For the smooth manifold $M$, we have the $\mb{R}$-vector space $H^k_{{\rm dR}}(M,\mb{R})$, the $k$-th de Rham cohomology group of $M$ for each $k\geq 0$. The direct sum \[H_{{\rm dR}}^*(M,\mathbb{R})=\bigoplus_{k\geq 0} H^{k}_{{\rm dR}}(M,\mathbb{R}),\] is turned into an $\mb{R}$-algebra, with the product being the wedge product of differential forms.

 For each $f\in \text{sym}_k(\mf{g})^G$, we define an element in $H^{2k}(M,\mb{R})$ giving an $\mb{R}$-linear map of vector spaces $\text{sym}_k(\mf{g})^G\ra H^{2k}_{{\rm dR}}(M,\mb{R})$, and extend it linearly to give an $\mb{R}$-algebra morphism $\text{sym}(\mf{g})^G\ra H^*_{{\rm dR}}(M,\mb{R})$. This map is called the Chern-Weil map associated to the principal bundle $P(M,G)$. The outline of the construction is as follows:
\begin{enumerate}
	\item Fix a connection $\omega\colon P\ra \Lambda^1_{\mf{g}}T^*P$ on $P(M,G)$. Let $\Omega\colon P\ra\Lambda^2_{\mf{g}}T^*P$ be the associated curvature form for the connection $\omega$.
	\item Given an element $f\in \text{sym}_k(\mf{g})^G$, we associate a differential $2k$-form $f(\Omega)$ on $P$. We have 
	\begin{equation}\label{definitionoffOmega}
		\begin{split}
			f(\Omega)(X_1,X_2,&\cdots,X_{2k-1},X_{2k})\\
			&=\frac{1}{(2k)!}\sum_{\sigma\in S_{2k}}\text{sgn}(\sigma)
			f(\Omega(X_{\sigma(1)},X_{\sigma(2)}),\cdots,\Omega(X_{\sigma(2k-1)},X_{\sigma(2k)})),
		\end{split}
	\end{equation}
	for vector fields $X_i$ on $P$ for $1\leq i\leq 2k$, 
	where, sum is taken over elements of $S_{2k}$; that is the set of permutations of the set $\{1,\cdots,2k\}$ and $\text{sgn}(\sigma)$ denotes the sign of the permutation $\sigma\in S_{2k}$.
	\item We generalize the construction in step $2$ to associate a differential $q_1+\cdots+q_k$-form on $P$ for each $f\in \text{sym}_k(\mf{g})^G$ and differential forms $\varphi_1,\cdots,\varphi_k$ on $P$ of degree $q_1,\cdots,q_k$ respectively. We have
	\begin{equation}\label{definitionoffvarphi}
		\begin{split}
			f(\varphi_1,&\cdots,\varphi_k)(X_1,X_2,\cdots,X_{q_1+\cdots+q_k-1},
			X_{q_1+\cdots+q_k})\\
			&=\frac{1}{(q_1+\cdots+q_k)!}\sum_{\sigma\in S_{2k}}\text{sgn}(\sigma)
			f(\varphi_1(X_{\sigma(1)},\cdots, X_{\sigma(q_1)}),\\
			&~~~~~~~~~~~~~~~~~~~~~~~~~~\cdots,\varphi_k
			(X_{\sigma(q_1+\cdots+q_{k-1})},\cdots,X_{\sigma(q_1+\cdots+q_k)})),
		\end{split}
	\end{equation}
	for vector fields $X_i$ on $P$ for $1\leq i\leq 2k$, where, sum is taken over elements of $S_{q_1+\cdots+q_k}$; that is the set of permutations of the set $\{1,\cdots,q_1+\cdots+q_k\}$ and $\text{sgn}(\sigma)$ denotes the sign of the permutation $\sigma\in S_{q_1+\cdots+q_k}$.
	\item Given $f\in \text{sym}_k(\mf{g})^G$, the $2k$-form $f(\Omega)$ on $P$ constructed in step $2$ is seen to be projected to a unique $2k$-form $\widetilde{f(\Omega)}$ on $M$; that is, $\pi^*\widetilde{f(\Omega)}=f(\Omega)$. It is also observed that, this $\widetilde{f(\Omega)}$ is a closed form on $M$, giving an element in the 
	de Rham cohomology group $H^{2k}_{{\rm dR}}(M,\mb{R})$.
\end{enumerate}
Thus, we have an $\mb{R}$-linear map $\text{sym}_k(\mf{g})^G\ra H^{2k}_{{\rm dR}}(M,\mb{R})$. As mentioned before, we extend this to get an $\mb{R}$-algebra morphism 
$W\colon \text{sym}(\mf{g})^G\ra H^*_{{\rm dR}}(M,\mb{R})$. This map is called the Chern-Weil homomorphism. It is a standard fact that, though we fix a connection $\omega$ on $P(M,G)$ when defining the map $W$, this map does not depend on the connection. In other words, for any two connections $\omega_0$ and $\omega_1$ on $P(M,G)$ with curvatures $\Omega_0$ and $\Omega_1$ respectively, the $2k$-forms $\widetilde{f(\Omega_0)}$ and 
$\widetilde{f(\Omega_1)}$ represent same element in the de Rham cohomology group $H^{2k}_{{\rm dR}}(M,\mb{R})$; that is, $\widetilde{f(\Omega_0)}-\widetilde{f(\Omega_1)}$
is an exact form on $M$. To prove this result, we recall the following standard result:
\begin{lemma}\label{useful}
	Let $\pi\colon P\rightarrow M$ be a surjective submersion. Let $\omega$ be a $k$-form on $M$ and $\tau$ a $k$-form on $P$ such that $\pi^*\omega=\tau$. If $\tau$ is exact, and  $\tilde{\tau}=\pi^*\tilde{\omega}$ for some $k-1$ form $\tilde{\omega}$ on $M$, then, $\omega$ is exact.
\end{lemma}
Observe, by above result, that, to prove exactness of $\widetilde{f(\Omega_0)}-\widetilde{f(\Omega_1)}$, it suffices to prove that $f(\Omega_0)-f(\Omega_1)=d\Phi$ for a $2k-1$ form $\Phi$ on $P$, and that $\Phi$ projects to a $2k-1$ form on $M$. The outline of this construction is as follows:
\begin{enumerate}
	\item Set $\alpha=\omega_1-\omega_0$ and $\omega_t=\omega_0+t\alpha$ for $0\leq t\leq 1$. Observe that these $\omega_t$ are connection forms on the principal bundle $P(M,G)$. Let $\Omega_t$ be the corresponding curvature forms.
	\item We need to get a $2k-1$ form on $P$. If we are given a differential $1$-form on $P$ and $k-1$ differential $2$-forms on $P$, we know (\Cref{definitionoffvarphi}) how to construct a $2k-1$ form. Considering the $1$-form $\alpha$ on $P$ and the $2$-form $\Omega_t$ for $k-1$ times, we get differential $2k-1$ form, $f(\alpha,\Omega_t,\cdots,\Omega_t)$ on $P$.
	\item We only needed one differential $2k-1$ form on $P$, but, we have a family of differential $2k-1$ forms on $P$, indexed by elements of $[0,1]$. Consider the integral 
	\begin{equation}\label{definitionofPhi}
		\Phi=k\int_0^1 f(\alpha,\Omega_t,\cdots,\Omega_t)dt.
	\end{equation}
	This $\Phi$ is a differential $2k-1$ form on $P$.
	\item Observe that this $\Phi$ projects to a $2k-1$ form on $M$ and that \[d\Phi=f(\Omega_1)-f(\Omega_0)=f(\Omega_1,\cdots,\Omega_1)-f(\Omega_0,\cdots,\Omega_0).\]
	\item Thus, as observed before, the difference $\widetilde{f(\Omega_1)}-\widetilde{f(\Omega_0)}$ is also an exact form on $M$. Thus, $\widetilde{f(\Omega_1)}$ and $\widetilde{f(\Omega_0)}$ gives same element in the $2k$-th de Rham cohomology group $H^{2k}_{{\rm dR}}(M,\mb{R})$.
\end{enumerate}
Thus, the Chern-Weil morphism $\text{sym}(\mf{g})^G\ra H^*_{{\rm dR}}(M,\mb{R})$ associated to the principal bundle $P(M,G)$ 
does not depend on the choice of connection on $P(M,G)$.



	\chapter{On two notions of a gerbe over a stack} \label{Chap.2}

In the \Cref{Section:LiegroupoidsMoritaequivalencedifferentiablestacks}, we have introduced the notions of Lie groupoids, differentiable stacks and their correspondence.
 In the \Cref{Section:associatingastackforaLiegroupoid}, we have assigned a differentiable stack $B\mc{G}$, for a Lie groupoid $\mc{G}$. 
In \Cref{Lemma:Liegroupoidrepresentingstack}, we saw that, for each differentiable stack $\mc{D}$, we can associate a Lie groupoid $\mc{G}$ with the property that $\mc{D}\cong B\mc{G}$. We have mentioned in \Cref{Lemma:LiegroupoidsareMoritaequivalent} that, any two Lie groupoids obtained from any two atlases of a differentiable stacks are Morita equivalent. Then, we observed in \Cref{Remark:MoritaequivalentLiegroupoidsgiveisomorphicstacks} that, for any two Lie groupoids that are Morita equivalent, the associated stacks $B\mc{G}$ and $B\mc{H}$ are isomorphic.  Thus, there is a one-one correspondence between the collection of Morita equivalent Lie groupoids and the collection of differentiable stacks.

The main aim of this chapter is to extend the above correspondence to a correspondence between a particular type of morphism of Lie groupoids and a particular type of morphism of differentiable stacks. More precisely, we study morphism of Lie groupoids that are Lie groupoid extensions, and morphism of differentiable stacks that are epimorphisms along with the condition that the associated diagonal morphism is an epimoprhism of stacks. In literature, these two seemingly different notions (Lie groupoid extension, morphism of stacks with extra conditions) goes by the name of a gerbe over a stack.
The notion of a gerbe over a stack as a (Morita equivalence class of a) Lie groupoid extension is studied by Laurent-Gengoux, Stienon, and Xu in \cite{MR2493616}. The notion of a gerbe over a stack as a morphism of differentiable stacks is studied by Behrend and Xu in \cite{MR2817778}, Ginot in \cite{Ginot}, and Heinloth in \cite{MR2206877}, to name a few.

This chapter, based on our paper \cite{MR4124773}, is split into three sections. In the first section, we will recall the two notions of a gerbe over a stack. In the second section, assuming an extra condition on a morphism of stacks that is a gerbe over a stack, we associate a Lie groupoid extension.
 We also prove that the Lie groupoid extension produced in this way is unique up to a Morita equivalence. In the last section, we prove that the morphism of stacks associated to a Lie groupoid extension is a gerbe over a stack. 

\section{Two notions of a gerbe over a stack}\label{Section:Twonotionsofgerbe}
In this section, we recall the existing notions of a gerbe over a stack from references \cite{MR2493616,MR2817778}. 
\begin{definition}\label{Definition:Liegroupidextension}
	Let $[\mc{H}_1\rra M]$ be a Lie groupoid. A \textit{Lie groupoid extension of $[\mc{H}_1\rra M]$} is given by a Lie groupoid $[\mc{G}_1\rra M]$, and a morphsim of Lie groupoids
	$(F,1_M):[\mc{G}_1\rra M]\ra [\mc{H}_1\rra M]$, such that $F:\mc{G}_1\ra \mc{H}_1$ is a surjective submersion.
\end{definition}

We denote a Lie groupoid extension by $F\colon X_1\ra Y_1\rightrightarrows M$.

In \cite{MR2493616} the authors assume $F:\mc{G}_1\ra \mc{H}_1$ to be a smooth fibration, but we restrict ourselves to the case of a surjective submersion. Note that every (smooth) fibration is a surjective submersion.

The notions of a Morita morphism of Lie groupoids (\Cref{Definition:MoritamorphismofLiegroupoids}) and Morita equivalent Lie groupoids (\Cref{Definition:MoritaequivalentLiegroupoids}) extends respectively to the notions of a Morita morphism of Lie groupoid extensions and Morita equivalent Lie groupoid extensions. 

\begin{definition}[Morphism of Lie groupoids extensions]
	Let $F:\mc{G}_1\ra \mc{H}_1\rra M$ and $F':\mc{G}_1'\ra \mc{H}_1'\rra M'$ be Lie groupoid extensions. A \textit{morphism of Lie groupoid extensions from $\mc{G}_1\ra \mc{H}_1\rra M$ to $\mc{G}_1'\ra \mc{H}_1'\rra M'$} is given by a triple  \[(\psi_{\mc{G}}:\mc{G}_1\ra \mc{G}_1', \psi_{\mc{H}}:\mc{H}_1\ra \mc{H}_1', f:M\ra M'),\] of smooth maps, such that, the following conditions are satisfied:
	\begin{enumerate}
		\item $F'\circ \psi_{\mc{G}}=\psi_{\mc{H}}\circ F$,
		\item $(\psi_{\mc{G}},f):[\mc{G}_1\rra M]\ra [\mc{G}'_1\rra M']$ is a morphism of Lie groupoids,
		\item $(\psi_{\mc{H}},f):[\mc{H}_1\rra M]\ra [\mc{H}'_1\rra M']$ is a morphism of Lie groupoids.
	\end{enumerate}
\end{definition}
We see a morphism of Lie groupoid extensions as the following diagram,
\begin{equation}
	\begin{tikzcd}
		\mc{G}_1 \arrow[dd, "\psi_{\mc{G}}"'] \arrow[rr,"F"] &  & \mc{H}_1 \arrow[rr,yshift=0.75ex,"s"]\arrow[rr,yshift=-0.75ex,"t"'] \arrow[dd, "\psi_{\mc{H}}"'] &  & M \arrow[dd, "f"'] \\
		&  &                                                  &  &                    \\
		\mc{G}_1' \arrow[rr,"F'"]                             &  & \mc{H}_1' \arrow[rr,yshift=0.75ex,"s"]\arrow[rr,yshift=-0.75ex,"t"']                              &  & M'                
	\end{tikzcd}.\end{equation}
We sometimes denote a morphism of Lie groupoid extensions as a pair $((\psi_{\mc{G}},f),(\psi_{\mc{H}},f))$ of morphisms of Lie groupoids instead of a triple $(\psi_{\mc{G}}, \psi_{\mc{H}}, f)$ of smooth maps. 
\begin{definition}[Morita morphism of Lie groupoid extensions]\label{Definition:MoritamorphismofGroupoidExtensions}
	Let $F:\mc{G}_1\ra \mc{H}_1\rra M$ and $F':\mc{G}_1'\ra \mc{H}_1'\rra M'$ be Lie groupoid extensions. A  morphism $((\psi_{\mc{G}},f),(\psi_{\mc{H}},f))$ of Lie groupoid extensions from $\mc{G}_1\ra \mc{H}_1\rra M$ to $\mc{G}_1'\ra \mc{H}_1'\rra M'$ is said to be \textit{a Morita morphism of Lie groupoid extensions}, if $(\psi_{\mc{G}},f)$, and $(\psi_{\mc{H}},f)$ are Morita morphisms of Lie groupoids. 
\end{definition}
\begin{definition}[Morita equivalent Lie groupoid extensions]\label{Definition:MoritaequivalentLiegroupoidextensions}
	Let $\phi'\colon X_1'\rightarrow Y_1'\rightrightarrows M'$ and 
	$\phi\colon X_1\rightarrow Y_1\rightrightarrows M$ be a pair of Lie groupoid extensions. We say that $\phi'\colon X_1'\rightarrow Y_1'\rightrightarrows M'$ and 
	$\phi\colon X_1\rightarrow Y_1\rightrightarrows M$ are \textit{Morita equivalent Lie groupoid extensions}, if there exists a third Lie groupoid extension $\phi''\colon X_1''\rightarrow Y''\rightrightarrows M''$ and a pair of Morita morphisms of Lie groupoid extensions
	\[(\phi''\colon X_1''\rightarrow Y''\rightrightarrows M'')\rightarrow (\phi\colon X_1\rightarrow Y_1\rightrightarrows M)\] and \[(\phi''\colon X_1''\rightarrow Y''\rightrightarrows M'')\rightarrow (\phi'\colon X_1'\rightarrow Y_1'\rightrightarrows M').\]
\end{definition}

\begin{definition}[Gerbe over a stack as a Lie groupoid extension \cite{MR2493616}] \label{Definition:GerbeoverstackVersion1} 
	A \textit{gerbe over a stack} is the Morita equivalence class of a Lie groupoid extension.
\end{definition}
This is all about the notion of a gerbe over a stack as a Morita equivalence class of a Lie groupoid extension. Now, we introduce the notion of a  gerbe over a satck as as morphism of differentiable stacks satisfying certain conditions. 
\begin{definition}[A Gerbe over a stack \cite{MR2817778}]
	\label{Definition:GerbeoverstackasinMR2817778} Let $\mc{C}$ be a differentiable stack. A morphism of stacks $F\colon \mc{D}\rightarrow \mc{C}$ is said to be \textit{a gerbe over the stack $\mc{C}$}, if the morphism $F\colon \mc{D}\rightarrow \mc{C}$ and the diagonal morphism $\Delta_F\colon \mc{D}\rightarrow \mc{D}\times_{\mc{C}}\mc{D}$ associated to $F\colon \mc{D}\ra \mc{C}$ are epimorphisms of stacks.
\end{definition} 

Using the $2$-Yoneda Lemma (\Cref{Lemma:2-yoneda}) one can give a following characterization of a gerbe over a stack.

Suppose that $F\colon \mc{D}\rightarrow \mc{C}$ is an epimorphism of stacks. Given a manifold $U$ and a morphism of stacks $q\colon \underline{U}\rightarrow \mc{C}$, there exists a cover $\{U_\alpha\rightarrow U\}$ of $U$ and a morphism of stacks $L_\alpha\colon \underline{U_\alpha}\ra \mc{D}$ with the following $2$-commutative diagram, 
\begin{equation}\label{Diagram:LalphaDtoC}\begin{tikzcd}
		\underline{U_\alpha} \arrow[dd,"L_\alpha"'] \arrow[rr] & & \underline{U} \arrow[dd, "q"] \\
		& & \\
		\mc{D} \arrow[Rightarrow, shorten >=20pt, shorten <=20pt, uurr] \arrow[rr, "F"] & & \mc{C}
	\end{tikzcd}.\end{equation}
As $\pi_{\mc{C}}\colon \mc{C}\ra \text{Man}$ and $\pi_{\mc{D}}\colon \mc{D}\ra \text{Man}$ are categories fibered in groupoids, we can use the $2$-Yoneda lemma. As $U$ is an object of $\text{Man}$, the morphism $q\colon \underline{U}\rightarrow \mc{C}$ corresponds 
to an object $a$ of $\mc{C}(U)$.
Similarly for $U_\alpha$, the morphism $L_\alpha\colon \underline{U}_\alpha\rightarrow \mc{D}$ corresponds 
to an object $x_\alpha$ of $\mc{D}(U_\alpha)$. The $2$-commutative Diagram \ref{Diagram:LalphaDtoC}, corresponds to an isomorphism $F(x_\alpha)\rightarrow a|_{U_\alpha}$ in $\mc{C}(U_\alpha)$ for each $\alpha$.

Thus, if a morphism of stacks $F\colon \mc{D}\rightarrow \mc{C}$ is an epimorphism, then given a manifold $U$ and an object $a$ of $\mc{C}(U)$, there exists an open cover $\{U_\alpha\rightarrow U\}$ of $U$ and objects $x_\alpha$ of $\mc{D}(U_\alpha)$ with an isomorphism $F(x_\alpha)\rightarrow a|_{U_\alpha}$ in $\mc{C}(U_\alpha)$ for each $\alpha$. It turns out that the converse is true as well.

Suppose that $F:\mc{D}\ra\mc{C}$ is a morphism of stacks with the following property: for a manifold $U$ and an object $a$ of $\mc{C}(U)$, there exists an open cover $\{U_\alpha\ra U\}$ of $U$ and objects $x_\alpha\in \mc{D}(U_\alpha)$ such that $F(x_\alpha)$ is isomorphic to $a|_{U_\alpha}$ for each $\alpha$. Then $F\colon \mc{D}\ra \mc{C}$ is an epimorphism of stacks.

Assume further that $\Delta_F\colon \mc{D}\rightarrow \mc{D}\times_{\mc{C}}\mc{D}$ is an epimorphism of stacks. Given a manifold $U$ and a morphism of stacks $q\colon \underline{U}\rightarrow \mc{D}\times_{\mc{C}}\mc{D}$, there exists a cover $\{U_\alpha\rightarrow U\}$ and a morphism of stacks $L_\alpha\colon \underline{U_\alpha}\rightarrow \mc{D}$ such that we have following $2$-commutative diagram,
\begin{equation} \begin{tikzcd}
		\underline{U_\alpha} \arrow[dd,"L_\alpha"'] \arrow[rr] & & \underline{U} \arrow[dd, "q"] \\
		& & \\
		\mc{D} \arrow[rr, "\Delta"] \arrow[Rightarrow, shorten >=30pt, shorten <=30pt, uurr] & & \mc{D}\times_{\mc{C}}\mc{D}
	\end{tikzcd}.\end{equation}
Then by $2$-Yoneda lemma,
the map $q\colon \underline{U}\rightarrow \mc{D}\times_{\mc{C}}\mc{D}$ corresponds to an object $\big(a,b,p\colon F(a)\rightarrow F(b)\big)$ in $(\mc{D}\times_{\mc{C}}\mc{D})(U)$; that is, $a\in \mc{D}(U)_0,b\in \mc{D}(U)_0 \text { and } p\in \mc{C}(U)_1$. 
The morphism of stacks $L_\alpha\colon \underline{U_\alpha}\rightarrow \mc{D}$ corresponds to an object $c$ of $\mc{D}(U_\alpha)$. We have $\Delta_F(c)=\big(c,c,\text{Id}\colon F(c)\rightarrow
F(c)\big)$. The $2$-commutative diagram yields an isomorphism \[\big(c,c,\text{Id}\colon F(c)\rightarrow F(c)\big)\rightarrow \big(a|_{U_\alpha},b|_{U_\alpha},p|_{U_\alpha}\colon F(a|_{U_\alpha})\rightarrow F(b|_{U_\alpha})\big).\] 
That is, there exists isomorphisms $a_{\alpha}\colon c\rightarrow a|_{U_\alpha},b_{\alpha}\colon c\rightarrow b|_{U_\alpha}$ in $\mc{D}(U_\alpha)$ satisfying the following commutative diagram,
\[\begin{tikzcd}
	F(c) \arrow[dd, "F(a_\alpha)"'] \arrow[rr, "F(\text{Id})"] & & F(c) \arrow[dd, "F(b_\alpha)"] \\
	& & \\
	F(a|_{U_\alpha}) \arrow[rr,"p|_{U_\alpha}"] & & F(b|_{U_\alpha})
\end{tikzcd}.\]
In other words, we have $p|_{U_\alpha}\circ F(a_\alpha)=F(b_\alpha)\circ F(\text{Id})=F(b_\alpha)$. As $a_\alpha$ is an isomorphism, we have 
$p|_{U_\alpha}=F(b_\alpha)\circ F(a_{\alpha}^{-1})=F(b_\alpha\circ a_{\alpha}^{-1})$; that is, $p|_{U_\alpha}\colon F(a|_{U_\alpha})\rightarrow F(b|_{U_\alpha})$ is equal to $F(\tau_\alpha)$ for some isomorphism $\tau_\alpha\colon a|_{U_\alpha}\rightarrow b|_{U_\alpha}$.

Thus, if the diagonal morphism $\Delta_F\colon \mc{D}\rightarrow\mc{D}\times_{\mc{C}}\mc{D}$ is an epimorphism, then given a manifold $U$ and an arrow $p\colon F(a)\rightarrow F(b)$ in $\mc{C}(U)$, there exists an open cover $\{U_\alpha\rightarrow U\}$ of $U$ and a family of isomorphisms $\{\tau_\alpha\colon a|_{U_\alpha}\rightarrow b|_{U_\alpha}\}$ such that $F(\tau_\alpha)=p|_{U_\alpha}$. 
Again the converse holds.
We summarize our last couple of observations in the following lemma.
\begin{lemma}\label{Lemma:equivalentnotionfordifferentiablestack}
	A morphism of stacks $F\colon \mc{D}\rightarrow \mc{C}$ is a gerbe over a stack if and only if the following two conditions hold:
	\begin{enumerate}
		\item Given a manifold $U$ and an object $a$ of $ \mc{C}(U)$, there exists an open cover $\{U_\alpha\rightarrow U\}$ of $U$ and objects $x_\alpha$ of $\mc{D}(U_\alpha)$ with an isomorphism $F(x_\alpha)\rightarrow a|_{U_\alpha}$ in $\mc{C}(U_\alpha)$ for each $\alpha$. 
		\item Given a manifold $U$ and an arrow $p\colon F(a)\rightarrow F(b)$ in $\mc{C}(U)$, there exists an open cover $\{U_\alpha\rightarrow U\}$ of $U$ and isomorphisms $\tau_\alpha\colon a|_{U_\alpha}\rightarrow b|_{U_\alpha}$ in $\mc{D}(U_\alpha)$ such that $F(\tau_\alpha)=p|_{U_\alpha}$ in $\mc{C}(U_\alpha)$ for each $\alpha$.
	\end{enumerate}
\end{lemma}

\begin{example}
	Let $X$ be a smooth manifold. Let $G$ be a Lie group acting on the smooth manifold $X$. Consider a central extension of Lie groups $1\ra S^1\ra \hat{G}\xra{\pi} G\ra 1$. Let $[X/G]$
	and $[X/\hat{G}]$ respectively be the quotient stacks  associated to the actions of $\hat{G}$ and $G$ on $X$. The morphism of Lie groups $\pi:\hat{G}\ra G$ defines a morphism of Lie groupoids $(X\times \hat{G}\rra X)\ra (X\times G\rra X)$, given by $x\mapsto x$ and $(x,\hat{g})\mapsto (x,\pi(\hat{g}))$. Then, as we will see in \Cref{Section:GoSassociatedtoALiegrupoidExtension}, this morphism of Lie groupoids $(X\times \hat{G}\rra X)\ra (X\times G\rra X)$ associates a morphism of stacks $[X/\hat{G}]\xra{\pi} [X/G]$. In fact, this morphism of stacks $[X/\hat{G}]\xra{\pi} [X/G]$ is a gerbe over the quotient stack $[X/G]$. 
\end{example}
\begin{example}
	Let $M,N$ be manifolds. Let $f:M\ra N$ be a diffeomorphism. Then, the associated morphism of stacks $F:\underline{M}\ra \underline{N}$ is a gerbe over the stack $\underline{N}$. More over, a morphism of stacks $G:\underline{M}\ra \underline{N}$ is a gerbe over the stack $\underline{N}$ implies that the associated map of manifolds $g:M\ra N$ is a diffeomorphism.
\end{example}
\begin{remark} 
	Let $\mc{D}\ra\mc{C}$ be a gerbe over the stack $\mc{C}$. When the stack $\mc{C}$ is representable by a manifold; that is $\mc{C}\cong \underline{M}$ for a manifold $M$, we recover the notion of \textit{a gerbe over a manifold $M$}. It is immediate from \Cref{Lemma:equivalentnotionfordifferentiablestack} that, a gerbe over a manifold $M$, associates a groupoid $\mc{G}(U)$ with each open set $U\subseteq M$, such that the following conditions are satisfied:
	\begin{itemize}
		\item  given $x\in M$ there is an open subset $U\subseteq M$ containing $x$ such that $\mc{G}(U)$ is non empty,
		\item  given $a,b\in\mc{G}(U)$ and $x\in U\subseteq M$, there exists an open subset $V$ of $U$ containing $x$ such that $a|_V$ is isomorphic to $b|_V$.
	\end{itemize}
	These two properties respectively, are called ``locally non empty'' and ``locally connected''. 
	For further details on this topic, we
	refer to the Section $3$ of \cite{Moerdijk2}
\end{remark}
\begin{example}
	Let $M$ be a manifold and $\mc{O}(M)$ be the category of open sets of the manifold $M$. Let $G$ be a Lie group. For an open set $U\subseteq M$, let ${\rm Tor}(G)|_U$ denote the groupoid of principal $G$ bundles over the manifold $U$. Then, the assignment $U\mapsto {\rm Tor}(G)|_U$ for the open set $U\subseteq M$ gives a gerbe over the manifold $M$.
\end{example}

\section{A Lie groupoid extension associated to a Gerbe over a stack}\label{Section:gerbegivingLiegroupoid}
Let $\pi_{\mc{D}}\colon \mc{D}\ra \text{Man}$ and $\pi_{\mc{C}}\colon \mc{C}\ra \text{Man}$ be differentiable stacks, and $F\colon \mc{D}\rightarrow \mc{C}$ a gerbe over a stack in the sense of \Cref{Definition:GerbeoverstackasinMR2817778}. Throughout this section, by a gerbe over a stack we mean an epimrorphism of stacks $F:\mc{D}\ra \mc{C}$  whose associated diagonal morphism $\Delta_F:\mc{D}\ra \mc{D}\times_{\mc{C}}\mc{D}$ is an epimorphism. Assume further that the diagonal morphism $\Delta_F\colon \mc{D}\ra \mc{D}\times_{\mc{C}}\mc{D}$ is a representable surjective submersion. In this section, we associate a Lie groupoid extension for this morphism of stacks $F:\mc{D}\ra \mc{C}$. 
The outline of the construction is as follows: 
\begin{enumerate}	\item We use the fact that $F\colon \mc{D}\ra \mc{C}$ is an epimorphism of stacks to prove that there exists an atlas $q\colon \underline{X}\ra \mc{C}$ for $\mc{C}$ and a morphism of stacks $p\colon \underline{X}\ra \mc{D}$ satisfying the following $2$-commutative diagram (\Cref{Lemma:existenceofatlasforC}),
	\[\begin{tikzcd}
		\underline{X}\arrow[dd, "p"']\arrow[rrdd, "q"{name=M}] & & \\
		& & \\
		\mc{D}\arrow[rr, "F", swap]
		\arrow[Rightarrow, shorten >=10pt, shorten <=10pt, to=M] & & \mc{C}
	\end{tikzcd}.\]
	By $2$ commutative diagram above, $F\circ p$ and $q$ can be identified up to a $2$-isomorphism.
	\item The fact that the diagonal morphism $\Delta_F\colon \mc{D}\ra \mc{D}\times_{\mc{C}}\mc{D}$ is an epimorphism of stacks implies that the morphism of stacks $p\colon \underline{X}\ra \mc{D}$ obtained in step $1$ is an epimorphism of stacks (\Cref{Lemma:pisanepimorphism}).
	\item Under the assumption that the diagonal morphism $\Delta_F\colon \mc{D}\ra \mc{D}\times_{\mc{C}}\mc{D}$ is a representable surjective submersion, we prove that, the morphism of stacks $p\colon \underline{X}\ra \mc{D}$ mentioned in step $2$ is an atlas for the stack $\pi_{\mc{D}}\colon \mc{D}\ra \text{Man}$ (\Cref{Lemma:extraassumptiononDiagonalmorphism}).
	\item For the choices made in step $(2)$ and step $(3)$ for atlases $p\colon \underline{X}\ra \mc{D}$ and $q\colon \underline{X}\ra\mc{C}$ we respectively obtain Lie groupoids $\mc{G}_p=(X\times_{\mc{D}}X\rightrightarrows X)$
	and $\mc{G}_q=(X\times_{\mc{C}}X\rightrightarrows X)$.
	We then prove that, the morphism of stacks $F\colon \mc{D}\ra \mc{C}$ gives a Lie groupoid extension $\mc{G}_p\ra \mc{G}_q$ (\Cref{Lemma:MoSgivingLiegroupoidExtension}).
	\item Finally we prove that the above construction does not depend on the choice of $q\colon \underline{X}\rightarrow \mc{C}$ (\Cref{Lemma:LiegroupoidExtnIsIndependentofq}).
\end{enumerate}

\subsection{Existence of an atlas for $\mc{C}$ with a special property} \label{Subsection:existenceofatlas}

Let $F:\mc{D}\ra \mc{C}$ be a gerbe over the stack $\mc{C}$. First we observe that the property of $F:\mc{D}\ra \mc{C}$ being an epimorphism implies existence of a special type of an atlas for $\mc{C}$.

Let $Y$ be a manifold and $\tilde{q}\colon \underline{Y}\ra \mc{C}$ an atlas for $\mc{C}$. As $F:\mc{D}\ra \mc{C}$ is an epimorphism, there exists a surjective submersion $G\colon \underline{X}\ra \underline{Y}$
and a morphism of stacks $p\colon \underline{X}\ra \mc{D}$ with the following $2$-commutative diagram,
\[\begin{tikzcd}
	\underline{X} \arrow[dd, "p"'] 
	\arrow[rr, "G"] & & \underline{Y} \arrow[dd, "\tilde{q}"] \\
	& & \\
	\mc{D} \arrow[rr, "F"] \arrow[Rightarrow, shorten >=20pt, shorten <=20pt, uurr] & & \mc{C}
\end{tikzcd}.\]
As $g\colon X\ra Y$ is a surjective submersion, the composition $q:=(\tilde{q}\circ G)\colon \underline{X}\ra \mc{C}$ is an atlas for the stack $\mc{C}$ (\Cref{Remark:compositionisatlas}). Thus, we have obtained an atlas $q\colon \underline{X}\ra \mc{C}$ for $\mc{C}$ and a morphism of stacks $p\colon \underline{X}\ra \mc{D}$ with the following $2$-commutative diagram,
\begin{equation}\label{Diagram:FpToq}
	\begin{tikzcd}
		\underline{X}\arrow[dd, "p"']\arrow[rrdd, "q"{name=M}] & & \\
		& & \\
		\mc{D}\arrow[rr, "F", swap]
		\arrow[Rightarrow, shorten >=10pt, shorten <=10pt, to=M] & & \mc{C}
	\end{tikzcd}.
\end{equation}

\begin{lemma}\label{Lemma:existenceofatlasforC} 
	Let $F\colon \mc{D}\ra \mc{C}$ be an epimorphism of stacks. Then, there exists an atlas $q\colon \underline{X}\ra \mc{C}$ for $\mc{C}$ and a morphism of stacks $p\colon \underline{X}\ra \mc{D}$ satisfying the $2$-commutative Diagram \ref{Diagram:FpToq}.
\end{lemma}
The significance of the \Cref{Lemma:existenceofatlasforC} is that in fact $p:\underline{X}\ra \mc{D}$ is an atlas for $\mc{D}$ as we see in the succeeding sections. 
\subsection{Proof that $p\colon \underline{X}\ra \mc{D}$ is an epimorphism of stacks}\label{Subsection:pIsAnEpimorphism}
As a next setp, we use the condition that $\Delta_F:\mc{D}\ra \mc{D}\times_{\mc{C}}\mc{D}$ is an epimorphism to conclude that the morphism of stacks $p:\underline{X}\ra \mc{D}$ obtained in \Cref{Lemma:existenceofatlasforC} is an epimorphism of stacks. 

Let $p\colon \underline{X}\ra \mc{D}$ and $q\colon \underline{X}\ra \mc{C}$ be as in \Cref{Lemma:existenceofatlasforC}. Let $M$ be a manifold and $r\colon \underline{M}\ra \mc{D}$ a morphism of stacks. Consider the following set up of morphism of stacks,
\[\begin{tikzcd}
	& & \underline{M} \arrow[dd, "r"] & & \\
	& & & & \\
	\underline{X} \arrow[rr, "p"] & & \mc{D} \arrow[rr, "F"] & & \mc{C}
\end{tikzcd}.\]
We prove that, there exists a surjective submersion $\underline{W}\ra \underline{M}$ and a morphism of stacks $\underline{W}\ra \underline{X}$ that fits in the $2$-commutative diagram with other sides as $p:\underline{X}\ra \mc{D}$, and $r:\underline{M}\ra \mc{D}$. 

Consider the composition morphisms $F\circ p=q:\underline{X}\ra \mc{C}$ and 
$F\circ r:\underline{M}\ra \mc{C}$. As $q\colon \underline{X}\ra \mc{C}$ is an atlas for $\mc{C}$, the $2$-fiber product $\underline{X}\times_{\mc{C}}\underline{M}$ in the following diagram, 
\begin{equation}\label{Diagram:pr2XMtoM}
	\begin{tikzcd}
		\underline{X}\times_{\mc{C}}\underline{M} \arrow[dd, "{\rm pr}_1"'] \arrow[rr, "{\rm pr}_2"] & & \underline{M} \arrow[dd, "F\circ r"] \\
		& & \\
		\underline{X} \arrow[Rightarrow, shorten >=30pt, shorten <=30pt, uurr] \arrow[rr, "F\circ p=q"] & & \mc{C} 
	\end{tikzcd},\end{equation}
 is representable by a manifold and the projection map ${\rm pr}_2\colon \underline{X}\times_{\mc{C}}\underline{M}\ra \underline{M}$ is a surjective submersion.     

Consider the morphism of stacks $(p,r)\colon \underline{X}\times_{\mc{C}}\underline{M}\ra \mc{D}\times_{\mc{C}}\mc{D}$. As the diagonal morphism $\Delta_F\colon \mc{D}\ra \mc{D}\times_{\mc{C}}\mc{D}$ is an epimorphism of stacks, there exists a surjective submersion $\Phi\colon \underline{W}\ra \underline{X}\times_{\mc{C}}\underline{M}$ and a morphism of stacks $\gamma\colon \underline{W}\ra \mc{D}$ producing the following $2$-commutative diagram,
\begin{equation}\label{Diagram:PhiWtoXCM}
	\begin{tikzcd}
		\underline{W} \arrow[dd, "\gamma"'] \arrow[rr, "\Phi"] & & \underline{X}\times_{\mc{C}}\underline{M} \arrow[dd, "{(p,r)}"] \\
		& & \\
		\mc{D} \arrow[Rightarrow, shorten >=30pt, shorten <=30pt, uurr] \arrow[rr, "\Delta_F"] & & \mc{D}\times_{\mc{C}}\mc{D}
	\end{tikzcd}.\end{equation}
Extending the Diagram \ref{Diagram:PhiWtoXCM} along the first projections ${\rm pr}_1\colon \underline{X}\times_{\mc{C}}\underline{M}\ra \underline{M}$ and 
${\rm pr}_1\colon \mc{D}\times_{\mc{C}}\mc{D}\ra \mc{D}$
we obtain the following diagram,
\begin{equation}\label{Diagram:extendalongpr1}
	\begin{tikzcd}
		\underline{W} 
		\arrow[dd, "\gamma"'] \arrow[rr, "\Phi"] & & \underline{X}\times_{\mc{C}}\underline{M}
		\arrow[rr, "{\rm pr}_1"] \arrow[dd, "{(p,r)}"] & & \underline{X} \arrow[dd, "p"] \\
		& & & & \\
		\mc{D} \arrow[Rightarrow, shorten >=30pt, shorten <=30pt, uurr] \arrow[rr, "\Delta_F"] & & \mc{D}\times_{\mc{C}}\mc{D} \arrow[Rightarrow, shorten >=30pt, shorten <=30pt, uurr]
		\arrow[rr, "{\rm pr}_1"] & & \mc{D} 
	\end{tikzcd}.
\end{equation}
Similarly, extending the Diagram \ref{Diagram:PhiWtoXCM} along the second projections ${\rm pr}_2\colon \underline{X}\times_{\mc{C}}\underline{M}\ra \underline{M}$ and 
${\rm pr}_2\colon \mc{D}\times_{\mc{C}}\mc{D}\ra \mc{D}$
we obtain the following diagram,
\begin{equation}\label{Diagram:extendalongpr2} 
	\begin{tikzcd}
		\underline{W} \arrow[dd, "\gamma"'] \arrow[rr, "\Phi"] & & \underline{X}\times_{\mc{C}}\underline{M} \arrow[rr, "{\rm pr}_2"] \arrow[dd, "{(p,r)}"] & & \underline{M} \arrow[dd, "r"] \\
		& & & & \\
		\mc{D} \arrow[Rightarrow, shorten >=30pt, shorten <=30pt, uurr] \arrow[rr, "\Delta_F"] & & \mc{D}\times_{\mc{C}}\mc{D} \arrow[Rightarrow, shorten >=30pt, shorten <=30pt, uurr] \arrow[rr, "{\rm pr}_2"] & & \mc{D} 
	\end{tikzcd}.
\end{equation}
For the convenience  of visualization we combine Diagrams \ref{Diagram:extendalongpr1}, and \ref{Diagram:extendalongpr2} to draw the following diagram,
\begin{equation}\label{Diagram:WtoXMtoMandX}\begin{tikzcd}
		\underline{W} \arrow[dd, "\gamma"'] \arrow[rr, "\Phi"] & & \underline{X}\times_{\mc{C}}\underline{M} \arrow[dd, "{(p,r)}"] \arrow[rr, "{\rm pr}_2"] \arrow[rrrd,"{\rm pr}_1"'] & & \underline{M} \arrow[dd, "r"] & \\
		& & & & & \underline{X} \arrow[dd, "p"] \\
		\mc{D} \arrow[rr, "\Delta_F"] & & \mc{D}\times_{\mc{C}}\mc{D} \arrow[rr, "{\rm pr}_2"] \arrow[rrrd, "{\rm pr}_1"'] & & \mc{D} & \\
		&& & & & \mc{D} 
	\end{tikzcd}.\end{equation}
Observe that the maps ${\rm pr}_2\colon \underline{X}\times_{\mc{C}}\underline{M}\ra \underline{M}$ (from Diagram \ref{Diagram:pr2XMtoM}) and $\Phi\colon \underline{W}\ra \underline{X}\times_{\mc{C}}\underline{M}$ (from Diagram \ref{Diagram:PhiWtoXCM}) are surjective submersions. Thus, the composition 
${\rm pr}_2\circ \Phi\colon \underline{W}\ra \underline{M}$ is a surjective submersion. The composition 
${\rm pr}_1\circ \Phi\colon \underline{W}\ra \underline{X}$  gives the following diagram of morphism of stacks,
\begin{equation}\label{Diagram:CommutativeWMXD}\begin{tikzcd}
		\underline{W} \arrow[dd, "{\rm pr}_1\circ \Phi"'] \arrow[rr, "{\rm pr}_2\circ \Phi"] & & \underline{M} \arrow[dd, "r"] \\
		& & \\
		\underline{X} \arrow[rr, "p"] & & \mc{D} 
	\end{tikzcd}.\end{equation}
We further note from Diagram \Cref{Diagram:WtoXMtoMandX} that, $r\circ {\rm pr}_2\circ \Phi={\rm pr}_2\circ \Delta_F\circ \gamma$ and  
$p\circ {\rm pr}_1\circ \Phi={\rm pr}_1\circ \Delta_F\circ \gamma$. 
As ${\rm pr}_1\circ \Delta={\rm pr}_2\circ \Delta$, we see that 
${\rm pr}_1\circ \Delta\circ \gamma={\rm pr}_2\circ \Delta\circ \gamma$.
Thus, $r\circ {\rm pr}_2\circ \Phi=p\circ {\rm pr}_1\circ \Phi$. So, the Diagram \Cref{Diagram:CommutativeWMXD} is a $2$-commutative diagram. Thus, given a morphism of stacks $r\colon \underline{M}\ra \mc{D}$, there exists a surjective submersion $\Gamma={\rm pr}_2\circ \Phi\colon \underline{W}\ra \underline{M}$ and a morphism of stacks $\Psi={\rm pr}_1\circ \Phi\colon \underline{W}\ra \underline{X}$ with following $2$-commutative diagram,
\begin{equation}\begin{tikzcd}
		\underline{W} \arrow[dd, "\Psi"'] \arrow[rr, "\Gamma"] & & \underline{M} \arrow[dd, "r"] \\
		& & \\
		\underline{X} \arrow[Rightarrow, shorten >=20pt, shorten <=20pt, uurr] \arrow[rr, "p"] & & \mc{D} 
	\end{tikzcd}.\end{equation}
Thus, we conclude that the morphism $p\colon \underline{X}\ra\mc{D}$ is an epimorphism of stacks.
\begin{lemma}\label{Lemma:pisanepimorphism} 
	The morphism of stacks $p\colon \underline{X}\ra \mc{D}$ mentioned in Diagram \Cref{Diagram:FpToq} is an epimorphism of stacks.
\end{lemma}
\subsection{$p\colon \underline{X}\ra \mc{D}$ is an atlas for $\mc{D}$}
\label{Subsection:pIsAnAtlas} Let $F\colon \mc{D}\ra \mc{C}$ be a gerbe over a stack. We further assume that the diagonal morphism $\Delta_F\colon \mc{D}\ra \mc{D}\times_{\mc{C}}\mc{D}$ is a representable surjective submersion. Let $p\colon \underline{X}\ra \mc{D}$ be as in \Cref{Lemma:pisanepimorphism}. In general $p\colon \underline{X}\ra \mc{D}$ is not an atlas for $\mc{D}$. We have proved that $p\colon \underline{X}\ra \mc{D}$ is an epimorphism of stacks. We use the following result to conclude that $p\colon \underline{X}\ra \mc{D}$ is in fact an atlas for $\mc{D}$.
The following proposition is a variant of Proposition $2.16$ in \cite{MR2817778}.
\begin{proposition}\label{Proposition:epimorphismbeinganatlas}
	A morphism of stacks $r\colon \underline{X}\rightarrow \mc{D}$ is an atlas for $\mc{D}$ if
	\begin{enumerate}
		\item the morphism $r\colon \underline{X}\rightarrow \mc{D}$ is an epimorphism of stacks, 
		\item the fibered product $\underline{X}\times_{\mc{D}}\underline{X}$ is representable by a manifold and that the projection maps
		${\rm pr}_1\colon X\times_{\mc{D}}X\rightarrow X$ and
		${\rm pr}_2\colon X\times_{\mc{D}}X\rightarrow X$ are submersions.
	\end{enumerate}
\end{proposition}
By above Proposition, to prove $p\colon \underline{X}\ra \mc{D}$ is an atlas for $\mc{D}$, it only remains to prove that $\underline{X}\times_{\mc{D}}\underline{X}$ is representable by a manifold and that the projection maps ${\rm pr}_1\colon X\times_{\mc{D}}X\ra X$ and 
${\rm pr}_2\colon X\times_{\mc{D}}X\ra X$ are submersions. We prove them below.

For $p\colon \underline{X}\ra \mc{D}$ and for $F\circ p\colon \underline{X}\ra \mc{C}$, we have following pull
back diagrams,
\[\begin{tikzcd}
	\underline{X}\times_{\mc{D}}\underline{X} \arrow[dd, "{\rm pr}_1^{\mc{D}}"'] \arrow[rr, "{\rm pr}_2^{\mc{D}}"] & & \underline{X} \arrow[dd, "p"] & \\
	& & & \\
	\underline{X} \arrow[Rightarrow, shorten >=30pt, shorten <=30pt, uurr] \arrow[rr, "p"] & & \mc{D} \arrow[rr, "F"] & & \mc{C}
\end{tikzcd} \begin{tikzcd}
	\underline{X}\times_{\mc{C}}\underline{X} \arrow[dd, "{\rm pr}_1^{\mc{C}}"'] \arrow[rr, "{\rm pr}_2^{\mc{C}}"] & & \underline{X} \arrow[dd, "F\circ p"] \\
	& & \\
	\underline{X} \arrow[Rightarrow, shorten >=30pt, shorten <=30pt, uurr] \arrow[rr, "F\circ p"] & & \mc{C} 
\end{tikzcd}.\]
By uniqueness of pull-back, there exists a unique morphism of stacks 
$\Psi\colon \underline{X}\times_{\mc{D}}\underline{X}\ra 
\underline{X}\times_{\mc{C}}\underline{X}$ with following $2$-commutative diagram,
\begin{equation}\label{Diagram:ExistenceofPsi}\begin{tikzcd}
		\underline{X}\times_{\mc{D}}\underline{X} \arrow[rddd, "{\rm pr}_1^{\mc{D}}"', bend right] \arrow[rrrd, "{\rm pr}_2^{\mc{D}}", bend left] \arrow[rd, "\Psi"] & & & \\
		& \underline{X}\times_{\mc{C}}\underline{X} \arrow[dd, "{\rm pr}_1^{\mc{C}}"'] \arrow[rr, "{\rm pr}_2^{\mc{C}}"] & & \underline{X} \arrow[dd, "F\circ p"] \\
		& & & \\
		& \underline{X} \arrow[rr, "F\circ p"] & & \mc{C} 
	\end{tikzcd}.\end{equation}
We have assumed that the diagonal morphism $\Delta_F\colon \mc{D}\ra \mc{D}\times_{\mc{C}}\mc{D}$ is a representable surjective submersion. Consider the morphism of stacks $(p,p)\colon \underline{X}\times_{\mc{C}}\underline{X}\ra \mc{D}\times_{\mc{C}}\mc{D}$. We have the following $2$-fiber product diagram,
\begin{equation}\label{Diagram:(p,p)diagonalmorphism}
	\begin{tikzcd}
		\mc{D}\times_{\mc{D}\times_{\mc{C}}\mc{D}}(\underline{X}\times_{\mc{C}}\underline{X}) \arrow[dd, "{\rm pr}_1"'] \arrow[rr, "{\rm pr}_2"] & & \underline{X}\times_{\mc{C}}\underline{X} \arrow[dd, "{(p,p)}"] \\
		& & \\
		\mc{D} \arrow[Rightarrow, shorten >=30pt, shorten <=30pt, uurr] \arrow[rr, "\Delta_F"] & & \mc{D}\times_{\mc{C}}\mc{D} 
	\end{tikzcd}.\end{equation}
We obtain an isomorphism of stacks 
$\mc{D}\times_{\mc{D}\times_{\mc{C}}\mc{D}}(\underline{X}\times_{\mc{C}}\underline{X})\cong\underline{X}\times_{\mc{D}}\underline{X}$ (\cite[Corollary $69$]{metzler2003topological}). Observe that the morphism of stacks ${\rm pr}_2\colon \mc{D}\times_{\mc{D}\times_{\mc{C}}\mc{D}}
(\underline{X}\times_{\mc{C}}\underline{X})\ra \underline{X}\times_{\mc{C}}\underline{X}$ in Diagram \ref{Diagram:(p,p)diagonalmorphism} is same as the morphism of stacks $\Psi\colon \underline{X}\times_{\mc{D}}\underline{X}\ra 
\underline{X}\times_{\mc{C}}\underline{X}$ in Diagram \ref{Diagram:ExistenceofPsi}.
Thus, the above $2$-fiber product diagram can be seen as the following $2$-commutative diagram,
\begin{equation}\label{Diagram:Psi(p,p)DeltaF}
	\begin{tikzcd}
		\underline{X}\times_{\mc{D}}\underline{X} \arrow[dd] \arrow[rr, "\Psi"] & & \underline{X}\times_{\mc{C}}\underline{X} \arrow[dd, "{(p,p)}"] \\
		& & \\
		\mc{D} \arrow[Rightarrow, shorten >=30pt, shorten <=30pt, uurr] \arrow[rr, "\Delta_F"] & & \mc{D}\times_{\mc{C}}\mc{D} 
	\end{tikzcd}.\end{equation} 
As the diagonal morphism $\Delta_F\colon \mc{D}\ra \mc{D}\times_{\mc{C}}\mc{D}$ is a representable surjective submersion, the $2$-fiber product $\mc{D}\times_{\mc{D}\times_{\mc{C}}\mc{D}}(\underline{X}\times_{\mc{C}}\underline{X})$ is representable by a manifold and the projection map ${\rm pr}_2\colon \mc{D}\times_{\mc{D}\times_{\mc{C}}\mc{D}}(\underline{X}\times_{\mc{C}}\underline{X})\ra \underline{X}\times_{\mc{C}}\underline{X}$ is a surjective submersion at the level of manifolds. Thus, we see that $\underline{X}\times_{\mc{D}}\underline{X}$ is representable by manifold and $\Psi\colon \underline{X}\times_{\mc{D}}\underline{X}\ra 
\underline{X}\times_{\mc{C}}\underline{X}$ is a surjective submersion at the level of manifolds.

Both being compositions of surjective submersions, we see that ${\rm pr}_1^{\mc{D}}={\rm pr}_1^{\mc{C}}\circ \Psi$ and 
${\rm pr}_2^{\mc{D}}={\rm pr}_2^{\mc{C}}\circ \Psi$ are surjective submersions at the level of manifolds as well. Thus, $p\colon \underline{X}\ra \mc{D}$ is an atlas for $\mc{D}$. So, we have shown the following:
\begin{lemma}\label{Lemma:extraassumptiononDiagonalmorphism}
	Let $F\colon \mc{D}\ra \mc{C}$ be a gerbe over a stack. Further assume that, the diagonal morphism $\Delta_F\colon \mc{D}\ra \mc{D}\times_{\mc{C}}\mc{D}$ is a representable surjective submersion. Then, the morphism of stacks $p\colon \underline{X}\ra \mc{D}$ mentioned in Diagram \ref{Diagram:FpToq} is an atlas for the stack $\pi_{\mc{D}}\colon \mc{D}\ra \text{Man}$. In particular, there exists an atlas $p\colon \underline{X}\ra \mc{D}$ for $\mc{D}$ and an atlas $q\colon \underline{X}\ra \mc{C}$ for $\mc{C}$ with a $2$-commutative diagram as in Diagram \ref{Diagram:FpToq}.
\end{lemma}

\subsection{A gerbe over a stack gives a Lie groupoid extension}
\label{Subsection:GoSgivingLiegroupoidextension}
Let $F\colon \mc{D}\ra \mc{C}$ be a gerbe over a stack. We further assume that the diagonal morphism $\Delta_F\colon \mc{D}\ra \mc{D}\times_{\mc{C}}\mc{D}$ is a representable surjective submersion.
 By \Cref{Lemma:extraassumptiononDiagonalmorphism}, there exists an atlas $p\colon \underline{X}\ra \mc{D}$ for $\mc{D}$ and an atlas $q\colon \underline{X}\ra \mc{C}$ with a $2$-commutative diagram as in Diagram \ref{Diagram:FpToq}. For atlases $p\colon \underline{X}\ra\mc{D}$ and $q\colon \underline{X}\ra\mc{C}$ for $\mc{C}$, we have respectively associated the Lie groupoids $\mc{G}_p=(X\times_{\mc{D}}X\rightrightarrows X)$ and $\mc{G}_q=(X\times_{\mc{C}}X\rightrightarrows Y)$. The morphism of stacks $F\colon \mc{D}\ra \mc{C}$ induces a morphism of stacks $\Psi\colon \underline{X}\times_{\mc{D}}\underline{X}\ra \underline{X}\times_{\mc{C}}\underline{X}$ as in Diagram Diagram \ref{Diagram:ExistenceofPsi}.
Explicitly, $\Psi\colon \underline{X}\times_{\mc{D}}\underline{X}\ra \underline{X}\times_{\mc{C}}\underline{X}$ is given at the level of objects by 
\[(y,z,\alpha\colon p(y)\rightarrow p(z))\mapsto (y,z,F(\alpha)\colon F(p(y))\rightarrow F(p(z))).\]
At the level of morphisms, an arrow \[(u,v)\colon \big(y,z,\alpha\colon p(y)\rightarrow p(z))\ra \big(y',z',\alpha'\colon p(y')\rightarrow p(z')\big)\]
in $\underline{X}\times_{\mc{D}}\underline{X}$ is mapped to the arrow \[(u,v)\colon \big(y,z,F(\alpha)\colon F(p(y))\rightarrow F(p(z)))\ra \big(y',z',F(\alpha')\colon F(p(y'))\rightarrow F(p(z'))\big)\]
in $\underline{X}\times_{\mc{C}}\underline{X}$.
This morphism of stacks $\Psi\colon \underline{X}\times_{\mc{D}}\underline{X}\rightarrow
\underline{X}\times_{\mc{C}}\underline{X}$ is compatible with projection morphisms ${\rm pr}_1,{\rm pr}_2\colon \underline{X}\times_{\mc{D}}\underline{X}
\rightarrow \underline{X}$ and ${\rm pr}_1,{\rm pr}_2\colon \underline{X}\times_{\mc{C}}\underline{X}
\rightarrow \underline{X}$ in the sense that following diagrams are commutative diagram of morphisms of stacks, 
\begin{equation}\label{Diagram:PsiXDXtoXCX}\begin{tikzcd}
		\underline{X}\times_{\mc{D}}\underline{X} \arrow[dd,xshift=0.75ex,"{\rm pr}_2"]\arrow[dd,xshift=-0.75ex,"{\rm pr}_1"'] \arrow[rr,"\Psi"] & & \underline{X}\times_{\mc{C}}\underline{X} \arrow[dd,xshift=0.75ex,"{\rm pr}_2"]\arrow[dd,xshift=-0.75ex,"{\rm pr}_1"'] \\
		& & \\
		\underline{X} \arrow[rr,"\text{Id}"] & & \underline{X}
	\end{tikzcd}.\end{equation}
Recall from  the \Cref{Subsection:Liegroupoidassociatedtoatlas} that, the morphisms of stacks ${\rm pr}_1\colon \underline{X}\times_{\mc{D}}\underline{X}\ra \underline{X}$ and ${\rm pr}_2\colon \underline{X}\times_{\mc{D}}\underline{X}\ra \underline{X}$ respectively corresponds to the source and target maps of the Lie groupoid $\mc{G}_p=(X\times_{\mc{D}}X\rightrightarrows X)$. Likewise for the Lie groupoid $\mc{G}_q=(X\times_{\mc{C}}X\rightrightarrows X)$. Let $\Theta\colon X\times_{\mc{D}}X\ra X\times_{\mc{C}}X$ be the map of manifolds associated to the morphism of stacks $\Psi\colon \underline{X}\times_{\mc{D}}\underline{X}\ra \underline{X}\times_{\mc{C}}\underline{X}$.
Then, the Diagram \ref{Diagram:PsiXDXtoXCX}, gives the following diagram of morphism of Lie groupoids,
\begin{equation}\label{Diagram:ThetaXDXtoXCX}\begin{tikzcd}
		X\times_{\mc{D}}X \arrow[dd,xshift=0.75ex,"t"]\arrow[dd,xshift=-0.75ex,"s"'] \arrow[rr,"\Theta"] & & X\times_{\mc{C}}X \arrow[dd,xshift=0.75ex,"t"]\arrow[dd,xshift=-0.75ex,"s"'] \\
		& & \\
		X \arrow[rr,"\text{Id}"] & & X
	\end{tikzcd}.\end{equation}
We have the following result.
\begin{proposition}\label{Proposition:GoSgivingMorphismofLiegroupoids}
	Let $F\colon\mc{D}\ra \mc{C}$ be a gerbe over a stack. Assume further that the diagonal morphism $\Delta_F\colon \mc{D}\rightarrow \mc{D}\times_{\mc{C}}\mc{D}$ is a representable surjective submersion. Then the morphism of stacks $F\colon \mc{D}\rightarrow \mc{C}$ gives a morphism of Lie groupoids \[\mc{G}_p\ra \mc{G}_q\colon (X\times_{\mc{D}}X\rightrightarrows X)\rightarrow(X\times_{\mc{C}}X\rightrightarrows X),\]
	where $p\colon \underline{X}\ra \mc{D}$ and $q\colon \underline{X}\ra \mc{C}$ are as in \Cref{Lemma:extraassumptiononDiagonalmorphism}. 
\end{proposition}
Observe that the morphism of stacks $\Psi\colon \underline{X}\times_{\mc{D}}\underline{X}
\rightarrow
\underline{X}\times_{\mc{C}}\underline{X}$ is a surjective submersion at the level of manifolds (discussion after Diagram \ref{Diagram:Psi(p,p)DeltaF}); that is, the map $\Theta\colon X\times_{\mc{D}}X\ra X\times_{\mc{C}}X$ is a surjective submersion. Thus, $(\Theta,\text{Id})\colon (X\times_{\mc{D}}X\rightrightarrows X)\rightarrow (X\times_{\mc{C}}X\rightrightarrows X)$ (Diagram \ref{Diagram:ThetaXDXtoXCX}) is a Lie groupoid extension.
\begin{lemma}\label{Lemma:MoSgivingLiegroupoidExtension}
	Let $F\colon \mc{D}\ra \mc{C}$ be a gerbe over a stack. Assume  that the diagonal morphism $\Delta_F\colon \mc{D}\rightarrow \mc{D}\times_{\mc{C}}\mc{D}$ is a representable surjective submersion. Then the morphism of stacks $F\colon \mc{D}\rightarrow \mc{C}$ gives a Lie groupoid extension \[\mc{G}_p\ra \mc{G}_q\colon (X\times_{\mc{D}}X\rightrightarrows X)\rightarrow(X\times_{\mc{C}}X\rightrightarrows X)\]
	where $p\colon \underline{X}\ra \mc{D}$ and $q\colon \underline{X}\ra \mc{C}$ are as in \Cref{Lemma:extraassumptiononDiagonalmorphism}.
\end{lemma}
\subsection{Uniqueness of a Lie groupoid extension associated to a gerbe over a stack}\label{Subsection:uniquenessOfExtension} Let $F\colon \mc{D}\rightarrow \mc{C}$ be a gerbe over a stack. We further assume that $\Delta_F\colon \mc{D}\rightarrow \mc{D}\times_{\mc{C}}\mc{D}$ is a representable surjective submersion. For atlases $p\colon \underline{X}\ra \mc{D}$ and $q\colon \underline{X}\ra \mc{C}$ mentioned in \Cref{Lemma:extraassumptiononDiagonalmorphism}, we have assigned a Lie groupoid extension 
\[ (X\times_{\mc{D}}X\rightrightarrows X)\rightarrow(X\times_{\mc{C}}X\rightrightarrows X)\] in  \Cref{Lemma:MoSgivingLiegroupoidExtension}. We prove in this subsection that,  up to a Morita equivalence, this Lie groupoid extension does not depend on the choice of atlases.

Let $q_Y\colon \underline{Y}\rightarrow \mc{C}$ be another atlas for $\mc{C}$ and $p_Y\colon \underline{Y}\rightarrow \mc{D}$ be the corresponding atlas for $\mc{D}$ as in \Cref{Lemma:extraassumptiononDiagonalmorphism}. These atlases $p_Y\colon \underline{Y}\ra \mc{D},q_Y\colon \underline{Y}\ra\mc{C}$ gives a Lie groupoid extension \[(Y\times_{\mc{D}}Y\rightrightarrows Y)\rightarrow 
(Y\times_{\mc{C}}Y\rightrightarrows Y)\] as in \Cref{Lemma:MoSgivingLiegroupoidExtension}.

We prove that 
$(X\times_{\mc{D}}X\rightrightarrows X)\rightarrow 
(X\times_{\mc{C}}X\rightrightarrows X)$ and $(Y\times_{\mc{D}}Y\rightrightarrows Y)\rightarrow 
(Y\times_{\mc{C}}Y\rightrightarrows Y)$ are Morita equivalent Lie groupoid extensions. 
It means that  there exists a Lie groupoid extension $(\mc{G}_1\rightrightarrows \mc{G}_0)\rightarrow 
(\mc{H}_1\rightrightarrows \mc{H}_0)$ and a pair of Morita morphisms of Lie groupoid extensions 
\begin{align*}
	\big((\mc{G}_1\rightrightarrows \mc{G}_0)\rightarrow (\mc{H}_1\rightrightarrows \mc{G}_0)\big)&\ra \big((X\times_{\mc{D}}X\rightrightarrows X)\rightarrow 
	(X\times_{\mc{C}}X\rightrightarrows X)\big),\\
	\big((\mc{G}_1\rightrightarrows \mc{G}_0)\rightarrow (\mc{H}_1\rightrightarrows \mc{G}_0)\big)&\ra \big((Y\times_{\mc{D}}Y\rightrightarrows Y)\rightarrow 
	(Y\times_{\mc{C}}Y\rightrightarrows Y)\big).
\end{align*}

Diagrammatically that means we have to find a Lie groupoid extension and a pair of Morita morphisms of Lie groupoid extensions, 
\begin{equation}\text{ from }\begin{tikzcd}
		\mc{G}_1 \arrow[dd,xshift=0.75ex,"t"]\arrow[dd,xshift=-0.75ex,"s"'] \arrow[rr] & & \mc{H}_1 \arrow[dd,xshift=0.75ex,"t"]\arrow[dd,xshift=-0.75ex,"s"'] \\
		& & \\
		\mc{G}_0 \arrow[rr] & & \mc{G}_0
	\end{tikzcd} ~ \text{to} \begin{tikzcd}
		X\times_{\mc{D}}X \arrow[dd,xshift=0.75ex,"t"]\arrow[dd,xshift=-0.75ex,"s"'] \arrow[rr] & & X\times_{\mc{C}}X \arrow[dd,xshift=0.75ex,"t"]\arrow[dd,xshift=-0.75ex,"s"'] \\
		& & \\
		X \arrow[rr] & & X
\end{tikzcd},\end{equation}
and 
\begin{equation} \text{ from }\begin{tikzcd}
		\mc{G}_1 \arrow[dd,xshift=0.75ex,"t"]\arrow[dd,xshift=-0.75ex,"s"'] \arrow[rr] & & \mc{H}_1 \arrow[dd,xshift=0.75ex,"t"]\arrow[dd,xshift=-0.75ex,"s"'] \\
		& & \\
		\mc{G}_0 \arrow[rr] & & \mc{G}_0
	\end{tikzcd} ~ \text{to} \begin{tikzcd}
		Y\times_{\mc{D}}Y \arrow[dd,xshift=0.75ex,"t"]\arrow[dd,xshift=-0.75ex,"s"'] \arrow[rr] & & Y\times_{\mc{C}}Y \arrow[dd,xshift=0.75ex,"t"]\arrow[dd,xshift=-0.75ex,"s"'] \\
		& & \\
		Y \arrow[rr] & & Y
	\end{tikzcd}.\end{equation}
For this, first, we need a smooth manifold $\mc{G}_0$ and a pair of smooth maps $\mc{G}_0\rightarrow X, \mc{G}_0\rightarrow Y$. 

For the morphisms of stacks $p_X\colon \underline{X}\rightarrow \mc{D}$ and $p_Y\colon \underline{Y}\rightarrow \mc{D}$, consider the following $2$-fiber product diagram,
\[\begin{tikzcd}
	\underline{X}\times_{\mc{D}}\underline{Y} \arrow[dd, "{\rm pr}_1"'] \arrow[rr, "{\rm pr}_2"] & & \underline{Y} \arrow[dd, "p_Y"] \\
	& & \\
	\underline{X} \arrow[Rightarrow, shorten >=20pt, shorten <=20pt, uurr] \arrow[rr, "p_X"'] & & \mc{D} 
\end{tikzcd}.\]

As $p_X\colon \underline{X}\rightarrow \mc{D}$ is an atlas, the morphism of stacks ${\rm pr}_2\colon \underline{X}\times_{\mc{D}}\underline{Y}\rightarrow \underline{Y}$ induces a surjective submersion, 
$g\colon X\times_{\mc{D}}Y\rightarrow Y$ at the level of manifolds.
Similarly 
${\rm pr}_1\colon \underline{X}\times_{\mc{D}}\underline{Y}\rightarrow \underline{X}$ induces a surjective submersion $f\colon X\times_{\mc{D}}Y\rightarrow X$ at the level of manifolds. 

We now construct a Lie groupoid extension 
of the form $(**\rightrightarrows X\times_{\mc{D}}Y)\rightarrow (***\rightrightarrows X\times_{\mc{D}}Y)$.
In order to describe an overall outline of the construction, for the time being we denote the respective morphism sets by $**$ and $***$.

Next, we find a Morita morphisms of Lie groupoid extensions, 
\begin{equation}\text{ from }
	\begin{tikzcd}
		** \arrow[dd,xshift=0.75ex,"t"]\arrow[dd,xshift=-0.75ex,"s"'] \arrow[rr] & & ***\arrow[dd,xshift=0.75ex,"t"]\arrow[dd,xshift=-0.75ex,"s"'] \\
		& & \\
		X\times_{\mc{D}}Y \arrow[rr] & & X\times_{\mc{D}}Y
	\end{tikzcd}\text{to} \begin{tikzcd}
		X\times_{\mc{D}}X \arrow[dd,xshift=0.75ex,"t"]\arrow[dd,xshift=-0.75ex,"s"'] \arrow[rr] & & X\times_{\mc{C}}X \arrow[dd,xshift=0.75ex,"t"]\arrow[dd,xshift=-0.75ex,"s"'] \\
		& & \\
		X \arrow[rr] & & X
	\end{tikzcd} 
\end{equation} 
and 
\begin{equation} \text{ from } \begin{tikzcd}
		** \arrow[dd,xshift=0.75ex,"t"]\arrow[dd,xshift=-0.75ex,"s"'] \arrow[rr] & & *** \arrow[dd,xshift=0.75ex,"t"]\arrow[dd,xshift=-0.75ex,"s"'] \\
		& & \\
		X\times_{\mc{D}}Y \arrow[rr] & & X\times_{\mc{D}}Y
	\end{tikzcd}\text{to} \begin{tikzcd}
		Y\times_{\mc{D}}Y \arrow[dd,xshift=0.75ex,"t"]\arrow[dd,xshift=-0.75ex,"s"'] \arrow[rr] & & Y\times_{\mc{C}}Y \arrow[dd,xshift=0.75ex,"t"]\arrow[dd,xshift=-0.75ex,"s"'] \\
		& & \\
		Y \arrow[rr] & & Y
\end{tikzcd}\end{equation} 
That means, we need a pair of Morita morphisms of Lie groupoids 
\begin{equation} \begin{tikzcd}
		** \arrow[dd,xshift=0.75ex,"t"]\arrow[dd,xshift=-0.75ex,"s"'] \arrow[rr] & & X\times_{\mc{D}}X \arrow[dd,xshift=0.75ex,"t"]\arrow[dd,xshift=-0.75ex,"s"'] \\
		& & \\
		X\times_{\mc{D}}Y \arrow[rr,"f"] & & X
\end{tikzcd}\end{equation} and \begin{equation}\begin{tikzcd}
		*** \arrow[dd,xshift=0.75ex,"t"]\arrow[dd,xshift=-0.75ex,"s"'] \arrow[rr] & & X\times_{\mc{C}}X \arrow[dd,xshift=0.75ex,"t"]\arrow[dd,xshift=-0.75ex,"s"'] \\
		& & \\
		X\times_{\mc{D}}Y \arrow[rr,"f"] & & X
\end{tikzcd}\end{equation}
which are compatible with maps $**\rightarrow ***$ and $ X\times_{\mc{D}}X\rightarrow X\times_{\mc{C}}X$. 
Similarly, we need a pair of Morita morphisms of Lie groupoids 
\begin{equation}
	\begin{tikzcd}
		** \arrow[dd,xshift=0.75ex,"t"]\arrow[dd,xshift=-0.75ex,"s"'] \arrow[rr] & & Y\times_{\mc{D}}Y \arrow[dd,xshift=0.75ex,"t"]\arrow[dd,xshift=-0.75ex,"s"'] \\
		& & \\
		X\times_{\mc{D}}Y \arrow[rr,"g"] & & Y
\end{tikzcd}\end{equation}
and
\begin{equation}
	\begin{tikzcd}
		*** \arrow[dd,xshift=0.75ex,"t"]\arrow[dd,xshift=-0.75ex,"s"'] \arrow[rr] & & Y\times_{\mc{C}}Y \arrow[dd,xshift=0.75ex,"t"]\arrow[dd,xshift=-0.75ex,"s"'] \\
		& & \\
		X\times_{\mc{D}}Y \arrow[rr,"g"] & & Y
\end{tikzcd}\end{equation} 
compatible with maps $**\rightarrow ***$ and $ Y\times_{\mc{D}}Y\rightarrow Y\times_{\mc{C}}Y$.

Our first task is to find $**$ and $***$. 
Recalling the set up of Morita morphisms of Lie groupoids (\Cref{Definition:MoritamorphismofLiegroupoids}), we see that given a surjective submersion $f\colon M\rightarrow N$ and a Lie groupoid $\mc{G}_1\rightrightarrows N$, the pull-back Lie groupoid as in below diagram,
\begin{equation} 
	\begin{tikzcd}
		f^*\mc{G}_1 \arrow[dd,xshift=0.75ex,"t"]\arrow[dd,xshift=-0.75ex,"s"'] \arrow[rr] & & \mc{G}_1 \arrow[dd,xshift=0.75ex,"t"]\arrow[dd,xshift=-0.75ex,"s"'] \\
		& & \\
		M \arrow[rr,"f"] & & N
\end{tikzcd},\end{equation}
gives a Morita morphism of Lie groupoids $(f^*\mc{G}_1\rightrightarrows M)\ra (\mc{G}_1\rightrightarrows N)$.

Let $(\mc{G}_1 \rightrightarrows X\times_{\mc{D}}Y)$ be the pull-back of the Lie groupoid $(X\times_{\mc{D}}X\rightrightarrows X)$ along $f\colon X\times_{\mc{D}}Y\rightarrow X$ and $(\mc{G}_1'\rightrightarrows X\times_{\mc{D}}Y) $ be the pull-back of the Lie groupoid $(Y\times_{\mc{D}}Y\rightrightarrows Y)$ along $g\colon X\times_{\mc{D}}Y\rightarrow Y$. 
However, as we will shortly see, we do not have to distinguish between these two pull-backs as they are isomorphic.
We have following diagrams representing the pull-back groupoids,
\begin{equation} \begin{tikzcd}
		\mc{G}_1 \arrow[dd,xshift=0.75ex,"t"]\arrow[dd,xshift=-0.75ex,"s"'] \arrow[rr] & & X\times_{\mc{D}}X \arrow[dd,xshift=0.75ex,"t"]\arrow[dd,xshift=-0.75ex,"s"'] \\
		& & \\
		X\times_{\mc{D}}Y \arrow[rr,"f"] & & X
	\end{tikzcd}
	\begin{tikzcd}
		\mc{G}_1' \arrow[dd,xshift=0.75ex,"t"]\arrow[dd,xshift=-0.75ex,"s"'] \arrow[rr] & & Y\times_{\mc{D}}Y \arrow[dd,xshift=0.75ex,"t"]\arrow[dd,xshift=-0.75ex,"s"'] \\
		& & \\
		X\times_{\mc{D}}Y \arrow[rr,"g"] & & Y
	\end{tikzcd}.
\end{equation} 
Similarly, we write $(\mc{H}_1 \rightrightarrows X\times_{\mc{D}}Y) $ for pull-back of the Lie groupoid $(X\times_{\mc{C}}X\rightrightarrows X)$ along $f\colon X\times_{\mc{D}}Y\rightarrow X$ and $(\mc{H}_1'\rightrightarrows X\times_{\mc{D}}Y) $ for pull-back of the Lie groupoid $(Y\times_{\mc{C}}Y\rightrightarrows Y)$ along $g\colon X\times_{\mc{D}}Y\rightarrow Y$. As before, the Lie groupoids 
$\mc{H}_1 \rightrightarrows X\times_{\mc{D}}Y $ and $\mc{H}_1'\rightrightarrows X\times_{\mc{D}}Y $ will be isomorphic. We have following diagrams representing the pull-back groupoids,
\begin{equation} \begin{tikzcd}
		\mc{H}_1 \arrow[dd,xshift=0.75ex,"t"]\arrow[dd,xshift=-0.75ex,"s"'] \arrow[rr] & & X\times_{\mc{C}}X \arrow[dd,xshift=0.75ex,"t"]\arrow[dd,xshift=-0.75ex,"s"'] \\
		& & \\
		X\times_{\mc{D}}Y \arrow[rr,"f"] & & X
	\end{tikzcd}
	\begin{tikzcd}
		\mc{H}_1' \arrow[dd,xshift=0.75ex,"t"]\arrow[dd,xshift=-0.75ex,"s"'] \arrow[rr] & & Y\times_{\mc{C}}Y \arrow[dd,xshift=0.75ex,"t"]\arrow[dd,xshift=-0.75ex,"s"'] \\
		& & \\
		X\times_{\mc{D}}Y \arrow[rr,"g"] & & Y
	\end{tikzcd}.
\end{equation} 
Now we give an isomorphism between $\mc{G}_1\rightrightarrows X\times_{\mc{D}}Y$ and $\mc{G}_1'\rightrightarrows X\times_{\mc{D}}Y$. The construction of isomorphism between $\mc{H}_1\rightrightarrows X\times_{\mc{D}}Y$ and $\mc{H}_1'\rightrightarrows X\times_{\mc{D}}Y$ is very much same. We define map $\mc{G}_1\rightarrow \mc{G}_1'$ by giving a morphism of stacks $\underline{\mc{G}_1}\rightarrow \underline{\mc{G}_1'}$, where
\[\underline{\mc{G}_1}=(\underline{X}\times_{\mc{D}}\underline{Y})\times_{f,X,s}(\underline{X}\times_{\mc{D}}\underline{X})\times_{f,X,t}(\underline{X}\times_{\mc{D}}\underline{Y}),\]
and 
\[\underline{\mc{G}_1'}=(\underline{X}\times_{\mc{D}}\underline{Y})\times_{g,X,s}(\underline{Y}\times_{\mc{D}}\underline{Y})\times_{g,X,t}(\underline{X}\times_{\mc{D}}\underline{Y}\big).\]
A typical element in the object set of $\underline{\mc{G}_1}$ is of the form 
\[\bigg(\big(m,n,\alpha\colon p(m)\rightarrow q(n)\big),\big(a,b,p(a)\rightarrow p(b)\big),\big(m',n',\alpha'\colon p(m')\rightarrow q(n')\big)\bigg)\]
such that $m=s\big(a,b,p(a)\rightarrow p(b)\big)=a$, $n=t(a,b,p(a)\rightarrow p(b))=b$ and $p(a)\rightarrow p(b)$ is just $p(m)\rightarrow p(m')$. So, this demands a typical element to be of the form 
\[\bigg(\big(m,n,\alpha\colon p(m)\rightarrow q(n)\big),\big(m,m',p(m)\rightarrow p(m')\big),\big(m',n',\alpha'\colon p(m')\rightarrow q(n')\big)\bigg).\]
The corresponding image in $\underline{\mc{G}_1}'$ is 
\[\bigg(\big(m,n,\alpha\colon p(m)\rightarrow q(n)\big),\big(n,n',q(n)\rightarrow q(n')\big),\big(m',n',\alpha'\colon p(m')\rightarrow q(n')\big)\bigg)\]
This gives a map of stacks $\underline{\mc{G}_1}\rightarrow \underline{\mc{G}_1'}$ at the level of objects. The map at the level of morphisms can be defined similarly. This gives an isomorphism of stacks $\underline{\mc{G}_1}\rightarrow \underline{\mc{G}_1'}$, which in turn induces an isomorphism of Lie groupoids $\mc{G}_1\rightrightarrows X\times_{\mc{D}}Y$ and $\mc{G}_1'\rightrightarrows X\times_{\mc{D}}Y$. Hence, the pull-backs are isomorphic. 
\begin{lemma}\label{Lemma:LiegroupoidExtnIsIndependentofq}
	Let $F\colon \mc{D}\rightarrow \mc{C}$ be a gerbe over a stack. 
	Assume that the diagonal morphism $\Delta_F\colon \mc{D} \rightarrow \mc{D} \times_{\mc{C}}\mc{D}$ is a representable surjective submersion. 
	Then, upto a Morita equivalence, the Lie groupoid extension in  \Cref{Lemma:MoSgivingLiegroupoidExtension} does not depend on the choice of $q\colon \underline{X}\rightarrow \mc{C}$.
\end{lemma}
Thus, using \Cref{Lemma:Liegroupoidrepresentingstack}, \Cref{Lemma:existenceofatlasforC}, \Cref{Lemma:pisanepimorphism}, \Cref{Lemma:extraassumptiononDiagonalmorphism}, \Cref{Lemma:MoSgivingLiegroupoidExtension} and \Cref{Lemma:LiegroupoidExtnIsIndependentofq} we have the following result.
\begin{theorem}\label{Theorem:GoSgivesLiegroupoidextension}
	Let $F\colon \mc{D}\rightarrow \mc{C}$ be a gerbe over a stack. Assume that the diagonal morphism $\Delta_F\colon \mc{D} \rightarrow \mc{D} \times_{\mc{C}}\mc{D}$ is a representable surjective submersion. Then there exists an atlas $p\colon \underline{X}\ra \mc{D}$ for $\mc{D}$ and an atlas $q\colon \underline{X}\ra \mc{C}$ for $\mc{C}$, as in \Cref{Lemma:extraassumptiononDiagonalmorphism}, producing a Lie groupoid extension $\phi\colon \mc{G}\rightarrow \mc{H}$, where 
	$\mc{G}=(X\times_{\mc{D}}X\rightrightarrows X)$ and 
	$\mc{H}=(X\times_{\mc{C}}X\rightrightarrows X)$. Explicitly, the morphism of stacks $\Phi\colon B\mc{G}\ra B\mc{H}$ associated to $\phi\colon \mc{G}\ra \mc{H}$   along with the morphism of stacks $F\colon \mc{D}\ra \mc{H}$ forms following $2$-commutative diagram,
	\[\begin{tikzcd}
		B\mc{G} \arrow[dd, "\cong"'] \arrow[rr,"\Phi"] & & B\mc{H} \arrow[dd, "\cong"] \\
		& & \\
		\mc{D} \arrow[Rightarrow, shorten >=20pt, shorten <=20pt, uurr] \arrow[rr, "F"] & & \mc{C} 
	\end{tikzcd}.\]
	Here the isomorphisms $\mc{D}\cong B\mc{G}$ and $\mc{C}\cong B\mc{H}$ are as mentioned in \Cref{Lemma:Liegroupoidrepresentingstack}. Further, if there exists another gerbe over the stack $F'\colon \mc{D}'\ra \mc{C}'$  isomorphic to the gerbe $F\colon \mc{D}\ra \mc{C}$, then the Lie groupoid extensions associated to $F'\colon \mc{D}'\ra \mc{C}'$ and $F\colon \mc{D}\ra \mc{C}$ are Morita equivalent. 
\end{theorem} 
\begin{remark}\label{Remark:FullCapacityofDeltaF}
	Observe that we have not made full use of the condition $\Delta_F$ being a surjective submersion. We have only used the following. The morphism of stacks $p\colon \underline{X}\ra \mc{D}$ obtained in \Cref{Lemma:extraassumptiononDiagonalmorphism} is such that, $\underline{X}\times_{\mc{D}}\underline{X}$ is representable by a manifold and the morphism of stacks 
	$\Psi\colon \underline{X}\times_{\mc{D}}\underline{X}
	\rightarrow
	\underline{X}\times_{\mc{C}}\underline{X}$ is a surjective submersion at the level of manifolds. 
\end{remark}

\section{A gerbe over a stack associated to a Lie groupoid extension}
\label{Section:GoSassociatedtoALiegrupoidExtension}
This section is divided into four parts. In the first part, we will associate a $\mc{G}-\mc{H}$-bibundle for a morphism of Lie groupoids $\mc{G}\ra \mc{H}$. In the second, we associate a morphism of stacks $BP:B\mc{G}\ra B\mc{H}$ for a $\mc{G}-\mc{H}$-bibundle. In the next, we recall the notion of composition of bibundles and use it to prove that for  Morita equivalent Lie groupoids, the associated differentiable stacks are isomorphic. In the last part, we prove that the morphism of stacks associated to a Lie groupoid extension is a gerbe over a stack. 

Before we proceed to these constructions, we recall the construction of a fiber bundle from a principal bundle $Q(M,G)$ and a manifold $F$ with an action of $G$ from the left side. This is a classical construction. More details about this can be found in \cite[Section $1.5$]{MR1393940}. 

Let $M$ be a manifold, $G,H$ are Lie groups, $\pi\colon Q\rightarrow M$ a principal $G$-bundle, $\mu:G\times H\ra H$ an action of the Lie group $G$ on the manifold $H$. Consider the action of $G$ on $Q\times H$, given by $(q,h)\cdot g=(qg,g^{-1}h)$. This action satisfies the conditions mentioned in \Cref{Proposition:quotientbundleproposition}, giving the principal $H$-bundle $(Q\times H)/G\ra M$, as in the following diagram,
\[\begin{tikzcd}
	Q\times H \arrow[dd] \arrow[rr] & & (Q\times H)/G \arrow[dd] \\
	& &   \\
	Q \arrow[rr, "\pi"]   & & M  \end{tikzcd}.
\]
Thus, given a principal bundle $Q(M, G)$, an action of $G$ on a Lie group $H$, we get a principal $H$-bundle $(Q\times H)/G\ra M$.  

Let us see the above construction in terms of Lie groupoid bundles and bibundles. The following diagram illustrates $\pi\colon Q\rightarrow M$ as a $[G\rra *]$-bundle over $M$, and 
$\mu:G\times H\ra H$ as a $[G\rra *]-[H\rra *]$-bibundle, 
\[
\begin{tikzcd}
	&                         & G \arrow[dd,xshift=0.75ex,"t"]\arrow[dd,xshift=-0.75ex,"s"'] &                         & H \arrow[dd,xshift=0.75ex,"t"]\arrow[dd,xshift=-0.75ex,"s"'] \\
	& Q \arrow[ld, "\pi"] \arrow[rd] &              & H \arrow[ld] \arrow[rd] &              \\
	M             &                         & *            &                         & *           
\end{tikzcd}.\]

The projection map ${\rm pr}_1\colon Q\times H\ra Q$ induces the map $\widetilde{{\rm pr}_1}\colon (Q\times H)/G\ra Q/G\cong M$. This produces a principal $H$-bundle $(Q\times H)/G\ra M$.
The following diagram illustrates the construction,
\begin{equation}\label{Diagram:BGBHLiegroupsmap}
	\begin{tikzcd}
		(Q\times H)/G \arrow[dd, "\widetilde{{\rm pr}_1}"'] & & Q\times H \arrow[ldd, "{\rm pr}_1"'] \arrow[rdd, "{\rm pr}_2"] \arrow[ll, "\text{quotient}"'] & & \\
		& & G \arrow[dd,xshift=0.75ex,"t"]\arrow[dd,xshift=-0.75ex,"s"'] & & H \arrow[dd,xshift=0.75ex,"t"]\arrow[dd,xshift=-0.75ex,"s"'] \\
		Q/G \arrow[d, "\cong","\widetilde{\pi}"'] & Q \arrow[ld] \arrow[rd] \arrow[ld, "\pi"'] \arrow[l, "\text{quotient}"'] & & H \arrow[ld] \arrow[rd] & \\
		M & & * & & * 
	\end{tikzcd}
\end{equation}
We interpret above diagram as,
\begin{equation} \begin{tikzcd}
		& & H \arrow[dd,xshift=0.75ex,"t"]
		\arrow[dd,xshift=-0.75ex,"s"'] \\
		& (Q\times H)/G \arrow[rd] \arrow[ld,"\widetilde{\pi}\circ \widetilde{{\rm pr}_1}"'] & \\
		M & & *
	\end{tikzcd}.\end{equation}
This $[H\rra *]$-bundle over $M$ is precisely the data of a principal $H$-bundle 
$(Q\times H)/G\ra M$ over $M$. Using the \Cref{Proposition:morphismofquotientbundles}, for each morphism of principal bundles $Q(M,G)\ra Q'(M',G)$, we have the 
morphism of principal bundles $(Q\times H)(M,G)\ra (Q'\times H)(M',G)$. 
Thus, we have a functor $BG\ra BH$.

The above construction assigns a morphism of stacks $BG\ra BH$, for Lie groups $G, H$, with an action of $G$ on $H$. 
In particular, given a morphism of Lie groups $G\ra H$, we have the associated morphism of stacks $BG\ra BH$. As Lie groups are special cases of Lie groupoids, we may aim for a similar construction in the case of Lie groupoids; that is, for a morphism of Lie groupoids $\mc{G}\ra \mc{H}$ to induce a morphism of stacks $B\mc{G}\ra B\mc{H}$. Before assigning a morphism of stacks $B\mc{G}\ra B\mc{H}$, for a morphism of Lie groupoids $\phi:\mc{G}\ra \mc{H}$, as promised in \Cref{Remark:generalizedmorphismofLiegroupoids}, we associate a $\mc{G}-\mc{H}$ bibundle $\left<\phi\right>\colon \mc{G}\ra \mc{H}$ (\cite[Remark $3.24, 3.27$]{MR2778793}).

\subsection{A morphism of Lie groupoids $\mc{G}\ra \mc{H}$ gives a $\mc{G}-\mc{H}$-bibundle}\label{Subsection:bibundleAssociatedtoMorphismofLiegroupoids} 

Let   $\phi\colon \mc{G}\rightarrow \mc{H}$ be a morphism of Lie groupoids. The target map $t\colon \mc{H}_1\ra \mc{H}_0$ is a principal $\mc{H}$-bundle (\Cref{Example:targetmapisLGbundle}). 

Consider the pull-back of the principal $\mc{H}$-bundle $t\colon \mc{H}_1\rightarrow \mc{H}_0$ along the map $\phi_0\colon \mc{G}_0\rightarrow \mc{H}_0$ to get a principal $\mc{H}$-bundle over $\mc{G}_0$ (\Cref{Subsection:PullbackofPrincipalLiegroupoidbundle}),
\begin{equation} \begin{tikzcd}
		& \mc{G}_0\times_{\mc{H}_0}\mc{H}_1 \arrow[ld, "{\rm pr}_1"'] \arrow[rd, "{\rm pr}_2"] & & \mc{H}_1 \arrow[dd,xshift=0.75ex,"t"]
		\arrow[dd,xshift=-0.75ex,"s"'] \\
		\mc{G}_0 \arrow[rd, "\phi_0"] & & \mc{H}_1 \arrow[ld, "t"'] \arrow[rd, "s"] & \\
		& \mc{H}_0 & & \mc{H}_0
	\end{tikzcd}.\end{equation}

The map $\mu\colon (\mc{G}_0\times_{\mc{H}_0}\mc{H}_1)
\times_{s\circ {\rm pr}_2,\mc{H}_0,t}\mc{H}_1\ra \mc{G}_0\times_{\mc{H}_0}\mc{H}_1$, $((u,h),\tilde{h})\mapsto (u,h\circ \tilde{h})$ gives a right action of $\mc{H}$ on $\mc{G}_0\times_{\mc{H}_0}\mc{H}_1$ (\Cref{Subsection:PullbackofPrincipalLiegroupoidbundle}).

The map $\tilde{\mu}\colon \mc{G}_1\times_{s,\mc{G}_0,{\rm pr}_1}
(\mc{G}_0\times_{\mc{H}_0}\mc{H}_1) \ra(\mc{G}_0\times_{\mc{H}_0}\mc{H}_1)$,
$(g,(u,h))\mapsto (t(g),\phi(g)\circ h)$ gives a left action of $\mc{G}$ on $\mc{G}_0\times_{\mc{H}_0}\mc{H}_1$.

Thus, the manifold $\mc{G}_0\times_{\mc{H}_0} \mc{H}_1$ along with maps ${\rm pr}_1\colon \mc{G}_0\times_{\mc{H}_0} \mc{H}_1\ra \mc{G}_0, s\circ {\rm pr}_2\colon \mc{G}_0\times_{\mc{H}_0} \mc{H}_1\ra \mc{H}_0$ produce a $\mc{G}-\mc{H}$ bibundle. 
This $\mc{G}-\mc{H}$ bibundle is described by the following diagram,
\begin{equation} \begin{tikzcd}
		\mc{G}_1 \arrow[dd,xshift=0.75ex,"t"]\arrow[dd,xshift=-0.75ex,"s"'] & & \mc{H}_1 \arrow[dd,xshift=0.75ex,"t"]\arrow[dd,xshift=-0.75ex,"s"'] \\
		& \mc{G}_0\times_{\mc{H}_0} \mc{H}_1 \arrow[rd, "s\circ {\rm pr}_2"] \arrow[ld, "{\rm pr}_1"'] & \\
		\mc{G}_0 & & \mc{H}_0
	\end{tikzcd}.\end{equation} 

Thus, given a morphism of Lie groupoids $\phi\colon \mc{G}\rightarrow \mc{H}$, the manifold $\mc{G}_0\times_{\mc{H}_0} \mc{H}_1$ along with the maps ${\rm pr}_1\colon \mc{G}_0\times_{\mc{H}_0} \mc{H}_1\ra \mc{G}_0, s\circ {\rm pr}_2\colon \mc{G}_0\times_{\mc{H}_0} \mc{H}_1\ra \mc{H}_0$ is a $\mc{G}-\mc{H}$ bibundle.
We denote the manifold $\mc{G}_0\times_{\mc{H}_0}\mc{H}_1$ by $\phi^*\mc{H}_1$ and the $\mc{G}-\mc{H}$ bibundle by $\left<\phi\right>\colon \mc{G}\ra\mc{H}$. 

As a Morita morphism of Lie groupoids is a special case of a morphism of Lie groupoids, it is natural to ask what properties does a $\mc{G}-\mc{H}$-bibundle 
$\left<\phi\right>\colon \mc{G}\ra \mc{H}$ associated to a Morita morphism of Lie groupoids $\phi \colon \mc{G}\ra \mc{H}$ would have. It turns out that the action of $\mc{G}$ on $\mc{G}_0\times_{\mc{H}_0}\mc{H}_1$ is such that the map $s\circ {\rm pr}_2$ is a principal $\mc{G}$-bundle. More precisely, we have the following result.
\begin{lemma}\cite[Lemma $3.34$]{MR2778793}
	\label{Lemma:<f>isGprincipalbibundle}
	A morphism of Lie groupoids $f\colon \mc{G}\rightarrow \mc{H}$ is a Morita morphism of Lie groupoids if and only if the corresponding $\mc{G}-\mc{H}$ bibundle $\left<f\right>\colon \mc{G}\ra \mc{H}$ is a $\mc{G}$-principal bibundle (\Cref{Remark:GPrincipalbibunde}).
\end{lemma} 
\subsection{A $\mc{G}-\mc{H}$-bibundle gives a morphism of stacks $B\mc{G}\ra B\mc{H}$}\label{SubSubsection:MorphismofStacksassociatedtobibundle}
In this section, we associate a morphism of stacks $B\mc{G}\ra B\mc{H}$ for a $\mc{G}-\mc{H}$-bibundle. As a morphism of Lie groupoids is a special case of a bibundle (\Cref{Subsection:bibundleAssociatedtoMorphismofLiegroupoids}), this construction would in particular associate a morphism of stacks for a morphism of Lie groupoids. 

Let $P:\mc{G}\ra \mc{H}$ be a $\mc{G}-\mc{H}$ bibundle, which we see as the following diagram,
\begin{equation}\label{Diagram:GHbibundleinSubSubsection:MorphismofStacksassociatedtobibundle}
	\begin{tikzcd}
		\mc{G}_1 \arrow[dd,xshift=0.75ex,"t"]\arrow[dd,xshift=-0.75ex,"s"'] & & \mc{H}_1 \arrow[dd,xshift=0.75ex,"t"]\arrow[dd,xshift=-0.75ex,"s"'] \\
		& P \arrow[rd] \arrow[ld] & \\
		\mc{G}_0 & & \mc{H}_0
	\end{tikzcd}.\end{equation}

To define a morphism of stacks $B\mc{G}\ra B\mc{H}$, we associate a functor $B\mc{G}(M)\ra B\mc{H}(M)$ for each smooth manifold $M$.  Let $Q(M,\mc{G})$ be an object of the category $B\mc{G}(M)$; that is  a principal $\mc{G}$-bundle $\pi:Q\ra M$. We see $Q(M,\mc{G})$ as the following diagram, 
\begin{equation}\label{Diagram:QMGinSubSubsection:MorphismofStacksassociatedtobibundle}
	\begin{tikzcd}
		& & \mc{G}_1 \arrow[dd,xshift=0.75ex,"t"]\arrow[dd,xshift=-0.75ex,"s"'] \\
		& Q \arrow[rd] \arrow[ld, "\pi"'] & \\
		M& & \mc{G}_0
	\end{tikzcd}.\end{equation}
Combining the Diagrams \ref{Diagram:GHbibundleinSubSubsection:MorphismofStacksassociatedtobibundle}, and \ref{Diagram:QMGinSubSubsection:MorphismofStacksassociatedtobibundle}, we get the following diagram
\[
\begin{tikzcd}
	&                         & \mc{G}_0\arrow[dd,xshift=0.75ex,"t"]\arrow[dd,xshift=-0.75ex,"s"'] &                         & \mc{H}_0 \arrow[dd,xshift=0.75ex,"t"]\arrow[dd,xshift=-0.75ex,"s"'] \\
	& Q \arrow[ld, "\pi"'] \arrow[rd] &              & P \arrow[ld] \arrow[rd] &              \\
	M             &                         & \mc{G}_0            &                         & \mc{H}_0           
\end{tikzcd}.\]
The same idea of considering the pullback and the quotient space, as in Diagram \ref{Diagram:BGBHLiegroupsmap}, gives the following diagram,
\[\begin{tikzcd}
	(Q\times_{\mc{G}_0}P)/\mc{G}_1 \arrow[dd, "\widetilde{{\rm pr}_1}"'] & & Q\times_{\mc{G}_0}P \arrow[ldd, "{\rm pr}_1"'] \arrow[rdd, "{\rm pr}_2"] \arrow[ll, "\text{quotient}"'] & & \\
	& & \mc{G}_1 \arrow[dd,xshift=0.75ex,"t"]\arrow[dd,xshift=-0.75ex,"s"'] & & \mc{H}_1 \arrow[dd,xshift=0.75ex,"t"]\arrow[dd,xshift=-0.75ex,"s"'] \\
	Q/\mc{G}_1 \arrow[d, "\cong","\widetilde{\pi}"'] & Q \arrow[ld] \arrow[rd] \arrow[ld, "\pi"'] \arrow[l, "\text{quotient}"'] & & P \arrow[ld,"a_{\mc{G}}"'] \arrow[rd] & \\
	M & & \mc{G}_0 & & \mc{H}_0 
\end{tikzcd}.\]
For our convenience, we interpret the above diagram as
\begin{equation} \label{Diagram:BGBHLiegroupoidsmap}\begin{tikzcd}
		& &  \mc{H}_1 \arrow[dd,xshift=0.75ex,"t"]
		\arrow[dd,xshift=-0.75ex,"s"'] \\
		&  (Q\times_{\mc{G}_0} P)/\mc{G}_1  \arrow[rd] \arrow[ld,"\widetilde{\pi}\circ \widetilde{{\rm pr}_1}"'] & \\
		M & & \mc{H}_0
	\end{tikzcd}.\end{equation}

At the level of objects, the morphism of stacks $BP\colon B\mc{G}\rightarrow B\mc{H}$ is defined as 
\begin{equation}\label{Equation:DefinitionofBP}
	BP(\pi\colon Q\rightarrow M)=(\widetilde{\pi}\circ \widetilde{{\rm pr}_1}\colon (Q\times_{\mc{G}_0}P)/\mc{G}_1\rightarrow M).
\end{equation}
At the level of morphisms, it is defined similarly  as in the case of $\mc{G}=(G\rightrightarrows *)$ and $\mc{H}=(H\rightrightarrows *)$.
Thus, given a $\mc{G}-\mc{H}$ bibundle $P\colon \mc{G}\rightarrow \mc{H}$ we have associated a morphism of stacks $BP\colon B\mc{G}\rightarrow B\mc{H}$. 

\begin{remark}\label{Remark:mapBGtoBHdeterminesGHbibundle}
	For a $\mc{G}-\mc{H}$-bibundle $P:\mc{G}\ra \mc{H}$ we have associated a morphism of stacks $BP:B\mc{G}\ra B\mc{H}$. We can ask if the converse is true; that is if we can associate a $\mc{G}-\mc{H}$-bibundle for an arbitrary morphism of stacks. The answer is yes, as we show below.
	
	Let $F:B\mc{G}\ra B\mc{H}$ be a morphism of stacks. To obtain a $\mc{G}-\mc{H}$-bibundle, we need a principal $\mc{H}$-bundle with $\mc{G}_0$ as the base space. We have the principal $\mc{G}$-bundle $t:\mc{G}_1\ra \mc{G}_0$ with $\mc{G}_0$ as the base space. As $F$ is a morphism of stacks, we have the functor $F:B\mc{G}(\mc{G}_0)\ra B\mc{H}(\mc{G}_0)$. The image $P=F(t:\mc{G}_1\ra \mc{G}_0)$ gives an object of $B\mc{H}(\mc{G}_0)$; that is, a principal $\mc{H}$-bundle with base space as $\mc{G}_0$. It turns out that this is a 
	$\mc{G}-\mc{H}$-bibundle, and that there is a natural isomophism of functors $F\Rightarrow BP:B\mc{G}\ra B\mc{H}$. Thus, giving a morphism of stacks $B\mc{G}\ra B\mc{H}$ is same as giving a $\mc{G}-\mc{H}$-bibundle. More details can be found in  \cite[Remark $4.18$]{MR2778793}.
\end{remark}

Combining the construction of associating a $\mc{G}-\mc{H}$-bibundle for a morphism of Lie groupoids $\phi:\mc{G}\ra \mc{H}$ and associating a morphism of stacks $B\mc{G}\ra B\mc{H}$ for a $\mc{G}-\mc{H}$-bibundle, we have the construction of a morphism of stacks $BP\colon B\mc{G}\rightarrow B\mc{H}$ for a morphism of Lie groupoids $f\colon \mc{G}\rightarrow \mc{H}$.

Now we are in a position to prove that, for a Lie groupoid extension 
$[\mc{G}_1\rra M]\ra [\mc{H}_1\rra M]$, the associated morphism of stacks $B\mc{G}\ra B\mc{H}$ is a gerbe over the stack $B\mc{H}$. To prove that the diagonal morphism $B\mc{G}\ra B\mc{G}\times_{B\mc{H}}B\mc{G}$ is an epimorphism of stacks  (\Cref{SubSubsection:diagonalisEpimorphism}), we use a trick. Instead of proving $B\mc{G}\ra B\mc{G}\times_{B\mc{H}}B\mc{G}$ is an epimorphism, we prove that there exists a Lie groupoid $\mc{K}$ such that $B\mc{G}\times_{B\mc{H}}B\mc{G}\cong B\mc{K}$, and that $B\mc{G}\ra B\mc{K}$ is an epimorphism. As $B\mc{G}\ra B\mc{K}$ is related to $B\mc{G}\ra B\mc{G}\times_{B\mc{H}}B\mc{G}$, the observation that $B\mc{G}\ra B\mc{K}$ is an epimorphism of stacks then imply that $B\mc{G}\ra B\mc{G}\times_{B\mc{H}}B\mc{G}$ is an epimorphism of stacks. In this, we use the result that if two Lie groupoids are Morita equivalent, then the associated differentiable stacks are isomorphic. To prove this, we need the terminology of the composition of bibundles, which we will describe below. 

\subsection{composition of bibundles} Observe that the Diagram \ref{Diagram:BGBHLiegroupoidsmap} associates a principal $\mc{H}$-bundle for each principal $\mc{G}$-bundle $Q(M,\mc{G})$ and a $\mc{G}-\mc{H}$-bibundle. As $BP(Q(M,\mc{G}))$ is a $[M\rra M]-\mc{H}$-bibundle, we can see the construction in Diagram \ref{Diagram:BGBHLiegroupoidsmap} as a construction that associates a $[M\rra M]-\mc{H}$-bibundle for each $[M\rra M]-\mc{G}$-bibundle and a $\mc{G}-\mc{H}$-bibundle. This construction can be generalized to the case of assigning a $\mc{G}-\mc{K}$-bibundle for each $\mc{G}-\mc{H}$-bibundle and a $\mc{H}-\mc{K}$-bibundle. 

Let $P\colon \mc{G}\ra \mc{H}$ be a $\mc{G}-\mc{H}$ bibundle and 
$Q\colon \mc{H}\ra\mc{K}$  a $\mc{H}-\mc{K}$ bibundle. We have the
following diagrams for bibundles,
\begin{equation} 
	\begin{tikzcd}
		\mc{G}_1 \arrow[dd,xshift=0.75ex,"t"]\arrow[dd,xshift=-0.75ex,"s"'] & & \mc{H}_1 \arrow[dd,xshift=0.75ex,"t"]\arrow[dd,xshift=-0.75ex,"s"'] & & \mc{H}_1'\arrow[dd,xshift=0.75ex,"t"]\arrow[dd,xshift=-0.75ex,"s"'] \\
		& P \arrow[ld] \arrow[rd] & & Q \arrow[ld] \arrow[rd ] & \\
		\mc{G}_0 & & \mc{H}_0 & & \mc{H}'_0 
	\end{tikzcd}.
\end{equation}
Ignoring the  action of $\mc{G}$ on $P$, we can consider $a_{\mc{G}}\colon P\ra \mc{G}_0$ as a principal $\mc{H}$-bundle,
\begin{equation} 
	\begin{tikzcd}
		& & \mc{H}_1 \arrow[dd,xshift=0.75ex,"t"]\arrow[dd,xshift=-0.75ex,"s"'] & & \mc{H}_1'\arrow[dd,xshift=0.75ex,"t"]\arrow[dd,xshift=-0.75ex,"s"'] \\
		& P \arrow[ld,"a_{\mc{G}}"'] \arrow[rd] & & Q \arrow[ld] \arrow[rd] & \\
		\mc{G}_0 & & \mc{H}_0 & & \mc{H}'_0 
	\end{tikzcd}.
\end{equation}
Given 
a principal $\mc{H}$-bundle and a $\mc{H}-\mc{H}'$ bibundle, we know (\Cref{Equation:DefinitionofBP}) how to associate a principal $\mc{H}'$-bundle. For the principal $\mc{H}$-bundle $a_{\mc{G}}\colon P\ra \mc{G}_0$, we associate the principal $\mc{H}'$-bundle $BQ(a_{\mc{G}})\colon (P\times_{\mc{H}_0}Q)/\mc{H}_1\ra \mc{G}_0$.  

Action of $\mc{G}$ on $P$ induces an action of $\mc{G}$ on $(P\times_{\mc{H}_0}Q)/\mc{H}_1$, producing the following $\mc{G}-\mc{H}'$ bibundle,
\begin{equation}
	\label{Diagram:compositionOfbibundles}\begin{tikzcd}
		\mc{G}_1 \arrow[dd,xshift=0.75ex,"t"]\arrow[dd,xshift=-0.75ex,"s"'] & & \mc{H}'_1 \arrow[dd,xshift=0.75ex,"t"]\arrow[dd,xshift=-0.75ex,"s"'] \\
		& (P\times_{\mc{H}_0} Q)/\mc{H}_1 \arrow[ld] \arrow[rd] & \\
		\mc{G}_0 & & \mc{H}'_0 
	\end{tikzcd}.\end{equation}
\begin{definition}\label{Definition:compositionOfbibundles}
	Let $P\colon \mc{G}\ra \mc{H}$ a $\mc{G}-\mc{H}$ bibundle and $Q\colon \mc{H}\ra \mc{H}'$ be a 
	$\mc{H}-\mc{H}'$ bibundle. We define \textit{the composition of $Q$ with $P$} to be the $\mc{G}-\mc{H}'$ bibundle (in the Diagram \ref{Diagram:compositionOfbibundles}),
	\begin{equation}\label{Equation:compositionOfbibundles}
		Q\circ P=(P\times_{\mc{H}_0} Q)/\mc{H}_1.
	\end{equation}
\end{definition}
Using the terminology of composition of bibundles, we now prove that, if $\mc{G}$ and $\mc{H}$ are Morita equivalent Lie groupoids, then the stacks $B\mc{G}$ and $B\mc{H}$ are isomorphic. 
\begin{proposition}
	\label{Proposition:Moritaequivalentimpliesisomorphicstacks}
	Let $\mc{H},\mc{H}'$ be Morita equivalent Lie groupoids, then, the stacks $B\mc{H}$ and $B\mc{H}'$ are isomorphic.
	\begin{proof} Let $\mc{H}$ and $\mc{H}'$ be Morita equivalent Lie groupoids; that is, there exists a Lie groupoid $\mc{G}$ and a pair of Morita morphisms of Lie groupoids $f\colon \mc{G}\rightarrow \mc{H}$ and $ g\colon \mc{G}\rightarrow \mc{H}'$. With this data, we produce an isomorphism of stacks $B\mc{H}\rightarrow B\mc{H}'$.
		
		As mentioned in \Cref{Remark:mapBGtoBHdeterminesGHbibundle}, giving a morphism of stacks $B\mc{H}\rightarrow B\mc{H}'$ is same as giving a $\mc{H}-\mc{H}'$ bibundle. 
		The morphism of Lie groupoids $f\colon \mc{G}\rightarrow \mc{H}$ gives the following $\mc{G}-\mc{H}$ bibundle,
		\begin{equation}
			 \begin{tikzcd}
				\mc{G}_1 \arrow[dd,xshift=0.75ex,"t"]\arrow[dd,xshift=-0.75ex,"s"'] & & \mc{H}_1 \arrow[dd,xshift=0.75ex,"t"]\arrow[dd,xshift=-0.75ex,"s"'] \\
				& \left<f\right> \arrow[rd, "s\circ {\rm pr}_2"] \arrow[ld, "{\rm pr}_1"'] & \\
				\mc{G}_0 & & \mc{H}_0
			\end{tikzcd}.
		\end{equation} 
			As $f\colon \mc{G}\rightarrow \mc{H}$ is a Morita morphism of Lie groupoids,  the bibundle $\left<f\right>\colon \mc{G}\ra \mc{H}$ is a $\mc{G}$-principal bibundle (\Cref{Lemma:<f>isGprincipalbibundle}). Thus, $\left<f\right>\colon \mc{G}\ra \mc{H}$ can be considered as a $\mc{H}-\mc{G}$ bibundle,
			\begin{equation}\label{Diagram:<f>asHGbibundle}
				\begin{tikzcd}
					\mc{H}_1 \arrow[dd,xshift=0.75ex,"t"]\arrow[dd,xshift=-0.75ex,"s"'] & & \mc{G}_1 \arrow[dd,xshift=0.75ex,"t"]\arrow[dd,xshift=-0.75ex,"s"'] \\
					& \left<f\right> \arrow[rd, "{\rm pr}_1"] \arrow[ld, "s\circ {\rm pr}_2"'] & \\
					\mc{H}_0 & & \mc{G}_0
					\end{tikzcd}.\end{equation} 
		The morphism of Lie groupoids $g\colon \mc{G}\rightarrow \mc{H}'$ gives the following $\mc{G}-\mc{H}'$ bibundle,
		\begin{equation}\label{Diagram:<g>asGH'bibundle} \begin{tikzcd}
				\mc{G}_1 \arrow[dd,xshift=0.75ex,"t"]\arrow[dd,xshift=-0.75ex,"s"'] & & \mc{H}'_1 \arrow[dd,xshift=0.75ex,"t"]\arrow[dd,xshift=-0.75ex,"s"'] \\
				& \left<g\right> \arrow[ld, "{\rm pr}_1"'] \arrow[rd, "s\circ {\rm pr}_2"] & \\
				\mc{G}_0 & & \mc{H}'_0
			\end{tikzcd}.\end{equation} 
		Composing the $\mc{H}-\mc{G}$ bibundle (Diagram \ref{Diagram:<f>asHGbibundle}) $\left<f\right>\colon \mc{H}\ra \mc{G}$ with the $\mc{G}-\mc{H}'$ bibundle (Diagram \ref{Diagram:<g>asGH'bibundle}) $\left<g\right>\colon \mc{G}\ra \mc{H}'$, 
		we get the $\mc{H}-\mc{H}'$ bibundle $\left<g\right>\circ \left<f\right>\colon \mc{H}\rightarrow \mc{H}'$ (\Cref{Equation:compositionOfbibundles}), as explained  in the following diagram,
		\begin{equation} \begin{tikzcd}
				& & \left<g\right>\circ \left<f\right> \arrow[ldd] \arrow[rdd] & & \\
				\mc{H}_1 \arrow[dd,xshift=0.75ex,"t"]\arrow[dd,xshift=-0.75ex,"s"'] & & \mc{G}_1 \arrow[dd,xshift=0.75ex,"t"]\arrow[dd,xshift=-0.75ex,"s"'] & & \mc{H}'_1 \arrow[dd,xshift=0.75ex,"t"]\arrow[dd,xshift=-0.75ex,"s"'] \\
				& \left<f\right> \arrow[ld] \arrow[rd] & & \left<g\right> \arrow[ld] \arrow[rd] & \\
				\mc{H}_0 & & \mc{G}_0 & & \mc{H}'_0
			\end{tikzcd}.\end{equation} 
		As $g\colon \mc{G}\rightarrow \mc{H}'$ is also a Morita morphism of Lie groupoids, interchanging $f$ and $g$ we obtain a $\mc{H}'-\mc{H}$ bibundle as follows,
		\begin{equation} \begin{tikzcd}
				& & \left<f\right>\circ \left<g\right> \arrow[ldd] \arrow[rdd] & & \\
				\mc{H}'_1 \arrow[dd,xshift=0.75ex,"t"]\arrow[dd,xshift=-0.75ex,"s"'] & & \mc{G}_1 \arrow[dd,xshift=0.75ex,"t"]\arrow[dd,xshift=-0.75ex,"s"'] & & \mc{H}_1 \arrow[dd,xshift=0.75ex,"t"]\arrow[dd,xshift=-0.75ex,"s"'] \\
				& \left<g\right> \arrow[ld] \arrow[rd] & & \left<f\right> \arrow[ld] \arrow[rd] & \\
				\mc{H}'_0 & & \mc{G}_0 & & \mc{H}_0
			\end{tikzcd}.\end{equation} 
		The $\mc{H}-\mc{H}'$ bibundle $\left<g\right>\circ \left<f\right>\colon \mc{H}\rightarrow \mc{H}'$ gives a morphism of stacks $B\mc{H}\rightarrow B\mc{H}'$ and the $\mc{H}'-\mc{H}$ bibundle $\left<f\right>\circ \left<g\right>\colon \mc{H}'\rightarrow \mc{H}$ gives a morphism of stacks $B\mc{H}'\rightarrow B\mc{H}$. It is easy to see that the maps $B\mc{H}'\rightarrow B\mc{H}$ and $B\mc{H}\rightarrow B\mc{H}'$ are inverses to each other, giving an isomorphism of stacks $B\mc{H}\rightarrow B\mc{H}'$. Thus, the stacks $B\mc{H}$ and $B\mc{H}'$ are isomorphic.
	\end{proof}
\end{proposition} 

\subsection{A Lie groupoid extension gives a gerbe over a stack}
\label{Subsection:LieGpdExtngivesGoS}
Let $\phi\colon (\mc{G}_1\rightrightarrows \mc{G}_0)\rightarrow (\mc{H}_1\rightrightarrows \mc{H}_0)$ be a Lie groupoid extension. 
For this morphism of Lie groupoids, we have a morphism of stacks $F\colon B\mc{G}\ra B\mc{H}$, which at the level of objects have the following description,
\begin{equation}\label{Equation:definitionofFBGBH}
	F(\pi\colon Q\ra M)=(\widetilde{\pi}\circ \widetilde{{\rm pr}_1}\colon (Q\times_{\mc{G}_0}\mc{H}_1)/\mc{G}_1\rightarrow M).
\end{equation} 
Now we will  prove that $F\colon B\mc{G}\ra B\mc{H}$ is a gerbe over a stack; that is, the morphism $F\colon B\mc{G}\ra B\mc{H}$ and the diagonal morphism $\Delta_F\colon B\mc{G}\ra B\mc{G}\times_{B\mc{H}}B\mc{G}$ are epimorphisms of stacks.

\subsubsection{Proof that $F\colon B\mc{G}\ra B\mc{H}$ is an epimorphism}
\label{SubSubsection:FisanEpimorphism} 
To prove that $F:B\mc{G}\ra B\mc{H}$ is an epimorphism of stacks, we prove that, for each manifold $U$ and a morphism of stacks $\underline{U}\ra B\mc{H}$, there exists an open cover $\{U_i\ra U\}$, and a morphism of stacks $\underline{U_i}\ra B\mc{G}$ such that the following diagram is $2$-commutative, 
\begin{equation}\label{Diagram:UiUBGBH}
	\begin{tikzcd}
		\underline{U_i} \arrow[dd, "l_i"'] \arrow[rr,"\Phi=\text{inclusion}"] & & \underline{U} \arrow[dd, "q"] \\
		& & \\
		B\mc{G} \arrow[Rightarrow, shorten >=20pt, shorten <=20pt, uurr] \arrow[rr, "F"] & & B\mc{H}
	\end{tikzcd}.
\end{equation}

Let $q:\underline{U}\ra B\mc{H}$ be a morphism of stacks. We see this as a morphism of stacks $B(U\rra U)\ra B\mc{H}$. Then, the \Cref{Remark:mapBGtoBHdeterminesGHbibundle} says that, there is a unique $[U\rra U]-\mc{H}$-bibundle representing the morphism $q:\underline{U}\ra B\mc{H}$. As we have mentioned before, an $[U\rra U]-\mc{H}$-bibundle
is same as a principal $\mc{H}$-bundle over the manifold $U$.

Let $\pi\colon P\ra U$ be the $\mc{H}$-bundle associated to the morphism of stacks $q\colon \underline{U}\ra B\mc{H}$. For similar reasons as in the case of a principal Lie group bundle, there exists an open cover $\{U_i\ra U\}$ of $U$, such that, the restriction is pull-back of the trivial $\mc{H}$-bundle $t:\mc{H}_1\ra \mc{H}_0$. As the trivial bundle is with base $\mc{H}_0$, we also have a collection of maps $r_i:U_i\ra \mc{H}_0$, such that the restriction $\pi^{-1}(U_i)\ra U_i$ is pull-back of the bundle $t:\mc{H}_1\ra \mc{H}_0$, as in the following diagram, 
\begin{equation}\label{Diagram:pullbackoftalongri}
	\begin{tikzcd}
		\pi^{-1}(U_i) \arrow[dd,"\pi|_{\pi^{-1}(U_i)}"'] \arrow[rr] & & \mc{H}_1 \arrow[dd, "t"] \\
		& & \\
		U_i \arrow[rr, "r_i"] & & \mc{H}_0=M 
	\end{tikzcd}. 
\end{equation}

As finding a morphism stacks $\underline{U_i}\ra B\mc{G}$ is same as finding a principal $\mc{G}$-bundle with base space $U_i$, we look for such a $\mc{G}$-bundle. As $\mc{G}_0=\mc{H}_0=M$, we have a distinguished principal $\mc{G}$-bundle $t:\mc{G}_1\ra \mc{G}_0$. The pull-back of $t\colon \mc{G}_1\ra \mc{G}_0=M$ along $r_i\colon U_i\ra \mc{G}_0$ gives a  principal $\mc{G}$-bundle $l_i\colon W_i\ra U_i$,
\begin{equation}\label{Diagram:WiG1UiM}
	~~~~~~~\begin{tikzcd}
		W_i \arrow[dd,"l_i"'] \arrow[rr] & & \mc{G}_1 \arrow[dd, "t"] \\
		& & \\
		U_i \arrow[rr, "r_i"] & & \mc{G}_0=M 
	\end{tikzcd}. 
\end{equation}
This principal $\mc{G}$-bundle $l_i\colon W_i\ra U_i$ gives a morphism of stacks $\underline{U_i} \ra B\mc{G}$, which we denote by $l_i$. So, we have the morphism of stacks $l_i\colon \underline{U_i}\ra B\mc{G}$ for each $i$. This gives a pair of  compositions of morphisms of stacks $q\circ \Phi\colon \underline{U}_i\ra\underline{U}\ra B\mc{H}$ and $F\circ l_i\colon \underline{U}_i\ra B\mc{G}\ra B\mc{H}$.
We prove that these two compositions give the $2$-commutative Diagram \ref{Diagram:UiUBGBH}, which would then imply that $F\colon B\mc{G}\ra B\mc{H}$ is an epimorphism of stacks. 

\begin{lemma}\label{Lemma:compatabilityofFwithpullback}
	Let $\pi\colon Q\ra U$ be the pull-back of the principal $\mc{G}$-bundle $t\colon \mc{G}_1\ra M$ along a smooth map $r\colon U\ra M$. Then $F(\pi\colon Q\ra U)$ is the pull-back of the principal $\mc{H}$-bundle $F(t\colon \mc{G}_1\ra M)$ along the smooth map $r\colon U\ra M$.
	\begin{proof}
		Consider the following pull-back diagram,
		\begin{equation}\label{Diagram:pullbackoftalongr}\begin{tikzcd}
				Q \arrow[dd, "\pi"'] \arrow[rr] & & \mc{G}_1 \arrow[dd, "t"] \\
				& & \\
				U \arrow[rr, "r"] & & M
			\end{tikzcd}.\end{equation}
		We have $F(t\colon \mc{G}_1\ra M)=(\widetilde{t_{\mc{G}}}\circ \widetilde{{\rm pr}_1}\colon (\mc{G}_1\times_M \mc{H}_1)/\mc{G}_1\ra M)$ (\Cref{Equation:definitionofFBGBH}) which can be expressed by the following diagram,
		\begin{equation}\begin{tikzcd}
				& & (\mc{G}_1\times_M \mc{H}_1)/\mc{G}_1 \arrow[ldd] \arrow[rdd] & & \\
				& & \mc{G}_1 \arrow[dd,xshift=0.75ex,"t"]\arrow[dd,xshift=-0.75ex,"s"'] & & \mc{H}_1 \arrow[dd,xshift=0.75ex,"t"]\arrow[dd,xshift=-0.75ex,"s"'] \\
				& \mc{G}_1 \arrow[ld, "t_{\mc{G}}"'] \arrow[rd, "s_{\mc{G}}"] & & \mc{H}_1 \arrow[ld, "t_{\mc{H}}"'] \arrow[rd] & \\
				M & & M & & M 
			\end{tikzcd}.\end{equation}
		Adjoining the pull-back diagram (Diagram \ref{Diagram:pullbackoftalongr}) with the above diagram, we have the following diagram,
		\begin{equation}\begin{tikzcd}
				& & & (\mc{G}_1\times_M \mc{H}_1)/\mc{G}_1 \arrow[ldd] \arrow[rdd] & & \\
				& Q \arrow[ld, "\pi"'] \arrow[rd] & & \mc{G}_1 \arrow[dd,xshift=0.75ex,"t"]\arrow[dd,xshift=-0.75ex,"s"'] & & \mc{H}_1 \arrow[dd,xshift=0.75ex,"t"]\arrow[dd,xshift=-0.75ex,"s"'] \\
				U \arrow[rd, "r"] & & \mc{G}_1 \arrow[ld, "t_{\mc{G}}"'] \arrow[rd, "s_{\mc{G}}"] & & \mc{H}_1 \arrow[ld, "t_{\mc{H}}"'] \arrow[rd] & \\
				& M & & M & & M 
			\end{tikzcd}.\end{equation}
		
		Thus, we have $F(\pi\colon Q\ra U)=((Q\times_M\mc{H}_1)/\mc{G}_1\ra U)$. 
		As $Q=U\times_M\mc{G}_1$, we have
		\begin{align*}
			(Q\times_M\mc{H}_1)/\mc{G}_1=
			(U\times_M\mc{G}_1\times_M\mc{H}_1)/\mc{G}_1=U\times_M (\mc{G}_1\times_M\mc{H}_1)/\mc{G}_1.
		\end{align*}
		Note that $U\times_M (\mc{G}_1\times_M\mc{H}_1)/\mc{G}_1$ is precisely the pull-back of $(\mc{G}_1\times_M \mc{H}_1)/\mc{G}_1$ along $r\colon U\ra M$. Thus, $F(\pi\colon Q\ra U)$ is the pull-back of $F(t\colon \mc{G}_1\ra M)$ along $r\colon U\ra M$.
	\end{proof}
\end{lemma}
As the Diagrams \ref{Diagram:pullbackoftalongri} and \ref{Diagram:WiG1UiM} are pull-back diagrams, observe that 
\[F(l_i:W_i\ra U_i)=(\pi|_{\pi^{-1}(U_i)}\colon \pi^{-1}(U_i)\ra U_i),\] and 
\[(\pi|_{\pi^{-1}(U_i)}\colon \pi^{-1}(U_i)\ra U_i)=
q(U_i\ra U)=(q\circ \Phi)(\text{Id}\colon U_i\ra U_i).\] Here, $(W_i\ra U_i)=l_i(\text{Id}\colon U_i\ra U_i)$.
So, $(F\circ l_i) (\text{Id}\colon U_i\ra U_i)$ is equal to $(q\circ \Phi)(\text{Id}\colon U_i\ra U_i)$. So, there is an isomorphism $(F\circ l_i)(\text{Id}\colon U_i\ra U_i)\rightarrow (q\circ \Phi)(\text{Id}\colon U_i\ra U_i)$. For similar reasons we see that there is an isomorphism 
$(F\circ l_i)(f\colon N\ra U_i)\rightarrow (q\circ \Phi)(f\colon N\ra U_i)$ for each $f\colon N\ra U$ in $\underline{U_i}$. Thus, we have the following $2$-commutative diagram, 
\begin{equation}
	\begin{tikzcd}
		\underline{U_i} \arrow[dd, "l_i"'] \arrow[rr,"\Phi=\text{inclusion}"] & & \underline{U} \arrow[dd, "q"] \\
		& & \\
		B\mc{G} \arrow[Rightarrow, shorten >=20pt, shorten <=20pt, uurr] \arrow[rr, "F"] & & B\mc{H}
	\end{tikzcd}.
\end{equation}
Thus, $F\colon B\mc{G}\ra B\mc{H}$ is an epimorphism of stacks. 
\begin{proposition}\label{Proposition:BGBHisanepimorphism}
	Given a Lie groupoid extension $f\colon (\mc{G}_1\rightrightarrows M)\ra (\mc{H}_1\rightrightarrows M)$, the corresponding morphism of stacks $F\colon B\mc{G}\ra B\mc{H}$ is an epimorphism of stacks.
\end{proposition}
\subsubsection{Proof that the diagonal morphism $\Delta_F\colon B\mc{G}\ra B\mc{G}\times_{B\mc{H}}B\mc{G}$ is an epimorphism.}
\label{SubSubsection:diagonalisEpimorphism} 

As $\phi\colon \mc{G}\ra \mc{H}$ is a Lie groupoid extension, the $2$-fiber product $\mc{G}\times_{\mc{H}}\mc{G}$ is a Lie groupoid. As stackification and Yoneda embedding preserves the $2$-fiber product (\cite[$\text{I}.2.4$]{Carchedi}, \cite[\href{https://stacks.math.columbia.edu/tag/04Y1}{Tag 04Y1}]{Johan}), we see that \[B\mc{G}\times_{B\mc{H}}B\mc{G}\cong B(\mc{G}\times_{\mc{H}}\mc{G}).\]
Further, the diagonal morphism of stacks $\Delta_F\colon B\mc{G}\ra B\mc{G}\times_{B\mc{H}}B\mc{G}$ is the morphism of stacks associated to the diagonal morphism of Lie groupoids $\Delta_{\phi}\colon \mc{G}\ra \mc{G}\times_{\mc{H}}\mc{G}$, given by $\Delta_{\phi}(a)=(a,\text{Id}\colon a\ra a,a)$ and $\Delta_{\phi}(g)=(g,g)$ for $a\in \mc{G}_0$ and $g\in \mc{G}_1$.
We have the following morphism of Lie groupoids, 
\[\begin{tikzcd}
	\mc{G}_1 \arrow[dd,xshift=0.75ex,"t"]\arrow[dd,xshift=-0.75ex,"s"'] \arrow[rr] & & (\mc{G}\times_{\mc{H}}\mc{G})_1 \arrow[dd,xshift=0.75ex,"t"]\arrow[dd,xshift=-0.75ex,"s"'] \\
	& & \\
	M \arrow[rr] & & (\mc{G}\times_{\mc{H}}\mc{G})_0 
\end{tikzcd}.\]
As the above morphism of Lie groupoids is not identity on base space, we can not use  \Cref{Proposition:BGBHisanepimorphism} to conclude that $\Delta_F\colon B\mc{G}\ra B\mc{G}\times_{B\mc{H}}B\mc{G}$ is an epimorphism of stacks. 
\begin{remark}\label{Remark:descriptionof(GHG)0}
	As $\phi\colon (\mc{G}_1\rightrightarrows \mc{G}_0)\ra (\mc{H}_1\rightrightarrows \mc{H}_0)$ is a Lie groupoid extension, the map $\mc{G}_0\ra \mc{H}_0$ is an identity map. Thus, we have \[(\mc{G}\times_{\mc{H}}\mc{G})_0= \mc{G}_0\times_{Id,\mc{H}_0,s}\mc{H}_1\times_{t\circ {\rm pr}_2,\mc{K}_0,Id}\mc{H}_0=\mc{H}_1.\]
\end{remark}
\begin{lemma}\label{Lemma:GHGisTransitiveLieGroupoid}
	The Lie groupoid $\mc{G}\times_{\mc{H}}\mc{G}$ is a transitive Lie groupoid . 
	\begin{proof}
		To prove $\mc{G}\times_{\mc{H}}\mc{G}$ is a transitive Lie groupoid, we prove that, for $h_1,h_2\in \mc{H}_1=(\mc{G}\times_{\mc{H}}\mc{G})_0$ (\Cref{Remark:descriptionof(GHG)0}) there exists 
		$(g_1,h,g_2)\in (\mc{G}\times_{\mc{H}}\mc{G})_1$ such that 
		$s(g_1,h,g_2)=h\circ \phi(g_1)=h_1$ and $t(g_1,h,g_2)=\phi(g_2)\circ h=h_2$.
		
		As $\phi\colon \mc{G}_1\ra \mc{H}_1$ is surjective, we can choose $g_1\in \mc{G}_1$ to be such that $\phi(g_1)=1_{s(h_1)}$. Choose such a $g_1\in \mc{G}_1$.
		Choose $h=h_1$ and $g_2\in \mc{G}_2$ such that 
		$\phi(g_2)=h_2\circ h_1^{-1}$. 
		So, given $h_1,h_2\in \mc{H}_1=(\mc{G}\times_{\mc{H}}\mc{G})_0$, there exists $(g_1,h,g_2)\in (\mc{G}\times_{\mc{H}}\mc{G})_1$ such that $s(g_1,h,g_2)=h\circ \phi(g_1)=h_1$ and $t(g_1,h,g_2)=\phi(g_2)\circ h=h_2$. Thus, $\mc{G}\times_{\mc{H}}\mc{G}$ is a transitive Lie groupoid.
	\end{proof}
\end{lemma}

Combining \Cref{Lemma:GHGisTransitiveLieGroupoid} and  \Cref{Lemma:TransitiveLiegroupoidisMEtoLiegroup} we see that $\mc{G}\times_{\mc{H}}\mc{G}$ is Morita equivalent to a Lie groupoid of the form $(K\rightrightarrows *)$. Thus, by  \Cref{Proposition:Moritaequivalentimpliesisomorphicstacks}, 
the stacks $B(\mc{G}\times_{\mc{H}}\mc{G})$ and $B(K\rightrightarrows *)$ are isomorphic. So, the morphism of stacks $\Delta_F\colon B\mc{G}\ra B\mc{G}\times_{B\mc{H}}B\mc{G}$ is isomorphic to the map $B\mc{G}\ra B(K\rightrightarrows *)$.

Using an argument similar to the proof of  \Cref{Proposition:BGBHisanepimorphism}, we conclude that for any morphism of Lie groupoids $(\mc{G}_1\rightrightarrows M)\ra (K\rightrightarrows *)$, the corresponding morphism of stacks $B\mc{G}\ra B\mc{K}$ is an epimorphism of stacks. Thus, $\Delta_F\colon B\mc{G}\ra B\mc{G}\times_{B\mc{H}}B\mc{G}$ is an epimorphism of stacks. Therefore we obtain the following:
\begin{proposition}\label{Proposition:DiagonalmapisanEpimorphism}
	Given a Lie groupoid extension $\phi\colon (\mc{G}_1\rightrightarrows M)\ra (\mc{H}_1\rightrightarrows M)$ the diagonal morphism of stacks $\Delta_F\colon B\mc{G}\ra B\mc{G}\times_{B\mc{H}}B\mc{G}$ is an epimorphism of stacks.
\end{proposition}
Finally we conclude the following.
\begin{theorem}\phantomsection
	\label{Theorem:BGBHisaGoS} 
	\begin{enumerate}
		\item{Given a Lie groupoid extension $\phi\colon (\mc{G}_1\rightrightarrows M)\ra (\mc{H}_1\rightrightarrows M)$, the corresponding morphism of stacks $F\colon B\mc{G}\ra B\mc{H}$ is a gerbe over the stack $B\mc{H}$.}
		\item{ Let $\phi\colon \mc{G}_1\ra \mc{H}_1\rightrightarrows M$ and $\phi''\colon \mc{G}''_1\ra \mc{H}''_1\rightrightarrows M''$ be Morita equivalent Lie groupoid extensions. Let $\Phi\colon B\mc{G}\ra B\mc{H}$ and $\Phi''\colon B\mc{G}''\ra B\mc{H}''$ be the respective morphism of stacks corresponding to $\phi$ and $\phi''$. Then, $\Phi$ and $\Phi''$ are isomorphic in the sense that, the following diagram is $2$-commutative,
			\[\begin{tikzcd}
				B\mc{G}'' \arrow[rr,"\Phi''"] \arrow[dd,"\cong"'] & & B\mc{H}'' \arrow[dd,"\cong"] \\
				& & \\
				B\mc{G} \arrow[rr,"\Phi"] \arrow[Rightarrow, shorten >=20pt, shorten <=20pt, uurr]& & B\mc{H} . 
			\end{tikzcd}\]}
	\end{enumerate}
	\begin{proof}
		\begin{enumerate}
			\item Immediate from  \Cref{Proposition:BGBHisanepimorphism} and \Cref{Proposition:DiagonalmapisanEpimorphism}.
			\item Let the Morita equivalence be given by the Lie groupoid extension $\phi' \colon\mc{G}'_1\ra \mc{H}'_1\rightrightarrows M'$. In particular, that means we have a Morita morphism from the Lie groupoid extension $\phi' \colon\mc{G}'_1\ra \mc{H}'_1\rightrightarrows M'$ to $\phi \colon\mc{G}_1\ra \mc{H}_1\rightrightarrows M$, expressed by the following diagram, 
			\begin{equation}
				\label{Diagram:MELiegroupoidextensions}
				\begin{tikzcd}
					\mc{G}_1 \arrow[dd,xshift=0.75ex,"t"]\arrow[dd,xshift=-0.75ex,"s"'] \arrow[rrrrrrr, "\phi_1", bend left] & & \mc{G}_1' \arrow[ll, "\psi_{\mc{G}}"'] \arrow[rrr, "\phi_1'"] \arrow[dd,xshift=0.75ex,"t"]\arrow[dd,xshift=-0.75ex,"s"'] & & & \mc{H}_1' \arrow[rr, "\psi_{\mc{H}}"] \arrow[dd,xshift=0.75ex,"t"]\arrow[dd,xshift=-0.75ex,"s"'] & & \mc{H}_1 \arrow[dd,xshift=0.75ex,"t"]\arrow[dd,xshift=-0.75ex,"s"'] \\
					& & & & & & & \\
					M \arrow[rrrrrrr, "\text{Id}", bend right] & & M' \arrow[rrr, "\text{Id}"] \arrow[ll, "f"'] & & & M' \arrow[rr, "f"] & & M
				\end{tikzcd}.
			\end{equation}
			Here, $(\psi_\mc{G},f)\colon (\mc{G}'_1\rightrightarrows M')\rightarrow (\mc{G}_1\rightrightarrows M)$ and 
			$(\psi_\mc{H},f)\colon (\mc{H}'_1\rightrightarrows M')\rightarrow (\mc{H}_1\rightrightarrows M)$ are Morita morphisms of Lie groupoids. Then, by  \Cref{Proposition:Moritaequivalentimpliesisomorphicstacks} $B\mc{G}'\cong B\mc{G}$ and $B\mc{H}'\cong B\mc{H}$, and commutativity of \Cref{Diagram:MELiegroupoidextensions} gives the following commutative diagram
			\[\begin{tikzcd}
				B\mc{G} \arrow[rrrrrr, bend left,"\Phi"] & & B\mc{G}' \arrow[ll,"\cong"'] \arrow[rr,"\Phi'"] & & B\mc{H}' \arrow[rr,"\cong"] & & B\mc{H}.
			\end{tikzcd}\]
			Reorganizing the above diagram we obtain, 
			\[\begin{tikzcd}
				B\mc{G}' \arrow[rr,"\Phi'"] \arrow[dd,"\cong"'] & & B\mc{H}' \arrow[dd,"\cong"] \\
				& & \\
				B\mc{G} \arrow[rr,"\Phi"]& & B\mc{H}. 
			\end{tikzcd}\] 
			Thus, the gerbe $\Phi\colon B\mc{G}\ra B\mc{H}$ is isomorphic to the gerbe $\Phi'\colon B\mc{G}'\ra B\mc{H}'$. Repeating the same argument for a Morita morphism from the Lie groupoid extension $\phi' \colon\mc{G}'_1\ra \mc{H}'_1\rightrightarrows M'$ to $\phi''\colon \mc{G}''_1\ra \mc{H}''_1\rightrightarrows M''$, we complete the proof.
		\end{enumerate}
	\end{proof}
\end{theorem}
\begin{remark}\label{Remark:PropertyofBGtoBH}
	Let $\mc{D}\ra \mc{C}$ be a gerbe over a stack. Assume further that $\mc{D}\ra \mc{D}\times_{\mc{C}}\mc{D}$ is a representable surjective submersion. In particular, this means there exists  atlases $\underline{X}\ra \mc{C}$ and $\underline{X}\ra \mc{D}$ respectively for the stacks $\mc{C}$ and $\mc{D}$ such that the smooth map $X\times_{\mc{D}}X\ra X\times_\mc{C}X$ is a surjective submersion (\Cref{Remark:FullCapacityofDeltaF}). 
	Now we make the following observation:
	
	Let $(\Phi,1_M):(\mc{G}_1\rightrightarrows M)\ra (\mc{H}_1\rightrightarrows M)$ be a Lie groupoid extension and 
	$B\mc{G}\ra B\mc{H}$ be the associated morphism of stacks. Then
	\begin{enumerate}
		\item the morphism of stacks $B\mc{G}\ra B\mc{H}$ is a gerbe over the stack $B\mc{H}$ (\Cref{Theorem:BGBHisaGoS}).
		\item there exists atlases $\underline{M}\ra B\mc{H}$ and $\underline{M}\ra B\mc{G}$ satisfying $\underline{M}\times_{B\mc{G}}\underline{M}=\mc{G}_1$ and $\underline{M}\times_{B\mc{H}}\underline{M}=\mc{H}_1$. Moreover, the smooth map $\Phi:\mc{G}_1\ra \mc{H}_1$ associated to the morphism of stacks $\underline{M}\times_{B\mc{G}}\underline{M}\ra \underline{M}\times_{B\mc{H}}\underline{M}$ is a surjective submersion.
	\end{enumerate}
\end{remark}

	\chapter{Chern-Weil theory for principal bundles over Lie groupoids}\label{Chap.3}
Let $G$ be a Lie group, $M$ a smooth manifold, and $P(M,G)$ a principal bundle. In \Cref{Section:Atiyahforsmoothmanifold}, we have given the construction of Atiyah sequence for $P(M,G)$. Then, we introduced the notion of connection (and curvature) for $P(M,G)$ using the Atiyah sequence. We have also mentioned the Chern-Weil map associated to the principal bundle $P(M,G)$. In this chapter, which is based on our paper \cite{biswas2020chern}, we have introduced the notion of Atiyah sequence and the Chern-Weil map for a principal bundle over Lie groupoid $[X_1\rra X_0]$.
\section{Principal bundle over Lie groupoids}\label{Section:principalbundnleoverLiegroupoids}
In this section we recall the notion of principal bundle over Lie groupoid. This notion of principal bundle over Lie groupoid can be found in \cite{MR2270285,MR3150770,MR1950948,MR2119241}.
\begin{definition}\label{Definition:principalG-bundleoverLiegroupoid}
	Let  $G$ be a Lie group, and $[X_1\rra X_0]$ a Lie groupoid. A \textit{principal $G$-bundle over $[X_1\rra X_0]$} consists of the following data:
	\begin{enumerate}
		\item a principal $G$-bundle $\pi:E_G\ra X_0$ over the manifold $X_0$,
		\item a smooth map $\mu:X_1\times_{s,X_0,\pi}E_G\ra E_G$ giving an action $(\pi,\mu)$ of the Lie groupoid $[X_1\rra X_0]$ on the manifold $E_G$,
	\end{enumerate}
	such that $(\gamma \cdot a)g=\gamma \cdot (ag)$ for each $(\gamma,a,g)\in X_1\times_{s,X_0,\pi}E_G\times G$.
\end{definition}

We denote a principal $G$-bundle over a Lie groupoid $[X_1\rra X_0]$ as
$(E_G\ra X_0,[X_1\rra X_0])$.

A principal $G$-bundle $(E_G\ra X_0,[X_1\rra X_0])$ 
gives a natural Lie groupoid structure on $[s^*E_G:=X_1\times_{X_0}E_G\rra E_G]$ with source, target maps respectively being $(\gamma,a)\mapsto a$ and $(\gamma,a)\mapsto \gamma a$, and composition $\big((\gamma,a),(\gamma',b)\big)\mapsto (\gamma'\circ \gamma,a)$.

This notion of a principal $G$-bundle over a Lie groupoid includes several classical cases of principal bundles, such as a principal bundle over a manifold, an equivariant principal bundle.

	\begin{example}
	Let $H$ be a Lie group, and $[H\rra *]$ the associated Lie groupoid.  Let $G$ be a Lie group. 
	A principal $G$-bundle over $[H\rra *]$ is the same as a
	left action of the Lie group $H$ on the manifold $G$ that commutes with the right translation action of $G$ on itself.
\end{example}

\begin{example}
For a smooth manifold $M$, a Lie group $G$, a principal $G$-bundle over $M$ is a principal $G$-bundle over the Lie groupoid $[M\rra M]$.
\end{example}

\begin{example}
Let $M$ be a smooth manifold and $G,H$ are Lie groups acting on $M$. Consider the Lie groupoid $[H\times M\rra M]$, associated to the action of $H$ on $M$. A principal $G$-bundle over $[H\times M\rra M]$ is given by a principal $G$-bundle over the manifold $M$, with an action of the Lie groupoid $[H\times M\rra M]$ on $P$. This action of $[H\times M\rra M]$ gives an action of $H$ on the manifold $P$. The compatibility condition between the actions of $[H\times M\rra M]$ and $G$ on $P$ means that the actions of $G,H$ on $P$ are compatible and that the map $\pi:P\ra M$ is $H$-equivariant. Thus, we have a $H$-equivariant principal $G$-bundle $\pi:P\ra M$. So, a principal $G$-bundle over $[H\times M\rra M]$ is precisely a $H$-equivariant principal $G$-bundle over $M$.
\end{example}

\begin{definition}
	Let $[X_1\rra X_0]$ be a Lie groupoid. A \textit{morphism of principal bundles from $(E_G\ra X_0,[X_1\rra X_0])$ to $(E_G'\ra X_0,[X_1\rra X_0])$} is given by a morphism $\psi:E_G\ra E_G'$ of principal bundles  over the manifold $X_0$, such that, the following diagram is commutative,
	\[
	\begin{tikzcd}
		X_1\times_{X_0}E_G \arrow[dd, "{(1,\psi)}"'] \arrow[rr, "\mu"] &  & E_G \arrow[dd, "\psi"] \\
		&  &                        \\
		X_1\times_{X_0}E_G' \arrow[rr, "\mu'"]                         &  & E_G'                  
	\end{tikzcd}.\]
\end{definition}

As mentioned in the introduction of this chapter, given a principal bundle $P(M,G)$, we have the short exact sequence of vector bundles over the manifold $M$ (Diagram \ref{Equation:AtiyahSequence}). Using similar ideas, we construct a short exact sequence of vector bundles over the Lie groupoid $[X_1\rra X_0]$, for a principal bundle $(E_G\ra X_0,[X_1\rra X_0])$. For this, firstly, we recall the notion of a vector bundle over a Lie groupoid $[X_1\rra X_0]$.

\begin{definition}[{\cite[Section $5.4$]{MR1950948}}]
	Let $[X_1\rra X_0]$ be a Lie groupoid. A \textit{vector bundle over $[X_1\rra X_0]$} consists of the following data:
	\begin{enumerate}
		\item a vector bundle $E\ra X_0$ over the manifold $X_0$,
		\item a smooth map $\mu:X_1\times_{s,X_0,\pi}E\ra E$ giving an action $(\pi,\mu)$ of the Lie groupoid $[X_1\rra X_0]$ on the manifold $E$,
	\end{enumerate}
	such that, for each $\gamma\in X_1$, the map $\mu(\gamma,-):E_{s(\gamma)}\ra E_{t(\gamma)}$ is a linear map.
\end{definition}

\begin{definition}
	Let $[X_1\rra X_0]$ be a Lie groupoid. A \textit{morphism of vector bundles from $(E\ra X_0,[X_1\rra X_0])$ to $(E'\ra X_0,[X_1\rra X_0])$} is given by a morphism $\psi:E\ra E'$ of vector bundles over the manifold $X_0$, such that the following diagram commutes,
	\[
	\begin{tikzcd}
		X_1\times_{X_0}E \arrow[dd, "{(1,\psi)}"'] \arrow[rr, "\mu"] &  & E \arrow[dd, "\psi"] \\
		&  &                      \\
		X_1\times_{X_0}E' \arrow[rr, "\mu'"]                         &  & E'                  
	\end{tikzcd}.\]
\end{definition}
 For a Lie groupoid $[X_1\rra X_0]$, the tangent bundle $TX_0\ra X_0$ does not always have a canonical structure of a vector bundle over $[X_1\rra X_0]$. We need an extra structure on the Lie groupoid $[X_1\rra X_0]$, the notion of a connection, which we define below.

Consider the vector bundle $TX_0\ra X_0$. For it to be a vector bundle over the Lie groupoid $[X_1\rra X_0]$, we need an action of the Lie groupoid $[X_1\rra X_0]$ on the manifold $TX_0$; that is a smooth map $\mu:X_1\times_{X_0}TX_0\ra TX_0$. As mentioned before, this would assign a map $\mu(\gamma,-):T_{s(\gamma)}X_0\ra T_{t(\gamma)}X_0$ for each $\gamma\in X_1$. This should remind the notion of a connection on the vector bundle $TX_0\ra X_0$, that assigns a linear map $T_{s(\gamma)}X_0\ra T_{t(\gamma)}X_0$ for each path $s(\gamma)\ra t(\gamma)$ in $X_0$. The case $\mu:X_1\times_{X_0}TX_0\ra TX_0$ is not exactly similar, but we would still call it ``a connection'', on the Lie groupoid $[X_1\rra X_0]$.

We introduce a connection on $[X_1\rra X_0]$ as a sub bundle $\mc{H}\subseteq TX_1$, that complements the kernel of the differential $s_*$, and satisfies certain compatibility conditions related to the other structure maps of the Lie groupoid $[X_1\rra X_0]$. To specify the compatibility with the structure maps, we would have to see it as a data that assigns $T_{s(\gamma)}X_0\ra T_{t(\gamma)}X_0$ for each $\gamma\in X_1$.

\begin{definition}
	Let $[X_1\rra X_0]$ be a Lie groupoid. A sub bundle $\mc{H}\subseteq TX_1$ is said to be \textit{complementing the kernel of $s_*$}, if, for each $\gamma\in X_1$, we have
	\[T_{\gamma}X_1=\mc{H}_\gamma X_1\oplus \ker(s_{*,\gamma}).\]
\end{definition}

Let $\mc{H}\subseteq TX_1$ be a sub bundle complementing the kernel of $s_*$. Let $\gamma\in X_1$. Consider the map $s_{*,\gamma}:T_\gamma X_1\ra T_{s(\gamma)}X_0$. As $s_{*,\gamma}$ is a surjective map, the restriction of $s_{*,\gamma}$ to the complement of $\ker(s_{*,\gamma})$ is an isomorphism; that is,
$s_{*,\gamma}|_{\mc{H}_{\gamma}X_1}:\mc{H}_{\gamma}X_1\ra T_{s(\gamma)}X_0$ is an isomorphism.

Consider the composition $T_{s(\gamma)}X_0\ra \mc{H}_\gamma X_1\ra T_{t(\gamma)}X_0$ of the maps $s^{-1}_{*,\gamma}|_{\mc{H}_{\gamma}X_1}$, and
$t_{*,\gamma}$. For our convenience, we denote this composition as $\theta_{\gamma}$. Now we are in a position to define the notion of a connection on a Lie groupoid as a sub bundle satisfying certain conditions.
\begin{definition}[{\cite[Definition $3.1.$]{MR3150770}}]\label{Definition:connectiononLiegroupoid}
	Let $[X_1\rra X_0]$ be a Lie groupoid. A \textit{connection on $[X_1\rra X_0]$} is a sub bundle $\mc{H}\subseteq TX_1$, that complements the kernel of $s_*$, such that, the following conditions holds:
	\begin{enumerate}
		\item $\theta_{\gamma\circ \gamma'}=\theta_\gamma\circ \theta_{\gamma'}$ for each $(\gamma,\gamma')\in X_1\times_{X_0}X_1$,
		\item $e_{*,a}(T_aX_0)=\mc{H}_{1_a}$ for all $a\in X_0$.
	\end{enumerate}
\end{definition}

Observe that, to construct $\theta_\gamma$, we have used the target map. The first condition in \Cref{Definition:connectiononLiegroupoid} is about the compatibility of the choice of $\mc{H}\subseteq TX_1$, with the composition map $m:X_1\times_{X_0}X_1\ra X_1$, and the second condition is about the compatibility with the identity map $e:X_0\ra X_1$. We have not mentioned about the compatibility with the inverse map $i:X_1\ra X_1$, because, it follows from the above conditions that, $i_{*,\gamma}(\mc{H}_\gamma X_1)=\mc{H}_{\gamma^{-1}}X_1$ for each $\gamma\in X_1$. 

Now we give some examples of connections on Lie groupoids.

\begin{example}
	Let $[M\rra M]$ be the Lie groupoid associated to a manifold $M$. The source map is the identity map. So, $\ker(s_{*,\gamma})=\{0\}$ for each $\gamma\in M$. The assignment $\gamma\mapsto \mc{H}_\gamma M=T_\gamma M$ for $\gamma\in M$ defines a connection on the Lie groupoid $[M\rra M]$.
\end{example}

\begin{example}
	Let $[X_1\rra X_0]$ be an \`etale Lie groupoid. As the source map is a local diffeomorphism, we have $\ker(s_{*,\gamma})=\{0\}$ for each $\gamma\in X_1$. The assignment $\gamma\mapsto \mc{H}_\gamma X_1=T_\gamma X_1$ for $\gamma\in X_1$ defines a connection on the Lie groupoid $[X_1\rra X_0]$.
\end{example}

\begin{example}
	Consider the Lie groupoid $[M\times G\rra M]$ associated to an action of a Lie group $G$ on a manifold $M$. In this case, the source map is the projection map. So, $\ker(s_{*,(m,g)})=\{0\}\times T_gG$ for each $(m,g)\in M\times G$. The assignment $(m,g)\mapsto \mc{H}_{(m,g)}(M\times G)=T_{m}M$ for $(m,g)\in M\times G$ defines a connection on the Lie groupoid $[M\times G\rra M]$.
\end{example}

\begin{remark}[Notation]\label{Remark:notationofhorizontalcomponent}
	Let $[X_1\rra X_0]$ be a Lie groupoid and $\mc{H}\subseteq TX_1$ a connection on $[X_1\rra X_0]$. Let $\gamma\in X_1$. We have the decomposition $\mc{H}_{\gamma}X_1\oplus \ker(s_{*,\gamma})=T_\gamma X_1$. For each $v\in T_\gamma X_1$, we have $\tilde{v}\in \ker (s_{*,\gamma})$ and $\mc{H}_{\gamma}(v)\in \mc{H}_\gamma X_1$ such that $v=\tilde{v}+\mc{H}_\gamma(v)$. We call $\mc{H}_\gamma(v)$ to be the \textit{horizontal component of $v$}. 
\end{remark}

\subsection{An explicit description of the map $\theta_{\gamma}$ for a connection $\mc{H}\subseteq TX_1$}\label{Section:explicitdescriptionoftheta}

For our future purpose, we will give an explicit description of the map $\theta_{\gamma}:T_{s(\gamma)}X_0\ra T_{t(\gamma)}X_0$ for a connection $\mc{H}\subseteq TX_1$ on the Lie groupoid $[X_1\rra X_0]$.

Let $v\in T_{s(\gamma)}X_0$. As $s_{*,\gamma}\colon T_{\gamma}X_1\ra T_{s(\gamma)}X_0$ is surjective, there exists $w\in T_{\gamma}X_1$ such that $s_{*,\gamma}(w)=v$. Suppose that there exists $w'\in T_{\gamma}X_1$ with $s_{*,\gamma}(w')=s_{*,\gamma}(w)=v$.
For $w,w'$ we have $w=\mc{H}_{\gamma}(w)+\widetilde{w}$ and $w'=\mc{H}_{\gamma}(w')+\widetilde{w}'$ with $\widetilde{w},\widetilde{w}'\in \ker(s_{*,\gamma})$. Thus, \[s_{*,\gamma}(\mc{H}_{\gamma}(w))
=s_{*,\gamma}(w)=s_{*,\gamma}(w')
=s_{*,\gamma}(\mc{H}_{\gamma}(w')).\]
Thus, $s_{*,\gamma}(\mc{H}_{\gamma}(w))
=s_{*,\gamma}(\mc{H}_{\gamma}(w'))$; that is, $\mc{H}_{\gamma}(w)-
\mc{H}_{\gamma}(w')\in \ker(s_{*,\gamma})$. As $\mc{H}_{\gamma}X_1\subseteq T_{\gamma}X_1$ is a vector subspace and $\mc{H}_{\gamma}(w),\mc{H}_{\gamma}(w')\in \mc{H}_{\gamma}X_1$, we have $\mc{H}_{\gamma}(w)-
\mc{H}_{\gamma}(w')\in \mc{H}_{\gamma}X_1$. As $\mc{H}_{\gamma}X_1\bigcap \ker (s_{*,\gamma})=\{0\}$, and $\mc{H}_{\gamma}(w)-
\mc{H}_{\gamma}(w')\in \mc{H}_{\gamma}X_1\bigcap \ker (s_{*,\gamma})$, we have $\mc{H}_{\gamma}(w)-
\mc{H}_{\gamma}(w')=0$; that is $\mc{H}_{\gamma}(w)=
\mc{H}_{\gamma}(w')$. So, for $w,w'\in T_\gamma X_1$ with $s_{*,\gamma}(w)=s_{*,\gamma}(w')$, we have $\mc{H}_{\gamma}(w)=
\mc{H}_{\gamma}(w')$. We have the description of $s_{*,\gamma}|_{\mc{H}_\gamma X_1}^{-1}\colon T_{s(\gamma)}X_0\ra \mc{H}_\gamma X_1$, given by $v\mapsto \mc{H}_{\gamma}(w)$ where $w\in T_\gamma X_1$ is such that $s_{*,\gamma}(w)=v$. Thus, the map $\theta_\gamma\colon T_{s(\gamma)}X_0\ra T_{t(\gamma)}X_0$ is given by $v\mapsto t_{*,\gamma}(\mc{H}_{\gamma}(w))$ where $w\in T_\gamma X_1$ is such that $s_{*,\gamma}(w)=v$.

\subsection{$TX_0\ra X_0$ as a vector bundle over Lie groupoid $[X_1\rra X_0]$}
\label{Section:TXXisavectorbundle}
Let $[X_1\rra X_0]$ be a Lie groupoid. Consider the vector bundle $TX_0\ra X_0$. Let $\mc{H}\subseteq TX_1$ be a connection on $[X_1\rra X_0]$. Consider the map $\mu:X_1\times_{X_0}TX_0\ra TX_0$, defined as $(\gamma,v)\mapsto \theta_{\gamma}(v)$ for each $(\gamma,v)\in X_1\times_{X_0}TX_0$. This map $\mu$ gives an action of  $[X_1\rra X_0]$ on $TX_0$. 

Observe that, for each $\gamma\in X_1$, the map $\mu(\gamma,-)=\theta_\gamma$, is by definition, a linear map, being composition of linear maps. Thus, with this action of $[X_1\rra X_0]$ on $TX_0$, the tangent bundle map $TX_0\ra X_0$ is considered as a vector bundle over the Lie groupoid $[X_1\rra X_0]$.

\subsection{Comparison of the notion of connections with existing notions in literature}

In this section, we will compare our notion of connection on Lie groupoids with the notion of connection on Lie groupoids from the papers \cite{MR3107517,MR2183389}. The notion of connection on Lie groupoids appear in other works 
\cite{MR3872622,MR2493616,MR2238946,MR3886162}.

\begin{definition}[{\cite[Definition $2.8$]{MR3107517}}]
Let $[X_1\rra X_0]$ be a Lie groupoid. An \textit{Ehressman connection on $[X_1\rra X_0]$} is a sub-bundle $\mc{H}\subseteq TX_1$ which is complementary to (kernel of) $s_*:TX_1\ra TX_0$ and has the property that $\mc{H}_x=T_xX_0$ for all $x\in X_0$.
 
\end{definition}

\begin{definition}[{\cite[Definition $2.1$]{MR2183389}}]
Let $[X_1\rra X_0]$ be a Lie groupoid. A \textit{connection on $[X_1\rra X_0]$} is a subbundle $E\subseteq TX_1$ such that $[E\rra TX_0]$ is a subgroupoid of $[TX_1\rra TX_0]$ and both diagrams 
\[
\begin{tikzcd}
	E \arrow[rr,yshift=0.75ex,"s"]\arrow[rr,yshift=-0.75ex,"t"'] \arrow[dd] &  & TX_0 \arrow[dd] \\
	&  &                 \\
	X_1 \arrow[rr,yshift=0.75ex,"s"]\arrow[rr,yshift=-0.75ex,"t"']         &  & X_0            
\end{tikzcd},\]
are cartesian.  
\end{definition}

Let $\mc{H}\subseteq TX_1$ be a connection on $[X_1\rra X_0]$ as in 
\Cref{Definition:connectiononLiegroupoid}. Consider the set 
\[s^*TX_0=\{(\gamma, a,v)\in X_1\times TX_0: s(\gamma)=a, v\in T_aX_0\}.\] This gives a Lie groupoid $[s^*TX_0\rra TX_0]$ with the following description:
\begin{itemize}
\item the source map is given by $(\gamma,v)\mapsto (s(\gamma),v)$,
\item the target map is given by $(\gamma,v)\mapsto (t(\gamma),\theta_{\gamma}(v))$.
\end{itemize}
The composition, inversion, unit maps are defined similarly. Consider the set \[t^*TX_0=\{(\gamma, a,v)\in X_1\times TX_0: t(\gamma)=a, v\in T_aX_0\}.\] This gives a Lie groupoid $[t^*TX_0\rra TX_0]$ with the following description:
\begin{itemize}
	\item the source map is given by $(\gamma,u)\mapsto (s(\gamma),\theta_{\gamma^{-1}}(v))$,
	\item the target map is given by $(\gamma,u)\mapsto (t(\gamma),u)$.
\end{itemize}
The composition, inversion, unit maps are defined similarly. The map $\theta$ gives a morphisms of Lie groupoids $[s^*TX_0\rra TX_0]\ra [t^*TX_0\rra TX_0]$ (\cite[Lemma $4.6$]{biswas2020atiyah}).

The compatibility of connection with composition of the Lie groupoid implies that, 
\[s_{*,\gamma_2\circ \gamma_1}|_{\mc{H}}^{-1}(v)-
m_{*,(\gamma_1,\gamma_2)}(s_{*,\gamma_2}|_{\mc{H}}^{-1}(\theta_{\gamma_1}(v)), s_{*,\gamma_1}|_{\mc{H}}^{-1}(v))\in \ker(t_{*,{\gamma_2\circ \gamma_1}}),\] for each $v\in T_{s(\gamma_1)}X_0$. Suppose that 
\[s_{*,\gamma_2\circ \gamma_1}|_{\mc{H}}^{-1}(v)-
m_{*,(\gamma_1,\gamma_2)}(s_{*,\gamma_2}|_{\mc{H}}^{-1}(\theta_{\gamma_1}(v)), s_{*,\gamma_1}|_{\mc{H}}^{-1}(v))=0\] for each $(\gamma_1,\gamma_2)\in X_1\times_{X_0}X_1$ and $v\in T_{s(\gamma_1)}X_0$, then, $[\mc{H}\subseteq TX_0]$ is a subgroupoid of $[TX_1\rra TX_0]$. It turns out that, our definition of connection on $[X_1\rra X_0]$ implies the definition of connection in \cite{MR3107517}. The isomorphism $[s^*TX_0\rra TX_0]\ra [t^*TX_0\rra TX_0]$ implies that $\mc{H}$ is a connection in terms of \cite{MR2183389}.
\section{Atiyah sequence and connection on principal bundle over Lie groupoid}\label{Section:AtiyahforbundleoverLiegroupoid}

Let $\mb{X}=[X_1\rra X_0]$ be a Lie groupoid, and $(E_G\ra X_0,[X_1\rra X_0])$ a principal bundle. As mentioned in the introduction, our next aim is to associate a sequence of vector bundles over $[X_1\rra X_0]$ for the principal bundle $(E_G\ra X_0,[X_1\rra X_0])$.

The notion of Atiyah sequence and connection on a principal bundle over a manifold is discussed in detail in \Cref{Section:Atiyahforsmoothmanifold}.

For the principal bundle $E_G\ra X_0$, we have the Atiyah sequence,
\[0\ra (E_G\times \mf{g})/G\ra (TE_G)/G\ra TX_0\ra 0,\]
of vector bundles over the manifold $X_0$. Our aim is to turn this into a short exact sequence of vector bundles over the Lie groupoid $[X_1\rra X_0]$. 
As mentioned in \Cref{Section:TXXisavectorbundle}, after fixing a connection on $[X_1\rra X_0]$, the tangent bundle $TX_0\ra X_0$ can be considered as  a vector bundle over $[X_1\rra X_0]$.

Let $\mc{H}\subseteq TX_1$ be a connection on $[X_1\rra X_0]$. This gives the tangent bundle $TX_0\ra X_0$ a vector bundle structure over $[X_1\rra X_0]$. We expect that this connection $\mc{H}\subseteq TX_1$ would turn $(TE_G)/G\ra X_0$ and $(E_G\times \mf{g})/G\ra X_0$ into vector bundles over $[X_1\rra X_0]$.


An action of $\mb{X}$ on $(TE_G)/G\ra X_0$ is given by  a map 
\[\mu\colon X_1\times_{X_0}(TE_G)/G\ra (TE_G)/G.\] For each $\gamma\in X_1$ with $(\gamma, -)\in X_1\times_{X_0}(TE_G)/G$ we need a map \[\mu_{\gamma}\colon ((TE_G)/G)_{s(\gamma)}\ra ((TE_G)/G)_{t(\gamma)}.\]
An arbitrary element of $X_1\times_{X_0}(TE_G)/G$ is of the form $(\gamma,[v])$ where $\gamma\in X_1$ and $v\in T_aE_G$ with $\pi(a)=s(\gamma)$. For  $(\gamma,a)\in s^*E_G$, we have $t(\gamma)=\pi(\gamma a)$.
That is, for each $(\gamma,a)\in s^*E_G$ we need a map \[\mu_{\gamma}\colon ((TE_G)/G)_{\pi(a)}\ra ((TE_G)/G)_{\pi(\gamma a)}.\]
The pull-back diagram
\[\begin{tikzcd}
	TE_G \arrow[dd] \arrow[rr, "q^{/G}"] &  & (TE_G)/G \arrow[dd] \\
	&  &                     \\
	E_G \arrow[rr, "\pi"]                &  & X_0                
\end{tikzcd}\]
yields following isomoprhism, $((TE_G)/G)_{\pi(a)}\cong T_aE_G$ and $((TE_G)/G)_{\pi(\gamma a)}\cong T_{\gamma a}E_G$.
That means giving  $\mu_{(\gamma,a)}$ is same as giving $\mu_{(\gamma,a)}:T_aE_G\ra T_{\gamma a}E_G$.

Recall after \Cref{Definition:principalG-bundleoverLiegroupoid} we have observed $[s^*E_G\rra E_G]$ is a Lie groupoid, and in fact, $\mu_{(\gamma,a)}$ is a connection on the Lie groupoid $[s^*E_G\rra E_G]$.

A principal bundle $(E_G\ra X_0,[X_1\rra X_0])$ can be viewed as a morphism of Lie groupoids
\[({\rm pr}_1,\pi ):[s^*E_G\rra E_G]\ra [X_1\rra X_0].\] We now see that for the differential of the projection map $d({\rm pr}_1)\colon Ts^*E_G\ra TX_1$, the pull-back $\widetilde{\mc{H}}=d({\rm pr}_1)^{-1}(\mc{H})$ is a connection on the Lie groupoid $[s^*E_G\rra E_G]$. 

\subsection{Pull-back connection on the Lie groupoid $[s^*E_G\rra E_G]$}\label{ConnectionOnsEGEG}
For $\widetilde{\mathcal{H}}=
(d({\rm pr}_1))^{-1}(\mathcal{H})\subseteq Ts^*E_G$ to form a connection on the Lie groupoid $[s^*E_G\rra E_G]$, we need 
\begin{equation}\label{connections^*E_Gcondition}\ker (({\rm pr}_2)_{*,(\gamma,a)})\oplus \widetilde{\mathcal{H}}_{(\gamma,a)}s^*E_G=T_{(\gamma,a)}s^*E_G\end{equation}
to hold for each $(\gamma,a)\in s^*E_G$. Here $\widetilde{\mathcal{H}}_{(\gamma,a)}s^*E_G=({\rm pr}_1)^{-1}_{*,(\gamma,a)}(\mathcal{H}_\gamma X_1)$.
As $\mathcal{H}\subseteq TX_1$ is a connection on $[X_1\rra X_0]$, we have 
\begin{equation}\label{connectioX_1condition}
	\ker(s_{*,\gamma})\oplus\mathcal{H}_{\gamma}X_1=T_\gamma X_1\end{equation} for each $\gamma\in X_1$. 

Let $(v,w)\in T_{(\gamma,a)}(s^*E_G)$; that is, $v\in T_{\gamma}X_1$ and $w\in T_aE_G$ such that $s_{*,\gamma}(v)=\pi_{*,a}(w)$. 
Suppose that $(v,w)\in \ker(({\rm pr}_2)_{*,(\gamma,a)})\bigcap  \widetilde{\mathcal{H}}_{(\gamma,a)}s^*E_G$.  As $(v,w)\in \ker(({\rm pr}_2)_{*,(\gamma,a)})$, we have $w=0$. 
Whereas, $(v,w)\in \widetilde{\mc{H}}_{(\gamma,a)}(s^*E_G)=({\rm pr}_1)_{*,(\gamma,a)}^{-1}(\mc{H}_{\gamma}X_1)$ means $v\in \mc{H}_{\gamma}X_1$.
The condition that $s_{*,\gamma}(v)=\pi_{*,a}(w)$ and $w=0$ implies that  $v\in \ker (s_{*,\gamma})$. But, $\mc{H}\subseteq TX_1$ is a connection, thus, the condition 
$v\in \ker(s_{*,\gamma})\bigcap \mathcal{H}_{\gamma}X_1$ implies that $v=0$.
Thus, we conclude 
\[\ker(({\rm pr}_2)_{*,(\gamma,a)})\bigcap  \widetilde{\mathcal{H}}_{(\gamma,a)}s^*E_G=\{0\}.\]

Let $(v,w)\in T_{(\gamma,a)}(s^*E_G)$. As $v\in T_{\gamma}X_1=\ker(s_{*,\gamma})\oplus \mc{H}_{\gamma}X_1$ we have $v=\widetilde{v}+\mc{H}_{\gamma}(v)$ where $\widetilde{v}$ is projection of $v$ onto $\ker(s_{*,\gamma})$ and $\mc{H}_{\gamma}(v)$ is the projection of $v$ onto $\mc{H}_{\gamma}X_1$. 
As $({\rm pr}_2)_{*,(\gamma,a)}(\widetilde{v},0)=0$, we obtain \[(\widetilde{v},0)\in \ker(({\rm pr}_2)_{*,(\gamma,a)}).\] Since $({\rm pr}_1)_{*,(\gamma,a)}(\mc{H}_{\gamma}(v),w)=\mc{H}_{\gamma}(v)\in \mc{H}_{\gamma}X_1$; that is, $(\mc{H}_{\gamma}(v),w)\in ({\rm pr}_1)_{*,(\gamma,a)}^{-1}(\mc{H}_{\gamma}X_1)$, we obtain, 
\[(\mc{H}_{\gamma}(v),w)\in \widetilde{\mathcal{H}}_{(\gamma,a)}s^*E_G.\]
So, given $(v,w)\in T_{\gamma}\in  T_{(\gamma,a)}(s^*E_G)$, we have 
\begin{align}\label{definitionofwidetildeH}
	(v,w)=(\widetilde{v},0)+(\mc{H}_{\gamma}(v),w)\in \ker(({\rm pr}_2)_{*,(\gamma,a)})+ \widetilde{\mathcal{H}}_{(\gamma,a)}s^*E_G,\end{align}
which along with the condition $\ker(({\rm pr}_2)_{*,(\gamma,a)})\bigcap  \widetilde{\mathcal{H}}_{(\gamma,a)}s^*E_G=\{0\}$
imply that 
\[\ker(({\rm pr}_2)_{*,(\gamma,a)})\oplus  \widetilde{\mathcal{H}}_{(\gamma,a)}s^*E_G=T_{(\gamma,a)}s^*E_G.\]

Explicitly we can write, for each $(\gamma,a)\in s^*E_G=X_1\times_{X_0}E_G$, we have 
\[\widetilde{\theta}_{(\gamma,a)}(w)
= t_{*,(\gamma,a)}(\widetilde{\mathcal{H}}_{(\gamma,a)}(v,w))\]
for some $(v,w)\in T_{(\gamma,a)}s^*E_G$. 
In case of $[s^*E_G\rra E_G]$, we have the target map $t=\mu$.
 So, given $w\in T_aE_G=T_{s(\gamma,a)}E_G$, we have  \begin{equation}\label{thetatilde}
	\widetilde{\theta}_{(\gamma,a)}(w)
	= \mu_{*,(\gamma,a)}
	(\mc{H}_{\gamma}(v),w)\end{equation}
for some $v\in T_{\gamma}X_1$ with 
$s_{\gamma}(v)=\pi_{*,a}(w)$. It is easy to see that the above mentioned $\widetilde{\mc{H}}\subseteq Ts^*E_G$ and the corresponding $\widetilde{\theta}_{(\gamma,a)}\colon T_{s(\gamma,a)}E_G\ra T_{t(\gamma,a)}E_G$ for each $(\gamma,a)\in s^*E_G$ 
is such that  $e_{*,a}(T_aE_G)=\widetilde{\mathcal{H}}_{(1_{\pi(a)},a)}s^*E_G$ and that $\widetilde{\theta}_{(\gamma,a)}\circ \widetilde{\theta}_{(\gamma',a')}=\widetilde{\theta}_{(\gamma\circ \gamma',a)}$ for all composable $(\gamma,a),(\gamma',a')$ in $s^*E_G$. Thus, the distribution $\widetilde{\mathcal{H}}=(d({\rm pr}_1))^{-1}(\mc{H})\subseteq Ts^*E_G$ defines a connection on the Lie groupoid $[s^*E_G\rra E_G]$.
\begin{lemma}\label{Lemma:pullbackofintegrableconnection}
	If $\mc{H}\subseteq TX_1$ is a connection on the Lie groupoid $[X_1\rra X_0]$, then, the pull-back $\widetilde{\mathcal{H}}\subseteq Ts^*E_G$ is a  connection on the Lie groupoid $[s^*E_G\rra E_G]$. Further, if $\mc{H}\subseteq TX_1$ is an integrable connection, then $\widetilde{\mathcal{H}}\subseteq Ts^*E_G$ is an integrable connection.
\end{lemma}

\subsection{{$(TE_G)/G\ra X_0$} as a vector bundle over Lie groupoid}\label{Section:TEGXisavectorbundle}
Now that we have a connection on the Lie groupoid $[s^*E_G\rra E_G]$, we can give an action of the Lie groupoid $[X_1\rra X_0]$ on the manifold $(TE_G)/G$.

Consider the map $X_1\times_{X_0}(TE_G)/G\ra (TE_G)/G$ given by $(\gamma,[(a,v)])\mapsto [\tilde{\theta}_{(\gamma,a)}(v)]$. This gives an action of the Lie groupoid $[X_1\rra X_0]$ on the manifold $(TE_G)/G$. With this action, we have the vector bundle $((TE_G)/G\ra X_0,[X_1\rra X_0])$.

\subsection{{$(E_G\times \mf{g})/G\ra X_0$} as a vector bundle over the Lie groupoid}\label{Section:EGgisavectorbundle} The action $X_1\times_{X_0}(TE_G)/G\ra (TE_G)/G$ mentioned above, restricts to smooth map $X_1\times_{X_0}(E_G\times \mf{g})/G\ra (E_G\times \mf{g})/G$, giving the vector bundle 
$((E_G\times \mf{g})/G\ra X_0,[X_1\rra X_0])$. To be more precise, we have $\gamma\cdot [(a,A)]=[\gamma\cdot a, A]$ for all 
$(\gamma,[(a,A)])\in X_1\times_{X_0} (E_G\times \mf{g})/G$.

Combining the constructions in  
\Cref{Section:TXXisavectorbundle}, \Cref{Section:TEGXisavectorbundle}, and 
\Cref{Section:EGgisavectorbundle}, we get a short exact sequence,
\begin{equation}\label{Equation:AtiyahsequencefprPBLG}
	0\ra (E_G\times \mf{g})/G\ra (TE_G)/G\ra TX_0\ra 0,
\end{equation}
of vector bundles over the Lie groupoid $[X_1\rra X_0]$.

\begin{definition}
	Let $G$ be a Lie group, $[X_1\rra X_0]$ a Lie groupoid, $(E_G\ra X_0,[X_1\rra X_0])$ a principal bundle. The short exact sequence in \Cref{Equation:AtiyahsequencefprPBLG} is called the \textit{Atiyah sequence of vector bundles over $[X_1\rra X_0]$ associated to $(E_G\ra X_0,[X_1\rra X_0])$}.
\end{definition}

\begin{definition}
	Let $G$ be a Lie group, $[X_1\rra X_0]$ a Lie groupoid, and  $(E_G\ra X_0,[X_1\rra X_0])$ a principal bundle. A \textit{connection on the principal bundle $(E_G\ra X_0,[X_1\rra X_0])$} is defined as a splitting of  the Atiyah sequence of vector bundles over $[X_1\rra X_0]$ associated to $(E_G\ra X_0,[X_1\rra X_0])$ (\Cref{Equation:AtiyahsequencefprPBLG}).
\end{definition}

\section{differential forms on a Lie groupoid equipped with an integrable connection}\label{Section:differentialformsassociatedtoconnectiononLiegroupoid}

Let $G$ be a Lie group, $M$ a smooth manifold, and $\pi:P\ra M$ a principal $G$-bundle. We have seen in  \Cref{Section:Atiyahforsmoothmanifold} that a connection on $P(M, G)$ defined as a splitting $\mc{D}$ of the Atiyah sequence 
(Diagram \ref{Equation:AtiyahSequence}) determines a differential $1$-form on $P$, which we called the connection $1$-form associated to $\mc{D}$. 
In this section, we introduce the notion of a differential form on a  Lie groupoid, so that we can associate a differential $1$-form with the splitting.

In literature, there is already a notion of differential forms/ de Rham cohomology on Lie groupoids (\cite[Section $3.1$]{MR2270285}). Here we introduce a new notion of a differential form on a Lie groupoid, equipped with a connection. So, we are actually defining a differential form on the pair $(\mb{X},\mc{H})$. 

\begin{definition}
	Let $[X_1\rra X_0]$ be a Lie groupoid, and $\mc{H}\subseteq TX_1$ a connection on $[X_1\rra X_0]$. Let $\omega:X_1\ra \Lambda^k T^*X_1$ be a differential form on $X_1$. The \textit{horizontal component of $\omega$}, denoted  by $\mc{H}(\omega)$, is a differential $k$-form on $X_1$, defined as 
	\[\mc{H}(\omega)(\gamma)(v_1,\cdots,v_k)=
	\omega(\gamma)(\mc{H}_{\gamma}(v_1),\cdots, \mc{H}_{\gamma}(v_k)),
	\]
	for all $\gamma\in X_1$, $v_1,\cdots,v_k$ in $T_\gamma X_1$. Here, $\mc{H}_\gamma(v)$ is the horizontal component of $v$ (\Cref{Remark:notationofhorizontalcomponent}).
\end{definition}

\begin{definition}[{\cite[Definition $5.1.$]{MR2270285}}]
	Let $[X_1\rra X_0]$ be a Lie groupoid, and  $\mc{H}\subseteq TX_1$ a connection on $[X_1\rra X_0]$. A \textit{differential $k$-form of $(\mb{X},\mc{H})$}, is a differential $k$-form $\omega$ on $X_0$, such that $\mc{H}(s^*\omega)=\mc{H}(t^*\omega)$. 
\end{definition}

Let $\mb{X}=[X_1\rra X_0]$ be a Lie groupoid, $\mc{H}\subseteq TX_1$ a connection on $[X_1\rra X_0]$. 
Let $\Omega^k(\mb{X},\mc{H})$ be the vector space of differential $k$-forms of $(\mb{X},\mc{H})$. Thus, we have a collection of vector spaces $\{\Omega^k(\mb{X},\mc{H})\}_{k=0}^{\dim X_0}$. Similar to the case of a cochain complex 
\[\Omega^0(M)\ra \Omega^1(M)\ra \Omega^2(M)\cdots,\]
which we use to define de Rham cohomology of $M$, we would like to construct a cochain complex 
\[\Omega^0(\mb{X},\mc{H})\ra\Omega^1(\mb{X},\mc{H})\ra \Omega^2(\mb{X},\mc{H})\cdots.\]
The usual exterior derivation operator $d:\Omega^k(X_0)\ra \Omega^{k+1}(X_0)$ do not always preserve the property of ``being a differential form on $(\mb{X},\mc{H})$''. As we will see below, an extra condition on $\mc{H}$, of being an integrable connection, ensures that  the exterior derivation operator preserve the property of ``being a differential form on $(\mb{X},\mc{H})$''.

The following result was mentioned in \cite[Lemma $5.2$]{MR3150770} without proof. In our paper \cite{biswas2020chern}, we have given a proof. Following is a more detailed version of the proof in \cite{biswas2020chern}. 
\begin{lemma}[{\cite[Lemma $5.2$]{MR3150770}}]\label{Lemma:differentialisdifferentialform}
	Let $[X_1\rra X_0]$ be a Lie groupoid and $\mc{H}\subseteq TX_1$ an integrable connection on the Lie groupoid $[X_1\rra X_0]$. If $\varphi$ is a  $k$-form on the Lie groupoid $[X_1\rra X_0]$, then, the exterior derivative $d\varphi$ of $\varphi$, when $\varphi$ is seen as a  $k$-form on the manifold $X_0$, is a  $k+1$-form on the Lie groupoid $[X_1\rra X_0]$.
	\begin{proof}
		We prove for the case $k=1$; that is, when $\varphi$ is a  $1$-form on the Lie groupoid $[X_1\rra X_0]$. The proof works almost verbatim for arbitrary $k$. As $\varphi\colon X_0\ra \Lambda^1T^*X_0$ is a  $1$-form on the Lie groupoid $[X_1\rra X_0]$, we have $\varphi$ a  $1$-form on $X_0$, satisfying the condition $\mc{H}(s^*\varphi)=\mc{H}(t^*\varphi)$. That is, $\mc{H}(s^*\varphi)(\gamma)(v)=\mc{H}(t^*\varphi)(\gamma)(v)$ for all $\gamma\in X_1$ and $v\in T_{\gamma}X_1$. By definition, $\mc{H}(s^*\varphi)(\gamma)(v)=(s^*\varphi)(\gamma)(\mc{H}_\gamma)(v)$ and $\mc{H}(t^*\varphi)(\gamma)(v)=(t^*\varphi)(\gamma)(\mc{H}_\gamma)(v)$. Thus, for the $1$-form $\varphi$ on the manifold $X_0$, we have 
		$(s^*\varphi)(\gamma)(\mc{H}_{\gamma}v)=
		(t^*\varphi)(\gamma)(\mc{H}_{\gamma}v)$ for all $\gamma\in X_1$ and $v\in T_{\gamma}X_1$.
		
		For such a $1$-form $\varphi$ on the manifold $X_0$, we prove that, the exterior derivative $d\varphi$, which is a $2$-form on the manifold $X_0$, is a $2$-form on the Lie groupoid $[X_1\rra X_0]$. For that, we need to prove $\mc{H}(s^*(d\varphi))=\mc{H}(t^*(d\varphi))$ as  $2$-forms on the manifold $X_1$; that is,  
		$\mc{H}(s^*(d\varphi))(\gamma)(v_1,v_2)=\mc{H}(t^*(d\varphi))(\gamma)(v_1,v_2)$ for all $\gamma\in X_1$ and $v_1,v_2\in T_{\gamma}X_1$. By definition, we have 
		\begin{align*}
			\mc{H}(s^*(d\varphi))(\gamma)(v_1,v_2)&=
			s^*(d\varphi)(\gamma)(\mc{H}_{\gamma}(v_1),\mc{H}_{\gamma}(v_2))=
			d(s^*\varphi)(\gamma)(\mc{H}_{\gamma}(v_1),\mc{H}_{\gamma}(v_2)),\\
			\mc{H}(t^*(d\varphi))(\gamma)(v_1,v_2)&=
			t^*(d\varphi)(\gamma)(\mc{H}_{\gamma}(v_1),\mc{H}_{\gamma}(v_2))=
			d(t^*\varphi)(\gamma)(\mc{H}_{\gamma}(v_1),\mc{H}_{\gamma}(v_2)).
		\end{align*} 
	Assuming $(s^*\varphi)(\gamma)(\mc{H}_{\gamma}v)=
		(t^*\varphi)(\gamma)(\mc{H}_{\gamma}v)$ for all $\gamma\in X_1$ and $v\in T_{\gamma}X_1$, we prove that $d(s^*\varphi)(\gamma)(\mc{H}_{\gamma}(v_1),\mc{H}_{\gamma}(v_2))
		=d(t^*\varphi)(\gamma)(\mc{H}_{\gamma}(v_1),\mc{H}_{\gamma}(v_2))$
		for all $\gamma\in X_1$ and $v_1,v_2\in T_{\gamma}X_1$.
		
		For notational convenience, we write $\varphi_1=s^*\varphi$ and $\varphi_2=t^*\varphi$. 
		
		Fix $\gamma\in X_1$. Let $v_1,v_2\in T_{\gamma}X_1$ and consider their horizontal components $\mc{H}_\gamma(v_1),\mc{H}_\gamma(v_2)$. Observe that, for these horizontal vectors $\mc{H}_{\gamma}(v_1),\mc{H}_{\gamma}(v_2)\in \mc{H}_{\gamma}X_1$, there exists horizontal vector fields $Z_1,Z_2\colon X_1\ra TX_1$ such that $Z_1(\gamma)=\mc{H}_{\gamma}(v_1)$ and 
		$Z_2(\gamma)=\mc{H}_{\gamma}(v_2)$. Then, we have 
		\begin{align*}
			(d\varphi_1)(\gamma)(\mc{H}_\gamma(v_1),\mc{H}_\gamma(v_2))=(d\varphi_1)(\gamma)(Z_1(\gamma),
			Z_2(\gamma))=(d\varphi_1)(Z_1,Z_2)(\gamma),\\
			(d\varphi_2)(\gamma)(\mc{H}_\gamma(v_1),\mc{H}_\gamma(v_2))
			=(d\varphi_2)(\gamma)(Z_1(\gamma),
			Z_2(\gamma))=(d\varphi_2)(Z_1,Z_2)(\gamma).
		\end{align*} 
		For vector fields $Z_1,Z_2\colon X_1\ra TX_1$, the exterior derivative $d\varphi_1$ is given by 
		\[(d\varphi_1)(Z_1,Z_2)=Z_1\varphi_1(Z_2)-Z_2(\varphi_1(Z_1))-
		\varphi_1([Z_1,Z_2]).\]
		
	Integrality of $\mc{H}\subseteq TX_1$ implies that, for horizontal vector fields $Y, Y'$ on $X_1$, the Lie bracket $[Y,Y']$ is a horizontal vector field on $X_1$. As $\varphi_1(\widetilde{\gamma})(\mc{H}_{\widetilde{\gamma}}(v))=\varphi_2(\widetilde{\gamma})(\mc{H}_{\widetilde{\gamma}}(v))$ for  $\widetilde{\gamma}\in X_1$ and $v\in T_{\widetilde{\gamma}}X_1$, we see that the differential forms $\varphi_1$ and $\varphi_2$ agree on horizontal vector fields; in particular, we have 
		$\varphi_1(Z_1)=\varphi_2(Z_1)$,  $\varphi_1(Z_2)=\varphi_2(Z_2)$ and $\varphi_1([Z_1,Z_2])=\varphi_2([Z_1,Z_2])$. So, we have 
		\begin{align*}
			(d\varphi_1)(Z_1,Z_2)&=Z_1\varphi_1(Z_2)-Z_2(\varphi_1(Z_1))
			-\varphi_1([Z_1,Z_2])\\
			&=Z_1\varphi_2(Z_2)-Z_2(\varphi_2(Z_1))-\varphi_2([Z_1,Z_2])\\
			&=(d\varphi_2)(Z_1,Z_2).\end{align*}
		For $\gamma\in X_1$, we have  $d\varphi_1(\gamma)(Z_1(\gamma),Z_2(\gamma))=
		d\varphi_2(\gamma)(Z_1(\gamma),Z_2(\gamma))$ and thus 
		\[(d\varphi_1)(\gamma)(\mc{H}_\gamma(v_1),\mc{H}_\gamma(v_2))=(d\varphi_2)(\gamma)(\mc{H}_\gamma(v_1),\mc{H}_\gamma(v_2)).\]
Since $\gamma, v_1, v_2$ were arbitrary, we have 	\[(d\varphi_1)(\gamma)(\mc{H}_\gamma(v_1),\mc{H}_\gamma(v_2))=(d\varphi_2)(\gamma)(\mc{H}_\gamma(v_1),\mc{H}_\gamma(v_2)),\]
		for all $\gamma\in X_1$ and $v_1,v_2\in T_\gamma X_1$.
		Thus, for a $1$-form $\varphi$ on the Lie groupoid $[X_1\rra X_0]$ with respect to an integrable connection $\mc{H}\subseteq TX_1$, the condition $\mc{H}(s^*\varphi)=\mc{H}(t^*\varphi)$ implies that $\mc{H}(s^*d\varphi)=\mc{H}(t^*d\varphi)$. Thus, $d\varphi$ is a $2$-form on the Lie groupoid $[X_1\rra X_0]$. Using similar arguments, it follows that, for a  $k$-form $\varphi$ on Lie groupoid $[X_1\rra X_0]$, the exterior derivative $d\varphi$ of $\varphi$, is a  $k+1$-form on the Lie groupoid $\mb{X}=[X_1\rra X_0]$.
	\end{proof}
\end{lemma} 

Now, our aim is to associate a $1$-form on the Lie groupoid $[s^*E_G\rra E_G]$ for a connection $\mc{D}$ on $(E_G\ra X_0,[X_1\rra X_0])$. Before we proceed to that, we will present some basic properties related to differential forms on a pair $(\mb{X},\mc{H})$ (\Cref{Proposition:pullbackisdifferrentialform}, \Cref{Proposition:pullbackofdifferrentialform}).

 Let $f:M\ra N$ be a smooth map and $\omega:N\ra \Lambda^kT^*N$ a  $k$-form on $N$. Then, we have the notion of the pull-back of $\omega$ along $f$, denoted by $f^*\omega$, is a  $k$-form on $M$, defined as,
 \[(f^*\omega)(m)(v_1,\cdots,v_k)=\omega(f(m))(f_{*,m}(v_1),\cdots,f_{*,m}(v_k)),\] for all $m\in M$, $v_1,\cdots, v_k\in T_mM$. In \Cref{Proposition:pullbackisdifferrentialform} we will see that for a special type of morphism $(F_1,F_0):(\mb{X},\mc{H}_X)\ra (\mb{Y},\mc{H}_Y)$ (\Cref{Definition:compatible with connections}) the pull-back of a $k$-form $\omega$ on $(\mb{Y},\mc{H}_Y)$, along the map $F_0:X_0\ra Y_0$ is a $k$-form on $(\mb{X},\mc{H}_X)$.
 
 \begin{definition}\label{Definition:compatible with connections}
Let $\mb{X}=[X_1\rra X_0]$ be a Lie groupoid equipped with an integrable connection $\mc{H}_X\subseteq TX_1$, and $\mb{Y}=[Y_1\rra Y_0]$ a Lie groupoid equipped with an integrable connection $\mc{H}_Y\subseteq TY_1$. We say that a morphism of Lie groupoids $(F_1,F_0):[X_1\rra X_0]\ra [Y_1\rra Y_0]$ is \textit{compatible with connections $\mc{H}_X,\mc{H}_Y$}, if, 
\[(F_1)_{*,\gamma}(\mc{H}_X(v))= \mc{H}_Y((F_1)_{*,\gamma}(v)),\]
for all $\gamma\in X_1$ and $v\in T_\gamma X_1$. 
 \end{definition}

 \begin{proposition}\label{Proposition:pullbackisdifferrentialform}
Let $\mb{X}=[X_1\rra X_0], \mb{Y}=[Y_1\rra Y_0]$ be Lie groupoids equipped with  integrable connections $\mc{H}_X\subseteq TX_1, \mc{H}_Y\subseteq TY_1$ respectively, and  $(F_1,F_0):[X_1\rra X_0]\ra [Y_1\rra Y_0]$ a morphism of Lie groupoids compatible with connections $\mc{H}_X$ and $\mc{H}_Y$. Let $\omega$ be a  $k$-form on $Y_0$, and $F_0^*\omega$ the pull-back of $\omega$ along the map $F_0:X_0\ra Y_0$. If $\omega$ is a  $k$-form on $(\mb{Y},\mc{H}_Y)$, then, $F_0^*\omega$ is a  $k$-form on $(\mb{X},\mc{H}_X)$.
\begin{proof}
We will prove this for $k=1$. 
The proof for an arbitrary $k$ is almost similar. 
 As $\omega$ is a $1$-form on $(\mb{Y},\mc{H}_Y)$, we have $\mc{H}_Y(s_Y^*\omega)=\mc{H}_Y(t_Y^*\omega)$. Consider the differential form $F_0^*\omega:X_0\ra \Lambda^1T^*X_0$. We prove that $F_0^*\omega$ is a $1$-form on $(\mb{X},\mc{H}_X)$; that is, $\mc{H}_X(s_X^*(F_0^*\omega))
=\mc{H}_X(t_X^*(F_0^*\omega))$.

Let $\gamma \in X_1$ and $v\in T_\gamma X_1$. We have 
\begin{align*}
\mc{H}_X(s_X^*(F_0^*\omega))(\gamma)(v)
&=(s_X^*(F_0^*\omega))(\gamma)(\mc{H}_\gamma(v))\\
&=((F_0\circ s_X)^*\omega)(\gamma)(\mc{H}_\gamma(v))\\
&=((s_Y\circ F_1)^*\omega)(\gamma)(\mc{H}_\gamma(v))\\
&=(F_1^*(s_Y^*\omega))(\gamma)(\mc{H}_\gamma(v))\\
&=(s_Y^*\omega)(F_1(\gamma))((F_1)_{*,\gamma}(\mc{H}_\gamma(v)))\\
&=(s_Y^*\omega)(F_1(\gamma))((\mc{H}_{F_1(\gamma)}((F_1)_{*,\gamma}(v))))\\
&=\mc{H}_Y(s_Y^*\omega)(F_1(\gamma))((F_1)_{*,\gamma}(v))\\
&=\mc{H}_Y(t_Y^*\omega)(F_1(\gamma))((F_1)_{*,\gamma}(v))\\
&=\mc{H}_X(t_X^*(F_0^*\omega))(\gamma)(v).
\end{align*}
Thus, we have $\mc{H}_X(s_X^*(F_0^*\omega))(\gamma)(v)=\mc{H}_X(t_X^*(F_0^*\omega))(\gamma)(v)$, for all $\gamma\in X_1$ and $v\in T_\gamma X_1$. So, $F_0^*\omega$ is a  $1$-form on $(\mb{X},\mc{H}_X)$.
\end{proof}
 \end{proposition}
  We would need a converse of the above result. Suppose that we have the setup as in \Cref{Proposition:pullbackisdifferrentialform}. One could ask the following question: does $F_0^*\omega$ being a differential form on $(\mb{X},\mc{H}_X)$ imply that $\omega$ is a differential form on $(\mb{Y},\mc{H}_Y)$. The answer to this question is, yes, if $F_1:X_1\ra Y_1$ is a surjective submersion.
  
   \begin{proposition}\label{Proposition:pullbackofdifferrentialform}
   	Let $\mb{X}=[X_1\rra X_0], \mb{Y}=[Y_1\rra Y_0]$ be Lie groupoids equipped with  integrable connections $\mc{H}_X\subseteq TX_1, \mc{H}_Y\subseteq TY_1$ respectively, and  $(F_1,F_0):[X_1\rra X_0]\ra [Y_1\rra Y_0]$ a morphism of Lie groupoids compatible with connections $\mc{H}_X$ and $\mc{H}_Y$. Let $\omega$ be a $k$-form on $Y_0$, and $F_0^*\omega$ the pull-back of $\omega$ along the map $F_0:X_0\ra Y_0$. Suppose further that $F_1:X_1\ra Y_1$ is a surjective submersion. If $F_0^*\omega$ is a $k$-form on $(\mb{X},\mc{H}_X)$, then, $\omega$ is a $k$-form on $(\mb{Y},\mc{H}_Y)$.
   	\begin{proof}
We will prove this for $k=1$. 
For an arbitrary $k$, the proof is  same. As $F_0^*\omega$ is a $1$-form on $(\mb{X},\mc{H}_X)$, we have $\mc{H}_X(s_X^*(F_0^*\omega))=
\mc{H}_X(t_X^*(F_0^*\omega))$. We prove that $\omega$ is a differential form on $(\mb{Y},\mc{H}_Y)$; that is, $\mc{H}_Y(s_Y^*\omega)=\mc{H}_Y(t_Y^*\omega)$.

Let $\gamma\in Y_1$, and $v\in T_\gamma Y_1$. As $F_1:X_1\ra Y_1$ is a surjective submersion, for $\gamma$ and $v$, there exists $\tilde{\gamma}\in X_1$ and 
$\tilde{v}\in T_{\tilde{\gamma}}X_1$ such that $F_1(\tilde{\gamma})=\gamma$ and 
$(F_1)_{*,\tilde{\gamma}}(\tilde{v})=v$. We have 
\begin{align*}
\mc{H}_Y(s_Y^*\omega)(\gamma)(v)=
&\mc{H}_Y(s_Y^*\omega)(F_1(\tilde{\gamma}))((F_1)_{*,\tilde{\gamma}}(\tilde{v}))\\
&=(s_Y^*\omega)(F_1(\tilde{\gamma}))(\mc{H}_{F_1(\tilde{\gamma})}((F_1)_{*,\tilde{\gamma}}(\tilde{v})))\\
&=(s_Y^*\omega)(F_1(\tilde{\gamma}))
((F_1)_{*,\tilde{\gamma}}(\mc{H}_{\tilde{\gamma}}(\tilde{v})))\\
&=(F_1^*\circ s_Y^*)(\omega)(\tilde{\gamma})(\mc{H}_{\tilde{\gamma}}(\tilde{v}))\\
&=(s_X^*\circ F_0^*)(\omega)(\tilde{\gamma})(\mc{H}_{\tilde{\gamma}}(\tilde{v}))\\
&=\mc{H}_X(s^*_X (F_0^*\omega))(\tilde{\gamma})(\tilde{v})\\
&=\mc{H}_X(t^*_X (F_0^*\omega))(\tilde{\gamma})(\tilde{v})\\
&=\mc{H}_Y(t_Y^*\omega)(\gamma)(v).
\end{align*}
Thus, for each $\gamma\in Y_1$ and $v\in T_\gamma Y_1$, we have $\mc{H}_Y(s_Y^*\omega)(\gamma)(v)=\mc{H}_Y(t_Y^*\omega)(\gamma)(v)$. So, $\omega$ is a  $1$-form on $(\mb{Y},\mc{H}_Y)$. 
   	\end{proof}
   \end{proposition}
Some other interesting properties of differential forms on $(\mb{X},\mc{H}_X)$ will be discussed in later sections. 
\section{$1$-form associated to a connection on $(E_G\ra X_0,[X_1\rra X_0])$}\label{Section:differentialformsassociatedtoconnection}
Let $[X_1\rra X_0] $ be a Lie groupoid equipped with an integrable connection $\mc{H}\subseteq TX_1$, $(E_G\ra X_0,[X_1\rra X_0])$ a principal bundle, and $\mc{D}:(TE_G)/G\ra (E_G\times \mf{g})/G$ a connection on $(E_G\ra X_0,[X_1\rra X_0])$.

Viewing $\mc{D}$ as a morphism of vector bundles over the manifold $X_0$, we get a differential form $\omega:E_G\ra \Lambda^1_{\mf{g}}T^*E_G$, with the property that
$\widetilde{\mc{D}}(v)=(a,\omega(a)(v))$ for all $a\in E_G, v\in T_aE_G$. Here $\widetilde{\mc{D}}:TE_G\ra E_G\times \mf{g}$ is such that $q^{P\times \mf{g}}\circ \widetilde{\mc{D}}=\mc{D}\circ q^{TP}$ (\Cref{Lemma:propertiesofomega}).

Considering $\mc{D}$ as a morphism of vector bundles over $[X_1\rra X_0]$, we have \[\mc{D}(\gamma\cdot [v])=\gamma \cdot \mc{D}([v])\] for all $(\gamma,[v])\in X_1\times (TE_G)/G$. We have \[
\mc{D}(\gamma\cdot [v])=\mc{D}([\tilde{\theta}_{(\gamma,a)}(v)])=\mc{D}(q^{TP}(\tilde{\theta}_{(\gamma,a)}(v)))=(\mc{D}\circ q^{TP})(\tilde{\theta}_{(\gamma,a)}(v))
=q^{P\times \mf{g}}(\widetilde{\mc{D}}(\tilde{\theta}_{(\gamma,a)}(v))),\] and 
\[\gamma\cdot \mc{D}([v])=\gamma\cdot (\mc{D}(q^{TP}(v)))=\gamma\cdot (\mc{D}\circ q^{TP})(v)
=\gamma\cdot q^{P\times \mf{g}}(\widetilde{\mc{D}}(v))
=q^{P\times \mf{g}}(\gamma\cdot \widetilde{\mc{D}}(v)).\]
So, we have $q^{P\times \mf{g}}(\widetilde{\mc{D}}(\tilde{\theta}_{(\gamma,a)}(v)))=
q^{P\times \mf{g}}(\gamma\cdot \widetilde{\mc{D}}(v))$. Rewriting this in terms of differential forms we obtain
\[q^{P\times \mf{g}}(\gamma a, \omega(\gamma a) (\tilde{\theta}_{(\gamma,a)}(v)))=
q^{P\times \mf{g}}(\gamma\cdot (a,\omega(a)(v)))=q^{P\times \mf{g}}(\gamma a,\omega(a)(v)).\]
Thus, there exists unique $g\in G$ such that, 
\[(\gamma a, \omega(\gamma a) (\tilde{\theta}_{(\gamma,a)}(v)))=
(\gamma a g, {\rm{ad}}(g^{-1})\omega(a)(v)).\]
Comparing the coordinates, and noting that the action of $G$ on $E_G$ is free, 
we see that $g=1$. This would then imply that 
\begin{equation}\label{Equation:propertyofomega}
\omega(\gamma a)(\tilde{\theta}_{(\gamma,a)}(v))=\omega(a)(v),
\end{equation} for all $(\gamma,a)\in X_1\times_{X_0}E_G$ and $v\in T_aE_G$. 
The condition in \Cref{Equation:propertyofomega} implies, that, $\omega$ is a $1$-form on the Lie groupoid $[s^*E_G\rra E_G]$ (with respect to the connection $\widetilde{\mc{H}}$).

\begin{proposition}\label{Proposition:connectionformonLieGroupoid}
	Let $[X_1\rra X_0] $ be a Lie groupoid equipped with an integrable connection $\mc{H}\subseteq TX_1$, $(E_G\ra X_0,[X_1\rra X_0])$ a principal bundle, and $\mc{D}:(TE_G)/G\ra (E_G\times \mf{g})/G$ a connection on $(E_G\ra X_0,[X_1\rra X_0])$. Let $\omega:E_G\ra \Lambda^1_{\mf{g}}T^*E_G$ be the differential $1$-form associated to $\mc{D}$. Then, $\omega$ is a differential $1$-form on the Lie groupoid $[s^*E_G\rra E_G]$ (with respect to the connection $\widetilde{H}\subseteq Ts^*E_G$ on $[s^*E_G\rra  E_G]$).
	\begin{proof}
We prove that $\omega$ is a $1$-form on the Lie groupoid $[s^*E_G\rra E_G]$; that is, $\widetilde{\mc{H}}
({\rm pr}_2^*\omega)=\widetilde{\mc{H}}(\mu^*\omega)$.

 Let $(\gamma,a)\in s^*E_G=X_1\times_{X_0}E_G$, and $(v,r)\in T_\gamma X_1\times_{TX_0}T_aE_G$. We have,
 \begin{align*}
\widetilde{\mc{H}}({\rm pr}_2^*\omega) (\gamma,a)(v,r)
&=({\rm pr}_2^*\omega) (\gamma,a)(\widetilde{\mc{H}}_{(\gamma,a)}(v,r))\\
&=({\rm pr}_2^*\omega) (\gamma,a)(\mc{H}_\gamma(v),r))\\
&=\omega(a)(r), 
\end{align*}
and 
 \begin{align*}
	\widetilde{\mc{H}}(\mu^*\omega) (\gamma,a)(v,r)
	&=(\mu^*\omega) (\gamma,a)(\widetilde{\mc{H}}_{(\gamma,a)}(v,r))\\
	&=(\mu^*\omega) (\gamma,a)(\mc{H}_\gamma(v),r)\\
	&=\omega(\gamma a)( \mu_{*,(\gamma,a)} (\mc{H}_\gamma(v),r))\\
	&=\omega(\gamma a)(\tilde{\theta}_{(\gamma,a)}(r)).
\end{align*}
As $\omega(\gamma a)(\tilde{\theta}_{(\gamma,a)}(v))=\omega(a)(v)$ (\Cref{Equation:propertyofomega}), we conclude that,
\[	\widetilde{\mc{H}}(\mu^*\omega) (\gamma,a)(v,r)=\widetilde{\mc{H}}({\rm pr}_2^*\omega) (\gamma,a)(v,r),\] for all $(\gamma,a)\in s^*E_G=X_1\times_{X_0}E_G$, and $(v,r)\in T_\gamma X_1\times_{TX_0}T_aE_G$. Thus, $\omega$ is a differential $1$-form on the Lie groupoid $[s^*E_G\rra E_G]$ (with respect to the connection $\widetilde{H}\subseteq Ts^*E_G$ on $[s^*E_G\rra  E_G]$).
	\end{proof}
\end{proposition}
Working out the converse of the above result is straightforward backward calculation.

\begin{proposition}\label{Lemma:conerseofconnectionformonLieGroupoid}
	Let $[X_1\rra X_0] $ be a Lie groupoid equipped with an integrable connection $\mc{H}\subseteq TX_1$, $(E_G\ra X_0,[X_1\rra X_0])$ a principal bundle. Let $\mc{D}:(TE_G)/G\ra (E_G\times \mf{g})/G$ be a morphism of vector bundles over $X_0$, which is  a connection on $E_G(X_0,G)$, and $\omega:E_G\ra \Lambda^1_{\mf{g}}T^*E_G$ the associated differential $1$-form. If $\omega$ is a differential $1$-form on the Lie groupoid $[s^*E_G\rra E_G]$ (with respect to the connection $\widetilde{H}\subseteq Ts^*E_G$ on $[s^*E_G\rra  E_G]$), then $\mc{D}$ is a morphism of vector bundles over $[X_1\rra X_0]$, thus a connection on $(E_G\ra X_0,[X_1\rra X_0])$.
	\end{proposition}

\section{$2$-form associated to a connection on $(E_G\ra X_0,[X_1\rra X_0])$}\label{Section:differential2formsassociatedtoconnection}
Let $[X_1\rra X_0] $ be a Lie groupoid equipped with an integrable connection $\mc{H}\subseteq TX_1$, $(E_G\ra X_0,[X_1\rra X_0])$ a principal bundle, and $\mc{D}:(TE_G)/G\ra (E_G\times \mf{g})/G$ a connection on $(E_G\ra X_0,[X_1\rra X_0])$. 

Viewing $\mc{D}$ as a morphism of vector bundles over the manifold $X_0$, we get a $\mf{g}$-valued differential $2$-form $\mc{K}_{\mc{D}}:X_0\ra \Lambda^2_{\mf{g}}T^*X_0$ on $X_0$, which we called the curvature form associated to the connection $\mc{D}$ (\Cref{Equation:curvatureKDforclassicalprincipalbundle}). As in the case of $\omega$, we prove that $\mc{K}_{\mc{D}}$ is a differential $2$-form on the Lie groupoid $[X_1\rra X_0]$ (with respect to connection $\mc{H}\subseteq TX_1$ on $[X_1\rra X_0]$). As mentioned in \Cref{Proposition:propertyofKD}, we have
\[\pi^*\mc{K}_{\mc{D}}=d\omega-[\omega,\omega].\]

We have already seen that $\omega$ is a differential $1$-form on the Lie groupoid $[s^*E_G\rra E_G]$ (\Cref{Proposition:connectionformonLieGroupoid}). As $\mc{H}\subseteq TX_1$ is an integrable connection on $[X_1\rra X_0]$, the \Cref{Lemma:pullbackofintegrableconnection} says that, $\widetilde{\mc{H}}\subseteq T s^*E_G$ is an integrable connection on $[s^*E_G\rra E_G]$. The \Cref{Lemma:differentialisdifferentialform} says that $d\omega$ is a differential $2$-form on $[s^*E_G\rra E_G]$. If we can prove that $[\omega,\omega]$ is also a differential form on $[s^*E_G\rra E_G]$, then, it turns out that $\pi^*\mc{K}_{\mc{D}}$ (being a linear combination of differential forms) is a differential form on $[s^*E_G\rra E_G]$.
We need a result to conclude that $\mc{K}_{\mc{D}}$ is a $2$-form on $[X_1\rra X_0]$ from the observation that $\pi^*\mc{K}_{\mc{D}}$ is a $2$-form on $[s^*E_G\rra E_G]$. As an application of \Cref{Proposition:pullbackisdifferrentialform} and \Cref{Proposition:pullbackofdifferrentialform}, we give the following result which implies $\mc{K}_{\mc{D}}$ is a $2$-form on $[X_1\rra X_0]$ when we know that $\pi^*\mc{K}_{\mc{D}}$ is a $2$-form on $[s^*E_G\rra E_G]$.
\begin{lemma}\label{pullbackdifferentialformonLiegroupoid}
	Let $[X_1\rra X_0]$ be a Lie groupoid and 
	$(E_G\ra X_0,[X_1\rra X_0])$ a principal bundle. Let $\mc{H}\subseteq TX_1$ be a connection on the Lie groupoid $[X_1\rra X_0]$ and $\widetilde{\mc{H}}\subseteq T(s^*E_G)$ the pull-back connection on the Lie groupoid $[s^*E_G\rra E_G]$. Let $\tau\colon X_0\ra \Lambda^kT^*X_0$ be a differential $k$-form on the manifold $X_0$ and $
	\pi^*\tau\colon E_G\ra \Lambda^kT^*E_G$ the pull-back form on the manifold $E_G$ along the map $\pi\colon E_G\ra X_0$. Then, $\tau$ is a differential form on the Lie
	groupoid $\mb{X}=[X_1\rra X_0]$, if and only if, $\pi^*\tau$ is a differential $k$-form on the Lie groupoid $[s^*E_G\rra E_G]$.
	\begin{proof}
		Consider the morphism of Lie groupoids 
		$({\rm pr}_1,\pi)\colon [s^*E_G\rra E_G]\ra [X_1\rra X_0]$. As $\pi\colon E_G\ra X_0$ is a surjective submersion, its pull-back ${\rm pr}_1\colon s^*E_G\ra X_1$ is a surjective submersion. Let $(\gamma,a)\in s^*E_G$ and $(v,w)\in T_{(\gamma,a)}s^*E_G$. We have 
		\begin{align*}
			({\rm pr}_1)_{*,(\gamma,a)}(\widetilde{\mc{H}}_{(\gamma,a)}(v,w))
			&=({\rm pr}_1)_{*,(\gamma,a)}(\mc{H}_{\gamma}(v),w)\\
			&=\mc{H}_{\gamma}(v)\\
			&=\mc{H}_{{\rm pr}_1(\gamma,a)}(({\rm pr}_1)_{*,(\gamma,a)}(v,w)).
		\end{align*}
		Thus, $({\rm pr}_1)_{*,(\gamma,a)}(\widetilde{\mc{H}}_{(\gamma,a)}(v,w))
		=\mc{H}_{{\rm pr}_1(\gamma,a)}
		(({\rm pr}_1)_{*,(\gamma,a)}(v,w))$ for all $(\gamma,a)\in s^*E_G$ and $(v,w)\in T_{(\gamma,a)}s^*E_G$. Thus, the morphism of Lie groupoids $({\rm pr}_1,\pi)\colon [s^*E_G\rra E_G]\ra [X_1\rra X_0]$ satisfies the following conditions:
		\begin{enumerate}
			\item the map ${\rm pr}_1\colon s^*E_G\ra X_1$ is a surjective submersion,
			\item the map ${\rm pr}_1\colon s^*E_G\ra X_1$ is compatible with the connections $\widetilde{\mc{H}}$ and $\mc{H}$ in the sense that, $({{\rm pr}}_1)_{*,(\gamma,a)}(\widetilde{\mc{H}}_{(\gamma,a)}(v,w))
			=\mc{H}_{{\rm pr}_1(\gamma,a)}(({\rm pr}_1)_{*,(\gamma,a)}(v,w))$ for all $(\gamma,a)\in s^*E_G$ and $(v,w)\in T_{(\gamma,a)}s^*E_G$.
		\end{enumerate}
		Thus, $\tau$ is a differential form on the Lie
		groupoid $[X_1\rra X_0]$, if and only if, $\pi^*\tau$ is a differential $k$-form on the Lie groupoid $[s^*E_G\rra E_G]$ (\Cref{Proposition:pullbackisdifferrentialform} and  \Cref{Proposition:pullbackofdifferrentialform}).
	\end{proof}
\end{lemma}

Thus, to prove $\mc{K}_{\mc{D}}$ is a $2$-form on $[X_1\rra X_0]$, all it remains to prove is that $[\omega,\omega]$ is a $2$-form on $[X_1\rra X_0]$.

Let us recall the notion of the Lie bracket of Lie algebra valued differential forms on a manifold. Let $\mf{g}$ be a Lie algebra. Let $\omega$ and  $\eta$ be $\mf{g}$-valued differential $k$ form and differential $l$-form on a manifold $M$. Then, the Lie bracket $[\omega,\eta]$ is a $\mf{g}$-valued differential $k+l$ form on the manifold $M$ defined as follows:
\begin{align}
	\begin{aligned}\label{Equation:Liebracketofdifferentialforms}
		[\omega,\eta](a)&(v_1,\cdots,v_{k+l})\\
		&=\frac{1}{(k+l)!}\sum_{\sigma\in S_{k+l}} \text{sgn}(\sigma)[\omega(a)(v_{\sigma(1)},\cdots,v_{\sigma(k)}),
		\eta(a)(v_{\sigma(k+1)},\cdots,v_{\sigma(k+l)})]
\end{aligned} \end{align}
for all $a\in M$ and  $v_i\in T_aM$ for $1\leq i\leq k+l$.
\begin{lemma}\label{Lemma:liebracketisdifferentialformonLiegroupoid}
	Let $[X_1\rra X_0]$ be a Lie groupoid and 
	$(E_G\ra X_0,[X_1\rra X_0])$ a principal bundle. 
	Let $\mc{H}\subseteq TX_1$ be a connection on $[X_1\rra X_0]$ and $\widetilde{\mc{H}}\subseteq T(s^*E_G)$ the pullback connection on  $[s^*E_G\rra E_G]$. Let $\mc{D}
	$ be a connection on the principal bundle $(E_G\ra X_0,[X_1\rra X_0])$ and  $\omega\colon E_G\ra \Lambda^1_{\mf{g}}T^*E_G$ the associated connection $1$-form. Then,  the  Lie bracket $[\omega,\omega]\colon E_G\ra \Lambda^2_{\mf{g}}T^*E_G$ (\Cref{Equation:Liebracketofdifferentialforms}), is a differential $2$-form on the Lie groupoid $[s^*E_G\rra E_G]$.
	\begin{proof}
		The differential form $[\omega,\omega]\colon E_G\ra \Lambda^2_{\mf{g}}T^*E_G$, on the manifold $E_G$, is defined as 
		\[[\omega,\omega](a)(v_1,v_2)=[\omega(a)(v_1),\omega(a)(v_2)],\]
		for all $a\in E_G$ and $v_1,v_2\in T_aE_G$, where $[\omega(a)(v_1),\omega(a)(v_2)]$ is the Lie bracket in $\mf{g}$.
		
		To prove $[\omega,\omega]\colon E_G\ra \Lambda^2_{\mf{g}}T^*E_G$  is a differential $2$-form on the Lie groupoid $[s^*E_G\rra E_G]$, we  prove that  $\widetilde{\mc{H}}({\rm pr}_2^*([\omega,
		\omega]))=
		\widetilde{\mc{H}}(\mu^*([\omega,
		\omega]))$ as differential $2$-forms on the manifold $s^*E_G$; that is, 
		\[\widetilde{\mc{H}}({\rm pr}_2^*([\omega,
		\omega]))(\gamma,a)((v_1,r_1),(v_2,r_2))=\widetilde{\mc{H}}(\mu^*([\omega,
		\omega]))(\gamma,a)((v_1,r_1),(v_2,r_2))\] for all 
		$(\gamma,a)\in s^*E_G$ and $(v_1,r_1),(v_2,r_2)\in T_{(\gamma,a)}(s^*E_G)$.
		
		For $(\gamma,a)\in s^*E_G$ and $(v_1,r_1),(v_2,r_2)\in T_{(\gamma,a)}(s^*E_G)$, observe that \begin{align*}
			\widetilde{\mc{H}}({\rm pr}_2^*([\omega,
			\omega]))(\gamma,a)&((v_1,r_1),(v_2,r_2))\\
			&
			={\rm pr}_2^*([\omega,
			\omega])(\gamma,a)(\widetilde{\mc{H}}_{(\gamma,a)}(v_1,r_1),
			\widetilde{\mc{H}}_{(\gamma,a)}(v_2,r_2))\\
			&={\rm pr}_2^*([\omega,
			\omega])(\gamma,a)((\mc{H}_{\gamma}(v_1),r_1),
			(\mc{H}_{\gamma}(v_2),r_2))\\
			&=[\omega,\omega]({\rm pr}_2(\gamma,a))(({\rm pr}_2)_{*,(\gamma,a)}
			(\mc{H}_{\gamma}(v_1),r_1),
			({\rm pr}_2)_{*,(\gamma,a)}
			(\mc{H}_{\gamma}(v_2),r_2))\\
			&=[\omega,\omega](a)(r_1,r_2)\\
			&=[\omega(a)(r_1),\omega(a)(r_2)],
		\end{align*} 
		
		and that, 
		\begin{align*}
			\widetilde{\mc{H}}(\mu^*([\omega,
			\omega]))(\gamma,a)&((v_1,r_1),(v_2,r_2))\\
			&
			=\mu^*([\omega,
			\omega])(\gamma,a)(\widetilde{\mc{H}}_{(\gamma,a)}(v_1,r_1),
			\widetilde{\mc{H}}_{(\gamma,a)}(v_2,r_2))\\
			&=\mu^*([\omega,
			\omega])(\gamma,a)((\mc{H}_{\gamma}(v_1),r_1),
			(\mc{H}_{\gamma}(v_2),r_2))\\
			&=[\omega,\omega](\mu(\gamma, a))(\mu_{*,(\gamma,a)}
			(\mc{H}_{\gamma}(v_1),r_1),
			\mu_{*,(\gamma,a)}
			(\mc{H}_{\gamma}(v_2),r_2))\\
			&=[\omega,\omega](\gamma a)(\widetilde{\theta}_{(\gamma,a)}(r_1),\widetilde{\theta}_{(\gamma,a)}(r_2))\\
			&=[\omega(\gamma a)(\widetilde{\theta}_{(\gamma,a)}(r_1)),\omega(\gamma a)(\widetilde{\theta}_{(\gamma,a)}(r_2))].
		\end{align*}
		From \Cref{Equation:propertyofomega} we see that  $\omega(\gamma a)(\widetilde{\theta}_{(\gamma,a)}(r))=\omega(a)(r)$ for all $(\gamma,a)\in s^*E_G$ and for  all $(v,r)\in T_{(\gamma,a)}(s^*E_G)$; in particular, 
		\begin{align*}
			\widetilde{\mc{H}}(\mu^*([\omega,
			\omega]))(\gamma,a)((v_1,r_1),(v_2,r_2))
			&=[\omega(\gamma a)(\widetilde{\theta}_{(\gamma,a)}(r_1)),\omega(\gamma a)(\widetilde{\theta}_{(\gamma,a)}(r_2))]\\
			&
			=[\omega(a)(r_1),\omega(a)(r_2)]\\
			&=\widetilde{\mc{H}}({\rm pr}_2^*([\omega,
			\omega]))(\gamma,a)((v_1,r_1),(v_2,r_2)).
		\end{align*}
		Thus, \[\widetilde{\mc{H}}(\mu^*([\omega,
		\omega]))(\gamma,a)((v_1,r_1),(v_2,r_2))=\widetilde{\mc{H}}({\rm pr}_2^*([\omega,
		\omega]))(\gamma,a)((v_1,r_1),(v_2,r_2)),\] for all 
		$(\gamma,a)\in s^*E_G$ and $(v_1,r_1),(v_2,r_2)\in T_{(\gamma,a)}(s^*E_G)$. Thus, $\widetilde{\mc{H}}(\mu^*([\omega,
		\omega]))=\widetilde{\mc{H}}({\rm pr}_2^*([\omega,
		\omega]))$  as a differential $2$-form on the manifold $s^*E_G$; that is, $[\omega,\omega]\colon E_G\ra \Lambda^2_{\mf{g}}T^*E_G$ is a differential $2$-form on the Lie groupoid $[s^*E_G\rra E_G]$.
	\end{proof}
\end{lemma}

\begin{lemma}\label{CurvatureformIsDifferentialformOnLiegroupoid}
Let $[X_1\rra X_0]$ be a Lie groupoid and 
$(E_G\xra{\pi} X_0,[X_1\rra X_0])$ a principal bundle. Let $\mc{H}\subseteq TX_1$ be a connection on the Lie groupoid $[X_1\rra X_0]$ and $\widetilde{\mc{H}}\subseteq T(s^*E_G)$ the pull-back connection on the Lie groupoid $[s^*E_G\rra E_G]$. 
Let  $\mc{D}$ be a connection on  
	$(E_G\ra X_0,[X_1\rra X_0])$, and $\mc{K}_{\mc{D}}$ the curvature $2$-form of the connection $\mc{D}$ on the underlying principal
	$G$-bundle $E_G \ra X_0$. Then $\mc{K}_{\mc{D}}$
	is a differential $2$-form on the Lie groupoid $[X_1\rra X_0]$. 
\end{lemma}

\section{Chern-Weil map associated to a principal bundle $(E_G\ra X_0,[X_1\rra X_0])$}\label{Section:ChernWeilforEGX0X1X0}
Let $[X_1\rra X_0]$ be a Lie groupoid and 
$(E_G\xra{\pi} X_0,[X_1\rra X_0])$ a principal bundle. Let $\mc{H}\subseteq TX_1$ be a connection on the Lie groupoid $[X_1\rra X_0]$ and $\widetilde{\mc{H}}\subseteq T(s^*E_G)$ the pull-back connection on the Lie groupoid $[s^*E_G\rra E_G]$. Now that we have the notion of curvature form $\mc{K}_{\mc{D}}$ associated to a connection $\mc{D}$ on principal bundle $(E_G\ra X_0,[X_1\rra X_0])$, we will introduce the notion of Chern-Weil map.

In \Cref{ChernWeiloversmooth manifold}, we have discussed the Chern-Weil map associated to a principal bundle over a manifold. Recall that, for a principal bundle $P(M,G)$, the Chern-Weil map is a morphism of $\mb{R}$-algebras $\text{sym}(\mf{g})^G\ra H^\bullet_{\rm{{\rm dR}}}(M;\mb{R})$. This is defined by fixing a connection $\omega$ on $P(M,G)$ and considering the associated curvature form $\Omega:P\ra \Lambda^2_{\mf{g}}T^*P$. We will follow a slight variant of the procedure. Instead of using the $2$-form $\Omega$ on $P$, we use the $2$-form $\mc{K}_{\mc{D}}$ on $M$ associated to a connection $\mc{D}$ (the splitting associated to the connection $1$-form $\omega$).

Let $f\in \rm{sym}_k(\mf{g})^G$. Then, we have the differential form $f(\Omega):P\ra \Lambda^2_{\mf{g}}T^*P$ (\Cref{definitionoffOmega}). We then see that, for this $2$-form $f(\Omega)$ on $P$, there exists a unique closed $2k$-form $\widetilde{f(\Omega)}$ such that $\pi^*(\widetilde{f(\Omega)})=f(\Omega)$. This defined the map ${\rm{sym}}_k(\mf{g})^G\ra H_{\rm{{\rm dR}}}^{2k}(M; \mb{R})$ as $f\mapsto [\widetilde{f(\Omega)}]$.  Using the Atiyah sequence approach, we have associated a $2$-form $\mc{K}_{\mc{D}}$ on $M$, whose pull-back $\pi^*\mc{K}_{\mc{D}}$ is the $2$-form on $P$, usually denoted by $\Omega$. 

Consider the $2k$-form $f(\mc{K}_{\mc{D}})$ on $M$. We have
\begin{align*}
	f(\Omega)&=\pi^*(\widetilde{f(\Omega)}),\\
	f(\pi^*\mc{K}_{\mc{D}})&=\pi^*(\widetilde{f(\Omega)}),\\
	\pi^*(f(\mc{K}_{\mc{D}}))&=\pi^*(\widetilde{f(\Omega)}).
\end{align*} 

As $\pi:P\ra M$ is a surjective submersion, and $\widetilde{f(\Omega)}$ is closed, the condition $\pi^*(f(\mc{K}_{\mc{D}}))=\pi^*(\widetilde{f(\Omega)})$  imply that $f(\mc{K}_{\mc{D}})$ is closed. As $\widetilde{f(\Omega)}$ is the unique closed $2k$-form on $M$ such that $\pi^*(\widetilde{f(\Omega)})=
f(\Omega)$, the condition that $f(\mc{K}_{\mc{D}})$   is closed and $\pi^*(f(\mc{K}_{\mc{D}}))=f(\Omega)$ imply that $\widetilde{f(\Omega)}=
f(\mc{K}_{\mc{D}})$. Thus, to define the Chern-Weil map, we can skip one step if we are using the notion of Atiyah sequence approach of connection forms and curvature forms. 

\subsection{de Rham cohomology of the pair $(\mb{X},\mc{H})$}\label{Section:deRhamcohomologyofLiegroupoidwithConnection} To define the notion of Chern-Weil map associated to a principal bundle $(E_G\ra X_0,[X_1\rra X_0])$, we need the notion of de Rham cohomology algebra associated to a pair $(\mb{X},\mc{H})$, which we will introduce in this section.

For a manifold $M$, after defining the notion of differential forms on $M$, we associate an $\mb{R}$-algebra, namely, the de Rham cohomology ring $H^*_{{\rm dR}}(M;\mb{R})=\bigoplus_{k=0}^{\dim (M)}H_{{\rm dR}}^k(M;\mb{R})$. In this subsection, for a Lie groupoid $\mb{X}=[X_1\rra X_0]$, equipped with an integrable connection  $\mc{H}\subseteq TX_1$, we associate, what we call, the de Rham cohomology of the pair $(\mb{X};\mc{H})$, denoted by $H^*_{{\rm dR}}(\mb{X};\mc{H})$. When there is no confusion regarding the connection $\mc{H}\subseteq TX_1$ on $\mb{X}$, we simply  call this the de Rham cohomology of the Lie groupoid $\mb{X}=[X_1\rra X_0]$.

Let $\mathcal{DF}^k_{{\rm dR}}(\mb{X},\mc{H})$ denote the vector space of differential $k$-forms on the Lie groupoid $\mb{X}$, with respect to the connection $\mc{H}\subseteq TX_1$. For a differential $k$-form $\varphi$ on the Lie groupoid $[X_1\rra X_0]$, we have seen that, the exterior derivative $d\varphi$ is a differential $k+1$-form on the Lie groupoid $[X_1\rra X_0]$. This gives a cochain complex,
\[\cdots\ra \mathcal{DF}^{k-1}_{{\rm dR}}(\mb{X},\mc{H})\ra \mathcal{DF}^k_{{\rm dR}}(\mb{X},\mc{H})\ra
\mathcal{DF}^{k+1}_{{\rm dR}}(\mb{X},\mc{H})\ra \cdots .\]
The $k^{th}$ cohomology group of this cochain complex, denoted by $H^k_{{\rm dR}}(\mb{X},\mc{H})$ is called the \textit{$k^{th}$-de Rham cohomology of the pair $(\mb{X},\mc{H})$}. Similar to that of the de Rham cohomology ring of a manifold, we can form the \textit{de Rham cohomology ring of the pair $(\mb{X},\mc{H})$}, denoted by $H^*_{{\rm dR}}(\mb{X},\mc{H})$ and defined as 
\begin{equation}\label{cohomologyofLiegroupoid}H_{{\rm dR}}^*(\mb{X},\mc{H})=\bigoplus_{k=0}^{\dim X_0} H^k_{{\rm dR}}(\mb{X},\mc{H}).\end{equation}
\begin{remark}\label{Remark:comparisionwithdeRhamcohomology}
	The paper \cite{MR2270285} defines the de Rham cohomology (groups/ring) of a Lie groupoid in a different way; using the notion of simplicial manifold associated to the Lie groupoid and the corresponding double complex $\Omega^{\bullet}(\Gamma_{\bullet})$. 
\end{remark}

\begin{example}
	Let $M$ be a manifold and $[M\rra M]$ the associated Lie groupoid. Fix the connection $\mc{H}_mM=T_mM$. Let $\varphi$ be a differential $k$-form on the manifold $M$ (this $M$ is seen as object set of the Lie groupoid $[M\rra M]$). As $s,t\colon M\ra M$ are identity maps, we have $s^*\varphi=\varphi$ and $t^*\varphi=\varphi$. As $\mc{H}=TM$, we have $\mc{H}(\varphi)=\varphi$. So, we have $\mc{H}(s^*\varphi)
	=\mc{H}(t^*\varphi)$ for all differential $k$-forms $\varphi$ on $M$. Thus, the $k$-th de Rham cohomology group $H^k_{{\rm dR}}(M;\mb{R})$ is equal to the $k$-th de Rham cohomology group of the pair $([M\rra M],\mc{H})$.
\end{example}

\begin{example}
	Let $G$ be a Lie group and $M$ a manifold. Let $\mu\colon M\times G\ra M$ be an action of the Lie group on the manifold $M$ and $[M\times G\rra M]$ the associated Lie groupoid. Fix the connection $\mc{H}_{(m,g)}=T_mM$. Then, differential $k$-forms on the Lie groupoid $[M\times G\rra M]$ are the $G$-invariant differential forms on $M$. So, $H^k_{{\rm dR}}(\mb{X};\mc{H})\cong H^k(\Omega(M)^G)$. 
\end{example}

\subsection{Chern-Weil map of the principal bundle $(E_G\ra X_0,[X_1\rra X_0])$}\label{Section:ChernWeilMap}
Now that we have the notion of de Rham cohomology algebra associated to a pair $(\mb{X},\mc{H})$, we can associate the Chern-Weil map as a morphism of $\mb{R}$-algebras $\text{sym}(\mf{g})^G\ra H^*_{\rm{{\rm dR}}}(\mb{X};\mc{H})$.
\begin{lemma}[{\cite[Theorem $5.3$]{MR3150770}}]\label{f(K_D)isadifferentialform}
	Let $\mb{X}=[X_1\rra X_0]$ be a Lie groupoid equipped with an integrable connection $\mc{H}\subseteq TX_1$. Let $(E_G\ra X_0,[X_1\rra X_0])$ be a principal bundle. Let $\mc{D}$ be a connection on the principal bundle $(E_G\ra X_0,[X_1 \rra X_0])$ and $\mc{K}_{\mc{D}}$ the associated curvature $2$-form on the Lie groupoid $[X_1\rra X_0]$. Then, $f(\mc{K}_{\mc{D}})$ is a differential $2k$-form on the Lie groupoid $[X_1\rra X_0]$. 
	\begin{proof}
		We prove this when $k=1$. The general result follows along the same lines of argument. We prove that, $\mc{H}(s^*(f(\mc{K}_{\mc{D}})))=
		\mc{H}(t^*(f(\mc{K}_{\mc{D}})))$.
		
		Let $\gamma\in X_1$ and $v_1,v_2\in T_\gamma X_1$. We have 
		\begin{align*}
			\mc{H}(s^*(f(\mc{K}_{\mc{D}})))(\gamma)(v_1,v_2)&=
			s^*(f(\mc{K}_{\mc{D}}))(\gamma)(\mc{H}_\gamma (v_1),\mc{H}_{\gamma}(v_2))\\
			&= (f(\mc{K}_{\mc{D}}))(s(\gamma))(s_{*,\gamma}(\mc{H}_\gamma (v_1)),s_{*,\gamma}(\mc{H}_{\gamma}(v_2))\\
			&=f(\mc{K}_{\mc{D}}(s(\gamma))(s_{*,\gamma}(\mc{H}_\gamma (v_1)),s_{*,\gamma}(\mc{H}_{\gamma}(v_2)))
		\end{align*}
		As $\mc{K}_{\mc{D}}\colon X_0\ra \Lambda^2_{\mf{g}}T^*X_0$ is a differential $2$-form on the Lie groupoid $[X_1\rra X_0]$, we see that $\mc{H}(s^*\mc{K}_{\mc{D}})=
		\mc{H}(t^*\mc{K}_{\mc{D}})$ as differential $2$-forms on the manifold $X_1$; that is, for each $\gamma\in X_1$ and $v_1,v_2\in T_{\gamma}X_1$, we have \[\mc{K}_{\mc{D}}(s(\gamma))(s_{*,\gamma}(\mc{H}_\gamma (v_1)),s_{*,\gamma}(\mc{H}_{\gamma}(v_2))
		=\mc{K}_{\mc{D}}(t(\gamma))(t_{*,\gamma}(\mc{H}_\gamma (v_1)),t_{*,\gamma}(\mc{H}_{\gamma}(v_2)).\]
		Thus, for $\mc{H}(s^*(f(\mc{K}_{\mc{D}})))\colon X_1\ra \Lambda^2_{\mf{g}}T^*X_1$, we have 
		\begin{align*}
			\mc{H}(s^*(f(\mc{K}_{\mc{D}})))(\gamma)(v_1,v_2)&=
			f(\mc{K}_{\mc{D}}(s(\gamma))(s_{*,\gamma}(\mc{H}_\gamma (v_1)),s_{*,\gamma}(\mc{H}_{\gamma}(v_2)))\\
			&=f(\mc{K}_{\mc{D}}(t(\gamma))(t_{*,\gamma}(\mc{H}_\gamma (v_1)),t_{*,\gamma}(\mc{H}_{\gamma}(v_2)))\\
			&=\mc{H}(t^*(f(\mc{K}_{\mc{D}})))(\gamma)(v_1,v_2).
		\end{align*}
		Thus, $\mc{H}(s^*(f(\mc{K}_{\mc{D}})))(\gamma)(v_1,v_2)=\mc{H}(t^*(f(\mc{K}_{\mc{D}})))(\gamma)(v_1,v_2)$ for  $\gamma\in X_1$ and $v_1,v_2\in T_{\gamma}X_1$; that is $f(\mc{K}_{\mc{D}})$ is a differential $2$-form on the Lie groupoid $[X_1\rra X_0]$. For an arbitrary $k\in \mb{N}$ and $f\in \text{sym}_k(\mf{g})^G$, the differential $2k$-form on the manifold $X_0$ is a differential $2k$-form on the Lie groupoid $[X_1\rra X_0]$. 
	\end{proof}
\end{lemma}
Thus, given a connection $\mc{D}$ on the  principal bundle $(E_G\ra X_0,[X_1\rra X_0])$, we have associated a morphism of $\mb{R}$-algebras 
$\text{sym}(\mf{g})^G\ra H^*_{\rm{{\rm dR}}}(\mb{X};\mc{H})$, with $f\in \text{sym}_k(\mf{g})^G$ mapping to $f(\mc{K}_{\mc{D}})$. Here,  $H^*_{\rm{{\rm dR}}}(\mb{X};\mc{H})$ is the de Rham cohomology of the Lie groupoid $\mb{X}$ equipped with the integrable connection $\mc{H}$ on the Lie groupoid $\mb{X}$.

Now, we see that, this map, $\text{sym}(\mf{g})^G\ra H^*_{\rm{{\rm dR}}}(\mb{X};\mc{H})$, does not depend on the choice of connection on the principal bundle $(E_G\ra X_0,[X_1\rra X_0])$. Let $\mc{D}$ and $\mc{D}'$ be connections on the  principal bundle $(E_G\ra X_0,[X_1\rra X_0])$ and $\mc{K}_{\mc{D}}$ and $\mc{K}_{\mc{D}'}$ be the associated curvature $2$-forms on the Lie groupoid $\mb{X}=[X_1\rra X_0]$. Given an element $f\in \text{sym}_{k}(\mf{g})^G$, we prove that $[f(\mc{K}_{\mc{D}})]=[f(\mc{K}_{\mc{D}'})]\in \mc{H}^{2k}(\mb{X};\mc{H})$; that is, there exists a $2k-1$ form $\widetilde{\Phi}$ on the Lie groupoid $\mb{X}=[X_1\rra X_0]$; such that $f(\mc{K}_{\mc{D}})-f(\mc{K}_{\mc{D}'})=d\widetilde{\Phi}$.

Let $\omega,\omega'$ be the associated connection $1$-forms for connections $\mc{D}$ and $\mc{D}'$ respectively.  Let $\alpha=\omega-\omega'$ and $\alpha_t=\omega'+t\alpha$ for $0\leq t\leq 1$. Observe that $\alpha_t$ is a connection $1$-form on the principal bundle $E_G\ra X_0$. Let $\Omega_t$ be the associated curvature $2$-forms. Consider the differential form 
\[\Phi=k\int_0^1f(\alpha_t,\Omega_t,\Omega_t,\cdots,\Omega_t)dt \]
as mentioned in \Cref{definitionofPhi}. It is known that, this $2k-1$ form $\Phi$ on $E_G$ projects uniquely, as $\Phi=\pi^*(\widetilde{\Phi})$, to a $2k-1$ form on $X_0$ and that $f(\Omega)-f(\Omega')=d\Phi$. So, we have 
\begin{align*}
	f(\Omega)-f(\Omega')&=d\Phi\\
	f(\pi^*\mc{K}_{\mc{D}})-f(\pi^*(\mc{K}_{\mc{D}'}))&=d(\pi^*(\widetilde{\Phi}))\\
	\pi^*(f(\mc{K}_{\mc{D}})-f(\mc{K}_{\mc{D}'}))&=\pi^*(d(\widetilde{\Phi}))
\end{align*} 
As $\pi\colon E_G\ra X_0$ is a surjective submersion, the condition $\pi^*(f(\mc{K}_{\mc{D}})-f(\mc{K}_{\mc{D}'}))=\pi^*(d(\widetilde{\Phi}))$ implies that $f(\mc{K}_{\mc{D}})-f(\mc{K}_{\mc{D}'})=d\widetilde{\Phi}$. To prove that $\widetilde{\Phi}$ is a differential form on the Lie groupoid $[X_1\rra X_0]$ it suffices (\Cref{pullbackdifferentialformonLiegroupoid}) to prove that $\Phi$ is a differential form on the Lie groupoid $[s^*E_G\rra E_G]$; that is, $\widetilde{\mc{H}}({\rm pr}_2^*\Phi)=\widetilde{\mc{H}}(\mu^*\Phi)$.
Observe that 
\begin{align*}
	{\rm pr}_2^*\Phi&=
	{\rm pr}_2^*\left(k\int_0^1f(\alpha_t,\Omega_t,\Omega_t,\cdots,\Omega_t)dt\right)\\
	&=\left(k\int_0^1f({\rm pr}_2^*(\alpha_t),{\rm pr}_2^*(\Omega_t),\cdots,{\rm pr}_2^*(\Omega_t))dt\right) 
\end{align*}
So, we have 
\begin{align*}
	\widetilde{\mc{H}}({\rm pr}_2^*\Phi)&=
	k\int_0^1 f(\mc{H}({\rm pr}_2^*\alpha_t),\cdots,\mc{H}({\rm pr}_2^*(\Omega_t)))dt.
\end{align*}
As $\alpha_t,\Omega_t$ are differential forms on the Lie groupoid $[s^*E_G\rra E_G]$, we have $\mc{H}({\rm pr}_2^*\alpha_t)=\mc{H}(\mu^*\alpha_t)$ and $\mc{H}({\rm pr}_2^*(\Omega_t))=\mc{H}(\mu^*(\Omega_t))$. So, we have 
\begin{align*}
	\widetilde{\mc{H}}({\rm pr}_2^*\Phi)&=
	k\int_0^1 f(\mc{H}({\rm pr}_2^*\alpha_t),\cdots,\mc{H}({\rm pr}_2^*(\Omega_t)))dt\\
	&=k\int_0^1 f(\mc{H}(\mu^*\alpha_t),\cdots,\mc{H}(\mu^*(\Omega_t)))dt\\
	&=\mc{H}\left(k\int_0^1 f((\mu^*\alpha_t),\cdots,(\mu^*(\Omega_t)))dt\right)\\
	&=\mc{H}\left(\mu^*\left(k\int_0^1 f( \alpha_t,\cdots, (\Omega_t))dt \right)\right)\\
	&=\mc{H}(\mu^*\Phi).
\end{align*}
Thus, $\Phi$ is a differential $2k-1$ form on the Lie groupoid $[s^*E_G\rra E_G]$; which would imply, as mentioned before that $\widetilde{\Phi}$ is a differential $2k-1$ form on the Lie groupoid $[X_1\rra X_0]$.  Thus, we have the following result:
\begin{lemma}\label{Lemma:independence}
	Let $\mb{X}=[X_1\rra X_0]$ be a Lie groupoid equipped with an integrable connection $\mc{H}\subseteq TX_1$. Let $(E_G\ra X_0,[X_1\rra X_0])$ be a principal bundle. The map $\text{sym}(\mf{g})^G\ra H^*_{\rm{{\rm dR}}}(\mb{X},\mc{H})$ defined as $f\mapsto f(\mc{K}_{\mc{D}})$ for a connection $\mc{D}$ on $(E_G\ra X_0,[X_1\rra X_0])$ does not depend on the connection $\mc{D}$.
\end{lemma}
\begin{definition}\label{Definition:Chern-Weilmap}
	Let $\mb{X}=[X_1\rra X_0]$ be a Lie groupoid equipped with an integrable connection $\mc{H}\subseteq TX_1$. For a principal bundle $(E_G\ra X_0,[X_1\rra X_0])$ the map $\text{sym}(\mf{g})^G\ra H^*_{\rm{{\rm dR}}}(\mb{X},\mc{H})$ is called \textit{the Chern-Weil map for the principal bundle $(E_G\ra X_0,[X_1\rra X_0])$ over a Lie groupoid $\mb{X}=[X_1\rra X_0]$ equipped with an integrable connection $\mc{H}\subseteq TX_1$}.
\end{definition}

	\chapter{Connections on principal bundles over differentiable stacks}\label{Chap.4}

In Chapter \ref{Chap.3}, we have studied the notion of principal 
$G$-bundle over a Lie groupoid $\mb{X}$ for a Lie group $G$. It is a well-known result that, if two Lie groupoids $\mathbb{X}$ and $\mathbb{Y}$ are Morita equivalent, then the associated categories of principal $G$-bundles over $\mathbb{X}$ is equivalent to the category of principal $G$-bundles over $\mathbb{Y}$ (\cite[Corollary $2.12$]{MR2270285}). This suggest the notion of a principal bundle over a differentiable stack, as a differentiable stack is nothing but a Morita equivalence class of a Lie groupoid (\Cref{Lemma:LiegroupoidsareMoritaequivalent} and \Cref{Proposition:Moritaequivalentimpliesisomorphicstacks}). 

In this chapter, we recall the notion of principal bundles over differentiable stacks and study connections on them. 

\section{Principal bundles over differentiable stacks}\label{Section:principalbundleoverdifferentiablestacks}

In literature, the notion of principal $G$-bundle over a differentiable stack, for a Lie group $G$, appears in different (but equivalent) forms. We will first recall some of these definitions.

Let $G$ be a Lie group, $M$ a manifold, and $\pi:P\ra M$ a principal $G$-bundle. Then, we have an action of the Lie groupoid $[G\rra *]$ on $P$, which is nothing but the action of $G$ on $P$ coming from the principal bundle structure. We also have an action of the 
Lie groupoid $[M\rra M]$ on the manifold $P$; which is precisely the map $\pi:P\ra M$. With these actions, the manifold $P$ is considered as a $[M\rra M]-[G\rra *]$-bibundle, as in the following diagram,
\[
\begin{tikzcd}
	M \arrow[dd,xshift=0.75ex,"t"]\arrow[dd,xshift=-0.75ex,"s"'] &                                 & G \arrow[dd,xshift=0.75ex,"t"]\arrow[dd,xshift=-0.75ex,"s"'] \\
	& P \arrow[ld, "\pi"'] \arrow[rd] &              \\
	M            &                                 & *           
\end{tikzcd}.\]

Conversely, suppose that $P$ is a $[M\rra M]-[G\rra *]$-bibundle. The anchor map $\pi:P\ra M$ is a principal $G$-bundle. The map $P\ra *$ being $[M\rra M]$-invariant is not an extra condition, as the map $P\ra *$ is trivially a $[M\rra M]$-equivariant map. Thus, an $[M\rra M]-[G\rra *]$-bibundle is same as a principal $G$-bundle. This we have mentioned as Example \ref{Example:P(M,G)isbibundle}. Similarly, for a Lie group $G$, and a Lie groupoid $[X_1\rra X_0]$a principal $G$-bundle over $[X_1\rra X_0]$ is same as a $[X_1\rra X_0]-[G\rra *]$-bibundle.

We have mentioned in \Cref{SubSubsection:MorphismofStacksassociatedtobibundle} that, a $\mc{G}-\mc{H}$-bibundle will assign a morphism of stacks $B\mathcal{G}\ra B\mc{H}$. Conversely, every morphism of stacks $B\mathcal{G}\ra B\mc{H}$ comes from a $\mc{G}-\mc{H}$-bibundle (\cite[Remark $4.18$]{MR2778793}). Thus, the collection of $\mc{G}-\mc{H}$-bibundles is in one-one correspondence with the collection of morphisms of stacks $B\mathcal{G}\ra B\mc{H}$. 
Combining these observations, we can see that the collection of principal $G$-bundles over a Lie groupoid $[X_1\rra X_0]$ is in bijective correspondence with morphisms of stacks $B\mb{X}\ra B[G\rra *]$. Thus, the following definition is a natural choice for a principal $G$-bundle over a differentiable stack.

\begin{definition}[{\cite[Definition 6.1.]{MR3521476}}]
	Let $G$ be a Lie group and $\mc{D}$  a differentiable stack. A \textit{principal $G$-bundle over $\mc{D}$} is defined to be a morphism of stacks $\mc{D}\ra BG$, where $BG$ is the differentiable stack associated to the Lie groupoid $[G\rra *]$. 
\end{definition}

Let $G$ be a Lie group, $\mc{D}$ a differentiable stack, and $\theta:\mc{D}\ra BG$ a principal $G$-bundle. As $\mc{D}$ is a differentiable stack, there exists manifold $M$ and an atlas $p:\underline{M}\ra \mc{D}$. Consider the composition $\underline{M}\ra \mc{D}\ra BG$. Ignoring $\mc{D}$, we see the composition as a morphism of stacks $B([M\rra M])\ra BG$. For the reasons mentioned before, the morphism of stacks $B([M\rra M])\ra BG$ is same as an $[M\rra M]\ra [G\rra *]$-bibundle; that is a principal $G$-bundle over the manifold $M$. Consider the following $2$-fiber product diagram,
\[
\begin{tikzcd}
	\underline{M}\times_{\mathcal{D}}\underline{M} \arrow[dd, "{\rm pr}_2"'] \arrow[rr, "{\rm pr}_1"] &  & \underline{M} \arrow[dd, "p"] \\
	&  &                               \\
	\underline{M} \arrow[rr, "p"]                                                                     &  & \mathcal{D}                  
\end{tikzcd}.\]
The $2$-isomorphism in the above diagram induces a $2$-isomoprhism $\theta\circ p\circ {\rm pr}_1\Rightarrow \theta\circ p\circ {\rm pr}_2$. Both of these are functors from $B(M\times_{\mc{D}}M)=\underline{M}\times_{\mc{D}}\underline{M}\ra BG$. That is, there is an isomorphism of principal $G$-bundles over $M\times_{\mc{D}}M$. Further, these bundles satisfy cocycle condition on $M\times_{\mc{D}}M\times_{\mc{D}}M$. Thus, given a principal $G$-bundle over $\mc{D}$, for an atlas $\underline{M}\ra \mc{D}$, we have a principal $G$-bundle over $M$, with an isomorphism of two principal bundles over $M\times_{\mc{D}}M$, such that there is a cocycle condition on the manifold $M\times_{\mc{D}}M\times_{\mc{D}}M$. This motivates the following alternate description of a principal $G$-bundle over a differentiable stack. 

\begin{definition}[{\cite[Definition $2.11$]{MR2206877}}]
	Let $G$ be a Lie group, and $\mc{D}$ a differentiable stack. A \textit{$G$-bundle over $\mc{D}$} is given by a $G$-bundle $\mc{P}_X\ra X$ over the manifold for an atlas $X\ra \mc{D}$ with an isomorphism of the two pull-backs of ${\rm pr}_1^*\mc{P}_X\ra {\rm pr}_2^*\mc{P}_X$ satisfying the cocycle condition on $X\times_{\mc{D}}X\times_{\mc{D}}X$.
\end{definition} 

Before we proceed to the notion of connection on principal bundle over a differentiable stack, we give another equivalent description of a principal $G$-bundle over a differentiable stack.

\begin{definition}[{\cite[Definition $4.1$]{MR2817778}}, {\cite[Definition $3.1$]{biswas2020atiyah}}]\label{Definition:S1bundleoverstack}
	Let $G$ be a Lie group, and $\mc{C}$ a differentiable stack. A \textit{$G$-bundle over $\mc{C}$} is a $2$-commutative diagram
	\[
	\begin{tikzcd}
		\mathcal{D}\times G \arrow[dd, "{\rm pr}_1"'] \arrow[rr, "\sigma"] &  & \mathcal{D} \arrow[dd, "\pi"] \\
		&  &                               \\
		\mathcal{D} \arrow[rr, "\pi"]                                        &  & \mathcal{C}                  
	\end{tikzcd},\]
	such that the pull-back via $U\ra \mc{C}$, for any submersion from a manifold $U$, defines an $G$-bundle over $U$.
\end{definition}

To further explore the \Cref{Definition:S1bundleoverstack}, we need the following result.
\begin{lemma}[{\cite[Lemma $2.11$]{MR2817778}}]\label{Lemma:amorphismarepresentablesubmersion}
	Let $f:\mc{D}\ra \mc{C}$ be a morphism of stacks. Suppose a manifold $U$ and an epimorphism $U\ra \mc{C}$ are given. If the fibered product $V=\mc{D}\times_{\mc{C}}U$ is representable and $V\ra U$ is a submersion, then $f$ is a representable submersion.
\end{lemma}

Let $\pi:\mc{D}\ra \mc{C}$ be an $S^1$-bundle. Then, for a submersion $U\ra \mc{C}$, the fibered product $\mc{D}\times_{\mc{C}}U$ is representable by a manifold, and that the projection map ${\rm pr}_2:\mc{D}\times_{\mc{C}}U\ra U$ is a principal $G$-bundle. Thus, by Lemma \ref{Lemma:amorphismarepresentablesubmersion}, we see that $\pi:\mc{D}\ra \mc{C}$ is a representable submersion.  

\begin{remark}\label{Remark:explainingdetailsofdefinitionofPBS}
	Let us reveal some more details about this definition. Let $u:U\ra \mc{C}$ be an atlas for $\mc{C}$. Then, we have the principal bundle $\mc{D}_u\ra U$. Let $v:V\ra \mc{C}$ be another atlas for $\mc{C}$ with a $2$-morphism $\alpha:u\ra v\circ \varphi$ as in the following $2$-commutative diagram, 
	\[
	\begin{tikzcd}
		U \arrow[rd, "u"'] \arrow[rr, "\varphi"] &             & V \arrow[ld, "v"] \\
		& \mathcal{C} &                  
	\end{tikzcd}.\]
	
	We can pull-back $\mc{D}_v\ra V$ along the morphism $\varphi:U\ra V$ to get the principal bundle $\varphi^*\mc{D}_v\ra U$. We now have two principal bundles $\mc{D}_u\ra U$ and $\varphi^*\mc{D}_v$.
	The $2$-commutative diagram $\alpha$ implies that, these principal bundles $\mc{D}_u\ra U$ and $\varphi^*\mc{D}_v$ are isomorphic with isomorphism $\theta_\alpha:\mc{D}_u\ra \varphi^*\mc{D}_v$, as in the following diagram,
	\[
	\begin{tikzcd}
		&  & \varphi^*\mathcal{D}_v \arrow[dd] \arrow[rr, "{\rm pr}_2"] &  & \mathcal{D}_v \arrow[dd] \\
		&  &                                              &  &                          \\
		\mathcal{D}_u \arrow[rr] \arrow[rruu, "\theta_\alpha"] &  & U \arrow[rr, "\varphi"]                      &  & V                       
	\end{tikzcd}.\]
	
	Let $u:U\ra \mc{C}, v:V\ra \mc{C}, w:W\ra \mc{C}$ be atlases of $\mc{C}$ with morphisms $\varphi:U\ra V$, and $\psi:V\ra W$ as in the following diagram,
	\[
	\begin{tikzcd}
		U \arrow[rr, "\varphi"] \arrow[rrd, "u"'] &  & V \arrow[rr, "\psi"] \arrow[d, "v"] &  & W \arrow[lld, "w"] \\
		&  & \mathcal{C}                         &  &                   
	\end{tikzcd},\]
	with isomorphisms $\alpha:u\Rightarrow v\circ \varphi$ and $\beta:v\Rightarrow w\circ \psi$. Thus, we have following diagrams of isomorphisms of principal bundles,
	\[
	\begin{tikzcd}
		\mathcal{D}_u \arrow[dd] \arrow[rr, "\theta_\alpha"] &  & \varphi^*\mathcal{D}_v \arrow[dd] &  & D_v \arrow[dd] \arrow[rr, "\theta_\beta"] &  & \psi^*\mathcal{D}_w \arrow[dd] &  & \mathcal{D}_u \arrow[dd] \arrow[rr, "\theta_{\beta\circ\alpha}"] &  & (\psi\circ\varphi)^*\mathcal{D}_w \arrow[dd] \\
		&  &                                   &  &                                           &  &                                &  &                                                                  &  &                                              \\
		U \arrow[rr]                                         &  & U                                 &  & V \arrow[rr]                              &  & V                              &  & U \arrow[rr]                                                     &  & U                                           
	\end{tikzcd}.\]
	The definition of principal bundle says that, these isomorphisms $\theta_\alpha, \theta_\beta,\theta_{\beta\circ\alpha}$ are related by the condition $\psi^*\theta_\beta\circ \theta_\alpha=\theta_{\beta\circ\alpha}$, as in the following diagram,
	\[
	\begin{tikzcd}
		\mathcal{D}_u \arrow[dd, "\theta_\alpha"'] \arrow[rr, "\theta_{\beta\circ\alpha}"] &  & (\psi\circ\varphi)^*\mathcal{D}_w \arrow[dd, "="] \\
		&  &                                                   \\
		\varphi^*\mathcal{D}_v \arrow[rr, "\varphi^*\theta_\beta"']                        &  & (\varphi^*\circ \psi^*) \mathcal{D}_w            
	\end{tikzcd}.\]
\end{remark}

\begin{definition}
	Let $G$ be a Lie group and $\mc{C}_1,\mc{C}_2$ differentiable stacks. A \textit{morphism of principal $G$-bundles from $\mc{D}_1\ra \mc{C}_1$ to $\mc{D}_2\ra \mc{C}_2$} is given by a pair of morphisms of differentiable stacks $( \Phi:\mc{D}_1\ra \mc{D}_2,\Psi:\mc{C}_1\ra \mc{C}_2)$, such that the following diagram commutes,
	
	\[
	\begin{tikzcd}
		\mathcal{D}_1\times G \arrow[dd, "{\rm pr}_1"'] \arrow[rr, "\sigma_1"] \arrow[rd, "{(\Phi,1)}"'] &                                                                       & \mathcal{D}_1 \arrow[dd, "\pi_1"] \arrow[rd, "\Phi"] &                                   \\
		& \mathcal{D}_2\times G \arrow[dd, "{\rm pr}_2"] \arrow[rr, "\sigma_2"] &                                                      & \mathcal{D}_2 \arrow[dd, "\pi_2"] \\
		\mathcal{D}_1 \arrow[rr, "\pi_1"] \arrow[rd, "\Phi"]                                             &                                                                       & \mathcal{C}_1 \arrow[rd, "\Psi"]                     &                                   \\
		& \mathcal{D}_2 \arrow[rr, "\pi_2"]                                     &                                                      & \mathcal{C}_2                    
	\end{tikzcd}.\] 
\end{definition}	

\section{connections on principal bundles over differentiable stacks}
\begin{definition}\label{Definition:principalbundleoverstacks}
	Let $G$ be a Lie group, $\mc{C}$ a differentiable stack and 
	$\mc{D}\ra\mc{C}$ a principal $G$-bundle over the stack $\mc{C}$.  A \textit{connection $\nabla$ on $\mc{D}\ra \mc{C}$} consists of a connection $\nabla_u$ on each principal $G$-bundle $\mc{D}_u\ra U$, where $u:U\ra \mc{C}$ is a smooth atlas for $\mc{C}$, which pulls back naturally with respect to each $2$-commutative diagram of the form 
	\[
	\begin{tikzcd}
		U \arrow[rd, "u"'] \arrow[rr, "\varphi"] &             & V \arrow[ld, "v"] \\
		& \mathcal{C} &                  
	\end{tikzcd}.\]
	A connection $\nabla$ on $\mc{D}$ is \textit{flat} or \textit{integrable} if it is in addition integrable on each principal $G$-bundle $\mc{D}_u\ra U$.
\end{definition}

Similar to the explanation in \Cref{Remark:explainingdetailsofdefinitionofPBS}, the description of connection can be given in terms of atlases with some compatibility conditions.

Let $\nabla$ be a connection on $\mc{D}\ra U$. Let $u,v$ be a pair of atlases of $\mc{C}$ as in the above diagram.
Let $\omega_{G,u}$ be the $\mf{g}$-valued $1$-form on $\mc{D}_u$ (associated to connection $\nabla_u$ on $U$) and $\omega_{G,v}$  the $\mf{g}$-valued $1$-form on $\mc{D}_v$ (associated to connection $\nabla_v$ on $V$). Consider the $1$-form ${\rm pr}_2^*\omega_{G,v}$, which is the pull-back of $\omega_{G,v}$ of along the smooth map ${\rm pr}_2:\varphi^*\mc{D}_v\ra \mc{D}_v$. The condition on connection says that, the differential forms ${\rm pr}_2^*\omega_{G,v}, \omega_{G,u}$ are related by the isomorphism $\theta_\alpha:\mc{D}_u\ra \varphi^*\mc{D}_v$ as \[\omega_{G,u}=\theta_\alpha^*({\rm pr}_2^*\omega_{G,v}).\]

\section{Atiyah sequence associated to principal bundle over differentiable stack}\label{Section:Atiyahsequenceforbundleoverstacks}
Let $G$ be a Lie group, and $\mc{D}\ra \mc{C}$ a principal $G$-bundle. For this principal bundle, we associate a short exact sequence of vector bundles over the stack $\mc{C}$. We have done the same for principal bundles over Lie groupoids in \Cref{Section:AtiyahforbundleoverLiegroupoid}. 

The notion of a tangent stack associated to a differentiable stack is a straightforward generalization of the tangent bundle for a manifold. But, for completeness, we will recall the definition here. 
\begin{definition}[Tangent stack {\cite[Definition/Lemma $4.12$]{MR2206877}}] Let $\mc{C}$ be a differentiable stack, and $X\ra \mc{C}$ an atlas of $\mc{C}$. Then, we can take the tangent spaces to the groupoid $[X\times_{\mc{C}}X\rra X]$, namely $[T(X\times_{\mc{C}}X)\rra TX]$, then the stack associated to $[T(X\times_{\mc{C}}X)\rra TX]$ is independent of choice of atlas $X\ra \mc{C}$, and is called $T\mc{C}$, the tangent stack to $\mc{C}$.
\end{definition}

Let $\mc{C}$ be a Lie groupoid, and $T\mc{C}$ be the associated tangent stack. We have seen that a morphism of Lie groupoids induce a morphism of associated differentiable stacks (\Cref{Subsection:bibundleAssociatedtoMorphismofLiegroupoids,SubSubsection:MorphismofStacksassociatedtobibundle}). In particular, the morphism of Lie groupoids $[TX_1\rra TX_0]\ra[X_1\rra X_0]$ induce a morphism of stacks $T\mc{C}\ra \mc{C}$. As $TX_1\ra X_1$ and $TX_0
\ra X_0$ are vector bundles; we can ask if the morphism of stacks $T\mc{C}\ra \mc{C}$ is a ``vector bundle over the stack $\mc{C}$''. It turns out that under some special conditions, the morphism of stacks $T\mc{C}\ra \mc{C}$ is a vector bundle. Firstly, we will recall the definition of a vector bundle over a differentiable stack. 
\begin{definition}[vector bundle over differentiable stack {\cite[Definition $5.1$]{Ginot}}] Let $\mc{C}$ be a differentiable stack. A \textit{vector bundle over $\mc{C}$} is a representable morphism of stacks $\mc{D}\ra \mc{C}$ such that, for every morphism of stacks $f:\underline{U}\ra \mc{C}$ for a manifold $U$, the pull-back $\mc{D}\times_{\mc{C}}U\ra U$ is endowed with structure of a vector bundle. We also require, for every $a:V\ra U$, that the natural isomorphism $\varphi_a:(f\circ a)^*\mc{D}\ra a^*(f^*\mc{D})$ is a morphism of vector bundles.
\end{definition}

\begin{remark}
	If $\mc{C}$ is an orbifold, then $T\mc{C}\ra \mc{C}$ is a vector bundle over $\mc{C}$ (\cite[Lemma $5.2$]{Ginot}). Moreover, if $\mc{C}$ is a Deligne-Mumford stack, then $T\mc{C}\ra \mc{C}$ is a vector bundle over $\mc{C}$. However, $T\mc{C}\ra \mc{C}$ is not a vector bundle over $\mc{C}$ in general.	
\end{remark}

From now, we work with principal bundles over Deligne-Mumford stacks.
Similar to the construction of the Atiyah bundle, and adjoint bundle for a principal bundle over a manifold and over a Lie groupoid, we can associate adjoint bundle and Atiyah bundle for a principal bundle over a differentiable stack.

Let $G$ be a Lie group, and $\mc{D}\ra \mc{C}$ a principal $G$-bundle over $\mc{C}$. We define the Atiyah bundle 
${\rm At} (\mc{D})$ as 
\[{\rm At} (\mc{D})_u={\rm At} (\mc{D}_u),\]
for any etale morphism $u:U\ra\mc{C}$. For a $2$-commutative diagram of the form 
\[\begin{tikzcd}
	U \arrow[rd, "u"'] \arrow[rr, "\varphi"] &             & V \arrow[ld, "v"] \\
	& \mathcal{C} &                  
\end{tikzcd},\]
where $u,v,\varphi$ are etale morphisms, we get an isomorphism 
${\rm At} (\mc{D})_u\cong \varphi^* \rm{At} (\mc{D})_v$. We define the adjoint bundle as ${\rm ad} (\mc{D})=\mc{D}\times_G \mf{g}$, where ${\rm ad}(\mc{D})_u=\rm{ad} (\mc{D}_u)$. We get the following diagram of morphism of vector bundles 
\begin{equation}\label{Equation:Atiyahforstacks}
	0\ra {\rm ad} (\mc{D})\ra {\rm At}(\mc{D})\ra T\mc{C}\ra 0
\end{equation} 
over the stack $\mc{C}$.

\begin{definition}
	Let $G$ be a Lie group, $\mc{C}$ a Deligne-Mumford stack, and $\mc{D}\ra \mc{C}$ a principal $G$-bundle. The short exact sequence of vector bundles over $\mc{C}$, mentioned in Equation \ref{Equation:Atiyahforstacks}, is called the \textit{Atiyah sequence for the principal bundle $\mc{D}\ra \mc{C}$}.
\end{definition}

In our paper \cite{biswas2020atiyah}, we have proved the following results:
\begin{proposition}[{\cite[Proposition $5.11$]{biswas2020atiyah}}]
	A principal $G$-bundle $\mc{D}\ra \mc{C}$ over a Deligne-Mumford stack $\mc{C}$ admits a connection if and only if its associated Atiyah sequence has splitting.
\end{proposition}

\begin{theorem}[{\cite[Theorem]{biswas2020atiyah}}]
	Giving a connection on a principal $G$-bundle $\mc{D}$ over a Deligne-Mumford stack $\mc{C}$ with etale atlas $x:X_0\ra \mc{C}$ is equivalent to giving a connection on the associated principal $G$-bundle over the Lie groupoid $[X_1\rra X_0]$.
\end{theorem}
	\chapter{Extensions of topological groupoids}\label{Chap.5}

In \Cref{Chap.Preliminaries}, we have 	mentioned about Lie groupoids, differentiable stacks, and correspondence between them. In \Cref{Chap.2}, we have given correspondence between a particular kind of morphism of Lie groupoids and a particular kind of differentiable stacks; that is, between Lie groupoid extensions and gerbes over differentiable stacks. In this Chapter, based on our paper \cite{chatterjee2021extension}, we introduce the notion of topological groupoid extensions and give a correspondence between topological groupoid extensions and gerbes over topological stacks.  Part of this chapter is analogous to  \Cref{Chap.2}. But, we have given enough details to make it self content.

\section{Topological groupoids and topological stacks}\label{Section:top groupoids and top stacks}
In this section, we will recall standard definitions and results about topological groupoids and topological stacks. A major part of this section is based on \cite{Carchedi, metzler2003topological, noohi2005foundations}.
\begin{definition}[topological groupoid]
	A \textit{topological groupoid} is a groupoid $\mc{G}=[\mc{G}_1\rra \mc{G}_0]$ whose object set $\mc{G}_0$, morphism set $\mc{G}_1$ are topological spaces, and the structure maps $s,t,m,e,i$ are continuous maps.
\end{definition}

\begin{example}
	Given a topological space $M$, the groupoid $[M\rra M]$ (mentioned in \Cref{Example:MMcategory} is a topological groupoid. 
\end{example}

\begin{example}
	For a topological group $G$, the groupoid $[G\rra *]$ (mentioned in \Cref{Example:G*category}) is a topological grouopoid. 
\end{example}	

\begin{example}
	For an action of a topological group $G$ on a topological space $M$, the action groupoid $[M\times G\rra M]$ is a topological groupoid. 
\end{example}

\begin{definition}
	Let $\mc{G},\mc{H}$ be topological groupoids. A \textit{morphism of topological groupoids from $\mc{G}$ to $\mc{H}$} is a functor $(F_1,F_0):[\mc{G}_1\rra \mc{G}_0]\ra [\mc{H}_1\rra \mc{H}_0]$, such that, the maps $F_1:\mc{G}_1\rra \mc{H}_1$ and $F_0:\mc{G}_0\ra \mc{H}_0$ are continuous maps. 
\end{definition}

Let $\mc{C}$ be a (locally-small) category, and $U$ an object of $\mc{C}$. Consider the functor of points (associated to the object $U$) $h_U:\mc{C}^{\op}\ra \text{Set}$. Seeing a set as a discrete category (\Cref{Example:MMcategory}), we can see $h_U$ as a pseudo-functor $\mc{C}^{\op}\ra \text{Cat}$. As mentioned before, we are interested in Grothendieck topologies, on categories, that are subcanonical; that is, $h_U:\mc{C}^{\rm op}\ra \text{Cat}$ is a stack for each object $U$ of $\mc{C}$.

Consider the  Grothendieck topology $\mc{J}$ on ${\rm Top}$ described as follows.  A covering of an object $U$ of $\text{Top}$ is given by a jointly surjective family $\{\sigma_\alpha\colon U_\alpha\ra U\}$, where each $\sigma_\alpha:U_\alpha \ra U$ is  such that $\sigma_\alpha(U_\alpha)$ is an open subset of $U$ and, $\sigma_\alpha$ is a homeomorphism onto its image. Here, by jointly surjective family, we mean $\bigcup_\alpha \sigma_\alpha(U_\alpha)=U$. We call this topology to be \textit{the open-cover topology on $\text{Top}$}.  Unless mentioned otherwise,  we will always assume $\rm {Top}$ to be equipped with the open-cover topology. 

In this chapter, we will be working in the category $\text{Top}$ of topological spaces, with Grothendieck topology $\mc{J}$ on $\text{Top}$ being the open cover topology. The stacks we are interested are over the site $(\text{Top},\mc{J})$.

For each object $M$ of $\text{Top}$ we have a stack $h_M:\text{Top}\ra \text{Cat}$ with the description at the level of objects being $h_M(N)=\text{Hom}_{\text{Top}}(N,M)$ for objects $N$ in $\text{Top}$. Using the notion of morphisms of topological groupoids, we can associate a pseudo-functor $\text{Top}\ra \text{Cat}$ for a topological groupoid $\mc{G}$. 

Let $\mc{G}$ be a topological groupoid. Similar to the description of $h_M$ for a topological space $M$, we describe a pseudo-functor $h_\mc{G}:\text{Top}\ra \text{Cat}$. For an object $N$ of $\text{Top}$, declare $h_{\mc{G}}(N)$ as
\[h_{\mc{G}}(N)=\text{Hom}_{\text{TopGpd}}((N\rra N),\mc{G}),\] for an object $N$ of $\text{Top}$. Here, $\text{Hom}_{\text{TopGpd}}((N\rra N),\mc{G})$ is the category whose objects are Lie groupoid morphisms and morphisms are natural transformations. This gives a pseudo-functor $h_{\mc{G}}:\text{Top}\ra \text{Gpd}$. The stackification of this pseudo functor, denoted by $B\mc{G}$, is called the stack associated to the topological groupoid $\mc{G}$. More details about this can be found  \cite[$I.2.4$]{Carchedi}. We do not follow this method of assigning a stack for a topological groupoid. We follow the method of using principal $\mc{G}$-bundles for a topological groupoid $\mc{G}$. For that purpose, we first introduce the notion of action of a topological groupoid on a topological space. 

\begin{definition}[action of topological groupoid {\cite[Definition $4.1$]{MR896907}}]
	Let $[\mc{G}_1\rra M]$ be a topological groupoid, and $\pi:P\ra M$  a continuous map. Let $s^*P$ denote the subspace 
	$\{(\gamma,p)\in \mc{G}_1\times P| s(\gamma)=\pi(p)\}$. An \textit{action of $\mc{G}$ on $(P,\pi,M)$} is given by a continuous map $\mu:s^*P\ra P$, with $(\gamma, p)\mapsto \gamma p$, such that 
	\begin{enumerate}
		\item $\pi(\gamma p)=t(\gamma)$ for all $(\gamma, p)\in s^*P$,
		\item $\gamma(\gamma'p)=(\gamma\circ \gamma')p$ for all $(\gamma,\gamma')\in \mc{G}_1\times_{t,\mc{G}_0,s}\mc{G}_1$ and $(\gamma',p)\in \mc{G}_1\times_{M}P$,
		\item $1_{\pi(p)}p=p$ for all $p\in P$.	
	\end{enumerate}
\end{definition}
Alternatively, we also refer to the above definition as ``action of $[\mc{G}_1\rra M]$ on a topological space $P$'' being given by a pair of maps $(\pi:P\ra M, \mu:\mc{G}_1\times_{\mc{G}_0}P\ra P)$. We call the above action to be a left action of the topological groupoid on the topological space.  

Similarly, a right action of a topological groupoid $\mc{G}$ on a topological space is given by a pair of continuous maps 
\begin{equation}
	(\pi:P\ra M, \mu:P\times_{\pi,\mc{G}_0,t}\mc{G}_1\ra P),
\end{equation}
satisfying the following conditions:
\begin{enumerate}
	\item $\pi(p\gamma)=s(\gamma)$ for all $(p,\gamma)\in P\times_{\mc{G}_0}\mc{G}_1$,
	\item $(p\gamma)\gamma'=p(\gamma\circ \gamma')$ for all $(p,\gamma)\in P\times_{\mc{G}_0}\mc{G}_1$ and $(\gamma,\gamma')\in \mc{G}_1\times_{\mc{G}_0}\mc{G}_1$,
	\item $p 1_{\pi(p)}=p$ for all $p\in P$.
\end{enumerate}

\begin{example}
	Let $[M\rra M]$ be the topological groupoid associated to a topological space $M$. An action of  $[M\rra M]$ on a topological space $N$ consists precisely of a continuous map $f:M\ra N$.
\end{example}

\begin{example}
	Let $[G\rra *]$ be the topological groupoid associated to a Lie group $G$. An action of $[G\rra *]$ on a topological space $M$ is given by an action $\mu:M\times G\ra M$ of the topological  group $G$ on the topological space $M$.
\end{example}	
\begin{example}
	Let $[\mc{G}_1\rra \mc{G}_0]$ be a topological groupoid. Consider the composition map
	\[m:\mc{G}_1\times_{t,\mc{G}_0,s} \mc{G}_1\ra \mc{G}_1.\] 
	The pair $(t:\mc{G}_1\ra \mc{G}_0, m:\mc{G}_1\times_{t,\mc{G}_0,s} \mc{G}_1\ra \mc{G}_1)$ gives a left action of $[\mc{G}_1\rra \mc{G}_0]$ on $\mc{G}_1$. The pair $(s:\mc{G}_1\ra \mc{G}_0, m:\mc{G}_1\times_{t,\mc{G}_0,s} \mc{G}_1\ra \mc{G}_1)$ gives a right action of $[\mc{G}_1\rra \mc{G}_0]$ on $\mc{G}_1$.
\end{example}

Now, for a topological groupoid $\mc{G}$, we introduce the notion of a principal $\mc{G}$-bundle over a topological space $M$. For that, we need to recall the notion of a map admitting local sections. 

\begin{definition}
	Let $f:X\ra Y$ be a continuous map. We say that $f:X\ra Y$ is a \textit{map of local sections} or that \textit{it admits local sections}, if there is an open cover ${U_\alpha}$ of $Y$ and sections $\sigma_\alpha:U_\alpha\ra X$ of $f:X\ra Y$ for each $\alpha$.
\end{definition}

\begin{definition}[{\cite[Definition $I.2.13$]{Carchedi}}]\label{Definition:G-bundle}
	Let $[\mc{G}_1\rra \mc{G}_0]$ a topological groupoid and $M$ a topological space. 
	A \textit{principal $\mc{G}$-bundle over $M$} consists of a topological space $P$, with a   continuous map  $\pi\colon P\ra M$ admitting local sections and a right action $(a_{\mc{G}}, \mu)$ of 	$\mc{G}$ on $P$ such that 		
	\begin{enumerate}
		\item $\pi(p\gamma)=\pi(p)$ for all $(p,\gamma)\in P\times_{\mc{G}_0}\mc{G}_1$,
		\item the map $P\times_{\mc{G}_0}\mc{G}_1\ra P\times_{M}P$, $(p,\gamma)\mapsto (p,p\cdot \gamma)$ is a homeomorphism.
	\end{enumerate}
\end{definition} 

\begin{definition}\label{Def:mapG-bundle}Let $\mc{G}=[\mc{G}_1\rra \mc{G}_0]$ be a topological groupoid and $(P,\pi,M), (P',\pi',M')$ principal $\mc{G}$-bundles. A \textit{morphism} from $P(\pi,M)$ to $P'(\pi',M')$ is given by a pair of continuous maps
	$(F\colon P\ra P', f\colon M\ra M')$, satisfying the following compatibility conditions
	\begin{enumerate}
		\item  $a_{\mc{G}}'\circ F=a_{\mc{G}}$, and 
		$F(p\cdot \gamma)=F(p)\cdot \gamma$ for all 
		$(p,\gamma)\in P\times_{\mc{G}_0}\mc{G}_1$,
		\item  $\pi'\circ F=f\circ \pi$.
	\end{enumerate}
\end{definition}

\begin{example}
	Let $\mc{H}=[\mc{H}_1\rra \mc{H}_0]$ be a topological groupoid. The target map $t\colon  \mc{H}_1\ra \mc{H}_0$ can be considered a principal $\mc{H}$-bundle, with action of $\mc{H}$  given by the pair, composition of arrows and the  source map $s\colon  \mc{H}_1\ra \mc{H}_0$.	
\end{example}

\subsection{stack associated to a topological groupoid} In the same way as we pull-back a topological structure on a topological space along a continuous map, we have the notion of pull-back of a principal $\mc{G}$-bundle over a topological space $M$ along a continuous map $f:M\ra M'$.

Let $\pi:P\ra M$ be a principal $\mc{G}$-bundle, and $f\colon M'\to M$ a continuous map. Consider the set $f^*P=\{(p,m')\in P\times M':f(m')=\pi(p)\}$. The action of $\mc{G}$ on $P$ induce an action of $\mc{G}$ as $((p,n),\gamma))\mapsto (\mu(p,\gamma),n)$. The projection map $\pr_2 :f^*P\ra M'$ is then considered as a principal $\mc{G}$-bundle. We call $\pr_2 :f^*P\ra M'$ to be the pull-back of $(P,\pi,M)$ along the map $f:M\ra M'$. 

We use the above notion of pull-back of principal $\mc{G}$-bundle along a continuous map to introduce the notion of a pseudo-functor associated to a topological groupoid. 

Let $\mc{G}$ be a topological groupoid. Let $M$ be an object of $\text{Top}$. Consider the category of principal $\mc{G}$-bundles over $M$, $B\mc{G}(M)$.

For the same reason as in the class of principal $G$-bundle for a Lie group $G$, any morphism of principal $\mc{G}$-bundles of the form $(F,1_M)$ is an isomorphism. So, $B\mc{G}(M)$ is a groupoid for each object $M$ of $\text{Top}$. Thus, we get a  pseudo-functor $B\mc{G}\colon \text{Top}^{\rm op}{\ra \text{Gpd}}\subset \text{Cat}$, with the following description:
\begin{itemize}
	\item  an object $M$ in $\rm {Top}$ to the category of $\mc{G}$-bundles over the topological space $M$,
	\item a continuous map $f\colon M'\to M$ is sent to a functor $B\mc{G}(f)\colon B\mc{G}(M)\to B\mc{G}(M')$ defined by the pull-back of principal bundles
\end{itemize}
For a morphism $(F\colon P_1\to P_2, {\rm Id}_M)$ in $B\mc{G}(M)$, the morphism	
$(f^*P_1\to f^*P_2, {\rm Id}_{M'})$ in $B\mc{G}(M')$ is the unique morphism defined by the universal property of the pull-back diagram with respect to the maps $f\colon M'\to M$  and $P_2\to M$.

\begin{example}
	Let $\mc{G}$ be a topological groupoid. The pseudo-functor $B\mc{G}: \text{Top}^{\op}\ra \text{Gpd}$  is a stack over the site $(\text{Top},\mc{J})$.
\end{example}

As in the smooth setup, we will be interested in stacks that are isomorphic to $B\mc{G}$ for some topological groupoid $\mc{G}$.
In the smooth case, we declare a stack $\mc{D}\ra \text{Man}$ to be a differentiable stack, if there is a morphism of stacks $\underline{X}\ra \mc{D}$ with an extra condition. In the continuous case, we define topological stacks to be stacks that come from topological groupoids.

\begin{definition}[topological stack 
	{\cite{Carchedi}}]
	\label{Definition:topologicalstack}
	Let $(\text{Top},\mc{J})$ be the site  of topological spaces with open cover topology. A stack $\mc{D}\ra \text{Top}$ is called a \textit{topological stack} if there exists a topological groupoid $\mc{G}=[\mc{G}_1\rra \mc{G}_0]$ and an isomorphism of  stacks $B\mc{G}\cong \mc{D}$.
\end{definition}
\begin{remark}
	In \cite[Definition $7.1$]{noohi2005foundations} the topological stack defined above has been called a pretopological stack. Whereas the topological stack in the same paper has been defined [Definition $13.8$.] as a pretopological stack with some additional conditions.
\end{remark}

\begin{example}
	Let $X, Y$ be topological space, and $f\colon   X\ra Y$ be a continuous map. Let $\underline{X}$ and $\underline{Y}$ be stacks associated to $X$ and $Y$ respectively. Then, the map $f\colon  X\ra Y$ induces a morphism of stacks 
	$F\colon  \underline{X}\ra \underline{Y}$  defined by compositions (both at the level of objects and at the level of morphisms)
\end{example}

Likewise, for the smooth set up, we want to characterize topological stacks in terms of morphisms of stacks from a topological space. First of all, let us recall the notions of representable morphism and epimorphism of stacks from \Cref{Chap.Preliminaries}.

Let $(\mc{S},\mc{J})$ be a subcanonical site. A stack $\mc{D}\ra \mc{S}$ over the site $(\mc{S},\mc{J})$ is said to be a \textit{representable stack} if there exists an object $U$ of $\mc{S}$ and an isomorphism of stacks $\mc{D}\cong \underline{U}$. 

A morphism of stacks  $F\colon  \mc{D}\ra \mc{C}$ a \textit{representable morphism},  if for each object $U$ of $\mc{S}$ and a morphism of stacks $q\colon  \underline{U}\ra \mc{C}$, the $2$-fiber product  stack $\mc{D}\times_{\mc{C}}\underline{U}$ 
is representable by an object of $\mc{S}$. 
A morphism of stacks  $F\colon  \mc{D}\ra \mc{C}$ is said to be \textit{an epimorphism of stacks}, if for each object $U$ of $\mc{S}$ and a morphism of stacks $q\colon  \underline{U}\ra \mc{C}$, there exists a covering $\{\sigma_\alpha\colon  U_\alpha\ra U\}$ of $U$ and a family of morphisms $\{q_\alpha\colon  \underline{U_\alpha}\ra \mc{D}\}$ such that the following diagram is $2$-commutative
\[
\begin{tikzcd} 
	\underline{U_{\alpha}} \arrow[dd, "q_{\alpha}"'] \arrow[rr, "\sigma_{\alpha}"] &  & \underline{U} \arrow[dd, "q"] \\
	&  &                               \\
	\mathcal{D} \arrow[rr,"F"] \arrow[Rightarrow, shorten >=10pt, shorten <=10pt, uurr]                                                   &  & \mathcal{C}                  
\end{tikzcd},\]
for each $\alpha$.

\begin{definition}
	Let $\mc{D}\ra \text{Top}$ be a stack. Let $X$ be a topological space. A morphism of stacks $f\colon \underline{X}\rightarrow \mc{D}$ is said to be \textit{an atlas for the stack $\mc{D}$} if, for every topological space $Y$ and a morphism of stacks $\underline{Y}\ra \mc{D}$, the $2$-fibered product $\underline{X}\times_{\mc{D}}\underline{Y}$ is representable by a topological space $X\times_{\mc{D}}Y$ 
	and the map of topological spaces $X\times_{\mc{D}}Y\rightarrow Y$ associated to the morphism of stacks ${\rm pr}_1\colon \underline{X}\times_{\mc{D}}\underline{Y}\rightarrow \underline{Y}$  is  a map of local sections.
\end{definition}

An epimorphism of stacks can be described in terms of local sections, which we discuss below.

\begin{remark}
	Let  $\mc{D}$ and  $\mc{C}$ be stacks over a subcanonical site $(\mc{S},\mc{J})$.
	Then it is easy to see that a morphism of stacks  $F\colon  \mc{D}\ra \mc{C}$ is an epimorphism of stacks,  if and only if  for each object $U$ of $\mc{S}$ and an object $b$ of $\mc{C}(U)$, there exists a cover $\{U_\alpha\ra U\}$ of $U$ and objects $a_\alpha\in \mc{D}(U_\alpha)$ such that $F(a_\alpha)\cong b|_{U_\alpha}$ in $\mc{C}(U_\alpha)$ for each $\alpha$. 
\end{remark}

As we will see in the next result, an epimorphism of topological stacks can be described in terms of maps that admit local sections. 
\begin{lemma}
	A morphism of topological stacks $F\colon  \mc{D}\ra \mc{C}$ is an epimorphism of topological stacks if and only if, for each topological space $U$ and a morphism of stacks $\underline{U}\ra \mc{C}$ there exists a topological space $V$, a map of local sections
	$\pi\colon  V\ra U$ and a morphism of stacks $\theta\colon \underline{V}\ra \mc{D}$ such that the following diagram is $2$-commutative,
	\[\begin{tikzcd}
		\underline{V} \arrow[dd, "\theta"'] \arrow[rr, "\pi"] &  & \underline{U} \arrow[dd, "q"] \\
		&  &                               \\
		\mathcal{D} \arrow[rr,"F"] \arrow[Rightarrow, shorten >=10pt, shorten <=10pt, uurr]                                                   &  & \mathcal{C}                  
	\end{tikzcd}{\color{red}.}\]
	\begin{proof}
		Suppose $F\colon  \mc{D}\ra \mc{C}$ is  an epimorphism. Let $U$ be a topological space and $q\colon  \underline{U}\ra \mc{C}$ be a morphism of topological stacks. Then, there exists an open cover $\{U_\alpha\ra U\}$ and a morphism of stacks $q_\alpha\colon  \underline{U_\alpha}\ra \mc{D}$ such that the following diagram is $2$-commutative, \[\begin{tikzcd}
			\underline{U_{\alpha}} \arrow[dd, "q_{\alpha}"'] \arrow[rr, "\sigma_{\alpha}"] &  & \underline{U} \arrow[dd, "q"] \\
			&  &                               \\
			\mathcal{D} \arrow[rr,"F"] \arrow[Rightarrow, shorten >=10pt, shorten <=10pt, uurr]                                                   &  & \mathcal{C}                  
		\end{tikzcd},\]
		for each $\alpha\in \Lambda$.
		
		Let $V$ denote the disjoint union $\bigsqcup_{\alpha\in \Lambda} U_\alpha$ and $\pi\colon  V\ra U$ be the map defined as $\pi(x)=x$, if $x\in U_\alpha$ for some $\alpha\in \Lambda$. This map $\pi\colon  V\ra U$ is a continuous map with local sections; since there exists a cover $\{U_\alpha\}$ of $U$ and sections $U_\alpha\ra V$ of $\pi\colon  V\ra U$. 
		
		Observe that, each morphism of stacks $\underline{U_\alpha}\ra \mc{D}$ assigns an object of $\mc{D}(U_\alpha)$ by $2$-Yoneda Lemma. So, for each $\alpha\in \Lambda$, we have an object $a_\alpha\in \mc{D}(U_\alpha)$. As $U_\alpha\cap U_\beta=\emptyset$ in $V$, the collection $\{a_\alpha\}_{\alpha\in \Lambda}$ trivially agree on the intersection. Now, as $\mc{D}$ is a stack, this compatible collection $\{a_\alpha\in \mc{D}(U_\alpha)\}$ glue together  to produce an object of $\mc{D}(V)$, which by Yoneda $2$-lemma gives a morphism of stacks $\underline{V}\ra \mc{D}$. So, for the morphism of stacks $\underline{U}\ra \mc{C}$, there exists a map of local sections $\pi\colon  V\ra U$ and a morphism of stacks $\theta\colon  \underline{V}\ra \mc{D}$ with the following $2$-commutative diagram,
		\[\begin{tikzcd} 
			\underline{V} \arrow[dd, "\theta"'] \arrow[rr, "\pi"] &  & \underline{U} \arrow[dd, "q"] \\
			&  &                               \\
			\mathcal{D} \arrow[rr,"F"] \arrow[Rightarrow, shorten >=10pt, shorten <=10pt, uurr]                                                   &  & \mathcal{C}                  
		\end{tikzcd}.\]
		
		Conversely, assume that for a morphism of stacks $F\colon  \mc{D}\ra \mc{C}$, and a  morphism of stacks $\underline{U}\ra \mc{C}$ for a topological space $U$, there exists a map of local sections $\pi\colon  V\ra U$ and a morphism of stacks $\theta\colon  \underline{V}\ra \mc{D}$ with a $2$-morphism $F\circ \theta\Rightarrow q\circ \pi$. 
		
		As $\pi\colon  V\ra U$ is a map of local sections, there exists an open cover $\{U_\alpha\}$ of $U$ and sections $\sigma_\alpha\colon  U_\alpha\ra V$ of $\pi\colon  V\ra U$. Thus, we have inclusion $i_\alpha=\pi\circ \sigma_\alpha\colon  \underline{U_\alpha}\ra \underline{V}\ra \underline{U}$ and morphism of stacks $q_\alpha=\theta\circ \sigma_\alpha\colon  \underline{U_\alpha}\ra \underline{V}\ra \mc{D}$. The $2$-morphism $F\circ \theta\Rightarrow q\circ \pi$ induce $2$-morphisms $F\circ q_\alpha\Rightarrow q\circ i_\alpha$ giving following $2$-commutative diagram 
		\[	\begin{tikzcd} 
			\underline{U_{\alpha}} \arrow[dd, "q_{\alpha}"'] \arrow[rr, "\sigma_{\alpha}"] &  & \underline{U} \arrow[dd, "q"] \\
			&  &                               \\
			\mathcal{D} \arrow[rr,"F"] \arrow[Rightarrow, shorten >=10pt, shorten <=10pt, uurr]                                                   &  & \mathcal{C}                  
		\end{tikzcd}.\]
	\end{proof}
\end{lemma}
\begin{proposition}[{\cite{Carchedi}}]
	A stack $\mc{D}$ over $\text{Top}$ with the open cover topology is a topological stack (\Cref{Definition:topologicalstack}) if and only if it has an atlas, that is, a representable epimorphism of stacks $\underline{X}\ra \mc{D}$ for an object $X$ of $\text{Top}$.
\end{proposition}

We end this section by making a note of the following properties of a topological stack.

\begin{lemma}\label{Lemma:allmapsarerepresentable}
	Let $\mc{X}$ be a topological stack. Then, for every topological space $X$, any morphism of stacks $X\ra \mc{X}$ is representable.
	\begin{proof}
		Let $X$ be a topological space, $\mc{X}$ a topological stack,  and $p:X\ra \mc{X}$  a morphism of topological stacks. As $\mc{X}$ be a topological stack, the diagonal morphism $\mc{X}\ra \mc{X}\times \mc{X}$ is a representable morphism (\cite[Proposition $7.2$]{noohi2005foundations}).  
		
		Let $Y$ be a topological space and $g: Y\ra \mc{X}$  a morphism of stacks. We show that $X\times_{\mc{X}}Y$ is a representable stack. 
		
		Consider the morphism of stacks $(p,g):X\times Y\ra \mc{X}\times \mc{X}$. As
		the diagonal morphism $\mc{X}\ra \mc{X}\times \mc{X}$ is representable, the $2$-fiber product $\mc{X}\times_{\mc{X}\times \mc{X}}(X\times Y)$ is representable (by a topological space). Moreover, by \cite[Corollary $69$]{metzler2003topological},  $\mc{X}\times_{\mc{X}\times \mc{X}}(X\times Y)\cong X\times_{\mc{X}}Y$, and hence, 
		$X\times_{\mc{X}}Y$ is representable. Thus, $p:X\ra \mc{X}$  is a representable morphism of topological stacks.
	\end{proof}
\end{lemma}	

As a consequence, we have the following result that produces an atlas from a map of local sections. We have a similar result in the case of smooth setup (\Cref{Remark:compositionisatlas}). We have omitted the proof in the smooth setup. Here we give proof for the continuous case. 
\begin{lemma}\label{Lemma:nonuniquenessofatlas}
	Let $\mc{D}$ be a topological stack and $X\ra \mc{D}$ an atlas for $\mc{D}$. Then, for any map of local sections $Y\ra X$, the composition $Y\ra X\ra \mc{D}$ is an atlas for $\mc{D}$.
	\begin{proof}
		As $\mc{D}$ is a topological stack, any morphism of stacks $\underline{P}\ra \mc{D}$ is a representable morphism. In particular,  $Y\ra \mc{D}$ is a representable morphism of stacks. 
		
		Let $M$ be a topological space and $M\ra \mc{D}$ a morphism of topological stacks. As $X\ra \mc{D}$ is an atlas, the $2$-fiber product $\underline{X}\times_{\mc{D}}\underline{M}$ is representable by a topological space and ${\rm pr}_2:\underline{X}\times_{\mc{D}}\underline{M}\ra \underline{M}$ is a map of local sections, as in the  $2$-commutative diagram,
		\[
		\begin{tikzcd}
			\underline{X}\times_{\mc{D}}\underline{M} \arrow[dd] \arrow[rr] &  & \underline{M} \arrow[dd] \\
			&  &              \\
			\underline{X} \arrow[rr]            &  & \mc{D}      
		\end{tikzcd}.\]
		Consider the pull-back of ${\rm pr}_1:\underline{X}\times_{\mc{D}}\underline{M}\ra \underline{X}$ along $\underline{Y}\ra \underline{X}$ to obtain the following diagram,
		\[
		\begin{tikzcd}
			\underline{Y}\times_{\underline{X}}\underline{X} \times_{\mc{D}}\underline{M}\arrow[rr] \arrow[dd] &  & \underline{X}\times_{\mc{D}}\underline{M} \arrow[dd] \\
			&  &              \\
			\underline{Y} \arrow[rr]                        &  & \underline{X}           
		\end{tikzcd}.\]
		As $\underline{Y},\underline{X},\underline{X}\times_{\mc{D}}\underline{M}$ are representable by topological spaces, the $2$-fiber product $\underline{Y}\times_{\underline{X}}\underline{X}\times_{\mc{D}}\underline{M}$ is representable by a topological space. 
		
		We identify $\underline{Y}\times_{\underline{X}}\underline{X}\times_{\mc{D}}\underline{M}$ with $\underline{Y}\times_{\mc{D}}\underline{M}$. As $Y\ra X$ is a map of local sections, its pull-back $\underline{Y}\times_{\mc{D}}\underline{M}\ra \underline{X}\times_{\mc{D}}\underline{M}$ will also be a map of local sections. In turn, we get the composition $\underline{Y}\times_{\mc{D}}\underline{M}\ra \underline{X}\times_{\mc{D}}\underline{M}\ra \underline{M}$ giving a map of local sections.

		Combining the above two diagrams, we have the following diagram,
		\[
		\begin{tikzcd}
			\underline{Y}\times_{\mc{D}}\underline{M} \arrow[rr] \arrow[dd] &  & \underline{M} \arrow[dd] \\
			&  &              \\
			\underline{Y} \arrow[rr]                        &  & \mc{D}      
		\end{tikzcd}.\]
		Thus, for any topological space $M$ and a morphism of topological stacks $\underline{M}\ra \mc{D}$, the $2$-fiber product $\underline{Y}\times_{\mc{D}}\underline{M}$ is representable by a topological space and the projection map $\underline{Y}\times_{\mc{D}}\underline{M}\ra \underline{M}$ is a map of local sections. Hence, $Y\ra \mc{D}$ is an atlas for $\mc{D}$.
	\end{proof}
\end{lemma}

\section{Morphism of topological groupoids and topological stacks}\label{Section:morphism of top groupoids and extensions}
\begin{definition}[{\cite{metzler2003topological}}]
	Let  $\mc{G}=[\mc{G}_1\rra \mc{G}_0]$ and $\mc{H}=[\mc{H}_1\rra \mc{H}_0]$ be a pair of topological groupoids. A \textit{$\mc{G}-\mc{H}$ bibundle} consists of \begin{itemize}
		\item a topological space $P$,
		\item a left action of $\mc{G}$ on $P$ given by the pair $(a_{\mc{G}}\colon  P\ra \mc{G}_0, \mu_{\mc{G}}\colon  \mc{G}_1\times_{s,\mc{G}_0,a_{\mc{G}}}P\ra P)$,
		\item a right action of $\mc{H}$ on $P$ given by the  pair $(a_{\mc{H}}\colon  P\ra \mc{H}_0, \mu_{\mc{H}}\colon  
		P\times_{a_{\mc{H}},\mc{H}_0,t}\mc{H}_1\ra P)$,
	\end{itemize}
	satisfying the following conditions:
	\begin{enumerate}
		\item the map $a_{\mc{G}}\colon  P\ra \mc{H}_0$ is a principal $\mc{H}$-bundle,
		\item the map $a_{\mc{H}}\colon  P\ra \mc{H}_0$ is $\mc{G}$-invariant, ; that is $a_{\mc{H}}(\gamma\cdot p)=a_{\mc{H}}(p)$ for all $(\gamma,p)\in \mc{G}_1\times_{\mc{G}_0}P$,
		\item the action of $\mc{G}$ on $P$
		is compatible with the action of $\mc{H}$ on $P$, that is 
		$\gamma\cdot(p\cdot\delta)=(\gamma\cdot p)\cdot\delta$ for each $(\gamma,p,\delta)\in \mc{G}_1\times_{s,\mc{G}_0,a_{\mc{G}}}P\times_{a_{\mc{H}},\mc{H}_0,t}\mc{H}_1$.
	\end{enumerate}
	We denote a $\mc{G}-\mc{H}$-bibundle by $P\colon\mc{G}\to{\mc{H}}$.
\end{definition}
\begin{example}\label{Example:principalbundleasbibundle}
	Let $[\mc{H}_1\rra \mc{H}_0]$ be a topological groupoid and $\pi\colon  P\ra M$ a principal $\mc{H}$-bundle. The map $\pi:P\ra M$ gives an action of the topological groupoid $[M\rra M]$ on the topological space $P$. Thus, $P$ is considered as a $[M\rra M]-\mc{H}$-bibundle. We call $P$ to be the $[M\rra M]-\mc{H}$-bibundle associated to the principal $\mc{H}$-bundle $\pi\colon  P\ra M$.
\end{example}
\begin{example}\label{Example:morphismofgroupoidsasbibundle}
	Let $(\phi_1,\phi_0)\colon  [\mc{G}_1\rra \mc{G}_0]\ra [\mc{H}_1\rra \mc{H}_0]$ be a morphism of topological groupoids. As the target map $t\colon  \mc{H}_1\ra \mc{H}_0$ is a principal $\mc{H}$-bundle, its pull-back  $\mc{G}_0\times_{\mc{H}_0}\mc{H}_1$  is given by the pull back of $t\colon  \mc{H}_1\ra \mc{H}_0$ along  $\phi_0\colon  \mc{G}_0\ra \mc{H}_0$. Then  the map $\text{pr}_1\colon  \mc{G}_0\times_{\phi_0,\mc{H}_0,t}\mc{H}_1\ra \mc{G}_0$ and  the map 
	\begin{eqnarray}
		&&\mc{G}_1\times_{s,\mc{G}_0}(\mc{G}_0\times_{\phi_0,\mc{H}_0,t}\mc{H}_1)\ra \mc{G}_0\times_{\phi_0,\mc{H}_0,t}\mc{H}_1,\nonumber\\
		&&(\gamma,(a,\delta))\mapsto (t(\gamma),\phi_1(\gamma)\circ \delta)\nonumber
	\end{eqnarray}
	define a left action of $\mc{G}$ on $\mc{G}_0\times_{\phi_0,\mc{H}_0,t}\mc{H}_1$. Whereas  $s\circ \text{pr}_2\colon  \mc{G}_0\times_{\phi_0,\mc{H}_0,t}\mc{H}_1\ra \mc{H}_0$ and the map  
	\begin{eqnarray}
		&&(\mc{G}_0\times_{\phi_0,\mc{H}_0,t}\mc{H}_1)\times_{\mc{H}_0,t}\mc{H}_1\ra \mc{G}_0\times_{\phi_0,\mc{H}_0,t}\mc{H}_1,\nonumber\\
		&&((a,\delta),\delta')\mapsto (a,\delta \circ \delta')\nonumber
	\end{eqnarray}
	define a right action of $\mc{H}$ on $\mc{G}_0\times_{\phi_0,\mc{H}_0,t}\mc{H}_1$. Then we have a $\mc{G}-\mc{H}$-bibundle  $\mc{G}_0\times_{\phi_0,\mc{H}_0,t}\mc{H}_1$. We call it the $\mc{G}-\mc{H}$-bibundle associated to the morphism of topological groupoids $(\phi_1,\phi_0)\colon  [\mc{G}_1\rra \mc{G}_0]\ra [\mc{H}_1\rra \mc{H}_0]$ and, denote it by  $\mc{G}_0\times_{\phi_0,\mc{H}_0,t}\mc{H}_1\colon  \mc{G}\ra \mc{H}$.
\end{example}

\subsection{morphism of topological stack $B\mc{G}\ra B\mc{H}$ associated to a $\mc{G}-\mc{H}$-bibundle}\label{subsection:morphismassocaitedtobibundle} Consider a $\mc{G}-\mc{H}$-bibundle  $Q\colon  \mc{G}\ra \mc{H}$. Let $M$ be a topological space, and  $\pi\colon  P\ra M$ be a principal $\mc{G}$-bundle. As mentioned in  \Cref{Example:principalbundleasbibundle}, we see $(P,\pi,M)$  as a $[M\rra M]-\mc{G}$-bibundle. The composition $P\colon  [M\rra M]\ra \mc{G}$ and $Q\colon  \mc{G}\ra \mc{H}$ then gives a $[M\rra M]-\mc{H}$ bibundle $P\circ Q\colon  [M\rra M]\ra \mc{H}$, which is nothing but a principal $\mc{H}$-bundle over a topological space $M$. This produces a morphism of stacks  $BQ\colon B\mc{G}\ra B\mc{H}$. The morphism of stacks $BQ\colon  B\mc{G}\ra B\mc{H}$ will be called the  morphism of stacks associated to the $\mc{G}-\mc{H}$-bibundle $Q\colon  \mc{G}\ra \mc{H}$. 

\subsection{a morphism of topological stacks associated to a morphism of topological groupoids} 
Let $(\phi_1,\phi_0)\colon  [\mc{G}_1\rra \mc{G}_0]\ra [\mc{H}_1\rra \mc{H}_0]$ be a morphism of topological groupoids. Combining  \Cref{Example:morphismofgroupoidsasbibundle} with  construction in  \Cref{subsection:morphismassocaitedtobibundle}, we have a   morphism of topological stacks $B\mc{G}\ra B\mc{H}$ associated to   $(\phi_1,\phi_0)\colon \mc{G}\ra \mc{H}$. This morphism  of stacks $B\mc{G}\ra B\mc{H}$ will be called the morphism of topological stacks associated to the morphism of topological groupoids $\mc{G}\ra \mc{H}$.

In this paper, we will particularly be interested in a   special type of morphism of stacks, known as a gerbe over topological stack.
\begin{definition}[a gerbe over a topological stack {\cite{noohi2005foundations}}] Let $\mc{C}$ be a topological stack. A morphism of topological stacks $F\colon  \mc{D}\ra \mc{C}$  is said to be \textit{a gerbe over the topological stack $\mc{C}$}, if 
	\begin{itemize}
		\item the morphism $F\colon  \mc{D}\ra \mc{C}$ is an epimorphism of topological stacks,
		\item the diagonal morphism $\Delta_F\colon  \mc{D}\ra \mc{D}\times_{\mc{C}}\mc{D}$ is an epimorphism of topological stacks.
	\end{itemize}
\end{definition}
With this, we conclude our review of standard material. 	
\section{Topological groupoid extension and the associated morphism of stacks}\label{Section:Topogrpdextensionmorstacks}
In the paper \cite{MR2493616}, the authors claimed without many details that there is a  possible relation between gerbes over differentiable stacks and Lie groupoid extensions. In \cite{MR4124773} we have given an explicit construction for the correspondence between differentiable gerbes and Lie groupoid extension, with some additional conditions imposed on the differentiable gerbe.  	
But, in the topological set-up, the notion of a groupoid extension has not been explored. So, in this paper, we introduce the notion of topological groupoid extensions and study the correspondence with gerbes over topological stacks.

\begin{definition}[topological groupoid extension]	
	Let $[\mc{H}_1\rra M]$ be a topological groupoid. A \textit{topological groupoid extension of $[\mc{H}_1\rra M]$} is a morphism of topological groupoids \[(F,1_M)\colon [\mc{G}_1\rra M]\ra [\mc{H}_1\rra M]\] such that $F\colon  \mc{G}_1\ra \mc{H}_1$ admits local sections. We often denote such a topological groupoid extension by $F\colon  \mc{G}_1\ra \mc{H}_1\rra M$.\end{definition}

Now, we introduce the notion of Morita morphism of Lie groupoid extensions and Morita equivalent topological groupoid extensions imitating the notion of Lie groupoids setup.

Let $[\Gamma_1\rra \Gamma_0]$ be a topological groupoid and $J:P_0\ra \Gamma_0$ be a map of local sections. Let $P_1$ denote the fiber product 
$P_0\times_{J,\Gamma_0,s}\Gamma_1\times_{t,\Gamma_0,J}P_0$. Then $[P_1\rra P_0]$ has a topological groupoid structure with source and target maps being the projection maps ${\rm pr}_1, {\rm pr}_3\colon P_1\ra P_0$ respectively. We call $[P_1\rra P_0]$ to be the \textit{pull-back groupoid} of $[\Gamma_1\rra \Gamma_0]$ along $J\,\colon\,P_0\ra \Gamma_0$.

A morphism of topological groupoids $[\Delta_1\rra \Delta_0]\ra [\Gamma_1\rra \Gamma_0]$ is said to be a \textit{Morita morphism of topological groupoids} if $\Delta_0\ra \Gamma_0$ is a map of local sections and $[\Delta_1\rra \Delta_0]$ is isomorphic to the pullback groupoid of $[\Gamma_1\rra \Gamma_0]$ along the morphism $\Delta_0\ra \Gamma_0$. 

\begin{definition}
	Let $[\Gamma_1\rra \Gamma_0]$ and $[\Delta_1\rra \Delta_0]$ be a pair of  topological groupoids. We say that $[\Gamma_1\rra \Gamma_0]$ and $[\Delta_1\rra \Delta_0]$ are \textit{Morita equivalent} if there exists a third topological groupoid $[P_1\rra P_0]$ with a pair of Morita morphisms 
	$[P_1\rra P_0]\ra [\Delta_1\rra \Delta_0]$ and $[P_1\rra P_0]\ra [\Gamma_1\rra \Gamma_0]$.
\end{definition}
\begin{definition}[Morita morphism of topological groupoid extensions]\label{Definition:moritamorphismtopologicalextension}
	Let $\phi'\colon X_1'\rightarrow  Y_1'\rightrightarrows M'$  and 
	$\phi\colon X_1\rightarrow  Y_1\rightrightarrows M$ be a pair of  topological groupoid extensions.
	A \textit{Morita morphism of topological groupoid extensions} from   $\phi'\colon X_1'\rightarrow  Y_1'\rightrightarrows M'$  to 
	$\phi\colon X_1\rightarrow  Y_1\rightrightarrows M$ is given by a pair of Morita morphisms of topological groupoids, 
	\begin{equation}\begin{tikzcd} 
			X_1' \arrow[dd,xshift=0.75ex,"t"]\arrow[dd,xshift=-0.75ex,"s"'] \arrow[rr, "\psi_X"] &  & X_1 \arrow[dd,xshift=0.75ex,"t"]\arrow[dd,xshift=-0.75ex,"s"'] \\
			&  &  \\
			M' \arrow[rr, "f"] &  & M
		\end{tikzcd}\text{ and }\begin{tikzcd} 
			Y_1' \arrow[dd,xshift=0.75ex,"t"]\arrow[dd,xshift=-0.75ex,"s"'] \arrow[rr, "\psi_Y"] &  & Y_1 \arrow[dd,xshift=0.75ex,"t"]\arrow[dd,xshift=-0.75ex,"s"'] \\
			&  &  \\
			M' \arrow[rr, "f"] &  & M
	\end{tikzcd}\end{equation}
	such that the following diagram
	\begin{equation}\begin{tikzcd} 
			X_1' \arrow[dd,"\psi_X"'] \arrow[rr, "\phi'"] &  & Y_1' 
			\arrow[dd,"\psi_Y"] \\
			&  &  \\
			X_1 \arrow[rr, "\phi"] &  & Y_1
	\end{tikzcd}\end{equation} is commutative.
\end{definition}
\begin{definition}[Morita equivalent topological groupoid extensions]\label{Definition:moritaequivalenttopologicalgroupoidextensions}
	Let $\phi'\colon X_1'\rightarrow  Y_1'\rightrightarrows M'$  and 
	$\phi\colon X_1\rightarrow  Y_1\rightrightarrows M$ be a pair of topological groupoid extensions. We say that  $\phi'\colon X_1'\rightarrow  Y_1'\rightrightarrows M'$  and 
	$\phi\colon X_1\rightarrow  Y_1\rightrightarrows M$ are  \textit{Morita equivalent topological groupoid extensions}, if there exists a third topological groupoid extension $\phi''\colon X_1''\rightarrow Y''\rightrightarrows M''$ and a pair of Morita morphisms of topological groupoid extensions
	\[(\phi''\colon X_1''\rightarrow Y''\rightrightarrows M'')\rightarrow (\phi\colon X_1\rightarrow Y_1\rightrightarrows M)\] and \[(\phi''\colon X_1''\rightarrow Y''\rightrightarrows M'')\rightarrow (\phi'\colon X_1'\rightarrow Y_1'\rightrightarrows M').\]
\end{definition}

Let $(F,1_M)\colon  [\mc{G}_1\rra M]\ra [\mc{H}_1\rra M]$ be a topological groupoid extension. This morphism of topological groupoids gives a morphism of stacks $F:B\mc{G}\to B\mc{H}$. 
We will see in the following result that this morphism $B\mc{G}\to B\mc{H}$ is a gerbe over a topological stack.

\begin{theorem}\label{Th:topgrpdextensiontogerbes}
	Let $F\colon  [\mc{G}_1\rra M]\ra [\mc{H}_1\rra M]$ be a topological groupoid extension. Then, the associated morphism of topological stacks $F\colon  B\mc{G}\ra B\mc{H}$ is a gerbe over the topological stack $B\mc{H}$.
	\begin{proof}
		{\bf Step1:} We first prove that $F\colon  B\mc{G}\ra B\mc{H}$ is an epimorphism of topological stacks. 
		
		Let $U$ be an object of $\text{Top}$ and $q\colon  \underline{U}\ra B\mc{H}$ be a morphism of topological stacks. Let $\pi:Q\ra U$ be the principal $\mc{H}$-bundle associated to the morphism of stacks $q\colon  \underline{U}\ra B\mc{H}$. 
		
		Since $\pi:Q\ra U$ admits local sections, there exists an open cover $\{U_\alpha\to U\}_{\alpha\in \Lambda}$ of $U$ and a family of maps $\{\sigma_{\alpha}:U_{\alpha}\to M\}$, such that $\pi|_{\pi^{-1}U_\alpha}\colon  \pi^{-1}(U_\alpha)\ra U_\alpha$ is pull-back of the trivial principal $\mc{H}$-bundle $t\colon  \mc{H}_1\ra M$.		
		
		Next, we pull-back the principal $\mc{G}$-bundle $t\colon  \mc{G}_1\ra M$ along the morphism $\sigma_\alpha\colon   U_\alpha\ra M$ to obtain  a principal $\mc{G}$-bundle over the topological space $U_\alpha$, for each $\alpha$. In turn we obtain  a morphism of topological stacks $q_\alpha\colon  \underline{U_\alpha}\ra B\mc{G}$. So,  given a morphism of stacks $q\colon  \underline{U}\ra B\mc{H}$, we have an open cover $\{U_\alpha\}_{\alpha\in \Lambda}$ and morphism of stacks $q_\alpha\colon  \underline{U_\alpha}\ra B\mc{G}$ forming the $2$-commutative diagram
		\begin{equation}
			\begin{tikzcd} 
				\underline{U_\alpha} \arrow[dd, "q_i"'] \arrow[rr,"\Phi=\text{inclusion}"] & & \underline{U} \arrow[dd, "q"] \\
				& & \\
				B\mc{G} \arrow[Rightarrow, shorten >=10pt, shorten <=10pt, uurr] \arrow[rr, "F"] & & B\mc{H}
			\end{tikzcd}.
		\end{equation}
		The $2$-commutativity of the diagram follows from the observation that the morphism of stacks $F\colon  B\mc{G}\ra B\mc{H}$ maps the pull-back of the principal $\mc{G}$-bundle along a map $\eta\colon  U\ra M$ to the pull-back of the principal $\mc{H}$-bundle along a map $\eta\colon  U\ra M$ (\cite[Lemma $5.5$.]{MR4124773}).
		
		We now prove that the diagonal morphism $\Delta_F\colon  B\mc{G}\ra B\mc{G}\times_{B\mc{H}}B\mc{G}$ is an epimorphism of stacks. Observe that, the $2$-fiber product $\mc{G}\times_{\mc{H}}\mc{G}$ is a topological groupoid. As the stackification and Yoneda embedding are preserved under $2$-fiber product, we see that $B(\mc{G}\times_{\mc{H}}\mc{G})\cong B\mc{G}\times_{B\mc{H}}B\mc{G}$. It turns out that the $2$-fibered product $\mc{G}\times_{\mc{H}}\mc{G}$ is a transitive topological groupoid and so is Morita equivalent to a topological groupoid $(K\rra *)$ for some topological group $K$ (\cite[Lemma $5.8$, Lemma $5.9$]{MR4124773}). So, the diagonal morphism $\Delta_F\colon  B\mc{G}\ra B\mc{G}\times_{B\mc{H}}B\mc{H}$ is equivalent to morphism of stacks $B\mc{G}\ra B(K\rra *)$. By similar justification as in the case of $B\mc{G}\ra B\mc{H}$, we can see that $B\mc{G}\ra B(K\rra *)$ is an epimorphism of stacks. Thus, the diagonal morphism $\Delta_F\colon  B\mc{G}\ra B\mc{G}\times_{B\mc{H}}B\mc{G}$ is an epimorphism of stacks. In conclusion, the morphism of topological stacks associated to a topological groupoid extension $\mc{G}\ra \mc{H}$ is a gerbe over the topological stack $B\mc{H}$.
	\end{proof}
\end{theorem}	

\begin{theorem}\label{th:gerbetotopgrpdexten}
	Let $\mc{D},\mc{C}$ be topological stacks and $F\colon  \mc{D}\ra \mc{C}$ a gerbe over the topological stack $\mc{C}$. Further assume that the diagonal morphism $\Delta_F:\mc{D}\ra \mc{D}\times_{\mc{C}}\mc{D}$ is a representable map of local sections. Then, there exists a topological groupoid extension $\mc{G}\ra \mc{H}$ inducing the morphism $\mc{D}\ra \mc{C}$. Further, this topological groupoid extension is unique up to Morita equivalence.
	\begin{proof}
		As $F\colon  \mc{D}\ra \mc{C}$ is an epimorphism of topological stacks, given an object $U$ of $\text{Top}$ and a morphism of topological stacks $\underline{U}\ra \mc{C}$, there exists a map of local sections $V\ra U$ and a morphism of stacks $\underline{V}\ra \underline{U}$ with the following $2$-commutative diagram, 
		\[\begin{tikzcd} 
			\underline{V} \arrow[dd, "p"'] \arrow[rr, "\pi"] &  & \underline{U} \arrow[dd, "\tilde{q}"] \\
			&  &                               \\
			\mathcal{D} \arrow[rr,"F"] \arrow[Rightarrow, shorten >=10pt, shorten <=10pt, uurr]                                                   &  & \mathcal{C}                  
		\end{tikzcd}.\]
		Suppose that $\underline{U}\ra \mc{C}$ is an atlas for the stack $\mc{C}$, then, as $V\ra U$ is a map of local sections, the composition $\underline{V}\ra \underline{U}\ra \mc{C}$ is an atlas for the stack $\mc{C}$ (\Cref{Lemma:nonuniquenessofatlas}). Thus, for the morphism of stacks $F\colon  \mc{D}\ra \mc{C}$ there exists an atlas $q\colon  \underline{V}\ra \mc{C}$ for $\mc{C}$ and a morphism of stacks $p\colon \underline{V}\ra\mc{D}$ with the following $2$-commutative diagram  
		\[
		\begin{tikzcd} 
			V \arrow[dd, "p"'] \arrow[rrdd, "q"] &  &        \\
			&  &        \\
			\mc{D} \arrow[rr, "F"]               &  & \mc{C}
		\end{tikzcd}.\]
		The condition that $\mc{D}\ra \mc{D}\times_{\mc{C}}\mc{D}$ is an epimorphism implies that the morphism $p\colon\underline{V}\ra \mc{D}$ is an epimorphism of stacks. On the other hand  the morphism $p:\underline{V}\ra \mc{D}$ is an atlas for the stack $\mc{D}$ follows from  the fact that $\mc{D}\ra \mc{D}\times_{\mc{C}}\mc{D}$ is a representable  map of local sections. The atlases $p:\underline{V}\ra \mc{D}$ and $q:\underline{V}\ra \mc{C}$  respectively produce the topological groupoids $(V\times_{\mc{D}}V\rra V)$ and $(V\times_{\mc{C}}V\rra V)$ representing the stacks $\mc{D}$ and $\mc{C}$ respectively. The morphism of stacks $F:\mc{D}\ra \mc{C}$ will then induce a morphism of topological groupoids $(f,1_V)\colon(V\times_{\mc{D}}V\rra V)\ra (V\times_{\mc{C}}V\rra V)$. The condition that $\Delta_F\colon\mc{D}\ra \mc{D}\times_{\mc{C}}\mc{D}$ is a representable morphism of local sections further impose the condition that $f\colon V\times_{\mc{D}}V\ra V\times_{\mc{C}}V$ is a map of local sections. Thus, we get a topological groupoid extension associated to a gerbe over a topological stack $\mc{D}\ra \mc{C}$ satisfying the extra condition that $\Delta_F$ is a representable morphism of local sections. The Morita invariance follows from the corresponding argument in \cite{MR4124773}.
	\end{proof}
\end{theorem}

\section{Serre, Hurewicz fibrations and gerbes}\label{Section:Serre-hurewiczfibrations}

In this section, we will recall the notions of Serre and Hurewicz stacks introduced in \cite{MR3144243} and study certain properties of a  gerbe over such stacks.

\begin{definition}[\cite{MR3144243}]A topological stack $\mc{D}$ is said to be a \textit{Hurewicz stack} (respectively, \textit{Serre stack}) if it admits a presentation $\mb{X}=[R\rra X]$ by a topological groupoid in which the source (hence also the target) map $s\colon  R\ra X$ is locally a Hurewicz fibration (respectively, Serre fibration).
\end{definition} 
Let $(\mc{G}_1\rra M)$ be a topological groupoid with the source map $s\colon  \mc{G}_1\ra M$ being a locally a Hurewicz fibration (respectively, Serre fibration). Then $B\mc{G}$ is a  Hurewicz stack (respectively, Serre stack).
\begin{definition}[\cite{MR3144243}]
	A morphism of topological stacks $F\colon  \mc{D}\ra \mc{C}$ is said to be a Hurewicz (respectively, Serre) morphism if for every object $U$ of $\text{Top}$ and a morphism of topological stacks $U\ra \mc{C}$, the fiber product $\mc{D}\times_{\mc{C}}U$ is a Hurewicz (respectively, Serre) stack. 
\end{definition}	
We say that a gerbe $F\colon\mc{D}\ra \mc{C}$ is a {\it{Hurewicz gerbe}} (respectively, a {\it{Serre gerbe}}), if the underlying morphism of stacks is Hurewicz (respectively, Serre).

A straightforward consequence of the results \Cref{Th:topgrpdextensiontogerbes} and { \cite[Lemma 2.4.]{MR3144243}} is the following proposition. 
\begin{proposition}
	Let $\mc{G}=(\mc{G}_1\rra M)$ be a topological groupoid with  the source map $s\colon  \mc{G}_1\ra M$ a locally a Hurewicz fibration (respectively, Serre fibration) and $F\colon  (\mc{G}_1\rra M)\ra (\mc{H}_1\rra M)$ a topological groupoid extension. Then the gerbe  $F\colon B\mc{G}\ra B{\mc{H}}$ constructed in \Cref{Th:topgrpdextensiontogerbes} is a Hurewicz gerbe (respectively, Serre gerbe).
\end{proposition}

We give an alternate proof of {\cite[Lemma 2.4.]{MR3144243}}.

\begin{lemma} {\cite{MR3144243}}
	Let $f:\mc{X}\ra \mc{Y}$ be a morphism of topological stacks. If $\mc{X}$ is a Serre or Hurewicz stack, then $f:\mc{X}\ra \mc{Y}$ is a Serre or Hurewicz morphism of stacks.
	\begin{proof} Let $\mc{X}$ be a Serre stack, that is, there exists a topological space $X$ and an atlas $X\ra \mc{X}$ such that the source and target maps of the associated topological groupoid $[X\times_{\mc{X}}X\rra X]$ are (locally) Serre fibrations of topological spaces.
		
		Let $Y$ be a topological space and $g: Y\ra \mc{Y}$ a morphism of topological stacks. We show that the fiber product $\mc{X}\times_\mc{Y}Y$ is representable by a topological groupoid $[M\rra N]$ whose source and target maps are (locally) Serre fibrations of topological spaces.
		
		Let $R$ denote the fiber product $X\times_{\mc{X}}X$, and $R\ra \mc{Y}$ be the composition $R\ra \mf{X}\ra \mc{Y}$. Let $X\ra \mc{Y}$ denote the composition $X\ra \mc{X}\ra \mc{Y}$.  Consider the following pull-back diagrams,
		\[
		\begin{tikzcd} 
			R\times_{\mc{Y}}Y \arrow[dd] \arrow[rr] &                                         & Y \arrow[dd] \\
			& X\times_{\mc{Y}}Y \arrow[dd] \arrow[ru] &              \\
			R \arrow[rr]                            &                                         & \mc{Y}       \\
			& X \arrow[ru]                            &             
		\end{tikzcd}.\]
		As $\mc{Y}$ is a topological stack, the map $R\ra \mc{Y}$ is representable (\Cref{Lemma:allmapsarerepresentable}), thus, the fiber product $R\times_{\mc{Y}}Y$ is representable by topological space. For the same reason, $X\ra \mc{Y}$ is representable, thus, the fiber product $X\times_{\mc{Y}}Y$ is representable by topological space.  
		
		By pulling back  $s:R\ra X$  along the map $X\times_{\mc{Y}}Y\ra X$ we obtain a map $R\times_{\mc{Y}}Y\ra X\times_{\mc{Y}}Y$. As $s:R\ra X$ is a Serre fibration, the pull-back  $R\times_{\mc{Y}}Y\ra X\times_{\mc{Y}}Y$ is a (locally) Serre fibration. Similarly, $t:R\ra X$ (a (locally) Serre fibration of spaces) is pulled back to give a morphism $R\times_{\mc{Y}}Y\ra X\times_{\mc{Y}}Y$ (which would be a (locally) Serre fibration of spaces). Similarly, the other structure maps on $[R\rra X]$ is pulled back to produce a topological groupoid $[R\times_{\mc{Y}}Y\rra X\times_{\mc{Y}}Y]$. It turns out that $\mc{X}\times_{\mc{Y}}Y$ is representable by $[R\times_{\mc{Y}}Y\rra X\times_{\mc{Y}}Y]$. Thus, $\mc{X}\times_{\mc{Y}}Y$ is a Serre stack. The same proof goes verbatim if we consider $\mc{X}$ to be a Hurewicz stack instead of Serre stack. 
	\end{proof}
\end{lemma}
An immediate outcome of the above lemma is the following. 
\begin{corollary}
	Any morphism of stacks $B[X\rra X]\ra B\mc{H}$ induced by a morphism of topological groupoids $[X\rra X]\ra \mc{H}$ is a Serre morphism as well as a Hurewicz morphism. 
\end{corollary}

	\chapter*{Publications arising out of the PhD thesis\hfill} \addcontentsline{toc}{chapter}{Publications arising
		out of the PhD thesis}
	\begin{enumerate}
		\item [{[1]}] (j/w Saikat Chatterjee) On two notions of a Gerbe over a stack \cite{MR4124773}.
		\item [{[2]}] (j/w Indranil Biswas, Saikat Chatterjee and Frank Neumann) Chern-{W}eil theory for principal bundles over {L}ie groupoids \cite{biswas2020chern}.
		\item [{[3]}] (j/w Indranil Biswas, Saikat Chatterjee and Frank Neumann) Atiyah sequences and connections on principal bundles over differentiable stacks \cite{biswas2020atiyah}.
		\item [{[4]}] (j/w Saikat Chatterjee) Extension of topological groupoids and Serre, Hurewicz morphisms \cite{chatterjee2021extension}.
	\end{enumerate}
	
	\bibliography{references}
	\bibliographystyle{plain}
	
	\cleardoublepage
	\singlespacing
	\printindex 
	
\end{document}